\documentclass[11pt,reqno]{article}

\usepackage{amsmath,amssymb,amsfonts,amsthm,subfloat,algorithm,url,float,epsf,psfrag,bm}
\usepackage{mathtools}

\usepackage[pdftex]{color,graphicx}
\usepackage[colorlinks,citecolor=blue,urlcolor=blue]{hyperref}

\usepackage[noabbrev,capitalize]{cleveref}

\usepackage{algorithm, algorithmic, verbatim, float}
\usepackage{marvosym}
\usepackage{titlesec}
\titlespacing*{\section}{0pt}{0.6ex plus 0.6ex minus .2ex}{0.6ex plus .2ex}
\titlespacing*{\subsection}{0pt}{0.6ex plus 0.6ex minus .2ex}{0.6ex plus .2ex}
\usepackage[pagewise]{lineno}

\usepackage[numbers]{natbib}
\usepackage{setspace,caption}
\usepackage{verbatim}
\usepackage{array, float}
\usepackage{xcolor}
\usepackage{fontenc}
\usepackage[toc,page]{appendix}
\usepackage[nottoc]{tocbibind}
\usepackage[english]{babel}
\usepackage{relsize}
\usepackage{cellspace}
\usepackage{booktabs}
\usepackage{blindtext}
\usepackage{subcaption}
\usepackage[utf8]{inputenc}

\usepackage{mathrsfs}
\usepackage{amsfonts}
\usepackage{enumerate}
\usepackage{latexsym}
\usepackage{multirow}
\allowdisplaybreaks
\newtheoremstyle{exampstyle}
{8pt} 
{8pt} 
{\it} 
{} 
{\bfseries} 
{.} 
{.5em} 
{} 

\theoremstyle{exampstyle}
\usepackage{bbm}
\newtheorem{theorem}{Theorem}
\newtheorem{example}{Example}
\newtheorem{lemma}{Lemma}
\newtheorem{corollary}{Corollary}
\newtheorem{remark}{Remark}

\newtheorem{prop}{Proposition}

\newtheorem{definition}{Definition}

\numberwithin{equation}{section}
\numberwithin{example}{section}
\numberwithin{theorem}{section}
\numberwithin{lemma}{section}
\numberwithin{corollary}{section}
\numberwithin{prop}{section}
\numberwithin{definition}{section}
\numberwithin{remark}{section}

\newcommand{\eat}[1]{}

\DeclareMathOperator*{\argmin}{\arg\!\min}
\DeclareMathOperator*{\argmax}{\arg\!\max}

\newcommand{\keywords}[1]{\textbf{\textit{Keywords:}} #1}
\renewcommand{\bar}[1]{\overline{#1}}
\renewcommand{\hat}[1]{\widehat{#1}}
\renewcommand{\tilde}[1]{\widetilde{#1}}

\newcommand{\R}{\mathbb{R}}

\newcommand{\Rdcov}{\mathrm{RdCov}}

\newcommand{\Rdcorr}{\mathrm{RdCorr}}
\newcommand{\Ren}{\mathrm{RE}}

\definecolor{LightCyan}{rgb}{0.88,1,1}
\definecolor{Gray}{gray}{0.9}

\usepackage{xr}
\usepackage{enumerate}

\newcommand{\blind}{1}

\begin{document}
	
	
	
	\if1\blind
	{
		\title{\LARGE \bf Multivariate Rank-based Distribution-free Nonparametric Testing using Measure Transportation}
		\author{Nabarun Deb\thanks{E-mail: nd2560@columbia.edu}\hspace{.2cm}\\
			\small{Department of Statistics, Columbia University}\\
			Bodhisattva Sen\thanks{Supported by NSF Grants DMS-17-12822 and AST-16-14743; e-mail: bodhi@stat.columbia.edu} \hspace{.2cm}\\
			\small{Department of Statistics, Columbia University}     
		}
		\maketitle
	} \fi
	
	\if0\blind
	{
		\bigskip
		\bigskip
		\bigskip
		\begin{center}
			{\LARGE\bf Multivariate Rank-based Distribution-free Nonparametric Testing using Measure Transportation}
		\end{center}
		\medskip
	} \fi
	\renewcommand{\baselinestretch}{1} 
	\bigskip
	\begin{abstract}
		In this paper, we propose a general framework for distribution-free nonparametric testing in multi-dimensions, based on a notion of multivariate ranks defined using the theory of measure transportation. Unlike other existing proposals in the literature, these multivariate ranks share a number of useful properties with the usual one-dimensional ranks; most importantly, these ranks are distribution-free. This crucial observation allows us to design nonparametric tests that are exactly distribution-free under the null hypothesis. We demonstrate the applicability of this approach by constructing exact distribution-free tests for two classical nonparametric problems: (I) testing for mutual independence between random vectors, and (II) testing for the equality of multivariate distributions. In particular, we propose (multivariate) rank versions of distance covariance (\citet{Gabor2007}) and energy statistic (\citet{Gabor2013}) for testing scenarios (I) and (II) respectively. In both these problems we derive the asymptotic null distribution of the proposed test statistics. We further show that our tests are consistent against all fixed alternatives. Moreover, the proposed tests are tuning-free, computationally feasible and are well-defined under minimal assumptions on the underlying distributions (e.g., they do not need any moment assumptions). We also demonstrate the efficacy of these procedures via extensive simulations. In the process of analyzing the theoretical properties of our procedures, we end up proving some new results in the theory of measure transportation and in the limit theory of permutation statistics using Stein's method for exchangeable pairs, which may be of independent interest.
	\end{abstract}
	
	
	\section{Introduction}\label{sec:intro}
	
	Let us consider the following two classical multivariate nonparametric hypothesis testing problems:
	
	\noindent \textbf{(I)} \textbf{Testing for mutual independence}: Given independent observations from a distribution $G$ on $\mathbb{R}^d$, $d=d_1+d_2$, $d_1, d_2\geq 1$, let $G_1$ and $G_2$ denote the marginals of $G$ corresponding to the first $d_1$ and last $d_2$ components respectively. Then, the  problem of \emph{mutual independence} testing reduces to $$\mathrm{H}_0: G=G_1\otimes G_2 \qquad \mathrm{versus} \qquad \mathrm{H}_1:G\neq G_1\otimes G_2$$ where by $G_1\otimes G_2$ we mean the product of the marginal distributions $G_1$ and $G_2$. A natural extension of this problem is to test for the mutual independence of $K$ marginals, with $K \ge  2$. The independence testing problem has found applications in a wide variety of disciplines such as in statistical genetics~\cite{liu2010versatile}, marketing and finance~\cite{grover1985probabilistic}, survival analysis~\cite{martin2005testing}, ecological risk assessment~\cite{dishion1999middle}, independent component analysis~\cite{lu2009financial}, etc., and has consequently inspired a long line of research over the past century (see e.g.,~\cite{puri1971nonparametric,gieser1997},~\cite[Chapters 1 and 8]{hollander2013nonparametric} and the references therein).
	
	\noindent \textbf{(II)} \textbf{Testing for equality of distributions}: Given independent observations from two multivariate distributions, say $F_1$ and $F_2$ on $\mathbb{R}^d$, $d\geq 1$, the nonparametric \emph{two-sample goodness-of-fit} testing problem can be formulated as $$\mathrm{H}_0:F_1=F_2  \qquad \mathrm{versus} \qquad \mathrm{H}_1:F_1\neq F_2.$$ The above problem can also be extended to the \emph{$K$-sample} setup ($K \ge 2$) when one observes independent samples from $K$ distributions and the goal is to nonparametrically test the equality of all the $K$ distributions. The two-sample (or $K$-sample) problem also has numerous applications, e.g., in pharmaceutical studies~\cite{farris1999between}, causal inference~\cite{folkes1987field}, remote sensing~\cite{conradsen2003test}, econometrics~\cite{mayer1975selecting}, etc., and has been studied extensively (see e.g.,~\cite{bickel1969distribution,weiss1960two},~\cite{hollander2013nonparametric} and the references therein).
	
	In this paper we mainly study the above two problems and develop \emph{nonparametric} testing procedures that are \emph{exactly distribution-free} (i.e., the null distributions of the test statistics are free of the underlying (unknown) data generating distributions, for all sample sizes), \emph{computationally feasible} and are \emph{consistent} against all fixed alternatives (i.e., the probability of rejecting the null, calculated under the alternative, converges to $1$ as the sample size increases). In fact, we develop a general framework for multivariate distribution-free nonparametric testing applicable much beyond the above two examples. To the best of our knowledge, the test proposed in this paper in the context of testing mutual independence is the first and only nonparametric test that guarantees the three aforementioned desirable properties. In the multivariate two-sample setting, the only other test with the above properties is due to Rosenbaum \cite{rosenbaum2005}; also see~\cite{munmun2016,bbbm2019}. 
	
	To construct our finite sample distribution-free tests we use a suitable notion of {\it multivariate ranks} (obtained from the theory of measure transportation, to be discussed below) which are themselves distribution-free. This is analogous to what is usually done in one-dimensional problems. Let us illustrate this principle in the context of testing for mutual independence (problem \textbf{(I)}). When $d_1=d_2 = 1$, the classical product-moment correlation --- which mainly captures linear dependence between the variables --- can be used to test this hypothesis. However, the exact distribution of the Pearson correlation coefficient, under $\mathrm{H}_0$, depends on the marginals $G_1$ and $G_2$. This gave way to Spearman's rank-correlation (another related measure is Kendall's $\tau$ coefficient; also see~\cite{Kendall1990,pearson1920notes,gibbons2011nonparametric}) which calculates the product-moment correlation between the one-dimensional ranks of the variables. Consequently the resulting test is distribution-free under the null hypothesis of mutual independence and can deal with non-linear (monotone) dependencies. Note that the use of ranks to obtain distribution-free tests is ubiquitous in one-dimensional problems in nonparametric statistics --- e.g., two-sample Kolmogorov-Smirnov test~\cite{smirnoff1939}, Wilcoxon signed-rank test~\cite{Wilcoxon1947}, Wald-Wolfowitz runs test~\cite{wald1940}, Mann-Whitney rank-sum test~\cite{mann1947}, Kruskal-Wallis test~\cite{kruskal1952}, Hoeffding's $D$-test~\cite{hoeffding1948non}, etc.

	
	In the $d$-dimensional Euclidean space, for $d \ge 2$, due to the absence of a canonical ordering, the existing extensions of concepts like ranks (such as component-wise ranks,  e.g.,~\cite{Bickel1965,Puri1965}; spatial ranks, e.g.,~\cite{Chaudhuri1996,Marden1999}; depth-based ranks, e.g.,~\cite{liu1993,Zuo2000}; and Mahalanobis ranks, e.g.,~\cite{Hallin2004,Davy2006}) and the corresponding rank-based tests no longer possess exact distribution-freeness. This raises a fundamental question: ``How do we define multivariate ranks that can lead to distribution-free testing procedures?". A major breakthrough in this regard was very recently made in the pioneering work of Marc Hallin and co-authors (\cite{del2018center,chernozhukov2017monge}) where they propose a notion of multivariate ranks, based on the theory of measure transportation, that possesses many of the desirable properties present in their one-dimensional counterparts. \par 
	
\begin{figure}[h]
		\begin{center}
			\includegraphics[height=6.5cm,width=7.5cm]{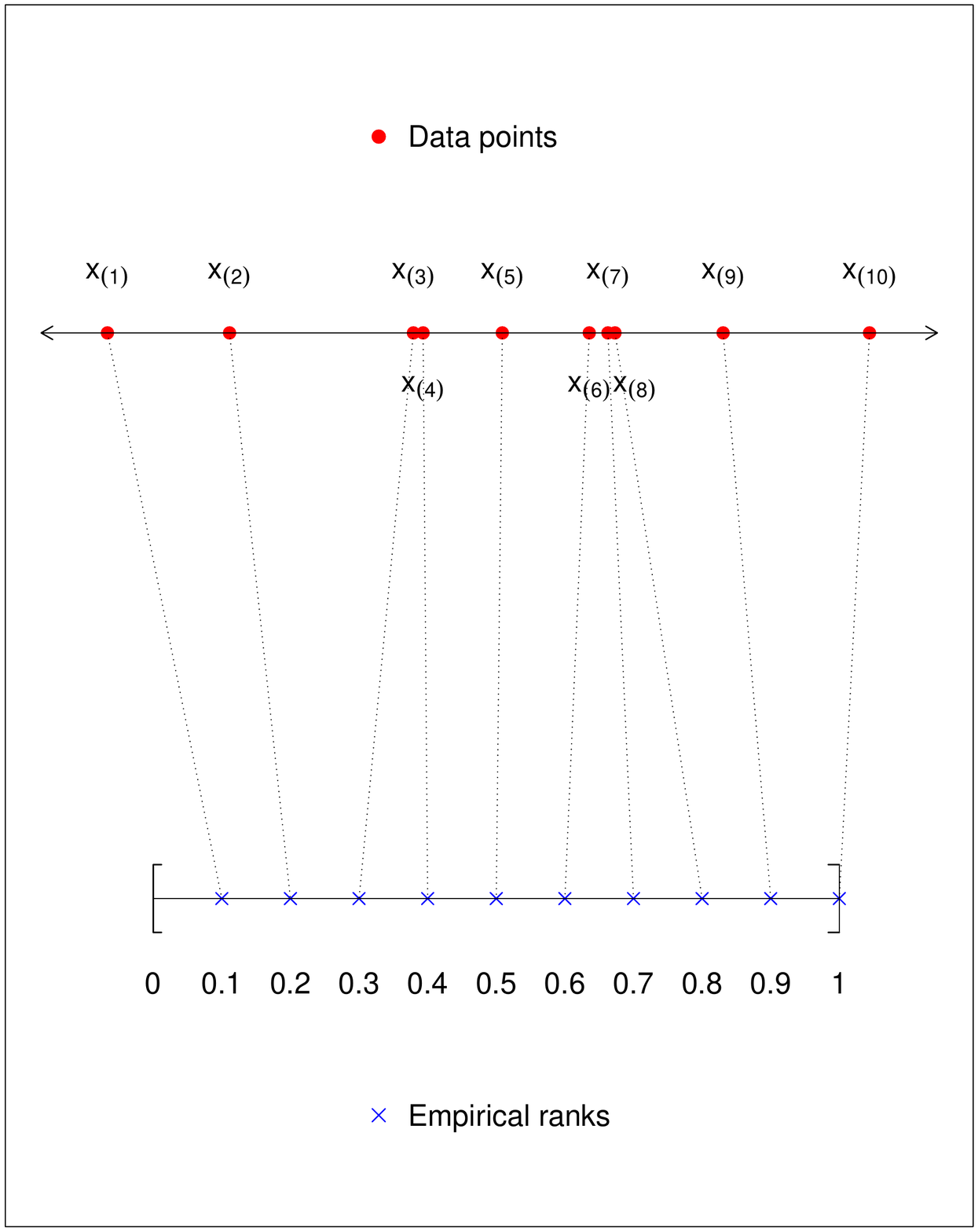}	
			\includegraphics[height=6.5cm,width=7.5cm]{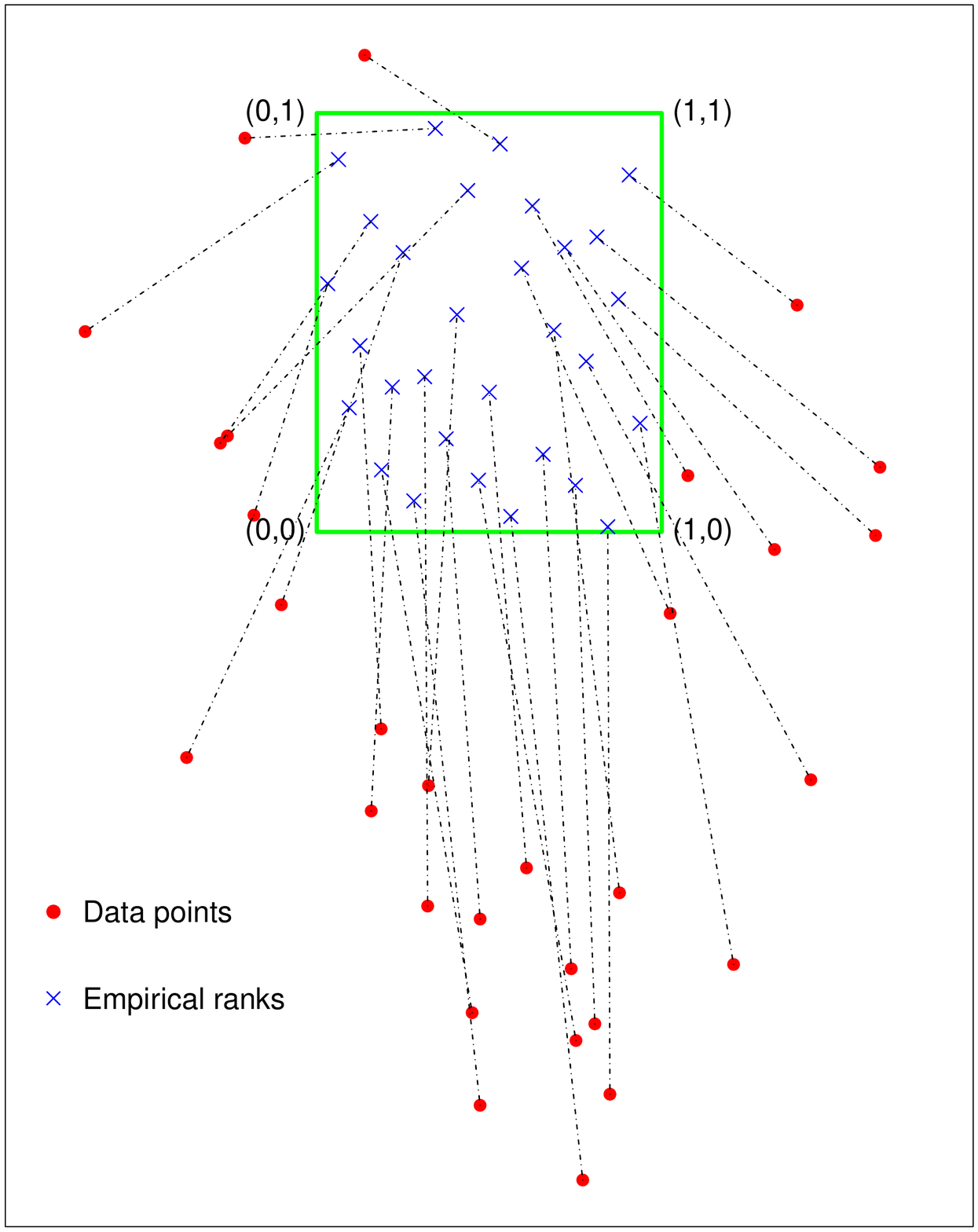}
		\end{center}	
		\caption{The left panel illustrates the correspondence between univariate data points and their ranks (which are the points $i/n$, for $i = 1,\ldots, n=10$). The right panel shows the similar correspondence between bivariate data points and their bivariate ranks which are now pseudo-random numbers in the unit square $[0,1]^2$. The rank of a data point (in solid red) is given by the blue cross at the other end of the dashed line joining them. Note that the points near the center of the data distribution are mapped close to $(1/2,1/2)$ whereas the points closer to the extremes of the data cloud are mapped to the corresponding extreme regions of the unit square, thereby giving rise to a natural bivariate ordering of the data points.}
		\label{fig:unimultiranknotion}
	\end{figure}

	In order to motivate this notion of multivariate ranks, let us start with the following interpretation of the one-dimensional ranks. Given a collection of $n$ i.i.d.~random variables $X_1,\ldots, X_n$ on $\R$ (having a continuous distribution) the {\it rank map} assigns these observations to elements of the set $\{1/n,2/n,\ldots, n/n\}$ (or $1,2,\ldots, n$, depending on interpretation) by solving the following optimization problem:
	\begin{equation}\label{eq:intro1}
	\hat{\sigma}\coloneqq \argmin_{\sigma = (\sigma(1), \ldots, \sigma(n))\,\in \, {S}_n} \sum_{i=1}^n \Big|X_{i}-\frac{\sigma(i)}{n} \Big|^2=\argmax_{\sigma = (\sigma(1), \ldots ,\sigma(n))\,\in\,{S}_n}\sum_{i=1}^n   \sigma(i) X_{i}
	\end{equation}
	where $S_n$ is the set of all permutations of $\{1,2,\ldots ,n\}$ (see~\cite[Chapter 1]{Villani2003}). It is not difficult to check (by using the rearrangement inequality, see e.g.,~\cite[Theorem 368]{Hardy1952}) that $\hat{\sigma}(i)/n$ (or simply $\hat{\sigma}(i)$) will equal the rank of $X_i$, for $i=1,\ldots, n$; see the left panel of~\cref{fig:unimultiranknotion}.
	
	Note that~\eqref{eq:intro1} can be readily extended to the multivariate setting where the discrete uniform numbers $\{i/n: 1\le i\leq n\}$ are replaced by the set of multivariate rank vectors $\{\mathbf{c}_1,\ldots ,\mathbf{c}_n\} \subset [0,1]^d$ --- a sequence of ``uniform-like" points in $[0,1]^d$ (see Section~\ref{sec:compare} for other choices of reference distributions; also see~\cite{del2018center, chernozhukov2017monge, boeckel2018multivariate}). In this paper we consider $\{\mathbf{c}_i: 1\le i\leq n\}$ as a quasi-Monte Carlo sequence --- in particular, we advocate the use of Halton sequences and employ it in our simulation experiments; other natural choices like the equally-spaced $d$-dimensional lattice are also possible (see Section~\ref{sec:revQMC} for a detailed discussion). Specifically, given i.i.d.~random vectors $\mathbf{X}_1,\ldots ,\mathbf{X}_n$ on $\R^d$, we consider the following optimization problem:
	\begin{equation}\label{eq:intro1-1}
	\hat{\sigma}\coloneqq \argmin_{\sigma = (\sigma(1), \ldots, \sigma(n))\,\in\, {S}_n} \sum_{i=1}^n \|\mathbf{X}_{i}- \mathbf{c}_{\sigma(i)}\|^2
	\end{equation}
	where, as before, the optimization is over $S_n$, the set of all permutations of $\{1,2,\ldots ,n\}$, and $\|\cdot\|$ denotes the usual Euclidean norm in $\R^d$. Note that~\eqref{eq:intro1-1} can be viewed as an \emph{assignment problem} (see e.g.,~\cite{munkres1957,bertsekas1988}) for which algorithms with worst case complexity $\mathcal{O}(n^3)$ are available in the literature (see~\cref{sec:compute} for a discussion). Based on~\eqref{eq:intro1-1}, one can then define the multivariate rank of $\mathbf{X}_{i}$ as $\mathbf{c}_{\hat{\sigma}(i)}$. This is illustrated in the right panel of~\cref{fig:unimultiranknotion} where the dashed lines join the data points (in red) with the corresponding rank vectors (indicated by blue crosses). 
	
	The above optimization problem (see~\eqref{eq:intro1-1}) indeed results in a distribution-free notion of empirical multivariate ranks as we demonstrate in~\cref{prop:finproperties} (also see~\cite[Proposition 1.6.1]{del2018center}). Note that~\eqref{eq:intro1-1} is connected to the theory of optimal measure transportation as we are ``transporting'' the empirical distribution of the $\mathbf{X}_i$'s to the empirical distribution of $\mathbf{c}_i$'s. 
	We review this literature and build on the work of~\cite{del2018center} in Sections~\ref{sec:revOT} and~\ref{sec:defn}.

	Having defined a suitable notion of multivariate ranks, the next natural question becomes: ``How does one use these multivariate ranks for nonparametric testing?''. In this regard we have a general yet powerful recipe: Given a set of multivariate observations for a nonparametric testing problem (e.g., \textbf{(I)} or \textbf{(II)}), define their multivariate ranks in such a way (depending on the problem) so that the distribution of these ranks is exactly universal (free of the data generating distribution(s)) under $\mathrm{H}_0$.  Next, take a ``good" test statistic for the corresponding nonparametric testing problem (which may not be distribution-free under $\mathrm{H}_0$). Then form a new test by evaluating the original test statistic on these obtained multivariate ranks instead of the data points themselves. Clearly, this will result in a distribution-free test statistic. We believe that this approach is quite general and can consequently be used in a variety of multivariate nonparametric inference problems, much beyond the two problems \textbf{(I)} and \textbf{(II)} discussed above (see~\cref{sec:otherp} for more on this). As we have observed before, this prescription indeed yields the Spearman's rank correlation coefficient when applied to the usual product-moment correlation for testing mutual independence when $d_1=d_2=1$.

	Let us now illustrate the above idea through a concrete application, namely, the problem of testing for mutual independence (i.e., problem \textbf{(I)}). Over the last 2-3 decades a plethora of nonparametric testing procedures have been proposed for this problem in the multivariate setting; see e.g.,~\cite{Gabor2007,berrett2017nonparametric,gretton2008kernel,heller2013,munmun2016,gieser1997,taskinen2005,oja2010,friedman1983} and the references therein. One particular testing procedure, namely {\it distance covariance} (introduced in~\cite{Gabor2007}; also see~\cite{bakirov2006}), has received much attention recently, mainly due to its simplicity and good power properties. Let us briefly describe this procedure. As the name suggests, it simply computes the covariance between pairwise distances. In particular, given the random sample $\{(\mathbf{X}_i,\mathbf{Y}_i)\}_{i=1}^n$ where $\mathbf{X}_i\in\mathbb{R}^{d_1}$ and $Y_i\in\mathbb{R}^{d_2}$, we compute the Euclidean
	distance matrices $(a_{kl})_{k,l=1}^n := \{\lVert \mathbf{X}_k - \mathbf{X}_l\rVert\}$ and $(b_{kl})_{k,l=1}^n := \{\lVert \mathbf{Y}_k - \mathbf{Y}_l\rVert\}$. We further define the (double) centered version of the $a_{kl}$'s as $A_{kl} \coloneqq a_{kl} - \bar a_{k\cdot} - \bar a_{\cdot l} + \bar a_{\cdot\cdot}$, for $k,l = 1,\ldots, n$, where $\bar a_{k\cdot} \coloneqq n^{-1}\sum_{l=1}^n a_{k l}$, $\bar a_{\cdot l} \coloneqq n^{-1} \sum_{k=1}^n  a_{k l}$, and $\bar a_{\cdot\cdot} \coloneqq n^{-2}\sum_{k,l=1}^n a_{k l}$. Similarly, we define $B_{kl} : = b_{kl} - \bar b_{k\cdot} - \bar b_{
		\cdot l} + \bar b_{\cdot \cdot}$, for $k,l = 1,\ldots,n$. Then the sample \emph{distance covariance} is defined as 
	\begin{equation}\label{eq:firstdcov}\mathrm{dCov}_n^2(\mathbf{X},\mathbf{Y}) \coloneqq \frac{1}{n^2} \sum_{k,l=1}^n A_{k l} B_{k l}.
	\end{equation} 
	The distance covariance based test has many appealing properties (see e.g.,~\cite{Mori2019}), is computationally simple  and is consistent against all alternatives that have a finite mean. 
	However, the test based on $\mathrm{dCov}_n$ is not distribution-free. In fact, as far as we are aware, till date there are no distribution-free testing procedures for \textbf{(I)} that guarantee consistency against all alternatives, under minimal assumptions.

	We introduce and study {\it rank distance covariance} in Section~\ref{sec:Indep}, where we replace the $\mathbf{X}_i$'s and $\mathbf{Y}_i$'s above by their marginal multivariate ranks (defined using~\eqref{eq:intro1-1}). This automatically yields a distribution-free nonparametric test for mutual independence. We also introduce a ``population" version of the rank-based distance covariance in~\cref{sec:depmes} and explore its connection with Spearman's rank correlation when $d=1$ in~\cref{prop:srcc}. In~\cref{lem:proprdcov} we demonstrate some basic desirable properties of this measure which has interesting connections to the properties of usual distance covariance (as shown in~\cite{Mori2019}); in particular we show that rank distance covariance also characterizes independence.
	
	In~\cref{rem:indepdistfree} we show that our proposed rank distance covariance test is exactly distribution-free as soon as the two marginal distributions are absolutely continuous. In fact, when $d_1=d_2=1$, we show in~\cref{lem:dcovequiv} that our proposed test is exactly equivalent to a modification of the celebrated Hoeffding's $D$-statistic (\cite{hoeffding1948non}) --- one of the first nonparametric tests for mutual independence. We further demonstrate, in~\cref{theo:indepconsis}, that our proposed test is consistent (i.e., has asymptotic power 1) as soon as the two marginals are absolutely continuous. In fact, we do not even need the underlying distributions to have finite means for this result (cf.~with usual distance covariance). We also go a step further and obtain the asymptotic distributional limit of our test statistic, under $\mathrm{H}_0$, in~\cref{theo:indepasdistn}. This result further demonstrates that the asymptotic limit of our test statistic does not depend on the underlying data generating distribution and is invariant to the choice of the sequence $\{\mathbf{c}_n\}_{n \ge 1}$ --- the multivariate ranks. 
	
	In~\cref{sec:twogofp} we study the problem of testing for the equality of two multivariate distributions (i.e., problem \textbf{(II)}) and propose a test for this goodness-of-fit using the {\it rank energy statistic} which is based on the usual energy statistic (as in~\cite{Gabor2013}, also see~\cite{Baringhaus2004,rizzo2005} for definitions and motivation). Similar to distance covariance, the energy statistic is also based on pairwise distances and is extremely easy to compute. The energy statistic equals $0$ if and only if the two underlying distributions are the same, as long as the two distributions have finite means. The energy test has also attracted a lot of attention recently in a variety of applications, see e.g., in robust statistics~\cite{klebanov2002class}, microarray data analysis~\cite{xiao2004multivariate}, material structure analysis~\cite{benevs2009statistical}, etc. 
	
	We demonstrate the distribution-free nature of our proposed test statistic for problem \textbf{(II)} in~\cref{rem:twosamdistfree}. An interesting property of this proposed statistic is that it is exactly equivalent to the famous two-sample Cram\'{e}r-von Mises statistic (see e.g.,~\cite{AndersonCVM1962}) when $d=1$. We explain this connection in~\cref{lem:gofeqcvm}. We further prove the consistency and derive the asymptotic distribution (under $\mathrm{H}_0$) of our proposed rank-based energy test statistic in Theorems~\ref{theo:twosamconsis} and~\ref{theo:twosamasdistn} respectively. The population version of this rank-based energy statistic exhibits several interesting and desirable properties which we highlight in~\cref{lem:propren}. 
	
	To derive the asymptotic distributional limits of our test statistics (see Theorems~\ref{theo:indepasdistn} and~\ref{theo:twosamasdistn}) we develop some new general results involving certain permutation-based statistics (see Lemmas~\ref{lem:rankdcov} and~\ref{lem:rankenergy}), based on Stein's method for exchangeable pairs (see e.g.,~\cite{Cha2011}), which could be of independent interest. Further, to prove that the above proposed tests are consistent, under all fixed alternatives, we needed a new result on the convergence of the multivariate rank maps, which we state in~\cref{lem:L2}. The specialty of~\cref{lem:L2} is that it is sufficient for proving consistency and proceeds under minimal assumptions on the underlying measures, as opposed to the much stronger conditions usually required to show uniform convergence of multivariate ranks (see~\cite{del2018center,ghosal2019multivariate}). This result is of independent interest in the theory of measure transportation (see Section~\ref{sec:proprankquan} for more details). Moreover, we extend both the above tests to their multi-sample versions in~\cref{sec:multext}; the corresponding theoretical results are presented in Propositions~\ref{prop:multextind} and~\ref{prop:multextgof}. 
	
	We also carry out extensive simulation experiments to study the power behavior of the proposed tests (see~\cref{sec:sim}). These simulations show that our proposed procedures for mutual independence testing and two-sample goodness-of-fit testing perform well under a variety of alternatives, often outperforming competing methods. In general these distribution-free tests have good efficiency, are more powerful for distributions with heavy tails and are more robust to outliers and contaminations. In~\cref{sec:realdata}, we demonstrate practical advantages of our proposals over other competing methods via the analysis of two benchmark data sets.\par  

	In the following, we encapsulate some of the main contributions of this paper, all the while comparing our procedures to existing approaches from the statistics literature.
	
	\noindent \textbf{(i)} \textbf{Exact distribution-freeness}: As mentioned before, our proposals are all exactly distribution-free in finite samples. This is a particularly desirable property as it avoids the need to estimate any nuisance parameters, or use resampling/permutation ideas, or conservative asymptotic approximations, for determining rejection thresholds. Moreover, distribution-free procedures can help reduce computational burden in statistical problems --- a very practical concern in this era of big data; see e.g.,~\cite[Section 7]{heller2012} for an interesting discussion on this topic. As far as we are aware, the only distribution-free methods available in the literature for tackling the above discussed problems \textbf{(I)} and \textbf{(II)}  are:~\cite{rosenbaum2005,boeckel2018multivariate,munmun2014} for the multivariate two-sample problem and~\cite{munmun2016,heller2012} for the mutual independence testing problem. 
	
	\noindent \textbf{(ii)} \textbf{Completely nonparametric and computationally feasible}: Being based on multivariate ranks, our proposal is completely nonparametric. Moreover, computing our proposed test statistics is computationally feasible under all dimensions and sample sizes, and further, it does not involve any tuning parameters. This is in sharp contrast to approaches based on estimating functionals of underlying densities --- such as mutual information (see~\cite{berrett2017nonparametric}) --- or tests based on arbitrary partitions of the sample space (such as~\cite{gretton2008}) that are reliant on the choice of tuning parameter(s). In~\cref{sec:compute}, we explain how our proposed test statistics can be computed in a few simple steps using readily available \texttt{R} packages. Although exactly distribution-free graph based tests for mutual independence and two-sample goodness-of-fit testing were proposed in~\cite{munmun2016} and~\cite{munmun2014} respectively, these tests are extremely expensive to compute and possibly not applicable even for moderate sample sizes. 

	\noindent \textbf{(iii)} \textbf{Consistency under absolute continuity}: The only condition we need on the underlying distributions for the consistency of our tests is that they be absolutely continuous (no moment conditions are necessary). This enables their direct usage for nonparametric inference under heavy-tailed data-generating distributions such as stable laws~\cite{yang2012} and Pareto distributions~\cite{Rizzo2009}, and also sets them apart from popular methods such as usual distance covariance and energy statistic.  
	To the best of our knowledge, there are only two computationally efficient exactly distribution-free multivariate mutual independence testing procedures in literature, both based on a similar graph-based framework and proposed simultaneously in~\cite{heller2012}. 
	However, the authors in that paper do not provide any results that guarantee consistency of their tests against fixed alternatives. 
	
	\noindent \textbf{(iv)} \textbf{Broader scope of applications in multivariate nonparametric testing}: As described before, our approach is holistic. Based on our ideas, one can easily construct multivariate rank-based distribution-free tests for mutual independence using other statistics, such as Hilbert-Schmidt Independence Criteria (\cite{gretton2008kernel}) or HHG (\cite{heller2013}), instead of distance covariance; same goes for the goodness-of-fit testing problem. Note that, although we delve deep into these two particular nonparametric problems, we essentially describe a general principle to construct distribution-free tests in multivariate nonparametric settings that can be used in a variety of other contexts; e.g., in tests of symmetry~\cite{Gabor2013}, hierarchical clustering~\cite{rizzo2005}, change point analysis~\cite{szekely2009}, etc. We provide a concrete example of this in~\cref{sec:otherp}, where we present a distribution-free test for multivariate symmetry.

	The rest of the paper is organized as follows. In~\cref{sec:defrank}, we start with a brief overview of measure transportation (\cref{sec:revOT}), followed by a description of our proposed multivariate ranks and their properties (Sections~\ref{sec:defn} and~\ref{sec:proprankquan}).~\cref{sec:nrbnt} introduces new measures of multivariate association and goodness-of-fit and also discusses some interesting properties of these measures that make them desirable. Our proposed procedures for testing mutual independence and equality of distributions are introduced in~\cref{sec:mtp} (along with their multi-sample extensions). In that section we also discuss interesting/useful properties of our test statistics and provide theoretical guarantees with regards to distribution-freeness, consistency and asymptotic null distribution. In~\cref{sec:otherp} we develop a distribution-free test for testing multivariate symmetry. Section~\ref{sec:sim},~\cref{sec:addsim} and~\cref{sec:realdata} illustrate the usefulness of our proposed methods via simulation experiments and real data analysis. We conclude the main paper with a brief discussion in~\cref{sec:discuss}. In~\cref{sec:compute}, we explain how the proposed test statistics can be computed using standard software packages (in~\texttt{R}).~\cref{sec:revQMC} is aimed at providing a very brief introduction to the field of quasi-Monte Carlo methods which plays a tangential role in our approach. Finally, in  Appendices~\ref{sec:appen} and~\ref{sec:permappen} we provide the proofs of our main results, while in~\cref{sec:auxi}, we discuss some existing results on convex analysis and Stein's method of exchangeable pairs, which are used in the proofs of our main results. \par

	All the methods described in this paper have been implemented using the \texttt{R} software. The relevant codes, including our simulation experiments, are available in the first author's \texttt{GitHub} page: \url{https://github.com/NabarunD/MultiDistFree}. 
	
After the first version of this paper was posted on \texttt{arxiv}, we were made aware of the very recent paper~\cite{shi2019distribution} (uploaded after our first submission on \texttt{arxiv}). The paper~\cite{shi2019distribution} considers distribution-free mutual independence testing of two random vectors (i.e., problem (I)) using multivariate ranks as described in~\cite{del2018center}. Their paper also shows the distribution-freeness and consistency of the same test-statistic as in~\cref{sec:Indep} of this paper. However, the asymptotic consistency results in~\cite{shi2019distribution} are derived under more stringent conditions (e.g., nonvanishing Lebesgue probability densities). Note that in our paper, we develop a general framework for multivariate distribution-free nonparametric testing using optimal transportation, applicable much beyond problem (I); in particular, we also consider problem (II) and the problem of testing for multivariate symmetry (see~\cref{sec:otherp}).
	
	\section{Multivariate ranks and quantiles}\label{sec:defrank}
	In this section, we define ranks and quantiles for multivariate distributions (both population and empirical versions) using the theory of measure transportation; our approach is similar to that of~\cite{del2018center} and~\cite{boeckel2018multivariate}. This will serve a pivotal role in defining the test statistics that appear later in the paper.  
	
	\subsection{Preliminaries: Overview of measure transportation}\label{sec:revOT}
	Let us introduce some notation for the rest of the paper. We will use $\|\cdot\|$ and $\langle\cdot,\cdot\rangle$ to denote the standard Euclidean norm and inner product on a suitable finite dimensional Euclidean space (say $\mathbb{R}^d$) respectively.  Weak convergence of distributions will be denoted by $\stackrel{w}{\to}$ while $\overset{d}{=}$ will denote equality in distribution. We will use $\mathcal{U}^d$ to denote the uniform distribution on $[0,1]^d$, and $S_n$ for the set of all permutations of $\{1,2,\ldots ,n\}$. Let $\delta_{\mathbf{a}}$ denote the Dirac measure that assigns probability $1$ to the point $\mathbf{a}$. Finally, let $\mathcal{P}(\mathbb{R}^d)$ and $\mathcal{P}_{ac}(\mathbb{R}^d)$ denote the families of all probability distributions and Lebesgue absolutely continuous probability measures on $\mathbb{R}^d$, respectively. 
	
	As the name suggests, ``measure transportation" (perhaps more commonly referred to as \textit{optimal transportation}) is the problem of finding ``nice" functions $F:\R^d\to \R^d$ such that $F$ pushes a given measure $\mu\in\mathcal{P}(\R^d)$ to $\nu\in\mathcal{P}(\R^d)$. Here, by $F$ pushes $\mu$ to $\nu$, usually written as $F\#\mu=\nu$, we mean that $F(\mathbf{X}) \sim \nu$ where $\mathbf{X} \sim \mu$. This rich area of mathematics was initiated by the work of Gaspard Monge in 1781 (see~\cite{monge1781memoire}). Based on already introduced notation, perhaps the simplest version of \textit{Monge's problem} is as follows:
	\begin{align}\label{eq:Mongeproblem}
	\inf_{F} \int  \lVert \mathbf{x}-F(\mathbf{x})\rVert^2 \,d\mu( \mathbf{x})\qquad \mbox{subject to}\quad F\#\mu=\nu;
	\end{align}
	this is technically a mis-characterization as {Monge} originally worked with the loss $\lVert\cdot\rVert$ instead of $\lVert\cdot\rVert^2$. A minimizer of~\eqref{eq:Mongeproblem}, if it exists, is referred to as an {\it optimal transport map}. 
	One of the most powerful results in this field came into being from Brenier's \textit{Polar Factorization Theorem} (see~\cite{Brenier1991}) which yields: 
	If $\mu,\nu\in\mathcal{P}_{ac}\big(\mathbb{R}^d\big)$ have finite second-order moments, then the corresponding Monge's problem admits a $\mu$-a.e.~unique solution which happens to be the gradient of a convex function.

	While the above approach addresses the problem of finding functions that push $\mu$ to $\nu$, the assumption on second-order moments (which is a basic requirement for Monge's problem to make sense) seems extraneous and inappropriate. Indeed, for $d=1$, if $F_{\mu}$ and $F_{\nu}$ are the distribution functions associated with $\mu$ and $\nu$ (assumed to be absolutely continuous) respectively, then $F_{\nu}^{-1}\circ F_{\mu}$ pushes $\mu$ to $\nu$ without any moment assumptions. A ground-breaking extension of this univariate property was proved by Robert McCann in 1995, where he took a geometric approach to the problem of measure transportation. His result is the defining tool we will need to make sense of the definitions in this section. Therefore, let us state \textit{McCann's theorem} in a form which will be useful to us; see e.g.,~\cite[Theorem 2.12 and Corollary 2.30]{Villani2003}.
	\begin{prop}[McCann's theorem~\cite{Mccann1995}]\label{prop:Mccan}
		Suppose that $\mu, \nu\in\mathcal{P}_{ac}(\R^d)$. Then there exists functions $R(\cdot)$ and $Q(\cdot)$ (hereafter referred to as ``transport maps"), both of which are gradients of (extended) real-valued $d$-variate convex functions (hereafter called ``transport potentials"), such that $R\# \mu=\nu$, $Q\# \nu=\mu$, $R$ and $Q$ are unique ($\mu$ and $\nu$ a.e.~respectively), $R\circ Q (\mathbf{x})=\mathbf{x}$ ($\mu$ a.e.) and $Q\circ R(\mathbf{y})=\mathbf{y}$ ($\nu$ a.e.). 
		
		Moreover, if $\mu$ and $\nu$ have finite second moments, $R(\cdot)$ is also the solution to Monge's problem in~\eqref{eq:Mongeproblem}.
	\end{prop}
	Observe that McCann's theorem does away with all moment assumptions and guarantees existence and (a.e.) uniqueness of \textit{transport maps} under minimal assumptions on $\mu$ and $\nu$. Note that any convex function on $\mathbb{R}^d$ is differentiable Lebesgue a.e., and consequently $\mu$ (or $\nu$) a.e.~by Alexandroff Theorem (see e.g.,~\citet{Alexandroff1939}). In~\cref{prop:Mccan}, by ``gradient of a convex function" we essentially refer to a function from $\mathbb{R}^d\to\mathbb{R}^d$ which is $\mu$ (or $\nu$) a.e.~equal to the gradient of some convex function.

	\subsection{Definitions of multivariate ranks}\label{sec:defn}
	\begin{definition}[Population multivariate ranks and quantiles]\label{def:popquanrank}
		Set $\nu=\mathcal{U}^d$. Given a measure $\mu\in\mathcal{P}_{ac}(\mathbb{R}^d)$, the corresponding population rank and quantile maps are defined as functions $R(\cdot)$ and $Q(\cdot)$ respectively (as in~\cref{prop:Mccan}). Note that these are unique only up to measure zero sets with respect to $\mu$ and $\nu$ respectively. 
	\end{definition}
	\begin{remark}\label{rem:smoothrank}
		The smoothness and regularity properties of the population rank and quantile maps as in~\cref{def:popquanrank} have been studied extensively over the past $30$ years or so. Since such discussions are beyond the scope of this paper, we would like to refer the interested reader to~\cite{Figalli2013, caffarelli1990} and~\cite[Chapter 12]{Villani2009}. 
	\end{remark}
	In standard statistical applications, the population rank map is not available to the practitioner. In fact, the only accessible information about the measure $\mu$ comes in the form of empirical observations $\mathbf{X}_1,\mathbf{X}_2,\ldots ,\mathbf{X}_n\overset{i.i.d.}{\sim}\mu\in\mathcal{P}_{ac}(\mathbb{R}^d)$. A natural question thus arises: ``How can we estimate population ranks from empirical observations?". In this direction, let 
	\begin{equation}\label{eq:H_n^d}
	\mathbf{\mathcal{H}}_n^d\coloneqq\{\mathbf{h}_1^d,\ldots ,\mathbf{h}_n^d\}
	\end{equation}
	denote the (fixed) set of sample {\it multivariate rank vectors} (analogous to $\mathbf{c}_i$'s in~\eqref{eq:intro1-1}). In practice, for $d\geq 2$ we may take $\mathbf{\mathcal{H}}_n^d$ to be the $d$-dimensional Halton sequence of size $n$ (as described in~\cref{sec:revQMC}), and the usual $\{i/n\}_{1\leq i\leq n}$ sequence when $d=1$. The empirical distribution on $\mathbf{\mathcal{H}}_n^d$ will serve as a discrete approximation of $\mathcal{U}^d$. Also, let $\mathbf{\mathcal{D}}_n^X :=\{\mathbf{X}_1,\ldots ,\mathbf{X}_n\}$ be the observed data. Let
	\begin{equation}\label{eq:mu_nu_n}
	\mu_n^{\mathbf{X}} := \frac{1}{n} \sum_{i=1}^n \delta_{\mathbf{X}_i} \qquad \;\; \mathrm{and} \qquad \;\;  \nu_n := \frac{1}{n} \sum_{i=1}^n \delta_{\mathbf{h}_i^d}
	\end{equation}
	denote the empirical distributions on $\mathcal{D}_n^X$ and $\mathcal{H}_n^d$ respectively. 
	
	\begin{definition}[Empirical rank map]\label{def:empquanrank}
		We define the empirical rank function $\hat{R}_n:\mathbf{\mathcal{D}}_n^X\to \mathbf{\mathcal{H}}_n^d$ as the optimal transport map which transports $\mu_n^{\mathbf{X}}$ (the empirical distribution on the data) to $\nu_n$ (the empirical distribution on $\mathcal{H}_n^d$), i.e., 
		\begin{equation}\label{eq:R_n}
		\hat{R}_n = \argmin_F \int \|\mathbf{x} -  F(\mathbf{x})\|^2\,d\mu_n^{\mathbf{X}}( \mathbf{x})\qquad \mathrm{subject \; to}\quad F\#\mu_n^{\mathbf{X}} = \nu_n
		\end{equation}
	\end{definition}
	Note that~\eqref{eq:R_n} can be thought of as the discrete analogue of~\eqref{eq:Mongeproblem} which defines the population rank function $R(\cdot)$ if $\mu$ has finite second moments. Further,~\eqref{eq:R_n} is equivalent to the following optimization problem:
	\begin{align}\label{eq:empopt}
	\mathbf{\hat{\sigma}}_n\coloneqq \argmin_{\mathbf{\sigma}\in S_n} \sum_{i=1}^n \lVert \mathbf{X}_{i}-\mathbf{h}_{{\sigma(i)}}^d\rVert^2 =\argmax_{\sigma \in S_n} \sum_{i=1}^n \langle \mathbf{X}_i,\mathbf{h}_{{\sigma(i)}}^d\rangle.
	\end{align}
	The equivalence between the two optimization problems in~\eqref{eq:empopt} can be easily established by writing out the norms in terms of standard inner products. Note that $\mathbf{\hat{\sigma}}_n$ is a.s.~uniquely defined (for each $n$). Now, based on~\eqref{eq:empopt}, observe  that the sample rank map $\hat{R}_n$ satisfies 
	\begin{equation}\label{eq:R_n_h}
	\hat{R}_n(\mathbf{X}_i)=\mathbf{h}_{\hat{\sigma}_n(i)}^d, \qquad \mathrm{for \;\;} i = 1,\ldots, n.
	\end{equation}
	\begin{remark}\label{rem:emcompute}
		The optimization problem in~\eqref{eq:empopt} is a combinatorial optimization problem. However it is known to be equivalent to a linear program and can consequently be solved by standard solvers. Moreover, the special structure of the above problem allows us to view it as an assignment problem (see~\cite{munkres1957,bertsekas1988}) for which algorithms with worst case complexity $\mathcal{O}(n^3)$ are available in the literature. We will discuss this in more detail  in~\cref{sec:compute}.
	\end{remark}

	\begin{remark}[Connection to usual ranks in one-dimension]\label{rem:dim1con}
		For $d=1$, if we use $\mathbf{\mathcal{H}}_n^1 = \{i/n\}_{i=1}^{n}$, then the empirical ranks $\hat{R}_n(\cdot)$ reduce to the usual notion of one-dimensional ranks.
	\end{remark}
	\subsection{Properties of multivariate ranks}\label{sec:proprankquan}
	When $d=1$, the notion of ranks has a number of desirable properties which have been useful in analyzing rank-based estimators and test statistics (see e.g.,~\cite[Part I]{del2018center} and the references therein). Below in~\cref{prop:finproperties} (proved in~\cref{proof:finproperties}), we reproduce some of these properties for the empirical multivariate ranks as in~\cref{def:empquanrank}.~\cref{prop:finproperties} is in fact very similar to~\cite[Proposition 1.6.1]{del2018center}, with some differences which we will elaborate in~\cref{sec:compare}.
	\begin{prop}\label{prop:finproperties}
		Suppose that $\mathbf{X}_1,\ldots ,\mathbf{X}_n\overset{i.i.d.}{\sim}\mu\in\mathcal{P}_{ac}(\mathbb{R}^d)$. We define an \textit{order statistic} $\mathbf{X}^{(n)}_{(\cdot)}$ of $\{\mathbf{X}_1,\ldots ,\mathbf{X}_n\}$ as any fixed, arbitrary ordered version of the same --- for example, $\mathbf{X}^{(n)}_{(\cdot)}=(\mathbf{X}_{(1)},\ldots,\mathbf{X}_{(n)})$ where $\mathbf{X}_{(i)}$ is such that the first coordinate of $\mathbf{X}_{(i)}$ is the $i^{th}$ order statistic of the $n$-tuple formed by the first coordinates of the $n$-vectors in $\mathbf{X}^{(n)}_{(\cdot)}$. Then:
		\begin{itemize}
			\item[(i)] The \textit{order statistic} $\mathbf{X}^{(n)}_{(\cdot)}$ is complete and sufficient.
			\item[(ii)] The vector $(\hat{R}_n(\mathbf{X}_1),\ldots ,\hat{R}_n(\mathbf{X}_n))$ is uniformly distributed over the $n!$ permutations of the fixed grid $\mathbf{\mathcal{H}}_n^d$ (see~\eqref{eq:H_n^d}).
			\item[(iii)] $(\hat{R}_n(\mathbf{X}_1),\ldots ,\hat{R}_n(\mathbf{X}_n))$ and $\mathbf{X}^{(n)}_{(\cdot)}$ are mutually independent.
		\end{itemize}
	\end{prop}
	\begin{remark}[On~\cref{prop:finproperties}]\label{rem:finprop}
		Property (ii) from~\cref{prop:finproperties} is an analogue of the distribution-freeness of one-dimensional ranks. Property (iii) may be interpreted as the independence between ranks and order statistics.   
	\end{remark}
	As we will see in Section~\ref{sec:mtp}, the distribution-free property of the empirical multivariate ranks (see $(ii)$) will lead to the distribution-freeness of the proposed test statistics. However, to guarantee  the consistency of the proposed tests, we need the sample rank maps to be well-behaved as the sample size grows. In fact, in the following theorem (proved in~\cref{proof:lemL-2}) we show that the sample rank map converges to its population counterpart (i.e., the population rank function $R(\cdot)$ as in~\cref{def:popquanrank}) in a suitable sense, under minimal assumptions.

	\begin{theorem}[\textit{$L^2$-convergence}]\label{lem:L2}
		Assume $\mathbf{X}_1,\ldots ,\mathbf{X}_n\overset{i.i.d.}{\sim}\mu\in \mathcal{P}_{ac}(\mathbb{R}^d)$. Suppose that $\nu_n \overset{w}{\longrightarrow} \mathcal{U}^d$;  see~\eqref{eq:mu_nu_n}. Then
		\begin{align*}
		\frac{1}{n}\sum_{i=1}^n \lVert \hat{R}_n(\mathbf{X}_i)-R(\mathbf{X}_i)\rVert \overset{a.s.}{\longrightarrow}0.
		\end{align*}
	\end{theorem}
	\begin{remark}[$L^p$-convergence]\label{rem:L2-1}
		As $\hat{R}_n(\cdot)$ and $R(\cdot)$ are uniformly bounded,~\cref{lem:L2} implies convergence with respect to any $L^p$-norm, for $1\leq p<\infty$.
	\end{remark}
	\begin{remark}[On absolute continuity of $\mu$]\label{rem:L2-2}
		It is easy to see that \textit{$L^2$-convergence}, as presented in~\cref{lem:L2}, is weaker than \textit{uniform convergence}  (see~\cite{ghosal2019multivariate,chernozhukov2017monge,del2018center}). However, compared to the above references, the assumptions in~\cref{lem:L2} are minimal and do not, in general, guarantee uniform convergence of $\hat{R}_n(\cdot)$. In~\cref{sec:mtp} we will highlight specifically how and why~\cref{lem:L2} provides a more useful notion of convergence necessary for the results in this paper. For the time being, it is perhaps instructive to note that, even for ranks in one-dimension, the distribution-free property does not hold if the data generating measure is not continuous. Since distribution-free inference is the main goal of this paper, it seems reasonable to assume absolute continuity of $\mu$.
	\end{remark}
	\section{New multivariate rank-based measures for nonparametric testing}\label{sec:nrbnt}
	We introduce new multivariate rank-based measures of dependence and goodness-of-fit in this section and study the properties of these population quantities.
	\subsection{Rank-based dependence measure}\label{sec:depmes}
	In order to motivate our proposal, let us start with $d=1$. Suppose that ${Z}_1$ and ${Z}_2$ are real-valued absolutely continuous random variables with distribution functions $G_1(\cdot)$ and $G_2(\cdot)$. It is a simple probability exercise to show that ${Z}_1$ and ${Z}_2$ are independent if and only if $G_1({Z}_1)$ and $G_2({Z}_2)$ are independent (a more general version of this result will be proved later in the paper, see~\cref{lem:proprdcov} (part (b)). Thus, ${Z}_1$ and ${Z}_2$ are independent if and only if the joint characteristic function of $(G_1(Z_1),G_2(Z_2))$ factors as the product of the marginal characteristic functions, i.e., for all $(t,s)\in\mathbb{R}^2$,
	\begin{equation*}
	\left|\mathbb{E}\exp\big(itG_1({Z}_1)+isG_2({Z}_2))-\mathbb{E}\exp\big(itG_1({Z}_1))\mathbb{E}\exp\big(isG_2({Z}_2))\right|^2=0.
	\end{equation*}
	This suggests the following natural measure of dependence:
	\begin{equation*}
	\mathcal{R}_w\coloneqq \int\int \left|\mathbb{E}\exp\big(iG_{t,s}(\mathbf{Z})\big)-\mathbb{E}\exp\big(itG_1({Z}_1))\mathbb{E}\exp\big(isG_2({Z}_2))\right|^2 w(t,s)\,dt\,ds
	\end{equation*}
	where $\mathbf{Z} = (Z_1,Z_2)$, $G_{t,s}(\mathbf{Z})=tG_1({Z}_1)+sG_2({Z}_2)$ and $w:\mathbb{R}\times\mathbb{R}\to[0,\infty)$ is a weight function such that $\mathcal{R}_w<+\infty$. The following proposition (proved in~\cref{proof:srcc}) draws a connection between $\mathcal{R}_w$ and the classical Spearman's rank correlation, which we think has not been observed before.
	\begin{prop}\label{prop:srcc}
		Consider the notation introduced above. Set $$f_{{Z}_1,{Z}_2}(t,s)\coloneqq\mathbb{E}\exp\big(itG_1({Z}_1)+isG_2({Z}_2))-\mathbb{E}\exp\big(itG_1({Z}_1))\mathbb{E}\exp\big(isG_2({Z}_2)).$$ Then,
		\begin{equation}\label{eq:srccmain1}
		\lim\limits_{t,s\to 0,|t|/|s|\to c} \frac{\big|f_{{Z}_1,{Z}_2}(t,s)\big|^2}{\big|f_{{Z}_1,{Z}_1}(t,s)\big|\big|f_{{Z}_2,{Z}_2}(t,s)\big|}=\rho^2\big(G_1({Z}_1),G_2({Z}_2)\big)
		\end{equation}
		where $\rho^2\big(G_1({Z}_1),G_2({Z}_2)\big)$ denotes the usual correlation between $G_1({Z}_1)$ and $G_2({Z}_2)$. In the above display, $c>0$ is finite, and ensures that $s$ and $t$ do not converge to $0$ at ``different rates".
	\end{prop}
	The right side of~\eqref{eq:srccmain1} may be interpreted as the population analogue of the classical \emph{Spearman's rank correlation}. Note that applications of Spearman's rank correlation as a measure of association have been extensively studied in the statistics literature (see e.g.,~\cite{hauke2011comparison,iman1982distribution,mukaka2012guide}).

	\begin{remark}\label{rem:consrc2}
		\cref{prop:srcc} shows how the Spearman's rank correlation measure effectively looks at the difference between the joint and marginal characteristic functions for small (in magnitude) choices of $t\mathrm{ and }s$. Therefore, $\mathcal{R}_w$ (after rescaling) offers a very natural extension to Spearman's rank correlation, but it can capture all kinds of departures from independence.
	\end{remark}
	
	\begin{remark}\label{rem:consrc1}
		It is easy to see that the right hand side of~\eqref{eq:srccmain1} being $0$ does not imply the independence of ${Z}_1$ and ${Z}_2$. For example, say ${Z}_1\sim\mathcal{U}^1$ and ${Z}_2={Z}_1$ if ${Z}_1\in [1/4,3/4]$, ${Z}_2=1-{Z}_1$ if ${Z}_1\in (0,1/4)\cup (3/4,1)$. Then ${Z}_2\sim\mathcal{U}^1$ and both $G_1(\cdot)$, $G_2(\cdot)$ are identity functions on $(0,1)$. Therefore, $\mathbb{E}[G_1({Z}_1)G_2({Z}_2)]-\mathbb{E}[G_1({Z}_1)] \mathbb{E}[G_2({Z}_2)]=\mathbb{E}[{Z}_1{Z}_2]-1/4=0$.
	\end{remark}
	
	The above discussion now raises the following two questions: ``Can we extend $\mathcal{R}_w$ beyond $d=1$? Also, how do we choose the weight function $w(\cdot,\cdot)$?". For the first question, we will proceed by replacing $G_1(\cdot)$ and $G_2(\cdot)$ with the notion of population multivariate ranks as introduced in~\cref{def:popquanrank}. For the second question, we will borrow the weight function from the seminal paper~\cite{Gabor2007} where the authors introduce the notion of \emph{distance covariance}. As in~\cite{Gabor2007}, we do not make any claims on the optimality of our proposed weight function except that it ensures simple, applicable empirical formulae and an exact equivalence between $\mathcal{R}_w$ and the independence between ${Z}_1$ and ${Z}_2$. We are now in a position to formally define the new rank-based multivariate measure of dependence.
	\begin{definition}[Rank distance covariance]\label{def:rdcov}
		Suppose that $\mathbf{Z}_1\sim\mu_1$ and $\mathbf{Z}_2\sim\mu_2$ (not necessarily independent) such that $\mu_1\in \mathcal{P}_{ac}(\mathbb{R}^{d_1})$ and $\mu_2\in \mathcal{P}_{ac}(\mathbb{R}^{d_2})$. Let $R_1(\cdot)$ and $R_2(\cdot)$ denote the corresponding population rank maps (\cref{def:popquanrank}). The rank distance covariance $(\Rdcov^2)$ between $\mathbf{Z}_1$ and $\mathbf{Z}_2$ is defined as the usual distance covariance between $R_1(\mathbf{Z}_1)$ and $R_2(\mathbf{Z}_2)$, i.e.,  
		\begin{equation}\label{eq:rdcov1}
		\hspace{-0.10in}	\Rdcov^2(\mathbf{Z}_1,\mathbf{Z}_2)\coloneqq\int\limits_{\mathbb{R}^{d_1+d_2}} \frac{\left|\mathbb{E}\exp\big(iR_{\mathbf{t},\mathbf{s}}(\mathbf{Z}))-\mathbb{E}\exp\big(i\mathbf{t}^{\top}R_1(\mathbf{Z}_1))\mathbb{E}\exp\big(i\mathbf{s}^{\top}R_2(\mathbf{Z}_2))\right|^2}{c(d_1)c(d_2)\lVert \mathbf{t}\rVert^{1+d_1}\lVert \mathbf{s}\rVert^{1+d_2}} \,d\mathbf{t}\,d\mathbf{s}    	
		\end{equation}
		where $\mathbf{Z} := (\mathbf{Z}_1,\mathbf{Z}_2)$,  $R_{\mathbf{t},\mathbf{s}}(\mathbf{Z}) :=\mathbf{t}^{\top}R_1(\mathbf{Z}_1)+\mathbf{s}^{\top}R_2(\mathbf{Z}_2)$ and $c(d) :=\pi^{(1+d)/2}\big(\Gamma((1+d)/2)\big)^{-1}$. 
	\end{definition}
	
	Now let us look into some of the properties of $\Rdcov$ that make it a desirable measure of dependence. The proof of the following lemma is given in~\cref{proof:proprdcov}.
	\begin{lemma}\label{lem:proprdcov}
		Under the same assumptions as in~\cref{def:rdcov}, we have:
		\begin{itemize}
			\item[(a)] Suppose that $(\mathbf{Z}_1^1,\mathbf{Z}^1_2),(\mathbf{Z}^2_1,\mathbf{Z}_2^2),(\mathbf{Z}^3_1,\mathbf{Z}^3_2)$ are independent observations having the same distribution as $(\mathbf{Z}_1,\mathbf{Z}_2)$. Then,
			\begin{align}\label{eq:proprdcovequiv}
			\Rdcov^2(\mathbf{Z}_1,\mathbf{Z}_2)&=\mathbb{E}\big[\lVert R_1(\mathbf{Z}_1^1)-R_1(\mathbf{Z}_1^2)\rVert\lVert R_2(\mathbf{Z}_2^1)-R_2(\mathbf{Z}_2^2)\rVert\big]\nonumber \\ &\qquad+\mathbb{E}\big[\lVert R_1(\mathbf{Z}_1^1)-R_1(\mathbf{Z}_1^2)\rVert\big] \mathbb{E}\big[\lVert R_2(\mathbf{Z}_2^1)-R_2(\mathbf{Z}_2^2)\rVert\big]\nonumber \\ &\qquad-2\mathbb{E}\big[\lVert R_1(\mathbf{Z}_1^1)-R_1(\mathbf{Z}_1^2)\rVert\lVert R_2(\mathbf{Z}^1_2)-R_2(\mathbf{Z}^3_2)\rVert\big].
			\end{align}
			\item[(b)] $\Rdcov(\mathbf{Z}_1,\mathbf{Z}_2)=0$ if and only if $\mathbf{Z}_1$ and $\mathbf{Z}_2$ are independent.
			\item[(c)] $\Rdcov(\mathbf{Z}_1,\mathbf{Z}_1)>0$.
			\item[(d)] (Invariance) Suppose $\mathbf{a}_1\in\mathbb{R}^{d_1}$, $\mathbf{a}_2\in\mathbb{R}^{d_2}$ and  $b_1,b_2>0$. Then $\Rdcorr(\mathbf{Z}_1,\mathbf{Z}_2)=\Rdcorr(\mathbf{a}_1 + b_1\mathbf{Z}_1,\mathbf{a}_2 +b_2\mathbf{Z}_2)$.
			\item[(e)] Suppose that $(\mathbf{Z}_1^n,\mathbf{Z}_2^n)\in\mathbb{R}^{d_1}\times\mathbb{R}^{d_2}$ is a sequence of random vectors that converge weakly to $(\mathbf{Z}_1,\mathbf{Z}_2)$; here we assume that $\mathbf{Z}_1^n$ and $\mathbf{Z}_2^n$ have absolutely continuous distributions for all $n$. Then, $\Rdcov^2(\mathbf{Z}_1^n,\mathbf{Z}_2^n)\longrightarrow\Rdcov^2(\mathbf{Z}_1,\mathbf{Z}_2)$ as $n\to\infty$.
		\end{itemize}
	\end{lemma}
	We would like to refer the interested reader to~\cite{Mori2019} for an elaborate discussion on the importance of these properties in a dependence measure.
	\begin{remark}\label{rem:dcovdiff}
		Unlike distance covariance (see~\cite{Gabor2007,Mori2019}),~\cref{lem:proprdcov} does not require any moment assumptions on $\mathbf{Z}_1$ and $\mathbf{Z}_2$. However, we do need absolute continuity of the underlying measures $\mu_1$ and $\mu_2$, an assumption which has been justified in~\cref{rem:L2-2}.
	\end{remark}
	
	\begin{remark}\label{rem:conbivnorm}
		In~\cite[Theorem 7]{Gabor2007}, a closed form expression for distance covariance when $(X,Y)$ has a bivariate normal distribution, parametrized by  correlation $\rho$,  is derived. Although for rank distance covariance (as defined in~\eqref{eq:rdcov1}) such a closed form expression is not easy to obtain, we can readily approximate it using Monte Carlo. In~\cref{sec:bivdcovcon}, we demonstrate that, in this bivariate normal setting, the population distance covariance and population rank distance covariance are both monotone in $|\rho|$ and essentially indistinguishable as functions of $\rho$.
	\end{remark}
	
	\subsection{Rank-based measure for two-sample goodness-of-fit}\label{sec:gofmes}
	We can use a similar approach as in~\cref{sec:depmes} to come up with a measure for multivariate two-sample goodness-of-fit testing. Define $\mathcal{S}^{d-1} :=\{\mathbf{x}\in\mathbb{R}^d:\lVert \mathbf{x}\rVert=1\}$ and let $\kappa(\cdot)$ denote the uniform measure on $\mathcal{S}^{d-1}$. Further, assume $\mathbf{Z}_1\sim\mu_1$ and $\mathbf{Z}_2\sim\mu_2$ are independent where $\mu_1,\mu_2\in\mathcal{P}_{ac}(\mathbb{R}^d)$. Then, using the continuity and uniqueness of characteristic functions, it is rather straightforward to check that $\mu_1=\mu_2$ if and only if $\mathbf{a}^{\top}\mathbf{Z}_1\overset{d}{=} \mathbf{a}^{\top}\mathbf{Z}_2$ for $\kappa$ a.e. $\mathbf{a}$ (for more details see~\cite[Theorem 2.1]{Baringhaus2004}). Therefore, a natural way to measure equality of distributions $\mu_1 = \mu_2$ would be to compare $\mathbb{P}(\mathbf{a}^{\top}\mathbf{Z}_1\leq t)$ and $\mathbb{P}(\mathbf{a}^{\top}\mathbf{Z}_2\leq t)$ for all $\mathbf{a}\in \mathcal{S}^{d-1}$ and all $t\in\mathbb{R}$. This provides the main motivation behind the {\it energy measure} for two-sample goodness-of-fit (see~\cite{Baringhaus2004,Gabor2013}), which is defined as:
	\begin{equation*}
	\mathrm{En}(\mathbf{Z}_1,\mathbf{Z}_2) :=\gamma_d\int_{\mathbb{R}}\int_{\mathcal{S}^{d-1}} \left(\mathbb{P}(\mathbf{a}^{\top}\mathbf{Z}_1\leq t)-\mathbb{P}(\mathbf{a}^{\top}\mathbf{Z}_2\leq t)\right)^2\,d\kappa(\mathbf{a})\,dt
	\end{equation*} 
	where $\gamma_d\coloneqq \big(2\Gamma(d/2))^{-1}\sqrt{\pi}(d-1)\Gamma\big((d-1)/2\big)$ for $d>1$ and $\gamma_d:=1$ for $d=1$. It can be shown that $\mathrm{En}(\mathbf{Z}_1,\mathbf{Z}_2)$ is well-defined if $\mathbf{Z}_1$ and $\mathbf{Z}_2$ have finite first moments (see~\cite[Lemma 2.3]{Baringhaus2004}). With the above discussion in mind, we are now in a position to define the rank-based version of the energy measure.
	\begin{definition}[Rank energy]\label{def:ren}
		Suppose that $\mathbf{Z}_1\sim\mu_1$ and $\mathbf{Z}_2\sim\mu_2$ are independent and $\mu_1,\mu_2\in \mathcal{P}_{ac}(\mathbb{R}^{d})$. Fix some $\lambda\in (0,1)$ (prespecified). Also let $R_\lambda(\cdot)$ denote the population rank map (see~\cref{def:popquanrank}) corresponding to the mixture distribution $\lambda \mu_1+(1-\lambda)\mu_2$. Then the rank energy ($\Ren_{\lambda}^2$) between $\mathbf{Z}_1$ and $\mathbf{Z}_2$ is defined as:
		\begin{equation}\label{eq:ren1}
		\Ren_{\lambda}^2(\mathbf{Z}_1,\mathbf{Z}_2) :=\gamma_d\int_{\mathbb{R}}\int_{\mathcal{S}^{d-1}} \left[\mathbb{P}(\mathbf{a}^{\top}R_\lambda(\mathbf{Z}_1)\leq t)-\mathbb{P}(\mathbf{a}^{\top}R_\lambda(\mathbf{Z}_2)\leq t)\right]^2\,d\kappa(\mathbf{a})\,dt.   	
		\end{equation}
		In other words, the rank energy between $\mathbf{Z}_1$ and $\mathbf{Z}_2$ is exactly equal to the usual energy measure between $R_\lambda(\mathbf{Z}_1)$ and $R_\lambda(\mathbf{Z}_2)$. Note that~\eqref{eq:ren1} is well-defined without any moment assumptions.
	\end{definition}
	\begin{remark}\label{rem:lamclarify}
		The choice of $\lambda\in (0,1)$ in~\cref{def:ren} may seem subjective.  However, in the kind of applications we are interested in, we will see that a natural choice of $\lambda$ will surface from the context of the problem itself.
	\end{remark}
	Now let us inspect the properties of $\Ren_{\lambda}^2$ which make it a desirable candidate for measuring two-sample goodness-of-fit. The proof of the following result is given in~\cref{proof:propren}.
	\begin{lemma}\label{lem:propren}
		Under the same assumptions as in~\cref{def:ren}, we have:
		\begin{itemize}
			\item[(a)] Suppose that $\mathbf{Z}_1^1,\mathbf{Z}_1^2$ are i.i.d.~with the same distribution as $\mathbf{Z}_1$, and $\mathbf{Z}_2^1,\mathbf{Z}_2^2$ are i.i.d.~with the same distribution as $\mathbf{Z}_2$. Then, 
			\begin{equation*}
			\hspace{-0.1in} \Ren_{\lambda}^2(\mathbf{Z}_1,\mathbf{Z}_2) = 2\mathbb{E}\lVert R_\lambda(\mathbf{Z}_1^1) -R_\lambda(\mathbf{Z}_2^1)\rVert - \mathbb{E}\lVert R_\lambda(\mathbf{Z}_1^1) - R_\lambda(\mathbf{Z}_1^2)\rVert - \mathbb{E}\lVert R_\lambda(\mathbf{Z}_2^1) - R_\lambda(\mathbf{Z}_2^2)\rVert.
			\end{equation*}
			\item[(b)] $\Ren_{\lambda}^2(\mathbf{Z}_1,\mathbf{Z}_2)=0$ if and only if  $\mathbf{Z}_1\overset{d}{=}\mathbf{Z}_2$.
			\item[(c)] (Invariance) Suppose that $\mathbf{a}\in\mathbb{R}^{d}$ and  $b>0$. Then $\Ren_{\lambda}^2(\mathbf{Z}_1,\mathbf{Z}_2)=\Ren_{\lambda}^2(\mathbf{a} +b\mathbf{Z}_1,\mathbf{a} +b\mathbf{Z}_2)$.
			\item[(d)] Suppose that $\mathbf{Z}_1^n$ and $\mathbf{Z}_2^n$ are two independent sequences of random vectors having absolutely continuous distributions such that $\mathbf{Z}_1^n\overset{w}{\longrightarrow}\mathbf{Z}_1$ and $\mathbf{Z}_2^n\overset{w}{\longrightarrow}\mathbf{Z}_2$ as $n\to\infty$. Then, $\Ren_{\lambda}^2(\mathbf{Z}_1^n,\mathbf{Z}_2^n)\longrightarrow\Ren_{\lambda}^2(\mathbf{Z}_1,\mathbf{Z}_2)$ as $n\to\infty$.
		\end{itemize}
	\end{lemma} 
	
	\section{Distribution-free multivariate independence and equality of distributions testing}\label{sec:mtp}
	This section is devoted to developing the new multivariate rank-based distribution-free testing procedures for the nonparametric problems discussed in the Introduction.
	\subsection{Distribution-free mutual independence testing}\label{sec:Indep}
	Suppose that $(\mathbf{X}_1,\mathbf{Y}_1),\ldots ,(\mathbf{X}_n,\mathbf{Y}_n)$ are i.i.d.~observations from some probability distribution $\mu\in\mathcal{P}(\mathbb{R}^{d_1+d_2})$ (here $d_1,d_2 \ge 1$) with marginals  $\mu_{\mathbf{X}}$ and $\mu_{\mathbf{Y}}$. In this subsection we assume that $$ \textbf{(AP1)}:\qquad\mu_{\mathbf{X}}\in\mathcal{P}_{ac}(\mathbb{R}^{d_1}) \quad \qquad \mathrm{and} \qquad \quad \mu_{\mathbf{Y}}\in\mathcal{P}_{ac}(\mathbb{R}^{d_2}).$$ 
	
	We are interested in testing the hypothesis:
	\begin{align*}
	\mathrm{H}_0: \mu=\mu_{\mathbf{X}}\otimes\mu_{\mathbf{Y}}\qquad \qquad  \mathrm{versus} \qquad\qquad \mathrm{H}_1:\mu\neq\mu_{\mathbf{X}}\otimes\mu_{\mathbf{Y}}.
	\end{align*}
	The above is certainly a classical problem in statistics and has received widespread attention across many decades. One of the earliest approaches, for $d_1=d_2=1$, was through the introduction of Pearson's correlation (see e.g.,~\cite{pearson1920notes}), which was later modified into rank-based correlation measures such as Spearman's rank correlation (see~\cite{spearman1904proof}) and Kendall's $\tau$ (see~\cite{kendall1938new,Kendall1990}). For an overview of other parametric approaches to the above problem see~\cite{wilks1938large} and~\cite{Pillai1967} and the references therein. However, nonparametric testing procedures soon replaced parametric ones as they do not require strong modeling assumptions and are consequently more robust and generally applicable.
	
	One of the first nonparametric approaches to the above problem, when $d_1=d_2=1$, was by Hoeffding~\cite{hoeffding1948non}, where the author proposed a test based on empirical distribution functions; also see~\cite{blum1961distribution}. A ``quadrant" based procedure was introduced in the late 1950s by Mosteller (see~\cite{Mosteller1946}) and later analyzed in~\cite{blomqvist1950measure}; also see~\cite{gieser1997}. A density estimation based approach to independence testing was proposed in~\cite{Rosenblatt1975}. In~\cite{Bregsma2014}, the authors introduced a consistent test of independence using a signed covariance measure that can be viewed as a modification of Kendall's $\tau$. This test was further analyzed and extended to a test of independence between multiple (more than $2$) random variables in~\cite{drton2018high}. When either $d_1>1$ or $d_2>1$, perhaps the most common approach historically used coordinate-wise or spatial ranks and signs (see e.g.,~\cite{puri1971nonparametric,oja2010,Oja2004} and the references therein). Such coordinate-wise rank-based extensions to Spearman's rank correlation, Kendall's $\tau$ and the quadrant statistic (mentioned above) for testing independence, when $d_1>1$ or $d_2>1$, were proposed in~\cite{taskinen2003,taskinen2005}. In~\cite{friedman1983}, the authors present a graph-based test of independence. A density based approach, involving the estimation of mutual information has been used in~\cite{berrett2017nonparametric}. Other proposals include the use of a maximal (or total) information coefficient (see~\cite{Reshef2016,Reshef2018}), empirical copula processes (see~\cite{kojadinovic2009,Quessy2010}), ranks of pairwise distances  (see~\cite{heller2013}), etc. A kernel based method, namely the Hilbert-Schmidt Independence criteria, which perhaps dates back to 1959 (see~\cite{Renyi1959}) has also been recently studied in great detail by Gretton and co-authors  (see e.g.,~\cite{gretton2008kernel,GrettonHSIC2005,Gretton2005}; also see e.g.,~\cite{Sen2014,ramdas2015decreasing}). Given the huge body of work in this area, we refer the reader to~\cite{Mari2001,Josse2016} for a survey on other testing procedures existing in the literature. While some of the tests discussed above guarantee consistency against fixed alternatives, a recurrent problem with all these approaches is that they lack the exact distribution-free property when either $d_1>1$ or $d_2>1$.
	
	The only distribution-free test in the context of mutual independence testing was proposed in~\cite{heller2012}; also see~\cite{munmun2016} for an extension. However, none of these tests come with any result that guarantees consistency against all fixed alternatives. In~\cite{heller2016multivariate}, the authors suggested testing whether two multivariate random vectors are independent, by testing whether the distance of one random vector from some arbitrary reference point is independent of the distance of the other random vector from another arbitrary reference point. Although this test is distribution-free, once the arbitrary	reference points are fixed, the test will not be consistent against all alternatives.

	Over the past $40$ years or so, multivariate tests of independence based on empirical characteristic functions have gained some prominence, thanks to early works in~\cite{Kankainen1995,csorgo1985,feuerverger1993consistent} and most significantly due to the seminal work by Szekely and co-authors (see~\cite{bakirov2006,Gabor2007,szekely2009}), where the notion of distance covariance was introduced; recall that in Section~\ref{sec:depmes} we have already encountered the population version this measure. Interestingly, distance covariance, which we have already seen can be viewed as the covariance between pairwise distances among the observed data points (see~\eqref{eq:firstdcov}), can also be interpreted as a weighted integral in terms of the difference between the joint empirical characteristic function and the product of marginal characteristic functions (see~\cite{Gabor2007}). Distance covariance also has interesting connections to kernel based methods; see e.g.,~\cite{sejdinovic2013}. On account of being simple to implement, easily explainable and providing consistency against any fixed alternatives (under suitable moment assumptions), this testing procedure has attracted a lot of attention, has inspired many applications, and is still a subject  of active research.

	In this section, we introduce a distribution-free multivariate rank-based version of the distance covariance test (see e.g.,~\cite{Gabor2007}) and demonstrate its appealing properties. We describe our method below. Let $\mu_n^{\mathbf{X}}$ and $\mu_n^{\mathbf{Y}}$ denote the empirical distributions on $\mathcal{D}_n^{\mathbf{X}}\coloneqq \{\mathbf{X}_1,\ldots ,\mathbf{X}_n\}$ and $\mathcal{D}_n^{\mathbf{Y}} :=\{\mathbf{Y}_1,\ldots ,\mathbf{Y}_n\}$ respectively. Moreover, let $\mathcal{H}^{d_1}_n\coloneqq \{\mathbf{h}_1^{d_1},\ldots ,\mathbf{h}_n^{d_1}\}$ and $\mathcal{H}_n^{d_2}\coloneqq \{\mathbf{h}_1^{d_2},\ldots ,\mathbf{h}_n^{d_2}\}$ denote the (fixed) sample of $d_1$ and $d_2$-dimensional ranks (analogous to $\mathbf{c}_i$'s in~\eqref{eq:intro1-1}). For $i=1,2$, as in~\cref{sec:defn} we recommend the use of the standard $d_i$-dimensional Halton sequence (see~\cref{sec:revQMC} for a discussion) when $d_i>1$ and the standard $\{i/n\}_{i\leq n}$ grid when $d_i=1$. We will work under the following assumption on $\mathcal{H}_n^{d_1}$ and $\mathcal{H}_n^{d_2}$:\par 
	\textbf{(AP2)}: The empirical distributions on $\mathcal{H}_n^{d_1}$ and $\mathcal{H}_n^{d_2}$ converge weakly to $\mathcal{U}^{d_1}$  and $\mathcal{U}^{d_2}$ respectively.

	Finally, we shall use $\hat{R}_n^{\mathbf{X}}(\cdot)$ and $\hat{R}_n^{\mathbf{Y}}(\cdot)$ to denote the empirical rank maps (see~\cref{def:empquanrank}) corresponding to the transportation of $\mu_n^{\mathbf{X}}$ and $\mu_n^{\mathbf{Y}}$ to the empirical distributions on $\mathcal{H}_n^{d_1}$ and $\mathcal{H}_n^{d_2}$ respectively (see~\eqref{eq:R_n}). Next, we define
	\begin{equation}\label{eq:indepstat}
	\Rdcov_n^2\coloneqq S_1+S_2-2S_3
	\end{equation}
	where 
	\begin{align*}
	S_1&\coloneqq \frac{1}{n^2}\sum_{k,l=1}^n \lVert \hat{R}_n^{\mathbf{X}}(\mathbf{X}_k)-\hat{R}_n^{\mathbf{X}}(\mathbf{X}_l)\rVert\lVert \hat{R}_n^{\mathbf{Y}}(\mathbf{Y}_k)-\hat{R}_n^{\mathbf{Y}}(\mathbf{Y}_l)\rVert,\\ S_2&\coloneqq\left(\frac{1}{n^2}\sum_{k,l=1}^n \lVert \hat{R}_n^{\mathbf{X}}(\mathbf{X}_k)-\hat{R}_n^{\mathbf{X}}(\mathbf{X}_l)\rVert\right)\times\left(\frac{1}{n^2}\sum_{k,l=1}^n \lVert \hat{R}_n^{\mathbf{Y}}(\mathbf{Y}_k)-\hat{R}_n^{\mathbf{Y}}(\mathbf{Y}_l)\rVert\right),\\ S_3&\coloneqq\frac{1}{n^3}\sum_{k,l,m=1}^n \lVert \hat{R}_n^{\mathbf{X}}(\mathbf{X}_k)-\hat{R}_n^{\mathbf{X}}(\mathbf{X}_l)\rVert\lVert \hat{R}_n^{\mathbf{Y}}(\mathbf{Y}_k)-\hat{R}_n^{\mathbf{Y}}(\mathbf{Y}_m)\rVert.
	\end{align*}
	Observe that the right side of~\eqref{eq:indepstat} can be viewed as an empirical version of population $\Rdcov$ (see~\eqref{eq:rdcov1}) through its alternate expression as in~\cref{lem:proprdcov} (part (a)). $\Rdcov_n^2$ can also be viewed as a rank-transformed version of the empirical distance covariance measure as introduced in~\cite[Equations (2.9) and (2.18)]{Gabor2007}. Another way to look at this is that $\Rdcov_n^2$ equals the covariance between pairwise distances of the multivariate rank vectors (analogous to~\eqref{eq:firstdcov} as motivated in the Introduction). By~\cite[Theorem 1]{Gabor2007}, it is easy to see that the right side of~\eqref{eq:indepstat} is always nonnegative. Moreover, note that, given the ranks, $\Rdcov_n^2$ can be computed in $\mathcal{O}(n^2 (d_1+d_2))$ steps (see~\cite{Huo2016}). In the following lemma we demonstrate the distribution-free property of $\Rdcov_n^2$ (see~\cref{proof:indepdistfree} for a proof).
	\begin{lemma}\label{rem:indepdistfree} 
		Under assumption \textbf{(AP1)} and $\mathrm{H}_0$, the distribution of $\Rdcov_n^2$, as defined in~\eqref{eq:indepstat}, is free of $\mu_{\mathbf{X}}$ and $\mu_{\mathbf{Y}}$.
	\end{lemma}

	\noindent \textbf{Distribution-free independence testing procedure}: Given a (prespecified) type-I error level $\alpha \in (0,1)$, let $$c_n\coloneqq \inf\{c >0: \mathbb{P}_{\textrm{H}_0}(n\Rdcov_n^2\geq c)\leq \alpha\}.$$ Note that, under $\mathrm{H}_0$, $\Rdcov_n^2$ is distribution-free (by~\cref{rem:indepdistfree}) and therefore, so is $c_n$. In other words, $c_n$ depends only on $n, d_1, d_2, \mathcal{H}^{d_1}_n, \mathcal{H}^{d_2}_n$ and $\alpha$, and can consequently be determined even before the data is observed. Moreover, we show in~\cref{theo:indepasdistn} below that, if assumption \textbf{(AP2)} is satisfied then asymptotically $c_n$ does not even depend on the particular choice of $\mathcal{H}^{d_1}_n$ and $\mathcal{H}^{d_2}_n$. Given $c_n$, our proposed testing procedure rejects $\mathrm{H}_0$ if $n\Rdcov_n^2\geq c_n$ and accepts $\mathrm{H}_0$ otherwise. By definition of $c_n$, this is clearly a level $\alpha$ test. \par
	\begin{remark}\label{rem:prevrdcov}
		The notion of rank-based distance covariance has attracted some interest in the literature. For $d_1=d_2=1$, it has been discussed in~\cite{szekely2009}, although to the best of our knowledge, its theoretical properties haven't been analyzed. In the discussion~\cite{remillard2009} based on~\cite{szekely2009}, the author proposed using distance covariance based on the vectors of component-wise ranks (for general $d_1$, $d_2$). This idea also has connections with existing copula based approaches; see e.g.,~\cite{kojadinovic2009} for details. This approach however does not yield a distribution-free test (if either $d_1$ or $d_2$ is $>1$), neither for finite $n$ nor asymptotically. In that sense, our proposal provides the ``correct" version of rank-based distance covariance.
	\end{remark}

	\noindent One of the interesting features of our proposed statistic, i.e.,~$\Rdcov_n^2$, is that it has a close connection with the celebrated Hoeffding's D-statistic (see~\cite{hoeffding1948non}) --- one of the earliest nonparametric approaches to testing for mutual independence when $d_1=d_2=1$. In fact,~$\Rdcov_n^2$ is exactly equivalent to the statistic proposed in~\cite{weihs2018} (also see the right sides of~\eqref{eq:dcovequivmain1} and~\eqref{eq:dcovequivmain2} below for the population and the empirical versions respectively), which in turn is a modified version of Hoeffding's D-statistic. The following lemma (see~\cref{proof:dcovequiv} for a proof) makes this connection precise.\par
	\begin{lemma}\label{lem:dcovequiv}
		Suppose that $(X,Y)\in\mathbb{R}^2$ with bivariate distribution function (DF) $F^{X,Y}(\cdot)$, and corresponding marginal DFs, $F^{X}$ and $F^Y$. Assume that $F^{X}(\cdot)$ and $F^{Y}(\cdot)$ are absolutely continuous. Also suppose that random samples $(X_1,Y_1),\ldots , (X_n,Y_n)$ are drawn according to the same distribution as $(X,Y)$. Further, we will use $F_n^{X,Y}(\cdot)$, $F_n^X(\cdot)$ and $F_n^Y(\cdot)$ to denote the joint and marginal empirical DFs of $X_i$'s and $Y_i$'s respectively. Then the following holds:
		\begin{equation}\label{eq:dcovequivmain1}
		\frac{1}{4}\Rdcov^2(X,Y)=\int_{\mathbb{R}^2}\left(F^{X,Y}(x,y)-F^X(x)F^Y(y)\right)^2\,dF^{X}(x)\,dF^{Y}(y) \qquad \mathrm{and},
		\end{equation}
		\begin{equation}\label{eq:dcovequivmain2}
		\frac{1}{4}\Rdcov_n^2=\int \left(F^{X,Y}_n(x,y)-F^X_n(x)F^Y_n(y)\right)^2\,dF^X_n(x)\,dF^Y_n(y).
		\end{equation}
	\end{lemma}
	\noindent We are now interested in two fundamental questions about our proposed test: (a) ``What is the limiting distribution of our test statistic?''; (b) ``Is our test consistent against all fixed alternatives, as the sample size grows?''. We investigate these two questions in Theorems~\ref{theo:indepasdistn} and~\ref{theo:indepconsis} respectively (see~\cref{proof:indepasdistn} and~\cref{proof:indepconsis} for the proofs).
	
	\begin{theorem}\label{theo:indepasdistn}
		Under assumptions \textbf{(AP1)}, \textbf{(AP2)} and under $\mathrm{H}_0$, there exists universal nonnegative constants $(\eta_1,\eta_2,\ldots )$ such that
		\begin{equation*}
		n\Rdcov_n^2\overset{w}{\longrightarrow} \sum_{j=1}^{\infty}\eta_j Z_j^2\qquad \mathrm{as } \;\; n\to\infty
		\end{equation*}
		where $Z_1,Z_2,\ldots $ are i.i.d.~standard Gaussian random variables. In fact, $\eta_j$'s do not depend on the specific choice of $\mathcal{H}_n^{d_1}$ or $\mathcal{H}_n^{d_2}$ as long as \textbf{(AP2)} is satisfied.
	\end{theorem}
	
	\begin{remark}[Limiting distribution]\label{rem:connectdcov}
		The limiting distribution in~\cref{theo:indepasdistn} is exactly the same as that of usual distance covariance (under $\mathrm{H}_0$) when $\mu_{\mathbf{X}}=\mathcal{U}^{d_1}$ and $\mu_{\mathbf{Y}}=\mathcal{U}^{d_2}$ (see~\cite[Theorem 5]{Gabor2007}).\end{remark}
	\begin{remark}[Distribution-freeness]\label{rem:asdistfreegood}
		Note that the asymptotic distribution of the usual distance covariance statistic, given in~\cite[Theorem 5]{Gabor2007}, depends on $\mu_{\mathbf{X}}$ and $\mu_{\mathbf{Y}}$, which are unknown. As a result, even for large $n$, in practice, one usually has to resort to resampling/permutation techniques or further worst case approximations (see~\cite[Theorem 6]{Gabor2007}) to determine the critical value of the test. Having finite sample (and asymptotic) distribution-freeness avoids the need for such approximation techniques (for small as well as large $n$). In~\cref{sec:univcutas}, we discuss (computationally) how large $n$ should be (depending on $d_1$ and $d_2$) so as to use quantiles from the asymptotic distribution of $n\Rdcov_n^2$ to approximate thresholds for our testing procedure. In~\cref{table:thresholdind} (in~\cref{sec:univcutas}), we provide the universal asymptotic 0.95-quantiles as $d_1,d_2$ varies (for $d_1,d_2\leq 8$).
	\end{remark}
	
	\begin{remark}[Our proof technique]\label{rem:diffanalysis}
		Observe that, contrary to the study of the usual distance covariance~\cite{Gabor2007} which can be analyzed using standard techniques from empirical process theory (as in~\cite[Theorem 5]{Gabor2007}) or results from degenerate V-statistics (as used in~\cite[Theorem 2.7]{Russell2013}), the study of $\Rdcov_n^2$ is more complicated as it involves dependent multivariate ranks. Our main technique for proving~\cref{theo:indepasdistn} is to use Hoeffding's Combinatorial Central Limit Theorem (see e.g.,~\cite{Chen2015}  or~\cref{prop:HCLT}). In the process, we prove some results on permutation statistics (see~\cref{lem:rankdcov}) which may be of independent interest. 
	\end{remark}

	The following result (proved in~\cref{proof:indepconsis}) shows that our proposed testing procedure yields a consistent sequence of tests under fixed alternatives (i.e., the power of our test converges to 1, as the sample size increases, for any fixed alternative).    
	\begin{theorem}\label{theo:indepconsis}
		Under assumptions \textbf{(AP1)} and \textbf{(AP2)}, $$\Rdcov_n^2\overset{a.s.}{\longrightarrow}\Rdcov^2(\mathbf{X},\mathbf{Y}) \quad \mathrm{as }\;\; n\to\infty$$ where $(\mathbf{X},\mathbf{Y})\sim \mu$. Moreover, $\mathbb{P}(n\Rdcov_n^2>c_n)\overset{}{\longrightarrow}1$, as $n\to\infty$, provided $\mu\neq \mu_{\mathbf{X}}\otimes \mu_{\mathbf{Y}}$. 
	\end{theorem}
	
	\begin{remark}[Minimal assumptions]\label{rem:lesassume}
		The proof of~\cref{theo:indepconsis}  reveals that only the $L^2$-convergence of empirical transport maps (see~\cref{lem:L2}) is necessary. Therefore, by resorting to a weaker form of convergence (as compared to the $L^{\infty}$-convergence as in~\cite{chernozhukov2017monge,ghosal2019multivariate,del2018center}) we have effectively reduced the set of assumptions needed on $\mu_{\mathbf{X}}$ and $\mu_{\mathbf{Y}}$ for getting a consistent sequence of tests (contrary to~\cite{ghosal2019multivariate}). Moreover we are able to establish consistency without any moment assumptions (contrary to~\cite[Theorem 2]{Gabor2007}).  
	\end{remark}
	
	\begin{remark}[Halton sequence]\label{rem:haltonprop}
		\cref{cor:weakcon} ensures that assumption \textbf{(AP2)} is satisfied for the Halton sequence (see~\cref{sec:revQMC} for details). The same is true for other pseudo-random sequences (see~\cref{sec:revQMC} for examples).
	\end{remark}
	
	\begin{remark}[Invariance under coordinate-wise monotone transformations]\label{rem:ranktransone}
		An alternate approach to testing mutual independence would be to transform the observed data into their marginal one-dimensional ranks first and then construct the multivariate ranks based on this transformed data. Let us elaborate on this briefly. Let us write $\mathbf{X}_{i}=(X_{i1},X_{i2},\ldots ,X_{id_1})$ in terms of its univariate components, for $1\leq i\leq n$. For $1\leq j\leq d_1$, construct $\tilde{\mathbf{X}}_i$ such that $\tilde{X}_{ij}$ equals the usual one-dimensional rank of $X_{ij}$ among $X_{1j},\ldots ,X_{nj}$. Repeat the same exercise with the $\mathbf{Y}_i$'s to form $\tilde{\mathbf{Y}}_i$'s. Now, consider $\{(\tilde{\mathbf{X}}_i,\tilde{\mathbf{Y}}_i)\}_{i=1}^{n}$ and obtain multivariate ranks of $\tilde{\mathbf{X}}_i$'s and $\tilde{\mathbf{Y}}_i$'s using measure transportation as described above (see~\eqref{eq:empopt}). Finally, calculate a suitable test statistic for independence (such as $\Rdcov_n^2$) based on these ranks. This approach has natural connections to copula based methods (see~\cite{kojadinovic2009}) and ensures that the constructed tests will be invariant under coordinate-wise monotone transformations of the data (cf.~\cref{lem:proprdcov}, part (d)). We believe that an analogous theoretical analysis can be carried out for this modified procedure. 
	\end{remark}
	
	\subsection{Distribution-free multivariate two-sample testing}\label{sec:twogofp}
	Here we shall consider the two-sample goodness-of-fit testing problem in a multivariate setting. Suppose $\mathbf{X}_1,\ldots ,\mathbf{X}_m\overset{i.i.d.}{\sim}\mu_{\mathbf{X}}$ and $\mathbf{Y}_1,\ldots ,\mathbf{Y}_n\overset{i.i.d.}{\sim}\mu_{\mathbf{Y}}$ (independent of the $\mathbf{X}_i$'s), where we assume that 
	$$\textbf{(AP3)}:\qquad\qquad\mu_{\mathbf{X}},\mu_{\mathbf{Y}}\in\mathcal{P}_{ac}(\mathbb{R}^d).$$
	
	We are interested in testing the hypothesis:
	\begin{align*}
	\mathrm{H}_0: \mu_{\mathbf{X}}=\mu_{\mathbf{Y}}\qquad \qquad \mathrm{versus}\qquad\qquad \mathrm{H}_1:\mu_{\mathbf{X}}\neq\mu_{\mathbf{Y}}.
	\end{align*}
	The two-sample problem (or its multi-sample extension) has been studied in great detail over the years. In this context, rank and data-depth based methods have mostly been restricted to testing against location-scale alternatives, see e.g.,~\cite{hettmansperger1998,randles1990,motto1995}. Distribution-free depth-based tests which are consistent if restricted to the above class of alternatives are discussed in~\cite{liu2010versatile,rousson2002}. An alternative route for testing against general alternatives includes graph based tests such as in~\cite{friedman1979}, where the authors construct a test based on the minimum spanning tree of a graph with the data points as its vertices and pairwise distances as edge weights. Various interesting modifications and extensions to this test have been proposed in literature, see e.g.,~\cite{chen2017,henze1988,schilling1986,petrie2016}. Theoretical properties of all these tests can be studied under a unified framework as shown in~\cite{bbb2019}. 
	
	As mentioned in the Introduction, the only other multivariate nonparametric distribution-free two-sample goodness-of-fit test that has the same guarantees as our approach, can be attributed to~\cite{rosenbaum2005} (also see~\cite{bbbm2019,Castro2016} for subsequent theoretical analysis). In~\cite{rosenbaum2005}, Rosenbaum constructed his proposed test statistic from a minimum non-bipartite matching (see~\cite{bo2011}) of the pooled sample of observations. It has motivated numerous extensions and applications in real-life problems (see e.g.,~\cite{petrie2016,heller2010,hellerrosenbaum2010using}). This test has been recently extended to a $K$-sample version in~\cite{bbbm2019}. Another graph-based distribution-free test for this problem was proposed in~\cite{munmun2014}; however, the test becomes computationally infeasible even for moderate sample sizes.
	
	Yet another class of pairwise-distance based tests use ideas from Reproducing Kernel Hilbert Spaces (RKHS), see e.g.,~\cite{gretton2009fast,gretton2012}. The principle idea here is to embed probability distributions in RKHSs through what are called \textit{mean embeddings} and measure goodness-of-fit between two distributions by the Hilbert-Schmidt norm between the corresponding mean embeddings. These kernel based measures can alternatively be expressed as \emph{probability integral metrics} which equal $0$ if and only if the underlying distributions are exactly the same. In fact, the energy statistic (see~\cite{Gabor2013,Baringhaus2004}) --- a popular and powerful goodness-of-fit measure ---  can also be viewed as a special case of kernel based methods (see~\cite{sejdinovic2013}). Due to its simplicity, the energy distance has been studied and applied extensively over the past decade, as we have already highlighted in the Introduction. However, note that a common disadvantage of these kernel based methods (including the usual energy statistic) is that they are not exactly distribution-free.

	In this subsection, we propose the rank energy statistic --- a distribution-free goodness-of-fit measure based on the energy distance --- for testing the equality of two multivariate distributions. We describe our method below. We will use $\mu_m^{\mathbf{X}}$ and $\mu_n^{\mathbf{Y}}$ to denote the empirical distributions on $\mathcal{D}_m^{\mathbf{X}}\coloneqq \{\mathbf{X}_1,\ldots ,\mathbf{X}_m\}$ and $\mathcal{D}_n^{\mathbf{Y}} :=\{\mathbf{Y}_1,\ldots ,\mathbf{Y}_n\}$ respectively. Let $$\mu_{m,n}^{\mathbf{X},\mathbf{Y}}:=(m+n)^{-1}(m\mu_n^{\mathbf{X}}+n\mu_n^{\mathbf{Y}})$$ and let $\mathcal{H}^d_{m+n}\coloneqq \{\mathbf{h}_1^d,\ldots ,\mathbf{h}_{m+n}^d\} \subset [0,1]^d$ denote the (fixed) sample multivariate ranks. We will further work under the following assumption on $\mathcal{H}_{m+n}^d$:

	\textbf{(AP4)} The empirical distribution on $\mathcal{H}_{m+n}^d$ converges weakly to $\mathcal{U}^d$ as $\min{(m,n)}\to\infty$.
	Note that choosing $\mathcal{H}_{m+n}^d$ to be the $d$-dimensional Halton sequence for $d\geq 2$, and $\{i/(m+n):1\leq i\leq m+n\}$ for $d=1$, ensures that \textbf{(AP4)} is satisfied (see~\cref{cor:weakcon} for details).   
	
	Finally, we shall use $\hat{R}_{m,n}^{\mathbf{X},\mathbf{Y}}(\cdot)$ to denote the joint empirical rank map (see~\cref{def:empquanrank}) corresponding to the transportation of $\mu_{m,n}^{\mathbf{X},\mathbf{Y}}$ to the empirical distribution on $\mathcal{H}_{m+n}^d$. The rank energy statistic is defined as:
	\begin{align}\label{eq:twogofp}
	\Ren_{m,n}^2&\coloneqq \frac{2}{mn}\sum\limits_{i=1}^m \sum\limits_{j=1}^n \lVert \hat{R}_{m,n}^{\mathbf{X},\mathbf{Y}}(\mathbf{X}_i)-\hat{R}_{m,n}^{\mathbf{X},\mathbf{Y}}(\mathbf{Y}_j)\rVert-\frac{1}{m^2}\sum\limits_{i,j=1}^m \lVert \hat{R}_{m,n}^{\mathbf{X},\mathbf{Y}}(\mathbf{X}_i)-\hat{R}_{m,n}^{\mathbf{X},\mathbf{Y}}(\mathbf{X}_j)\rVert\nonumber \\&\qquad\qquad-\frac{1}{n^2}\sum\limits_{i,j=1}^n \lVert\hat{R}_{m,n}^{\mathbf{X},\mathbf{Y}}(\mathbf{Y}_i)-\hat{R}_{m,n}^{\mathbf{X},\mathbf{Y}}(\mathbf{Y}_j)\rVert.
	\end{align}
	Observe that the right hand side of~\eqref{eq:twogofp} can be viewed as an empirical version of $\Ren$ (see~\eqref{eq:ren1}) through its alternate expression as in~\cref{lem:propren} (part (a)). $\Ren_{m,n}^2$ can also be viewed as a rank-transformed version of the empirical energy measure as in~\cite[Equation (6.1)]{Gabor2013}. Due to space constraints, we will refer the interested reader to~\cite{Gabor2013} for further motivation of the energy statistic. By~\cite[Equation (5)]{Baringhaus2004}, it is easy to see that the right side of~\eqref{eq:twogofp} is always nonnegative. Just as $\Rdcov_n^2$ (in~\eqref{eq:indepstat}), $\Ren_{m,n}^2$ above can also be computed in $\mathcal{O}(m n d)$ steps (see~\cite{Zhao2015}) given the corresponding vector of multivariate ranks. In the following lemma we illustrate the distribution-free property of $\Ren_{m,n}^2$ (see~\cref{proof:twosamdistfree} for a proof).
	\begin{lemma}\label{rem:twosamdistfree}
		Under assumption \textbf{(AP3)} and under $\mathrm{H}_0$, the distribution of $\Ren_{m,n}^2$, as defined in~\eqref{eq:twogofp}, is free of $\mu_{\mathbf{X}}  \equiv \mu_{\mathbf{Y}}$.
	\end{lemma}
	\noindent \textbf{Distribution-free two-sample testing procedure}: Given a (prespecified) type-I error level $\alpha \in (0,1)$, let $$c_{m,n}\coloneqq \inf\{c >0: \mathbb{P}_{\textrm{H}_0}(mn(m+n)^{-1}\Ren_{m,n}^2\geq c)\leq \alpha\}.$$ As $\Ren_{m,n}^2$ is distribution-free under $\mbox{H}_0$ (by~\cref{rem:twosamdistfree}), so is $c_{m,n}$. Given $c_{m,n}$, our proposed testing procedure rejects $\mathrm{H}_0$ if $mn(m+n)^{-1}\Ren_{m,n}^2\geq c_{m,n}$ and accepts $\mathrm{H}_0$ otherwise. Clearly, this results in a level $\alpha$ test. \par
	An interesting feature of our proposed statistic $\Ren_{m,n}^2$ is its equivalence with the celebrated Cram\'{e}r-von Mises statistic for two sample equality of distributions testing (see e.g.,~\cite{AndersonCVM1962} and the right side of~\eqref{eq:gofeqcvmmain1} below) when $d=1$. The following lemma (see~\cref{proof:gofeqcvm} for a proof) makes this connection precise.
	\begin{lemma}\label{lem:gofeqcvm}
		For $d=1$, let $F_m^{X}$, $G_n^{Y}$ and $H_{m+n}^{X,Y}$ denote the empirical distribution functions on $\{X_1,\ldots ,X_m\}$, $\{Y_1,\ldots ,Y_n\}$ and the pooled sample respectively. Then,
		\begin{equation}\label{eq:gofeqcvmmain1}
		\frac{1}{2}\,\Ren_{m,n}^2=\int \bigg(F_m^X(t)-G_n^Y(t)\bigg)^2\,dH_{m+n}^{X,Y}(t).
		\end{equation}
		The right hand side of~\eqref{eq:gofeqcvmmain1} is the exact Cram\'{e}r-von Mises statistic as in~\cite{AndersonCVM1962}. At the population level, fix any $\lambda\in (0,1)$ and let $F^X$, $G^Y$ and $H^{X,Y}_{\lambda}$ be the distribution functions associated with the probability measures $\mu_{X}$, $\mu_{Y}$ and $\lambda \mu_X+(1-\lambda)\mu_Y$. Assume also that $F^X$ and $G^Y$ are absolutely continuous. Then,
		\begin{equation}\label{eq:gofeqcvmmain2}
		\frac{1}{2}\Ren_{\lambda}^2(X,Y)=\int_{-\infty}^{\infty} \bigg(F^X(t)-G^Y(t)\bigg)^2\,dH^{X,Y}_{\lambda}(t).
		\end{equation}
	\end{lemma}
	Next we find the asymptotic distribution of $\Ren_{m,n}^2$ in Theorem~\ref{theo:twosamasdistn} and prove the consistency of our proposed procedure in Theorem~\ref{theo:twosamconsis}; see~\cref{proof:indepasdistn} and~\cref{proof:twosamconsis} for their proofs.
	\begin{theorem}\label{theo:twosamasdistn}
		Suppose that $\min{(m,n)}\to\infty$. Under assumptions \textbf{(AP3)}, \textbf{(AP4)} and under $\textrm{H}_0$, we have:
		\begin{equation*}
		\frac{mn}{m+n}\Ren_{m,n}^2\overset{w}{\longrightarrow} \sum_{j=1}^{\infty} \tau_j Z_j^2\qquad \mathrm{as } \;\;n\to\infty
		\end{equation*}
		where $Z_1,Z_2,\ldots $ are i.i.d.~standard normals and $\tau_j$'s are fixed nonnegative constants. In fact, $\tau_j$'s do not depend on the specific choice of $\mathcal{H}_{m +n}^{d}$ as long as \textbf{(AP4)} is satisfied.
	\end{theorem}
	\begin{remark}[Limiting distribution]\label{rem:connectenergy}
		The limiting distribution in~\cref{theo:twosamasdistn} is exactly the same as that of the usual energy statistic (under $\mathrm{H}_0$) when $\mu_{\mathbf{X}}=\mu_{\mathbf{Y}}=\mathcal{U}^{d}$ (\cite[Theorem 2.3]{Baringhaus2004}).\end{remark}

	\begin{theorem}\label{theo:twosamconsis}
		Suppose that $m/(m+n)\longrightarrow\lambda\in (0,1)$. Then, under assumptions \textbf{(AP3)} and \textbf{(AP4)}, $$\Ren_{m,n}^2\overset{a.s.}{\longrightarrow}\Ren_{\lambda}^2(\mathbf{X},\mathbf{Y}) \qquad \mathrm{as } \; n\to\infty,$$ where $\mathbf{X}\sim\mu_{\mathbf{X}}$ and $\mathbf{Y}\sim \mu_{\mathbf{Y}}$ (note the connection with~\cref{rem:lamclarify}). Moreover, $\mathbb{P}(mn(m+n)^{-1}\Ren_{m,n}^2>c_n)\overset{}{\longrightarrow} 1$, as $n\to\infty$, provided $\mu_{\mathbf{X}}\neq\mu_{\mathbf{Y}}$.
	\end{theorem}

	The multivariate two-sample testing procedure described above bears all the useful properties of our independence testing procedure from~\cref{sec:Indep}. In particular, the proposed test is distribution-free for each fixed $m$ and $n$ and also in an asymptotic sense. In~\cref{sec:univcutas}, we study, using simulations, how large $m,n$ should be (depending on $d$) so as to reasonably use quantiles from the asymptotic distribution of $mn(m+n)^{-1}\Ren_{m,n}^2$ to determine thresholds for our testing procedure. In~\cref{table:twosamascutoff} (in~\cref{sec:univcutas}), we provide universal asymptotic quantiles ($5\%$) up to $d\leq 8$. 
	
	Our proposed test is also consistent against fixed alternatives without any moment assumptions, as opposed to the usual test based on the energy statistic (see~\cite{Gabor2013,Baringhaus2004}). Moreover, we are also able to reduce the smoothness assumptions on the underlying measures $\mu_{\mathbf{X}}$ and $\mu_{\mathbf{Y}}$ necessary for consistency (cf.~\cite[Proposition 5.2]{ghosal2019multivariate} and~\cite[Theorem 3.1]{boeckel2018multivariate}).

	
	\subsection{Extensions to the $K$-sample problem}\label{sec:tempksam}
	The methods we discussed in Sections~\ref{sec:Indep} and~\ref{sec:twogofp} have natural extensions to the $K$-sample setting; namely, testing for mutual independence of $K$ random vectors, and multivariate goodness-of-fit testing for $K$ populations (as mentioned in the Introduction). Using the same principles as above, we can again construct exact distribution-free tests for the above problems that will be consistent against all fixed alternatives. Due to space constraints, we relegate a detailed  discussion of this to~\cref{sec:addis}; in particular, see Propositions~\ref{prop:multextind} and~\ref{prop:multextgof}.

	\section{Constructing other distribution-free tests --- testing for symmetry}\label{sec:otherp}
	As we have discussed in the Introduction, our proposed recipe --- of using the notion of multivariate ranks obtained from the theory of measure transportation to define a suitable test statistic --- can also be used to construct distribution-free tests in other nonparametric testing problems (besides those discussed in~\cref{sec:mtp}). Let us illustrate this by constructing a distribution-free nonparametric test of multivariate symmetry. 
	
	The notion of symmetry in one-dimensional distributions is quite unambiguous. We say $X\sim \mu$ is symmetric if and only if  $X\overset{d}{=}-X$. This notion makes perfect sense even in dimensions larger than $1$, although there are various other notions of symmetry that might also be of interest (see~\cite{serfling2014multivariate}) when $d>1$. Nevertheless, we will focus on providing a distribution-free test based on the above notion in the multivariate setting. A comprehensive review of the literature on tests of multivariate symmetry is beyond the scope of this paper; we therefore refer the interested reader to~\cite{aki1993,beran1979,heathcote1995} and the references therein.

	So our problem may be stated as follows: Given i.i.d.~data $\mathbf{X}_1,\ldots ,\mathbf{X}_n\sim \mu\in\mathcal{P}_{ac}(\mathbb{R}^d)$, we want to test the hypothesis:
	\begin{equation}\label{eq:otherpsym}
	\textrm{H}_0: \mathbf{X}_1\overset{d}{=}-\mathbf{X}_1\qquad\qquad \textrm{versus} \qquad\qquad\textrm{H}_1: \textrm{ not H}_0.
	\end{equation}
	Observe that~\eqref{eq:otherpsym} may be interpreted as a two-sample equality of distributions testing problem, except that the collection of $\mathbf{X}_i$'s is not independent of the collection of $-\mathbf{X}_i$'s  --- a crucial difference with the two-sample problem setting discussed in~\cref{sec:twogofp}. In fact, if we pool the $\mathbf{X}_i$'s and $-\mathbf{X}_i$'s, then the pooled sample ranks are no longer uniformly distributed over the set of all $(2n)!$ permutations of the elements of the set $\mathcal{H}_{2n}^d$ (even for $d=1$). As a result, we need to be more careful while defining the joint multivariate ranks ($\hat{R}_n(\cdot)$) for this problem.

	We will propose a distribution-free test for this problem in the following three steps:\par 
	\noindent \textbf{(I)} Set $\mathbf{Z}_i :=(\mathbf{X}_i,-\mathbf{X}_i)$ for $1\leq i\leq n$. Consider the $2d$-dimensional (fixed) sample multivariate ranks $\mathcal{H}_n^{2d}$ and let $\tilde{R}_n(\cdot)$ denote the corresponding empirical transport map for the $\mathbf{Z}_i$'s (obtained by solving~\eqref{eq:empopt}). Set $\tilde{R}_n(\mathbf{Z}_{i})\coloneqq\big(\tilde{R}_n(\mathbf{Z}_{i})^{d:},\tilde{R}_n(\mathbf{Z}_{i})^{:d}\big)$ for $1\leq i\leq n$, where $\tilde{R}_n(\mathbf{Z}_{i})^{d:}$ and $\tilde{R}_n(\mathbf{Z}_{i})^{:d}$ denote the first and the last $d$ components of $\tilde{R}_n(\mathbf{Z}_{i})$ respectively.\par 
	\noindent\textbf{(II)} For all $1\leq i\leq n$, solve another empirical transportation problem between $\{\mathbf{X}_{i},-\mathbf{X}_{i}\}$ and $\big\{\tilde{R}_n(\mathbf{Z}_{i})^{d:},\tilde{R}_n(\mathbf{Z}_{i})^{:d}\big\}$ to get the map $R^*_{n,i}(\cdot):\{\mathbf{X}_{i},-\mathbf{X}_{i_0}\}\to \big\{\tilde{R}_n(\mathbf{Z}_{i})^{d:},\tilde{R}_n(\mathbf{Z}_{i})^{:d}\big\}$. To conclude, define $\hat{R}_n(\mathbf{X}_{i}) :=R_{n,i}^*(\mathbf{X}_{i})$ and $\hat{R}_n(-\mathbf{X}_{i}) :=R_{n,i}^*(-\mathbf{X}_{i})$ for $1\leq i\leq n$.\par 
	\noindent\textbf{(III)} Now, we can use any two-sample goodness-of-fit test statistic for testing~\eqref{eq:otherpsym}, e.g., the energy statistic from~\cite{Gabor2013} to test for equality of distributions between the ranks ($\hat{R}_n(\cdot)$) corresponding to $\mathbf{X}_i$'s and those corresponding to $-\mathbf{X}_i$'s. Let us call this statistic $T_n$.\par 
	\begin{lemma}\label{lem:multsym}
		If $\mu\in\mathcal{P}_{ac}(\mathbb{R}^d)$, then the distribution of $T_n$ (as in \textbf{(III)} above) is universal (free of $\mu$) under $\mathrm{H}_0$.
	\end{lemma}
	Note that~\cref{lem:multsym} (proved in~\cref{proof:multsym}) demonstrates the exact distribution-free nature of the test of multivariate symmetry proposed above. We leave the detailed theoretical analysis of this test as a subject for future research.

	While discussing the full scope of our proposal for constructing distribution-free tests is beyond the scope of this paper, we would like to refer the interested reader to the various applications of distance correlation and the energy statistic as elucidated in~\cite{Gabor2007,szekely2009} and~\cite{Gabor2013}, such as hierarchical clustering, detecting ``influential" observations, testing for non-linear dependence, change point analysis, etc. Our methodology suggests that one can possibly design distribution-free procedures for the above inference problems by using our ideas.

	\section{Numerical experiments}\label{sec:sim}
	In this section, we will discuss the empirical performance of our proposed tests in a wide variety of settings. Due to space constraints, we will restrict to $\Rdcov_n^2$ here; the performance of $\Ren^2_{m,n}$ will be discussed in detail in~\cref{sec:syngof}.
	
	\subsection{Synthetic data experiments for mutual independence testing}\label{sec:synind}
	We illustrate the empirical performance of $\Rdcov_n^2$ (hereafter referred to as $\Rdcov$) for testing mutual independence of two random vectors, based on synthetic data. Throughout our simulation settings, we fix $n=200$, $d_1=3$ and $d_2=3$. In \cref{table:indepcompare}, we compare the performance of our method with the following standard methods already existing in literature and already available in the \texttt{R} software: Pearson's correlation (P) (\cite{pearson1920notes}; computed using the \texttt{stats} package), distance covariance (DCoV) (\cite{Gabor2007}; implemented in \texttt{energy} package), Hilbert-Schmidt Independence Criteria (HSIC) (\cite{gretton2008kernel}) from the \texttt{dHSIC} package, mutual information (MINT) (\cite{berrett2017nonparametric}) from the \texttt{IndepTest} package, and Heller's graph-based test (HHG) (\cite{heller2013}) from the \texttt{HHG} package. 
	
	As a general rule, unless otherwise specified, we construct $\mathbf{X} \in \R^{3}$ (and $\mathbf{Y} \in \R^3$) from 3 independent random variables drawn from $X$ (and $Y$) according to the following settings:
	\begin{itemize}
		\item[(V1)] $A\sim \mathcal{N}(0,1)$, $X\sim 0.2\times \mbox{Cauchy}(0,1)+A$ and $Y\sim 0.2\times \mbox{Cauchy}(0,1)+A$.
		\item[(V2)] $X\sim \mathcal{U}(-1,1)$ and $Y=(X^2+\mathcal{U}(0,1))/2$.
		\item[(V3)] $X\sim \mathcal{N}(0,2)$, $E\sim \mbox{Ber}(0.04)$, $V\sim\mathcal{N}(0,2)$ and $Y=(1-E) V+E X$. $E$ and $V$ are independent.
		\item[(V4)] $W\sim\mathcal{U}(-1,1)$, $W_1\sim\mathcal{U}(0,1)$, $W_2\sim\mathcal{U}(0,1)$, $V_1=W+W_1/3$ and $V_2=4\times (W^2-0.5)^2+W_2$. Finally, $X=V_1$ and $Y+A \times \mathcal{N}(5,1)+(1-A) V_2$. $W$, $W_1$ and $W_2$ are independent.
		\item[(V5)] $(U_1,U_2,U_3,V_1,V_2,V_3)\sim\mathcal{N}_6(\mathbf{0},\Sigma)$ is independent of $(W_1,W_2,W_3,Z_1,Z_2,Z_3)\sim\mathcal{N}_6(\mathbf{1},\Sigma/2)$, where $\Sigma_{ii}=1$ and $\Sigma_{ij}=0.3$ if $i\leq 3,j>3$ or $i>3,j\leq 3$. Finally, $(X_1,X_2,X_3)\sim (1-A_1)(U_1,U_2,U_3)+A_1(W_1,W_2,W_3)$ and $(Y_1,Y_2,Y_3)\sim (1-A_2) (V_1,V_2,V_3)+A_2 (Z_1,Z_2,Z_3)$, $A_1\sim \mbox{Ber}(0.5)$ and $A_2\sim \mbox{Ber}(0.3)$. 
		\item[(V6)] $A\sim \mathcal{N}(0,1)$, $X\sim \mbox{Pareto}(1,2)^2+A$ and $Y\sim \mbox{Pareto}(1,1)^2+A$. 
		\item[(V7)] $\epsilon\sim\mathcal{N}(0,5)$, $X\sim\mathcal{U}(0,1)$, $Y=X^{1/4}+\epsilon$. Here $X$ and $\epsilon$ are independent.
		\item[(V8)] $X\sim\mathcal{N}(0,1)$ and $Y=\log{(4X^2)}$.
		\item[(V9)] $A\sim \mathcal{N}(0,1)$, $X\sim |A+\mbox{Pareto}(1,1)|^{1.5}$ and $Y\sim |A+\mbox{Pareto}(1,1)|^{1.5}$.
		\item[(V10)] Same setup as in (V5), but with $A_2=0$.
	\end{itemize}
	Our simulation settings are similar to a variety of settings popular in the mutual independence testing literature, see e.g.,~\cite{heller2013,munmun2016,Gabor2007}. For example, settings (V2) and (V4) are from~\cite{Newton2009} and have been used later in~\cite{heller2013,munmun2016}. Using mixture distributions to distinguish between different methods of independence testing is also common in literature (e.g.,~\cite{heller2013} and~\cite[Settings 3 and 5 from Table 1]{munmun2016}). This is the motivation behind settings (V3), (V5) and (V10). Note that (V7) is a simple regression model; (V8) was first used in~\cite{Gabor2007}. Finally, settings (V1), (V6) and (V9) have been chosen to illustrate the superior performance of $\Rdcov$ when dealing with heavy-tailed distributions, as has been highlighted in the Introduction. In \cref{table:indepcompare} we present our findings. The two columns corresponding to each method represents the rejection probabilities (estimated from $1000$ independent replications) at nominal levels $0.05$ and $0.1$.
	
	\begin{table}[h]
		\centering
		\begin{tabular}{| *{13}{Sc|}}
			\hline
			& \multicolumn{2}{c|}{(P)} & \multicolumn{2}{c|}{(DCoV)} & \multicolumn{2}{c|}{(HSIC)} & \multicolumn{2}{c|}{(MINT)} & \multicolumn{2}{c|}{(HHG)} & \multicolumn{2}{c|}{($\Rdcov$)}\\ 
			\hline
			V1 & \textbf{1.00} & \textbf{1.00} & 0.81 & 0.60 & \textbf{1.00} & \textbf{1.00} & 0.93 & 0.84 & \textbf{1.00} & \textbf{1.00} & \textbf{1.00} & \textbf{1.00} \\ 
			\hline
			V2 & 0.14 & 0.09 & \textbf{1.00} & \textbf{1.00} & \textbf{1.00} & \textbf{1.00} & \textbf{1.00} & \textbf{1.00} & \textbf{1.00} & \textbf{1.00} & \textbf{1.00} & \textbf{1.00} \\ 
			\hline
			V3 & 0.15 & 0.08 & 0.15 & 0.08 & 0.15 & 0.09 & \textbf{0.25} & \textbf{0.16} & 0.14 & 0.09 & 0.14 & 0.07 \\ 
			\hline
			V4 & 0.09 & 0.05 & 0.50 & 0.32 & \textbf{0.95} & \textbf{0.88} & 0.20 & 0.13 & 0.92 & 0.82 & 0.72 & 0.52 \\ 
			\hline
			V5 & 0.12 & 0.06 & 0.18 & 0.09 & 0.21 & 0.10 & 0.95 & 0.89 & 1.00 & 1.00 & \textbf{0.77} & \textbf{0.63} \\ 
			\hline
			V6 & 0.60 & 0.48 & 0.08 & 0.04 & 0.38 & 0.24 & 0.11 & 0.05 & 0.27 & 0.14 & \textbf{0.98} & \textbf{0.95} \\ 
			\hline
			V7 & 0.25 & 0.16 & \textbf{0.26} & \textbf{0.16} & 0.22 & 0.13 & 0.11 & 0.06 & 0.19 & 0.10 & 0.21 & 0.13 \\ 
			\hline
			V8 & 0.15 & 0.09 & \textbf{1.00} & \textbf{1.00} & \textbf{1.00} & \textbf{1.00} & 0.75 & 0.62 & \textbf{1.00} & \textbf{1.00} & \textbf{1.00} & \textbf{1.00} \\ 
			\hline
			V9 & 0.38 & 0.28 & 0.10 & 0.05 & 0.16 & 0.08 & 0.11 & 0.05 & 0.21 & 0.10 & \textbf{0.93} & \textbf{0.89} \\ 
			\hline
			V10 & \textbf{0.92} & \textbf{0.86} & 0.90 & 0.83 & 0.80 & 0.68 & 0.32 & 0.21 & 0.58 & 0.44 & 0.82 & 0.69 \\ 
			\hline
		\end{tabular}
		\vspace{0.1in}
		\caption{Proportion of times the null hypothesis was rejected across $10$ settings. Here $n=200$, $d_1=d_2=3$. The tests with the best empirical performances have been highlighted in bold.}
		\label{table:indepcompare}
	\end{table}
	Below we discuss the performance of $\Rdcov_n^2$ along with the competing procedures.
	
	\noindent\textbf{(P)}: As one may expect, most of the methods (including ours) outperform the Pearson's correlation based testing procedure consistently. The best performances of Pearson's correlation can be observed in settings where the dependence arises largely due to a linear relationship between the variables, such as (V1), (V6) and (V10). Note that $\Rdcov$ performs just as well in (V1), a lot better in (V6) and slightly worse in (V10).\par 
	\noindent \textbf{(DCoV)}: As distance covariance requires finite moments for consistency, it is expected to perform worse in the heavy-tailed settings such as (V1), (V6) and (V9). This is clearly seen in~\cref{table:indepcompare} where $\Rdcov$ outperforms usual distance covariance convincingly. Surprisingly, in (V5) where both $\mathbf{X}$ and $\mathbf{Y}$ are Gaussian mixtures, DCoV performs rather poorly whereas $\Rdcov$ performs significantly better. This phenomenon can also be seen in another setting (V4) where the coordinates of $\mathbf{Y}$ arise out of a mixture distribution. Across all the other settings, namely (V2), (V3), (V7), (V8) and (V10), observe that $\Rdcov$ and DCoV have almost identical performance. This leads us to believe that $\Rdcov$ can be expected to perform similar to DCoV except that $\Rdcov$ is a lot more robust to the presence of heavy-tailed distributions.

	\noindent \textbf{(HSIC)}: An interesting observation is that, in settings (V6) and (V9), where we have introduced a Pareto noise to Gaussian data, $\Rdcov$ performs better than the test based on the HSIC (with a bounded Gaussian kernel). Note that with bounded kernels, HSIC  does not need finite moment assumptions for consistency. In settings (V4) and (V5), both of which involve mixture distributions, the performance of HSIC fluctuates heavily, outperforming $\Rdcov$ in (V4) and underperforming in (V5). In all the other settings, its performance is more or less similar to $\Rdcov$.\par 
	\noindent \textbf{(MINT)}: Table~\ref{table:indepcompare} reveals that in almost all the settings $\Rdcov$ outperforms MINT, except in (V3) (where all methods perform poorly) and in (V5). This could be because the consistency type results in mutual information based tests (see~\cite[Theorem 4]{berrett2017nonparametric}) require a number of regularity conditions which may not hold for some of the settings above. Moreover, in (V6), (V7) and (V9), MINT performs only as well as a random guess (i.e., the trivial unbiased test).\par
	\noindent \textbf{(HHG)}: The superior performance of our test in the heavy-tailed settings compared to competing methods persists in the case of HHG as well. This leads us to believe, in general, that $\Rdcov$ is perhaps a {better choice when dealing with heavy-tailed distributions} than tests based on distances between data points. Note however that HHG (being based on ranks of distances) also does not require finite moment assumptions for consistency. Having said that however, the performance of HHG is definitely very competitive. Apart from the usual heavy-tailed distribution settings, the only other settings where $\Rdcov$ does slightly better would be (V7) and (V10). In fact, (V10) is a slightly modified version of a Gaussian mixture model, where the dependence is mostly linear (as is supported by the superior performance of Pearson's correlation). Our general sense is that, for near-linear dependence $\Rdcov$ is probably a better test than HHG. In settings (V4) and (V5) which are based on marginals having a mixture distribution, HHG convincingly outperforms $\Rdcov$ (and all the other competing tests). In the rest of the settings, HHG and $\Rdcov$ have similar performance. \par 
	Next we look at two very natural simulation settings based on correlated multivariate Gaussian random variables. Note that multivariate Gaussians are a very popular modeling choice in many practical scenarios. In~\cite{Gabor2007}, the authors use correlated multivariate Gaussians to compare the empirical performance of DCoV with likelihood ratio based tests.

	\noindent (IG) $(Z_1,\ldots ,Z_6)\sim\mathcal{N}_6(\mathbf{0},\Sigma)$ where $\Sigma_{i,j}=\rho\in [-1,1]$ for $i\geq 4, j\leq 3$, $i\leq 3, j\geq 4$ ,$\Sigma_{i,i}=1$ and $\Sigma_{i,j}=0$ otherwise. Set $(X_1,X_2,X_3)=(Z_1,Z_2,Z_3)$ and $(Y_1,Y_2,Y_3)=(Z_4,Z_5,Z_6)$.\par 
	\noindent (IGL) $(W_1,\ldots ,W_6)=(\exp{(Z_1)},\ldots ,\exp{(Z_6)})$ where $(Z_1,\ldots ,Z_6)\sim\mathcal{N}_6(\mathbf{0},\Sigma)$ and $\Sigma_{i,j}=\rho\in [-1,1]$ for $i\geq 4, j\leq 3$, $i\leq 3, j\geq 4$ ,$\Sigma_{i,i}=1$ and $\Sigma_{i,j}=0$ otherwise. Set $(X_1,X_2,X_3)=(W_1,W_2,W_3)$ and $(Y_1,Y_2,Y_3)=(W_4,W_5,W_6)$.\par
\begin{figure}[h]
		\begin{center}
			\includegraphics[height=7.5cm,width=7.5cm]{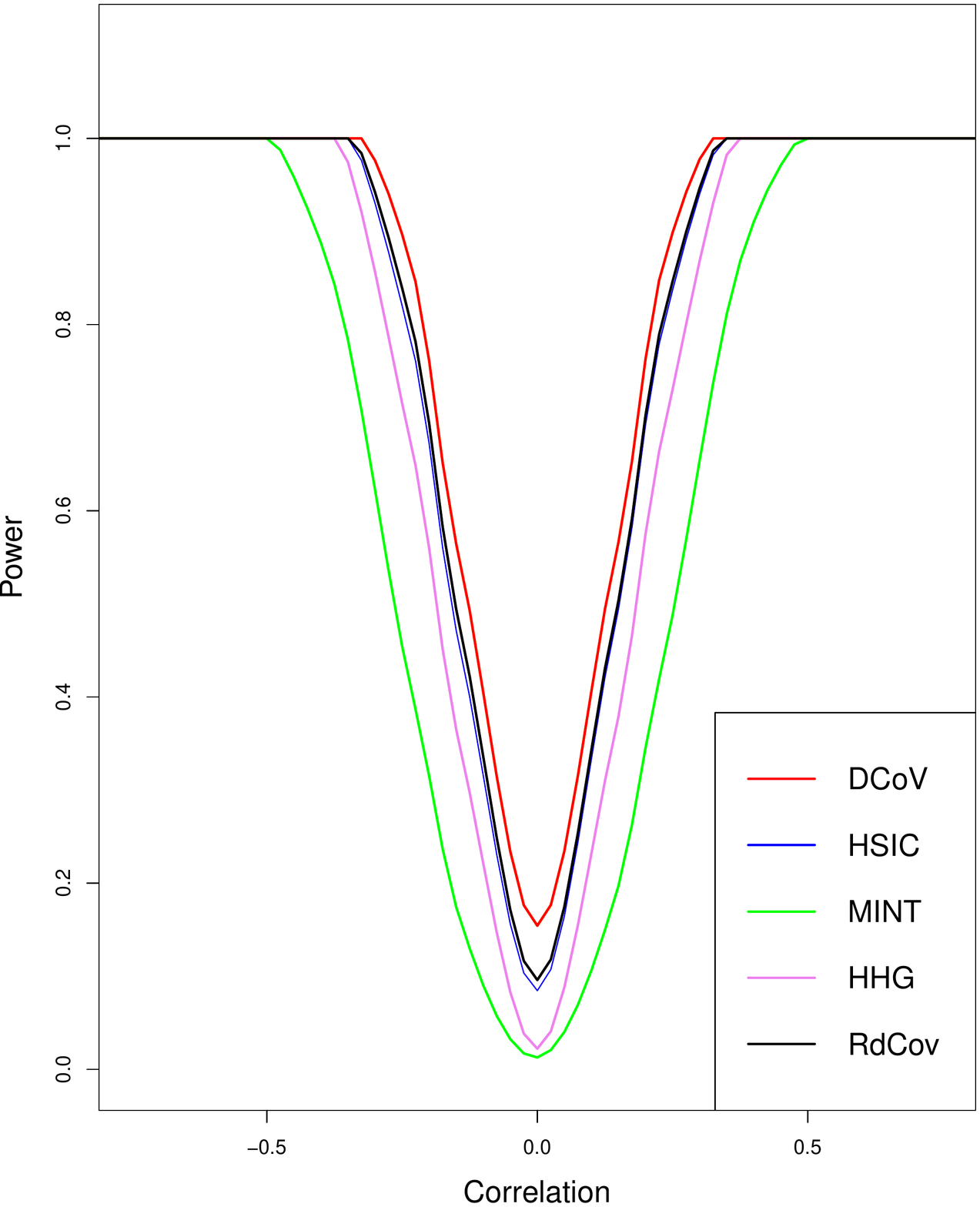}	
			\includegraphics[height=7.5cm,width=7.5cm]{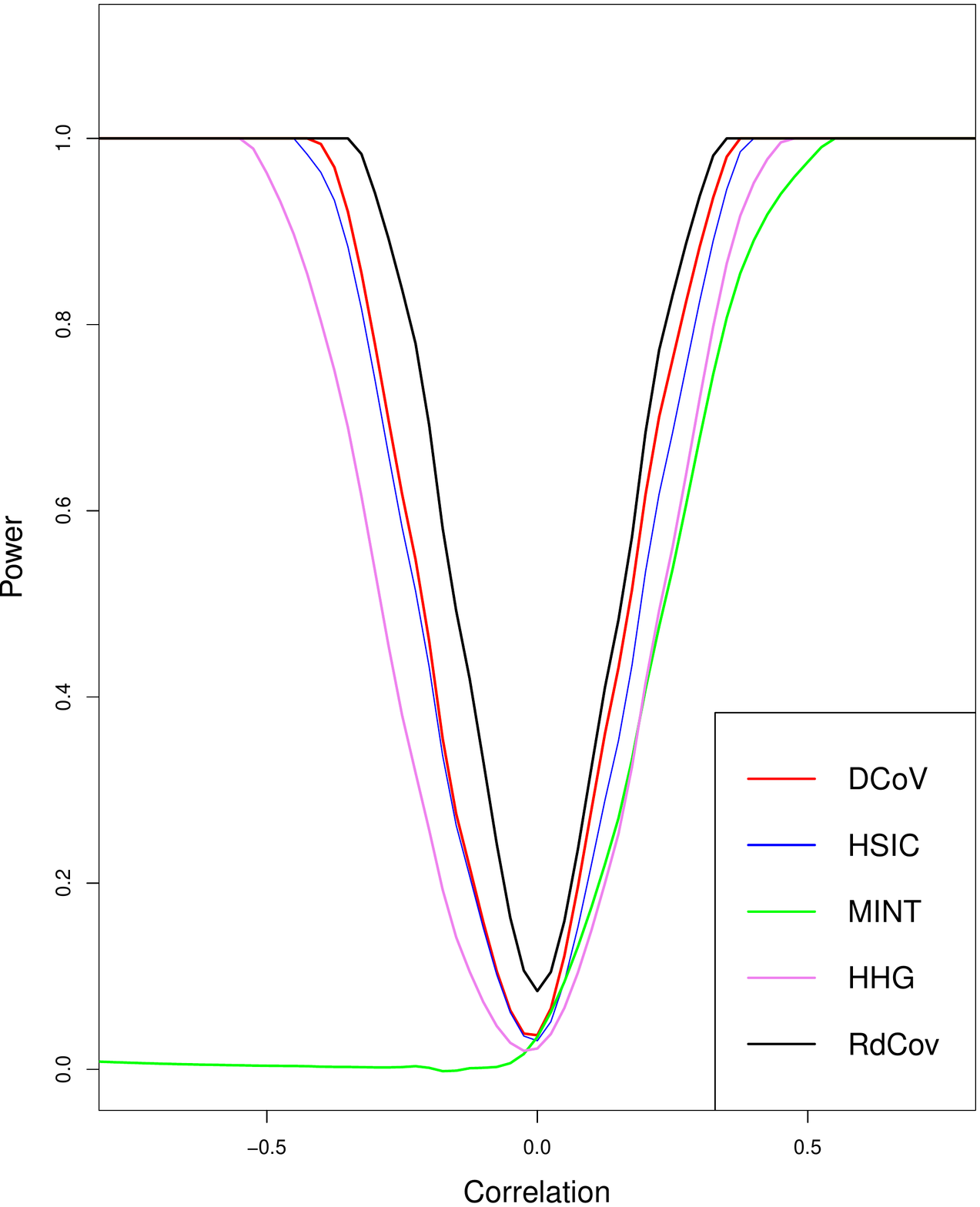}
		\end{center}	
		\caption{The left panel shows the power curves for (IG) at level $5\%$ for $\rho\in [-1,1]$. The right panel shows the same for (IGL). The above plots are obtained by estimating the powers of aforementioned tests on a grid (of size $20$ ranging from $-0.9$ to $0.9$) of possible correlation values, followed by polynomial smoothing using the \texttt{loess} function in \texttt{R}.}
		\label{fig:PowerIndcomp}
	\end{figure}

 
	In~\cref{fig:PowerIndcomp} we present plots corresponding to the power curves of different tests of independence as $\rho$ varies in $[-1,1]$.
	The left panel of~\cref{fig:PowerIndcomp} reveals that $\Rdcov$ significantly outperforms MINT and HHG for the setting (IG). This reinforces our belief that $\Rdcov$ is a better test under near-linear dependence relations. In fact, $\Rdcov$ also marginally outperforms HSIC. In setting (IG) DCoV has the best power curve, which is close to both $\Rdcov$ and HSIC. Next, for setting (IGL) $\Rdcov$ has the best performance as it convincingly outperforms all its competitors. In this setting, MINT does not appear to be consistent when $\rho<0$. This could perhaps be due to the fact that some regularity assumptions required for the consistency of MINT are not satisfied in this case. Note that (IGL) is a somewhat heavy-tailed setting, although the associated distribution has all moments finite (the exponential moment becomes infinite). The superior performance of $\Rdcov$ in this case further strengthens our belief that $\Rdcov$ can provide a better test of independence in heavy-tailed settings. For more simulations, please refer to~\cref{sec:addsim}.
	
	\section{Discussion}\label{sec:discuss}
	We have developed a framework for multivariate distribution-free
	nonparametric testing using the method of multivariate ranks defined using the theory of optimal transportation. We have illustrated our general approach through many problems: (I) testing for mutual
	independence of $K$ ($\ge 2$) random vectors, (II) goodness-of-fit
	testing for $K$ ($\ge 2$) multivariate distributions, (III) testing
	for multivariate symmetry of a random vector, etc. We show that our
	proposed tests are finite sample distribution-free, consistent against
	all alternatives (under minimal assumptions), and are computationally
	feasible. In fact, the proposed tests reduce to well-known
	one-dimensional tests for problems (I) and (II). We further derive the
	asymptotic weak limits of our test statistics, under the null
	hypotheses. In the process we also derive results on the asymptotic
	regularity of optimal transport maps (aka multivariate ranks) which is
	of independent interest. As far as we are aware, this is the first
	attempt to systematically develop distribution-free
	multivariate tests that are consistent against all alternatives and are computationally feasible.
	\par
	A natural future research direction is to theoretically investigate the power behavior of these
	proposed tests, in the sense of Pitman efficiency (see e.g.,~\cite{bbb2019}) or consistency against local alternatives (see e.g.,~\cite{gretton2012}). Further, to make our methods more scalable it would be interesting to explore approximate ``greedy" algorithms with lower computational complexity, that solve the assignment problem in~\eqref{eq:empopt}. Finally, we believe that our proposed general framework can be used to construct distribution-free tests in many other multivariate nonparametric testing problems beyond those discussed in this paper. We hope that more of such multivariate rank-based distribution-free tests will be studied in future.	
	\begin{spacing}{1.0}
		{\small{
				\bibliographystyle{chicago}
				\bibliography{template}}}

\begin{thebibliography}{}

\bibitem[\protect\citeauthoryear{Agarwal, Mukherjee, Bhattacharya, and
  Zhang}{Agarwal et~al.}{2019}]{bbbm2019}
Agarwal, D., S.~Mukherjee, B.~B. Bhattacharya, and N.~R. Zhang (2019).
\newblock Distribution-free multisample test based on optimal matching with
  applications to single cell genomics.
\newblock {\em arXiv preprint arXiv:1906.04776\/}.

\bibitem[\protect\citeauthoryear{Aki}{Aki}{1993}]{aki1993}
Aki, S. (1993).
\newblock On nonparametric tests for symmetry in {${\bf R}^m$}.
\newblock {\em Ann. Inst. Statist. Math.\/}~{\em 45\/}(4), 787--800.

\bibitem[\protect\citeauthoryear{Alexandroff}{Alexandroff}{1939}]{Alexandroff1939}
Alexandroff, A.~D. (1939).
\newblock Almost everywhere existence of the second differential of a convex
  function and some properties of convex surfaces connected with it.
\newblock {\em Leningrad State Univ. Annals [Uchenye Zapiski] Math.
  Ser.\/}~{\em 6}, 3--35.

\bibitem[\protect\citeauthoryear{Anderson}{Anderson}{1962}]{AndersonCVM1962}
Anderson, T.~W. (1962).
\newblock On the distribution of the two-sample {C}ram\'{e}r-von {M}ises
  criterion.
\newblock {\em Ann. Math. Statist.\/}~{\em 33}, 1148--1159.

\bibitem[\protect\citeauthoryear{Arias-Castro and Pelletier}{Arias-Castro and
  Pelletier}{2016}]{Castro2016}
Arias-Castro, E. and B.~Pelletier (2016).
\newblock On the consistency of the crossmatch test.
\newblock {\em J. Statist. Plann. Inference\/}~{\em 171}, 184--190.

\bibitem[\protect\citeauthoryear{Atkinson}{Atkinson}{1989}]{Atkinson1989}
Atkinson, K.~E. (1989).
\newblock {\em An introduction to numerical analysis\/} (Second ed.).
\newblock John Wiley \& Sons, Inc., New York.

\bibitem[\protect\citeauthoryear{Bakirov, Rizzo, and Sz\'{e}kely}{Bakirov
  et~al.}{2006}]{bakirov2006}
Bakirov, N.~K., M.~L. Rizzo, and G.~J. Sz\'{e}kely (2006).
\newblock A multivariate nonparametric test of independence.
\newblock {\em J. Multivariate Anal.\/}~{\em 97\/}(8), 1742--1756.

\bibitem[\protect\citeauthoryear{Baringhaus and Franz}{Baringhaus and
  Franz}{2004}]{Baringhaus2004}
Baringhaus, L. and C.~Franz (2004).
\newblock On a new multivariate two-sample test.
\newblock {\em J. Multivariate Anal.\/}~{\em 88\/}(1), 190--206.

\bibitem[\protect\citeauthoryear{Bene{\v{s}}, Lechnerov{\'a}, Klebanov,
  Sl{\'a}mov{\'a}, and Sl{\'a}ma}{Bene{\v{s}}
  et~al.}{2009}]{benevs2009statistical}
Bene{\v{s}}, V., R.~Lechnerov{\'a}, L.~Klebanov, M.~Sl{\'a}mov{\'a}, and
  P.~Sl{\'a}ma (2009).
\newblock Statistical comparison of the geometry of second-phase particles.
\newblock {\em Materials Characterization\/}~{\em 60\/}(10), 1076--1081.

\bibitem[\protect\citeauthoryear{Beran}{Beran}{1979}]{beran1979}
Beran, R. (1979).
\newblock Testing for ellipsoidal symmetry of a multivariate density.
\newblock {\em Ann. Statist.\/}~{\em 7\/}(1), 150--162.

\bibitem[\protect\citeauthoryear{Bergsma and Dassios}{Bergsma and
  Dassios}{2014}]{Bregsma2014}
Bergsma, W. and A.~Dassios (2014).
\newblock A consistent test of independence based on a sign covariance related
  to {K}endall's tau.
\newblock {\em Bernoulli\/}~{\em 20\/}(2), 1006--1028.

\bibitem[\protect\citeauthoryear{Berrett and Samworth}{Berrett and
  Samworth}{2019}]{berrett2017nonparametric}
Berrett, T.~B. and R.~J. Samworth (2019).
\newblock Nonparametric independence testing via mutual information.
\newblock {\em Biometrika\/}~{\em 106\/}(3), 547--566.

\bibitem[\protect\citeauthoryear{Bertsekas}{Bertsekas}{1988}]{bertsekas1988}
Bertsekas, D.~P. (1988).
\newblock The auction algorithm: a distributed relaxation method for the
  assignment problem.
\newblock {\em Ann. Oper. Res.\/}~{\em 14\/}(1-4), 105--123.

\bibitem[\protect\citeauthoryear{Bhattacharya}{Bhattacharya}{2019}]{bbb2019}
Bhattacharya, B.~B. (2019).
\newblock A general asymptotic framework for distribution-free graph-based
  two-sample tests.
\newblock {\em Journal of the Royal Statistical Society: Series B (Statistical
  Methodology)\/}~{\em 81\/}(3), 575--602.

\bibitem[\protect\citeauthoryear{Bickel}{Bickel}{1965}]{Bickel1965}
Bickel, P.~J. (1965).
\newblock On some asymptotically nonparametric competitors of {H}otelling's
  {$T^{2}$}.
\newblock {\em Ann. Math. Statist.\/}~{\em 36}, 160--173; correction, ibid.
  1583.

\bibitem[\protect\citeauthoryear{Bickel}{Bickel}{1968}]{bickel1969distribution}
Bickel, P.~J. (1968).
\newblock A distribution free version of the {S}mirnov two sample test in the
  {$p$}-variate case.
\newblock {\em Ann. Math. Statist.\/}~{\em 40}, 1--23.

\bibitem[\protect\citeauthoryear{Biswas, Mukhopadhyay, and Ghosh}{Biswas
  et~al.}{2014}]{munmun2014}
Biswas, M., M.~Mukhopadhyay, and A.~K. Ghosh (2014).
\newblock A distribution-free two-sample run test applicable to
  high-dimensional data.
\newblock {\em Biometrika\/}~{\em 101\/}(4), 913--926.

\bibitem[\protect\citeauthoryear{Biswas, Sarkar, and Ghosh}{Biswas
  et~al.}{2016}]{munmun2016}
Biswas, M., S.~Sarkar, and A.~K. Ghosh (2016).
\newblock On some exact distribution-free tests of independence between two
  random vectors of arbitrary dimensions.
\newblock {\em J. Statist. Plann. Inference\/}~{\em 175}, 78--86.

\bibitem[\protect\citeauthoryear{Blomqvist}{Blomqvist}{1950}]{blomqvist1950measure}
Blomqvist, N. (1950).
\newblock On a measure of dependence between two random variables.
\newblock {\em Ann. Math. Statistics\/}~{\em 21}, 593--600.

\bibitem[\protect\citeauthoryear{Blum, Kiefer, and Rosenblatt}{Blum
  et~al.}{1961}]{blum1961distribution}
Blum, J.~R., J.~Kiefer, and M.~Rosenblatt (1961).
\newblock Distribution free tests of independence based on the sample
  distribution function.
\newblock {\em Ann. Math. Statist.\/}~{\em 32}, 485--498.

\bibitem[\protect\citeauthoryear{Boeckel, Spokoiny, and Suvorikova}{Boeckel
  et~al.}{2018}]{boeckel2018multivariate}
Boeckel, M., V.~Spokoiny, and A.~Suvorikova (2018).
\newblock Multivariate brenier cumulative distribution functions and their
  application to non-parametric testing.
\newblock {\em arXiv preprint arXiv:1809.04090\/}.

\bibitem[\protect\citeauthoryear{Brenier}{Brenier}{1991}]{Brenier1991}
Brenier, Y. (1991).
\newblock Polar factorization and monotone rearrangement of vector-valued
  functions.
\newblock {\em Comm. Pure Appl. Math.\/}~{\em 44\/}(4), 375--417.

\bibitem[\protect\citeauthoryear{Caffarelli}{Caffarelli}{1990}]{caffarelli1990}
Caffarelli, L.~A. (1990).
\newblock Interior {$W^{2,p}$} estimates for solutions of the
  {M}onge-{A}mp\`ere equation.
\newblock {\em Ann. of Math. (2)\/}~{\em 131\/}(1), 135--150.

\bibitem[\protect\citeauthoryear{Chatterjee and Shao}{Chatterjee and
  Shao}{2011}]{Cha2011}
Chatterjee, S. and Q.-M. Shao (2011).
\newblock Nonnormal approximation by {S}tein's method of exchangeable pairs
  with application to the {C}urie-{W}eiss model.
\newblock {\em Ann. Appl. Probab.\/}~{\em 21\/}(2), 464--483.

\bibitem[\protect\citeauthoryear{Chaudhuri}{Chaudhuri}{1996}]{Chaudhuri1996}
Chaudhuri, P. (1996).
\newblock On a geometric notion of quantiles for multivariate data.
\newblock {\em J. Amer. Statist. Assoc.\/}~{\em 91\/}(434), 862--872.

\bibitem[\protect\citeauthoryear{Chen and Friedman}{Chen and
  Friedman}{2017}]{chen2017}
Chen, H. and J.~H. Friedman (2017).
\newblock A new graph-based two-sample test for multivariate and object data.
\newblock {\em J. Amer. Statist. Assoc.\/}~{\em 112\/}(517), 397--409.

\bibitem[\protect\citeauthoryear{Chen and Fang}{Chen and Fang}{2015}]{Chen2015}
Chen, L. H.~Y. and X.~Fang (2015).
\newblock On the error bound in a combinatorial central limit theorem.
\newblock {\em Bernoulli\/}~{\em 21\/}(1), 335--359.

\bibitem[\protect\citeauthoryear{Chernozhukov, Galichon, Hallin, Henry,
  et~al.}{Chernozhukov et~al.}{2017}]{chernozhukov2017monge}
Chernozhukov, V., A.~Galichon, M.~Hallin, M.~Henry, et~al. (2017).
\newblock Monge--kantorovich depth, quantiles, ranks and signs.
\newblock {\em The Annals of Statistics\/}~{\em 45\/}(1), 223--256.

\bibitem[\protect\citeauthoryear{Conradsen, Nielsen, Schou, and
  Skriver}{Conradsen et~al.}{2003}]{conradsen2003test}
Conradsen, K., A.~A. Nielsen, J.~Schou, and H.~Skriver (2003).
\newblock A test statistic in the complex wishart distribution and its
  application to change detection in polarimetric sar data.
\newblock {\em IEEE Transactions on Geoscience and Remote Sensing\/}~{\em
  41\/}(1), 4--19.

\bibitem[\protect\citeauthoryear{Cruz-Uribe and Neugebauer}{Cruz-Uribe and
  Neugebauer}{2002}]{Uribe2002}
Cruz-Uribe, D. and C.~J. Neugebauer (2002).
\newblock Sharp error bounds for the trapezoidal rule and {S}impson's rule.
\newblock {\em JIPAM. J. Inequal. Pure Appl. Math.\/}~{\em 3\/}(4), Article 49,
  22.

\bibitem[\protect\citeauthoryear{Cs\"{o}rg\H{o}}{Cs\"{o}rg\H{o}}{1985}]{csorgo1985}
Cs\"{o}rg\H{o}, S. (1985).
\newblock Testing for independence by the empirical characteristic function.
\newblock {\em J. Multivariate Anal.\/}~{\em 16\/}(3), 290--299.

\bibitem[\protect\citeauthoryear{De~Philippis and Figalli}{De~Philippis and
  Figalli}{2013}]{Figalli2013}
De~Philippis, G. and A.~Figalli (2013).
\newblock {$W^{2,1}$} regularity for solutions of the {M}onge-{A}mp\`ere
  equation.
\newblock {\em Invent. Math.\/}~{\em 192\/}(1), 55--69.

\bibitem[\protect\citeauthoryear{del Barrio, Cuesta-Albertos, Hallin, and
  Matr{\'a}n}{del Barrio et~al.}{2018}]{del2018center}
del Barrio, E., J.~A. Cuesta-Albertos, M.~Hallin, and C.~Matr{\'a}n (2018).
\newblock Center-outward distribution functions, quantiles, ranks, and signs in
  rd. arxiv e-prints, art.
\newblock {\em arXiv preprint arXiv:1806.01238\/}.

\bibitem[\protect\citeauthoryear{Dishion, Capaldi, and Yoerger}{Dishion
  et~al.}{1999}]{dishion1999middle}
Dishion, T.~J., D.~M. Capaldi, and K.~Yoerger (1999).
\newblock Middle childhood antecedents to progressions in male adolescent
  substance use: An ecological analysis of risk and protection.
\newblock {\em Journal of Adolescent Research\/}~{\em 14\/}(2), 175--205.

\bibitem[\protect\citeauthoryear{Drouet~Mari and Kotz}{Drouet~Mari and
  Kotz}{2001}]{Mari2001}
Drouet~Mari, D. and S.~Kotz (2001).
\newblock {\em Correlation and dependence}.
\newblock Imperial College Press, London; distributed by World Scientific
  Publishing Co., Inc., River Edge, NJ.

\bibitem[\protect\citeauthoryear{Drton, Han, and Shi}{Drton
  et~al.}{2018}]{drton2018high}
Drton, M., F.~Han, and H.~Shi (2018).
\newblock High dimensional independence testing with maxima of rank
  correlations.
\newblock {\em arXiv preprint arXiv:1812.06189\/}.

\bibitem[\protect\citeauthoryear{Dudley}{Dudley}{1978}]{Dudley1978}
Dudley, R.~M. (1978).
\newblock Central limit theorems for empirical measures.
\newblock {\em Ann. Probab.\/}~{\em 6\/}(6), 899--929 (1979).

\bibitem[\protect\citeauthoryear{Farris and Schopflocher}{Farris and
  Schopflocher}{1999}]{farris1999between}
Farris, K.~B. and D.~P. Schopflocher (1999).
\newblock Between intention and behavior: an application of community
  pharmacists' assessment of pharmaceutical care.
\newblock {\em Social science \& medicine\/}~{\em 49\/}(1), 55--66.

\bibitem[\protect\citeauthoryear{Feuerverger}{Feuerverger}{1993}]{feuerverger1993consistent}
Feuerverger, A. (1993).
\newblock A consistent test for bivariate dependence.
\newblock {\em International Statistical Review/Revue Internationale de
  Statistique\/}, 419--433.

\bibitem[\protect\citeauthoryear{Folkes, Koletsky, and Graham}{Folkes
  et~al.}{1987}]{folkes1987field}
Folkes, V.~S., S.~Koletsky, and J.~L. Graham (1987).
\newblock A field study of causal inferences and consumer reaction: the view
  from the airport.
\newblock {\em Journal of consumer research\/}~{\em 13\/}(4), 534--539.

\bibitem[\protect\citeauthoryear{Fox and Weisberg}{Fox and
  Weisberg}{2019}]{Fox2019}
Fox, J. and S.~Weisberg (2019).
\newblock {\em An {R} Companion to Applied Regression\/} (Third ed.).
\newblock Thousand Oaks {CA}: Sage.

\bibitem[\protect\citeauthoryear{Friedman and Rafsky}{Friedman and
  Rafsky}{1979}]{friedman1979}
Friedman, J.~H. and L.~C. Rafsky (1979).
\newblock Multivariate generalizations of the {W}ald-{W}olfowitz and {S}mirnov
  two-sample tests.
\newblock {\em Ann. Statist.\/}~{\em 7\/}(4), 697--717.

\bibitem[\protect\citeauthoryear{Friedman and Rafsky}{Friedman and
  Rafsky}{1983}]{friedman1983}
Friedman, J.~H. and L.~C. Rafsky (1983).
\newblock Graph-theoretic measures of multivariate association and prediction.
\newblock {\em Ann. Statist.\/}~{\em 11\/}(2), 377--391.

\bibitem[\protect\citeauthoryear{Ghosal and Sen}{Ghosal and
  Sen}{2019}]{ghosal2019multivariate}
Ghosal, P. and B.~Sen (2019).
\newblock Multivariate ranks and quantiles using optimal transportation and
  applications to goodness-of-fit testing.
\newblock {\em arXiv preprint arXiv:1905.05340\/}.

\bibitem[\protect\citeauthoryear{Gibbons and Chakraborti}{Gibbons and
  Chakraborti}{2011}]{gibbons2011nonparametric}
Gibbons, J.~D. and S.~Chakraborti (2011).
\newblock {\em Nonparametric statistical inference\/} (Fifth ed.).
\newblock Statistics: Textbooks and Monographs. CRC Press, Boca Raton, FL.

\bibitem[\protect\citeauthoryear{Gieser and Randles}{Gieser and
  Randles}{1997}]{gieser1997}
Gieser, P.~W. and R.~H. Randles (1997).
\newblock A nonparametric test of independence between two vectors.
\newblock {\em J. Amer. Statist. Assoc.\/}~{\em 92\/}(438), 561--567.

\bibitem[\protect\citeauthoryear{Gretton, Borgwardt, Rasch, Sch\"{o}lkopf, and
  Smola}{Gretton et~al.}{2012}]{gretton2012}
Gretton, A., K.~M. Borgwardt, M.~J. Rasch, B.~Sch\"{o}lkopf, and A.~Smola
  (2012).
\newblock A kernel two-sample test.
\newblock {\em J. Mach. Learn. Res.\/}~{\em 13}, 723--773.

\bibitem[\protect\citeauthoryear{Gretton, Bousquet, Smola, and
  Sch\"{o}lkopf}{Gretton et~al.}{2005}]{GrettonHSIC2005}
Gretton, A., O.~Bousquet, A.~Smola, and B.~Sch\"{o}lkopf (2005).
\newblock Measuring statistical dependence with {H}ilbert-{S}chmidt norms.
\newblock In {\em Algorithmic learning theory}, Volume 3734 of {\em Lecture
  Notes in Comput. Sci.}, pp.\  63--77. Springer, Berlin.

\bibitem[\protect\citeauthoryear{Gretton, Fukumizu, Harchaoui, and
  Sriperumbudur}{Gretton et~al.}{2009}]{gretton2009fast}
Gretton, A., K.~Fukumizu, Z.~Harchaoui, and B.~K. Sriperumbudur (2009).
\newblock A fast, consistent kernel two-sample test.
\newblock In {\em Advances in neural information processing systems}, pp.\
  673--681.

\bibitem[\protect\citeauthoryear{Gretton, Fukumizu, Teo, Song, Sch{\"o}lkopf,
  and Smola}{Gretton et~al.}{2008}]{gretton2008kernel}
Gretton, A., K.~Fukumizu, C.~H. Teo, L.~Song, B.~Sch{\"o}lkopf, and A.~J. Smola
  (2008).
\newblock A kernel statistical test of independence.
\newblock In {\em Advances in neural information processing systems}, pp.\
  585--592.

\bibitem[\protect\citeauthoryear{Gretton and Gy\"{o}rfi}{Gretton and
  Gy\"{o}rfi}{2008}]{gretton2008}
Gretton, A. and L.~Gy\"{o}rfi (2008).
\newblock Nonparametric independence tests: space partitioning and kernel
  approaches.
\newblock In {\em Algorithmic learning theory}, Volume 5254 of {\em Lecture
  Notes in Comput. Sci.}, pp.\  183--198. Springer, Berlin.

\bibitem[\protect\citeauthoryear{Gretton, Herbrich, Smola, Bousquet, and
  Sch\"{o}lkopf}{Gretton et~al.}{2005}]{Gretton2005}
Gretton, A., R.~Herbrich, A.~Smola, O.~Bousquet, and B.~Sch\"{o}lkopf (2005).
\newblock Kernel methods for measuring independence.
\newblock {\em J. Mach. Learn. Res.\/}~{\em 6}, 2075--2129.

\bibitem[\protect\citeauthoryear{Grover and Dillon}{Grover and
  Dillon}{1985}]{grover1985probabilistic}
Grover, R. and W.~R. Dillon (1985).
\newblock A probabilistic model for testing hypothesized hierarchical market
  structures.
\newblock {\em Marketing Science\/}~{\em 4\/}(4), 312--335.

\bibitem[\protect\citeauthoryear{Hallin and Paindaveine}{Hallin and
  Paindaveine}{2004}]{Hallin2004}
Hallin, M. and D.~Paindaveine (2004).
\newblock Rank-based optimal tests of the adequacy of an elliptic {VARMA}
  model.
\newblock {\em Ann. Statist.\/}~{\em 32\/}(6), 2642--2678.

\bibitem[\protect\citeauthoryear{Hallin and Paindaveine}{Hallin and
  Paindaveine}{2006}]{Davy2006}
Hallin, M. and D.~Paindaveine (2006).
\newblock Parametric and semiparametric inference for shape: the role of the
  scale functional.
\newblock {\em Statist. Decisions\/}~{\em 24\/}(3), 327--350.

\bibitem[\protect\citeauthoryear{Halton}{Halton}{1964}]{halton1964algorithm}
Halton, J.~H. (1964).
\newblock Algorithm 247: Radical-inverse quasi-random point sequence.
\newblock {\em Communications of the ACM\/}~{\em 7\/}(12), 701--702.

\bibitem[\protect\citeauthoryear{Hardy, Littlewood, and P\'{o}lya}{Hardy
  et~al.}{1952}]{Hardy1952}
Hardy, G.~H., J.~E. Littlewood, and G.~P\'{o}lya (1952).
\newblock {\em Inequalities}.
\newblock Cambridge, at the University Press.
\newblock 2d ed.

\bibitem[\protect\citeauthoryear{Hastings}{Hastings}{1970}]{hastings1970monte}
Hastings, W.~K. (1970).
\newblock Monte carlo sampling methods using markov chains and their
  applications.

\bibitem[\protect\citeauthoryear{Hauke and Kossowski}{Hauke and
  Kossowski}{2011}]{hauke2011comparison}
Hauke, J. and T.~Kossowski (2011).
\newblock Comparison of values of pearson's and spearman's correlation
  coefficients on the same sets of data.
\newblock {\em Quaestiones geographicae\/}~{\em 30\/}(2), 87--93.

\bibitem[\protect\citeauthoryear{Heathcote, Rachev, and Cheng}{Heathcote
  et~al.}{1995}]{heathcote1995}
Heathcote, C.~R., S.~T. Rachev, and B.~Cheng (1995).
\newblock Testing multivariate symmetry.
\newblock {\em J. Multivariate Anal.\/}~{\em 54\/}(1), 91--112.

\bibitem[\protect\citeauthoryear{Heller, Gorfine, and Heller}{Heller
  et~al.}{2012}]{heller2012}
Heller, R., M.~Gorfine, and Y.~Heller (2012).
\newblock A class of multivariate distribution-free tests of independence based
  on graphs.
\newblock {\em J. Statist. Plann. Inference\/}~{\em 142\/}(12), 3097--3106.

\bibitem[\protect\citeauthoryear{Heller and Heller}{Heller and
  Heller}{2016}]{heller2016multivariate}
Heller, R. and Y.~Heller (2016).
\newblock Multivariate tests of association based on univariate tests.
\newblock In {\em Advances in Neural Information Processing Systems}, pp.\
  208--216.

\bibitem[\protect\citeauthoryear{Heller, Heller, and Gorfine}{Heller
  et~al.}{2013}]{heller2013}
Heller, R., Y.~Heller, and M.~Gorfine (2013).
\newblock A consistent multivariate test of association based on ranks of
  distances.
\newblock {\em Biometrika\/}~{\em 100\/}(2), 503--510.

\bibitem[\protect\citeauthoryear{Heller, Jensen, Rosenbaum, and Small}{Heller
  et~al.}{2010}]{heller2010}
Heller, R., S.~T. Jensen, P.~R. Rosenbaum, and D.~S. Small (2010).
\newblock Sensitivity analysis for the cross-match test, with applications in
  genomics.
\newblock {\em J. Amer. Statist. Assoc.\/}~{\em 105\/}(491), 1005--1013.

\bibitem[\protect\citeauthoryear{Heller, Rosenbaum, and Small}{Heller
  et~al.}{2010}]{hellerrosenbaum2010using}
Heller, R., P.~R. Rosenbaum, and D.~S. Small (2010).
\newblock Using the cross-match test to appraise covariate balance in matched
  pairs.
\newblock {\em The American Statistician\/}~{\em 64\/}(4), 299--309.

\bibitem[\protect\citeauthoryear{Henze}{Henze}{1988}]{henze1988}
Henze, N. (1988).
\newblock A multivariate two-sample test based on the number of nearest
  neighbor type coincidences.
\newblock {\em Ann. Statist.\/}~{\em 16\/}(2), 772--783.

\bibitem[\protect\citeauthoryear{Hettmansperger, M\"{o}tt\"{o}nen, and
  Oja}{Hettmansperger et~al.}{1998}]{hettmansperger1998}
Hettmansperger, T.~P., J.~M\"{o}tt\"{o}nen, and H.~Oja (1998).
\newblock Affine invariant multivariate rank tests for several samples.
\newblock {\em Statist. Sinica\/}~{\em 8\/}(3), 785--800.

\bibitem[\protect\citeauthoryear{Hlawka}{Hlawka}{1961}]{Hlawka1961}
Hlawka, E. (1961).
\newblock Funktionen von beschr\"{a}nkter {V}ariation in der {T}heorie der
  {G}leichverteilung.
\newblock {\em Ann. Mat. Pura Appl. (4)\/}~{\em 54}, 325--333.

\bibitem[\protect\citeauthoryear{Hoeffding}{Hoeffding}{1948}]{hoeffding1948non}
Hoeffding, W. (1948).
\newblock A non-parametric test of independence.
\newblock {\em Ann. Math. Statistics\/}~{\em 19}, 546--557.

\bibitem[\protect\citeauthoryear{Hofer}{Hofer}{2009}]{Hofer2009}
Hofer, R. (2009).
\newblock On the distribution properties of {N}iederreiter-{H}alton sequences.
\newblock {\em J. Number Theory\/}~{\em 129\/}(2), 451--463.

\bibitem[\protect\citeauthoryear{Hofer and Larcher}{Hofer and
  Larcher}{2010}]{Hofer2010}
Hofer, R. and G.~Larcher (2010).
\newblock On existence and discrepancy of certain digital
  {N}iederreiter-{H}alton sequences.
\newblock {\em Acta Arith.\/}~{\em 141\/}(4), 369--394.

\bibitem[\protect\citeauthoryear{Hoffmann-J\o~rgensen}{Hoffmann-J\o~rgensen}{1991}]{Hoffman1991}
Hoffmann-J\o~rgensen, J. (1991).
\newblock {\em Stochastic processes on {P}olish spaces}, Volume~39 of {\em
  Various Publications Series (Aarhus)}.
\newblock Aarhus Universitet, Matematisk Institut, Aarhus.

\bibitem[\protect\citeauthoryear{Hollander, Wolfe, and Chicken}{Hollander
  et~al.}{2014}]{hollander2013nonparametric}
Hollander, M., D.~A. Wolfe, and E.~Chicken (2014).
\newblock {\em Nonparametric statistical methods\/} (Third ed.).
\newblock Wiley Series in Probability and Statistics. John Wiley \& Sons, Inc.,
  Hoboken, NJ.

\bibitem[\protect\citeauthoryear{Huo and Sz\'{e}kely}{Huo and
  Sz\'{e}kely}{2016}]{Huo2016}
Huo, X. and G.~J. Sz\'{e}kely (2016).
\newblock Fast computing for distance covariance.
\newblock {\em Technometrics\/}~{\em 58\/}(4), 435--447.

\bibitem[\protect\citeauthoryear{Iman and Conover}{Iman and
  Conover}{1982}]{iman1982distribution}
Iman, R.~L. and W.-J. Conover (1982).
\newblock A distribution-free approach to inducing rank correlation among input
  variables.
\newblock {\em Communications in Statistics-Simulation and Computation\/}~{\em
  11\/}(3), 311--334.

\bibitem[\protect\citeauthoryear{Jonker and Volgenant}{Jonker and
  Volgenant}{1987}]{jonker1987}
Jonker, R. and A.~Volgenant (1987).
\newblock A shortest augmenting path algorithm for dense and sparse linear
  assignment problems.
\newblock {\em Computing\/}~{\em 38\/}(4), 325--340.

\bibitem[\protect\citeauthoryear{Josse and Holmes}{Josse and
  Holmes}{2016}]{Josse2016}
Josse, J. and S.~Holmes (2016).
\newblock Measuring multivariate association and beyond.
\newblock {\em Stat. Surv.\/}~{\em 10}, 132--167.

\bibitem[\protect\citeauthoryear{Kankainen}{Kankainen}{1995}]{Kankainen1995}
Kankainen, A.-L. (1995).
\newblock {\em Consistent testing of total independence based on the empirical
  characteristic function}.
\newblock ProQuest LLC, Ann Arbor, MI.
\newblock Thesis (D.Phil.)--Jyvaskylan Yliopisto (Finland).

\bibitem[\protect\citeauthoryear{Karatzoglou, Smola, Hornik, and
  Zeileis}{Karatzoglou et~al.}{2004}]{kernlab}
Karatzoglou, A., A.~Smola, K.~Hornik, and A.~Zeileis (2004).
\newblock kernlab -- an {S4} package for kernel methods in {R}.
\newblock {\em Journal of Statistical Software\/}~{\em 11\/}(9), 1--20.

\bibitem[\protect\citeauthoryear{Kaufman, based in part on an earlier
  implementation~by Ruth~Heller, and Heller.}{Kaufman et~al.}{2019}]{HHGtwosam}
Kaufman, B. B. .~S., based in part on an earlier implementation~by Ruth~Heller,
  and Y.~Heller. (2019).
\newblock {\em HHG: Heller-Heller-Gorfine Tests of Independence and Equality of
  Distributions}.
\newblock R package version 2.3.2.

\bibitem[\protect\citeauthoryear{Kendall and Gibbons}{Kendall and
  Gibbons}{1990}]{Kendall1990}
Kendall, M. and J.~D. Gibbons (1990).
\newblock {\em Rank correlation methods\/} (Fifth ed.).
\newblock A Charles Griffin Title. Edward Arnold, London.

\bibitem[\protect\citeauthoryear{Kendall}{Kendall}{1938}]{kendall1938new}
Kendall, M.~G. (1938).
\newblock A new measure of rank correlation.
\newblock {\em Biometrika\/}~{\em 30\/}(1/2), 81--93.

\bibitem[\protect\citeauthoryear{Klebanov}{Klebanov}{2002}]{klebanov2002class}
Klebanov, L.~B. (2002).
\newblock A class of probability metrics and its statistical applications.
\newblock In {\em Statistical Data Analysis Based on the L1-Norm and Related
  Methods}, pp.\  241--252. Springer.

\bibitem[\protect\citeauthoryear{Kojadinovic and Holmes}{Kojadinovic and
  Holmes}{2009}]{kojadinovic2009}
Kojadinovic, I. and M.~Holmes (2009).
\newblock Tests of independence among continuous random vectors based on
  {C}ram\'{e}r-von {M}ises functionals of the empirical copula process.
\newblock {\em J. Multivariate Anal.\/}~{\em 100\/}(6), 1137--1154.

\bibitem[\protect\citeauthoryear{Kruskal}{Kruskal}{1952}]{kruskal1952}
Kruskal, W.~H. (1952).
\newblock A nonparametric test for the several sample problem.
\newblock {\em Ann. Math. Statistics\/}~{\em 23}, 525--540.

\bibitem[\protect\citeauthoryear{Kuipers and Niederreiter}{Kuipers and
  Niederreiter}{1974}]{Kuipers1974}
Kuipers, L. and H.~Niederreiter (1974).
\newblock {\em Uniform distribution of sequences}.
\newblock Wiley-Interscience [John Wiley \& Sons], New York-London-Sydney.
\newblock Pure and Applied Mathematics.

\bibitem[\protect\citeauthoryear{Kuo}{Kuo}{1975}]{Kuo1975}
Kuo, H.~H. (1975).
\newblock {\em Gaussian measures in {B}anach spaces}.
\newblock Lecture Notes in Mathematics, Vol. 463. Springer-Verlag, Berlin-New
  York.

\bibitem[\protect\citeauthoryear{Liu, Mcrae, Nyholt, Medland, Wray, Brown,
  Hayward, Montgomery, Visscher, Martin, et~al.}{Liu
  et~al.}{2010}]{liu2010versatile}
Liu, J.~Z., A.~F. Mcrae, D.~R. Nyholt, S.~E. Medland, N.~R. Wray, K.~M. Brown,
  N.~K. Hayward, G.~W. Montgomery, P.~M. Visscher, N.~G. Martin, et~al. (2010).
\newblock A versatile gene-based test for genome-wide association studies.
\newblock {\em The American Journal of Human Genetics\/}~{\em 87\/}(1),
  139--145.

\bibitem[\protect\citeauthoryear{Liu and Singh}{Liu and Singh}{1993}]{liu1993}
Liu, R.~Y. and K.~Singh (1993).
\newblock A quality index based on data depth and multivariate rank tests.
\newblock {\em J. Amer. Statist. Assoc.\/}~{\em 88\/}(421), 252--260.

\bibitem[\protect\citeauthoryear{Lu, Greevy, Xu, and Beck}{Lu
  et~al.}{2011}]{bo2011}
Lu, B., R.~Greevy, X.~Xu, and C.~Beck (2011).
\newblock Optimal nonbipartite matching and its statistical applications.
\newblock {\em Amer. Statist.\/}~{\em 65\/}(1), 21--30.

\bibitem[\protect\citeauthoryear{Lu, Lee, and Chiu}{Lu
  et~al.}{2009}]{lu2009financial}
Lu, C.-J., T.-S. Lee, and C.-C. Chiu (2009).
\newblock Financial time series forecasting using independent component
  analysis and support vector regression.
\newblock {\em Decision Support Systems\/}~{\em 47\/}(2), 115--125.

\bibitem[\protect\citeauthoryear{Lyons}{Lyons}{2013}]{Russell2013}
Lyons, R. (2013).
\newblock Distance covariance in metric spaces.
\newblock {\em Ann. Probab.\/}~{\em 41\/}(5), 3284--3305.

\bibitem[\protect\citeauthoryear{Mann and Whitney}{Mann and
  Whitney}{1947}]{mann1947}
Mann, H.~B. and D.~R. Whitney (1947).
\newblock On a test of whether one of two random variables is stochastically
  larger than the other.
\newblock {\em Ann. Math. Statistics\/}~{\em 18}, 50--60.

\bibitem[\protect\citeauthoryear{Marden}{Marden}{1999}]{Marden1999}
Marden, J.~I. (1999).
\newblock Multivariate rank tests.
\newblock In {\em Multivariate analysis, design of experiments, and survey
  sampling}, Volume 159 of {\em Statist. Textbooks Monogr.}, pp.\  401--432.
  Dekker, New York.

\bibitem[\protect\citeauthoryear{Martin and Betensky}{Martin and
  Betensky}{2005}]{martin2005testing}
Martin, E.~C. and R.~A. Betensky (2005).
\newblock Testing quasi-independence of failure and truncation times via
  conditional kendall's tau.
\newblock {\em Journal of the American Statistical Association\/}~{\em
  100\/}(470), 484--492.

\bibitem[\protect\citeauthoryear{Matteson and Tsay}{Matteson and
  Tsay}{2017}]{matteson2014}
Matteson, D.~S. and R.~S. Tsay (2017).
\newblock Independent component analysis via distance covariance.
\newblock {\em J. Amer. Statist. Assoc.\/}~{\em 112\/}(518), 623--637.

\bibitem[\protect\citeauthoryear{Mayer}{Mayer}{1975}]{mayer1975selecting}
Mayer, T. (1975).
\newblock Selecting economic hypotheses by goodness of fit.
\newblock {\em The Economic Journal\/}~{\em 85\/}(340), 877--883.

\bibitem[\protect\citeauthoryear{McCann}{McCann}{1995}]{Mccann1995}
McCann, R.~J. (1995).
\newblock Existence and uniqueness of monotone measure-preserving maps.
\newblock {\em Duke Math. J.\/}~{\em 80\/}(2), 309--323.

\bibitem[\protect\citeauthoryear{Monge}{Monge}{1781}]{monge1781memoire}
Monge, G. (1781).
\newblock M{\'e}moire sur la th{\'e}orie des d{\'e}blais et des remblais.
\newblock {\em M{\'e}moires Acad. Royale Sci. 1781\/}, 666--704.

\bibitem[\protect\citeauthoryear{M\'{o}ri and Sz\'{e}kely}{M\'{o}ri and
  Sz\'{e}kely}{2019}]{Mori2019}
M\'{o}ri, T.~F. and G.~J. Sz\'{e}kely (2019).
\newblock Four simple axioms of dependence measures.
\newblock {\em Metrika\/}~{\em 82\/}(1), 1--16.

\bibitem[\protect\citeauthoryear{Morokoff and Caflisch}{Morokoff and
  Caflisch}{1995}]{Morokoff1995}
Morokoff, W.~J. and R.~E. Caflisch (1995).
\newblock Quasi-{M}onte {C}arlo integration.
\newblock {\em J. Comput. Phys.\/}~{\em 122\/}(2), 218--230.

\bibitem[\protect\citeauthoryear{Mosteller}{Mosteller}{1946}]{Mosteller1946}
Mosteller, F. (1946).
\newblock On some useful ``inefficient'' statistics.
\newblock {\em Ann. Math. Statistics\/}~{\em 17}, 377--408.

\bibitem[\protect\citeauthoryear{M\"{o}tt\"{o}nen and Oja}{M\"{o}tt\"{o}nen and
  Oja}{1995}]{motto1995}
M\"{o}tt\"{o}nen, J. and H.~Oja (1995).
\newblock Multivariate spatial sign and rank methods.
\newblock {\em J. Nonparametr. Statist.\/}~{\em 5\/}(2), 201--213.

\bibitem[\protect\citeauthoryear{Mukaka}{Mukaka}{2012}]{mukaka2012guide}
Mukaka, M.~M. (2012).
\newblock A guide to appropriate use of correlation coefficient in medical
  research.
\newblock {\em Malawi Medical Journal\/}~{\em 24\/}(3), 69--71.

\bibitem[\protect\citeauthoryear{Munkres}{Munkres}{1957}]{munkres1957}
Munkres, J. (1957).
\newblock Algorithms for the assignment and transportation problems.
\newblock {\em J. Soc. Indust. Appl. Math.\/}~{\em 5}, 32--38.

\bibitem[\protect\citeauthoryear{Newman, Hettich, Blake, and Merz}{Newman
  et~al.}{1998}]{Sonardata}
Newman, D., S.~Hettich, C.~Blake, and C.~Merz (1998).
\newblock Uci repository of machine learning databases.

\bibitem[\protect\citeauthoryear{Newton}{Newton}{2009}]{Newton2009}
Newton, M.~A. (2009).
\newblock Introducing the discussion paper by {S}z\'{e}kely and {R}izzo
  [mr2752127].
\newblock {\em Ann. Appl. Stat.\/}~{\em 3\/}(4), 1233--1235.

\bibitem[\protect\citeauthoryear{Niederreiter}{Niederreiter}{1992}]{Niederreiter1992}
Niederreiter, H. (1992).
\newblock {\em Random number generation and quasi-{M}onte {C}arlo methods},
  Volume~63 of {\em CBMS-NSF Regional Conference Series in Applied
  Mathematics}.
\newblock Society for Industrial and Applied Mathematics (SIAM), Philadelphia,
  PA.

\bibitem[\protect\citeauthoryear{Oja}{Oja}{2010}]{oja2010}
Oja, H. (2010).
\newblock {\em Multivariate nonparametric methods with {R}}, Volume 199 of {\em
  Lecture Notes in Statistics}.
\newblock Springer, New York.
\newblock An approach based on spatial signs and ranks.

\bibitem[\protect\citeauthoryear{Oja and Randles}{Oja and
  Randles}{2004}]{Oja2004}
Oja, H. and R.~H. Randles (2004).
\newblock Multivariate nonparametric tests.
\newblock {\em Statist. Sci.\/}~{\em 19\/}(4), 598--605.

\bibitem[\protect\citeauthoryear{Pearson}{Pearson}{1920}]{pearson1920notes}
Pearson, K. (1920).
\newblock Notes on the history of correlation.
\newblock {\em Biometrika\/}~{\em 13\/}(1), 25--45.

\bibitem[\protect\citeauthoryear{Petrie}{Petrie}{2016}]{petrie2016}
Petrie, A. (2016).
\newblock Graph-theoretic multisample tests of equality in distribution for
  high dimensional data.
\newblock {\em Comput. Statist. Data Anal.\/}~{\em 96}, 145--158.

\bibitem[\protect\citeauthoryear{Pillai and Jayachandran}{Pillai and
  Jayachandran}{1967}]{Pillai1967}
Pillai, K. C.~S. and K.~Jayachandran (1967).
\newblock Power comparisons of tests of two multivariate hypotheses based on
  four criteria.
\newblock {\em Biometrika\/}~{\em 54}, 195--210.

\bibitem[\protect\citeauthoryear{Puri and Sen}{Puri and Sen}{1966}]{Puri1965}
Puri, M.~L. and P.~K. Sen (1966).
\newblock On a class of multivariate multisample rank-order tests.
\newblock {\em Sankhy\={a} Ser. A\/}~{\em 28}, 353--376.

\bibitem[\protect\citeauthoryear{Puri and Sen}{Puri and
  Sen}{1971}]{puri1971nonparametric}
Puri, M.~L. and P.~K. Sen (1971).
\newblock {\em Nonparametric methods in multivariate analysis}.
\newblock John Wiley\thinspace \&\thinspace Sons, Inc., New York-London-Sydney.

\bibitem[\protect\citeauthoryear{Quessy}{Quessy}{2010}]{Quessy2010}
Quessy, J.-F. (2010).
\newblock Applications and asymptotic power of marginal-free tests of
  stochastic vectorial independence.
\newblock {\em J. Statist. Plann. Inference\/}~{\em 140\/}(11), 3058--3075.

\bibitem[\protect\citeauthoryear{{R Core Team}}{{R Core Team}}{2019}]{Rsoft}
{R Core Team} (2019).
\newblock {\em R: A Language and Environment for Statistical Computing}.
\newblock Vienna, Austria: R Foundation for Statistical Computing.

\bibitem[\protect\citeauthoryear{Rabinowitz}{Rabinowitz}{1987}]{Philip1986}
Rabinowitz, P. (1987).
\newblock The convergence of noninterpolatory product integration rules.
\newblock In {\em Numerical integration ({H}alifax, {N}.{S}., 1986)}, Volume
  203 of {\em NATO Adv. Sci. Inst. Ser. C Math. Phys. Sci.}, pp.\  1--16.
  Reidel, Dordrecht.

\bibitem[\protect\citeauthoryear{Ramdas, Reddi, P{\'o}czos, Singh, and
  Wasserman}{Ramdas et~al.}{2015}]{ramdas2015decreasing}
Ramdas, A., S.~J. Reddi, B.~P{\'o}czos, A.~Singh, and L.~Wasserman (2015).
\newblock On the decreasing power of kernel and distance based nonparametric
  hypothesis tests in high dimensions.
\newblock In {\em Twenty-Ninth AAAI Conference on Artificial Intelligence}.

\bibitem[\protect\citeauthoryear{Randles and Peters}{Randles and
  Peters}{1990}]{randles1990}
Randles, R.~H. and D.~Peters (1990).
\newblock Multivariate rank tests for the two-sample location problem.
\newblock {\em Comm. Statist. Theory Methods\/}~{\em 19\/}(11), 4225--4238
  (1991).

\bibitem[\protect\citeauthoryear{R\'{e}millard}{R\'{e}millard}{2009}]{remillard2009}
R\'{e}millard, B. (2009).
\newblock Discussion of: {B}rownian distance covariance [mr2752127].
\newblock {\em Ann. Appl. Stat.\/}~{\em 3\/}(4), 1295--1298.

\bibitem[\protect\citeauthoryear{R\'{e}nyi}{R\'{e}nyi}{1959}]{Renyi1959}
R\'{e}nyi, A. (1959).
\newblock On measures of dependence.
\newblock {\em Acta Math. Acad. Sci. Hungar.\/}~{\em 10}, 441--451 (unbound
  insert).

\bibitem[\protect\citeauthoryear{Reshef, Reshef, Sabeti, and
  Mitzenmacher}{Reshef et~al.}{2018}]{Reshef2018}
Reshef, D.~N., Y.~A. Reshef, P.~C. Sabeti, and M.~Mitzenmacher (2018).
\newblock An empirical study of the maximal and total information coefficients
  and leading measures of dependence.
\newblock {\em Ann. Appl. Stat.\/}~{\em 12\/}(1), 123--155.

\bibitem[\protect\citeauthoryear{Reshef, Reshef, Finucane, Sabeti, and
  Mitzenmacher}{Reshef et~al.}{2016}]{Reshef2016}
Reshef, Y.~A., D.~N. Reshef, H.~K. Finucane, P.~C. Sabeti, and M.~Mitzenmacher
  (2016).
\newblock Measuring dependence powerfully and equitably.
\newblock {\em J. Mach. Learn. Res.\/}~{\em 17}, Paper No. 212, 63.

\bibitem[\protect\citeauthoryear{Rizzo}{Rizzo}{2009}]{Rizzo2009}
Rizzo, M.~L. (2009).
\newblock New goodness-of-fit tests for {P}areto distributions.
\newblock {\em Astin Bull.\/}~{\em 39\/}(2), 691--715.

\bibitem[\protect\citeauthoryear{Robert and Casella}{Robert and
  Casella}{2013}]{robert2013monte}
Robert, C. and G.~Casella (2013).
\newblock {\em Monte Carlo statistical methods}.
\newblock Springer Science \& Business Media.

\bibitem[\protect\citeauthoryear{Rockafellar}{Rockafellar}{1966}]{Rockafellar1966}
Rockafellar, R.~T. (1966).
\newblock Characterization of the subdifferentials of convex functions.
\newblock {\em Pacific J. Math.\/}~{\em 17}, 497--510.

\bibitem[\protect\citeauthoryear{Rosenbaum}{Rosenbaum}{2005}]{rosenbaum2005}
Rosenbaum, P.~R. (2005).
\newblock An exact distribution-free test comparing two multivariate
  distributions based on adjacency.
\newblock {\em J. R. Stat. Soc. Ser. B Stat. Methodol.\/}~{\em 67\/}(4),
  515--530.

\bibitem[\protect\citeauthoryear{Rosenblatt}{Rosenblatt}{1975}]{Rosenblatt1975}
Rosenblatt, M. (1975).
\newblock A quadratic measure of deviation of two-dimensional density estimates
  and a test of independence.
\newblock {\em Ann. Statist.\/}~{\em 3}, 1--14.

\bibitem[\protect\citeauthoryear{Rousson}{Rousson}{2002}]{rousson2002}
Rousson, V. (2002).
\newblock On distribution-free tests for the multivariate two-sample
  location-scale model.
\newblock {\em J. Multivariate Anal.\/}~{\em 80\/}(1), 43--57.

\bibitem[\protect\citeauthoryear{Schilling}{Schilling}{1986}]{schilling1986}
Schilling, M.~F. (1986).
\newblock Multivariate two-sample tests based on nearest neighbors.
\newblock {\em J. Amer. Statist. Assoc.\/}~{\em 81\/}(395), 799--806.

\bibitem[\protect\citeauthoryear{Sejdinovic, Sriperumbudur, Gretton, and
  Fukumizu}{Sejdinovic et~al.}{2013}]{sejdinovic2013}
Sejdinovic, D., B.~Sriperumbudur, A.~Gretton, and K.~Fukumizu (2013).
\newblock Equivalence of distance-based and {RKHS}-based statistics in
  hypothesis testing.
\newblock {\em Ann. Statist.\/}~{\em 41\/}(5), 2263--2291.

\bibitem[\protect\citeauthoryear{Sen and Sen}{Sen and Sen}{2014}]{Sen2014}
Sen, A. and B.~Sen (2014).
\newblock Testing independence and goodness-of-fit in linear models.
\newblock {\em Biometrika\/}~{\em 101\/}(4), 927--942.

\bibitem[\protect\citeauthoryear{Sen, Banerjee, and Woodroofe}{Sen
  et~al.}{2010}]{Sen2010}
Sen, B., M.~Banerjee, and M.~Woodroofe (2010).
\newblock Inconsistency of bootstrap: the {G}renander estimator.
\newblock {\em Ann. Statist.\/}~{\em 38\/}(4), 1953--1977.

\bibitem[\protect\citeauthoryear{Serfling}{Serfling}{2014}]{serfling2014multivariate}
Serfling, R.~J. (2014).
\newblock Multivariate symmetry and asymmetry.
\newblock {\em Wiley StatsRef: Statistics Reference Online\/}.

\bibitem[\protect\citeauthoryear{Shi, Drton, and Han}{Shi
  et~al.}{2019}]{shi2019distribution}
Shi, H., M.~Drton, and F.~Han (2019).
\newblock Distribution-free consistent independence tests via hallin's
  multivariate rank.
\newblock {\em arXiv preprint arXiv:1909.10024\/}.

\bibitem[\protect\citeauthoryear{Smirnoff}{Smirnoff}{1939}]{smirnoff1939}
Smirnoff, N. (1939).
\newblock On the estimation of the discrepancy between empirical curves of
  distribution for two independent samples.
\newblock {\em Bull. Math. Univ. Moscou\/}~{\em 2\/}(2), 16.

\bibitem[\protect\citeauthoryear{Spearman}{Spearman}{1904}]{spearman1904proof}
Spearman, C. (1904).
\newblock The proof and measurement of association between two things.
\newblock {\em American journal of Psychology\/}~{\em 15\/}(1), 72--101.

\bibitem[\protect\citeauthoryear{Sz\'{e}kely and Rizzo}{Sz\'{e}kely and
  Rizzo}{2005}]{rizzo2005}
Sz\'{e}kely, G.~J. and M.~L. Rizzo (2005).
\newblock Hierarchical clustering via joint between-within distances: extending
  {W}ard's minimum variance method.
\newblock {\em J. Classification\/}~{\em 22\/}(2), 151--183.

\bibitem[\protect\citeauthoryear{Sz\'{e}kely and Rizzo}{Sz\'{e}kely and
  Rizzo}{2009}]{szekely2009}
Sz\'{e}kely, G.~J. and M.~L. Rizzo (2009).
\newblock Brownian distance covariance.
\newblock {\em Ann. Appl. Stat.\/}~{\em 3\/}(4), 1236--1265.

\bibitem[\protect\citeauthoryear{Sz\'{e}kely and Rizzo}{Sz\'{e}kely and
  Rizzo}{2013}]{Gabor2013}
Sz\'{e}kely, G.~J. and M.~L. Rizzo (2013).
\newblock Energy statistics: a class of statistics based on distances.
\newblock {\em J. Statist. Plann. Inference\/}~{\em 143\/}(8), 1249--1272.

\bibitem[\protect\citeauthoryear{Sz\'{e}kely, Rizzo, and Bakirov}{Sz\'{e}kely
  et~al.}{2007}]{Gabor2007}
Sz\'{e}kely, G.~J., M.~L. Rizzo, and N.~K. Bakirov (2007).
\newblock Measuring and testing dependence by correlation of distances.
\newblock {\em Ann. Statist.\/}~{\em 35\/}(6), 2769--2794.

\bibitem[\protect\citeauthoryear{Taskinen, Kankainen, and Oja}{Taskinen
  et~al.}{2003}]{taskinen2003}
Taskinen, S., A.~Kankainen, and H.~Oja (2003).
\newblock Sign test of independence between two random vectors.
\newblock {\em Statist. Probab. Lett.\/}~{\em 62\/}(1), 9--21.

\bibitem[\protect\citeauthoryear{Taskinen, Oja, and Randles}{Taskinen
  et~al.}{2005}]{taskinen2005}
Taskinen, S., H.~Oja, and R.~H. Randles (2005).
\newblock Multivariate nonparametric tests of independence.
\newblock {\em J. Amer. Statist. Assoc.\/}~{\em 100\/}(471), 916--925.

\bibitem[\protect\citeauthoryear{van~de Geer}{van~de Geer}{2000}]{vande2000}
van~de Geer, S.~A. (2000).
\newblock {\em Applications of empirical process theory}, Volume~6 of {\em
  Cambridge Series in Statistical and Probabilistic Mathematics}.
\newblock Cambridge University Press, Cambridge.

\bibitem[\protect\citeauthoryear{van~der Vaart and Wellner}{van~der Vaart and
  Wellner}{1996}]{vaart1996}
van~der Vaart, A.~W. and J.~A. Wellner (1996).
\newblock {\em Weak convergence and empirical processes}.
\newblock Springer Series in Statistics. Springer-Verlag, New York.
\newblock With applications to statistics.

\bibitem[\protect\citeauthoryear{Varadarajan}{Varadarajan}{1958}]{Varadarajan1958}
Varadarajan, V.~S. (1958).
\newblock On the convergence of sample probability distributions.
\newblock {\em Sankhy\={a}\/}~{\em 19}, 23--26.

\bibitem[\protect\citeauthoryear{Villani}{Villani}{2003}]{Villani2003}
Villani, C. (2003).
\newblock {\em Topics in optimal transportation}, Volume~58 of {\em Graduate
  Studies in Mathematics}.
\newblock American Mathematical Society, Providence, RI.

\bibitem[\protect\citeauthoryear{Villani}{Villani}{2009}]{Villani2009}
Villani, C. (2009).
\newblock {\em Optimal transport}, Volume 338 of {\em Grundlehren der
  Mathematischen Wissenschaften [Fundamental Principles of Mathematical
  Sciences]}.
\newblock Springer-Verlag, Berlin.
\newblock Old and new.

\bibitem[\protect\citeauthoryear{Wald and Wolfowitz}{Wald and
  Wolfowitz}{1940}]{wald1940}
Wald, A. and J.~Wolfowitz (1940).
\newblock On a test whether two samples are from the same population.
\newblock {\em Ann. Math. Statistics\/}~{\em 11}, 147--162.

\bibitem[\protect\citeauthoryear{Weihs, Drton, and Meinshausen}{Weihs
  et~al.}{2018}]{weihs2018}
Weihs, L., M.~Drton, and N.~Meinshausen (2018).
\newblock Symmetric rank covariances: a generalized framework for nonparametric
  measures of dependence.
\newblock {\em Biometrika\/}~{\em 105\/}(3), 547--562.

\bibitem[\protect\citeauthoryear{Weiss}{Weiss}{1960}]{weiss1960two}
Weiss, L. (1960).
\newblock Two-sample tests for multivariate distributions.
\newblock {\em Ann. Math. Statist.\/}~{\em 31}, 159--164.

\bibitem[\protect\citeauthoryear{Wilcoxon}{Wilcoxon}{1947}]{Wilcoxon1947}
Wilcoxon, F. (1947).
\newblock Probability tables for individual comparisons by ranking methods.
\newblock {\em Biometrics\/}~{\em 3}, 119--122.

\bibitem[\protect\citeauthoryear{Wilks}{Wilks}{1938}]{wilks1938large}
Wilks, S.~S. (1938).
\newblock The large-sample distribution of the likelihood ratio for testing
  composite hypotheses.
\newblock {\em The Annals of Mathematical Statistics\/}~{\em 9\/}(1), 60--62.

\bibitem[\protect\citeauthoryear{Xiao, Frisina, Gordon, Klebanov, and
  Yakovlev}{Xiao et~al.}{2004}]{xiao2004multivariate}
Xiao, Y., R.~Frisina, A.~Gordon, L.~Klebanov, and A.~Yakovlev (2004).
\newblock Multivariate search for differentially expressed gene combinations.
\newblock {\em BMC bioinformatics\/}~{\em 5\/}(1), 164.

\bibitem[\protect\citeauthoryear{Yang}{Yang}{2012}]{yang2012}
Yang, G. (2012).
\newblock {\em The energy goodness-of-fit test for univariate stable
  distributions}.
\newblock ProQuest LLC, Ann Arbor, MI.
\newblock Thesis (Ph.D.)--Bowling Green State University.

\bibitem[\protect\citeauthoryear{Zhao and Meng}{Zhao and Meng}{2015}]{Zhao2015}
Zhao, J. and D.~Meng (2015).
\newblock Fast{MMD}: ensemble of circular discrepancy for efficient two-sample
  test.
\newblock {\em Neural Comput.\/}~{\em 27\/}(6), 1345--1372.

\bibitem[\protect\citeauthoryear{Zuo and Serfling}{Zuo and
  Serfling}{2000}]{Zuo2000}
Zuo, Y. and R.~Serfling (2000).
\newblock General notions of statistical depth function.
\newblock {\em Ann. Statist.\/}~{\em 28\/}(2), 461--482.

\end{thebibliography}
	\end{spacing}
     \newpage
     \appendix
     \section{Real data analysis}\label{sec:realdata}
     In this section, we will provide some real data examples where we shall compare the performance of $\Rdcov_n^2$ (henceforth referred to as $\Rdcov$) with usual distance covariance, and also Pearson's correlation. For the sake of transparency, we have chosen these real data examples from benchmark data sets which have been analyzed previously in the literature. In subsequent data analysis, we have implemented the distance correlation (DCR) based test using the \texttt{dcor.test} function in the \texttt{energy} package, and the Pearson's correlation test (P) using the \texttt{cor.test} function in the \texttt{stats} package in \texttt{R} (\cite{Rsoft}). \par 
     \begin{example}[US crimes data set]\label{ex:datasetcrime}
     	Consider the US crime data set from $110$ metropolitan areas with populations larger than $250,000$. The particular attributes in the data set include --- (i) {\it population} (in thousands, from $1968$), (ii) {\it nonwhite} (percentage of nonwhite population, $1960$), (iii) {\it density} (population per square mile, $1968$) and (iv) {\it crime} (crime rate per thousand, $1969$). This data set is available in~\cite{Fox2019}.  \par 
     	The data set contains missing values. In fact, complete data is available for $100$ out of the $110$ metropolitan areas, and we will restrict all subsequent data analysis to these $100$ metropolitan areas. This preprocessing is exactly the same as in~\cite{szekely2009}. Also,  as in~\cite{szekely2009}, we are interested in answering which two of the four attributes (mentioned above) are associated among each other. The following table describes our findings:
     	\begin{table}[h]
     		\begin{tabular}{| *{10}{Sc|}}
     			\hline
     			\multirow{2}{*}{Dataset} &
     			\multicolumn{3}{c|}{Nonwhite} &
     			\multicolumn{3}{c|}{Density} &
     			\multicolumn{3}{c|}{Crime} \\
     			& (P) & (DCR) & ($\Rdcov$) & (P) & (DCR) & $(\Rdcov$) & (P) & (DCR) & $(\Rdcov$)\\
     			\hline
     			Population & 0.49 & 0.01 & 0.38 & 0.00 & 0.00 & 0.00  & 0.00 & 0.00 & 0.00\\
     			\hline
     			Nonwhite &  &  &  & 0.98 & 0.29 & 0.02 & 0.00 & 0.00 & 0.00\\
     			\hline
     			Density &  &  &  &  &  &  & 0.27 & 0.03 & 0.03\\
     			\hline
     		\end{tabular}
     		\vspace{0.1in}
     		\caption{$P$-values for Pearson's correlation, distance correlation and rank distance covariance for the US crime data set.}
     		\label{table:datacrimetab}
     	\end{table}
     	
     	\noindent In~\cref{table:datacrimetab}, apart from population (i) versus nonwhite population (ii), and nonwhite population (ii) versus population density (iii), all the other possible combinations are shown to be associated (at least at the $5\%$ level) by all the three tests. So, let us inspect the above two possible combinations more carefully.
     	\noindent First let us look at the scatter plot between total population and nonwhite population (top left panel in~\cref{fig:UScrime}). It seems that there is no linear relationship between these variables, which is why Pearson's correlation does not detect any association between these variables. However, if we look at the scatter plot between the ranks of total and nonwhite populations, there does seem to be a linear relationship (see top right panel in~\cref{fig:UScrime}). This is borne out by a standard linear regression analysis with nonwhite population ranks as the response and total population ranks as the covariate, which leads to the following output:
     	\begin{verbatim}
     	> summary(fit)
     	Call: lm(formula = var2[, 2] ~ var1[, 2])
     	Coefficients:
     	Estimate Std. Error t value Pr(>|t|)    
     	(Intercept)          0.37175    0.05668   6.559 2.56e-09 ***
     	var1[, 2]  0.26386    0.09744   2.708  0.00799 ** 
     	Signif. codes:  0 ‘***’ 0.001 ‘**’ 0.01 ‘*’ 0.05 ‘.’ 0.1 ‘ ’ 1
     	\end{verbatim}
     	\noindent This suggests a clear association between these variables which is picked up by both distance correlation (at $5\%$) and rank distance covariance (at $1\%$).\par 
     	\begin{figure}[H]
     		\begin{center}
     			\includegraphics[height=6.5cm,width=7.5cm]{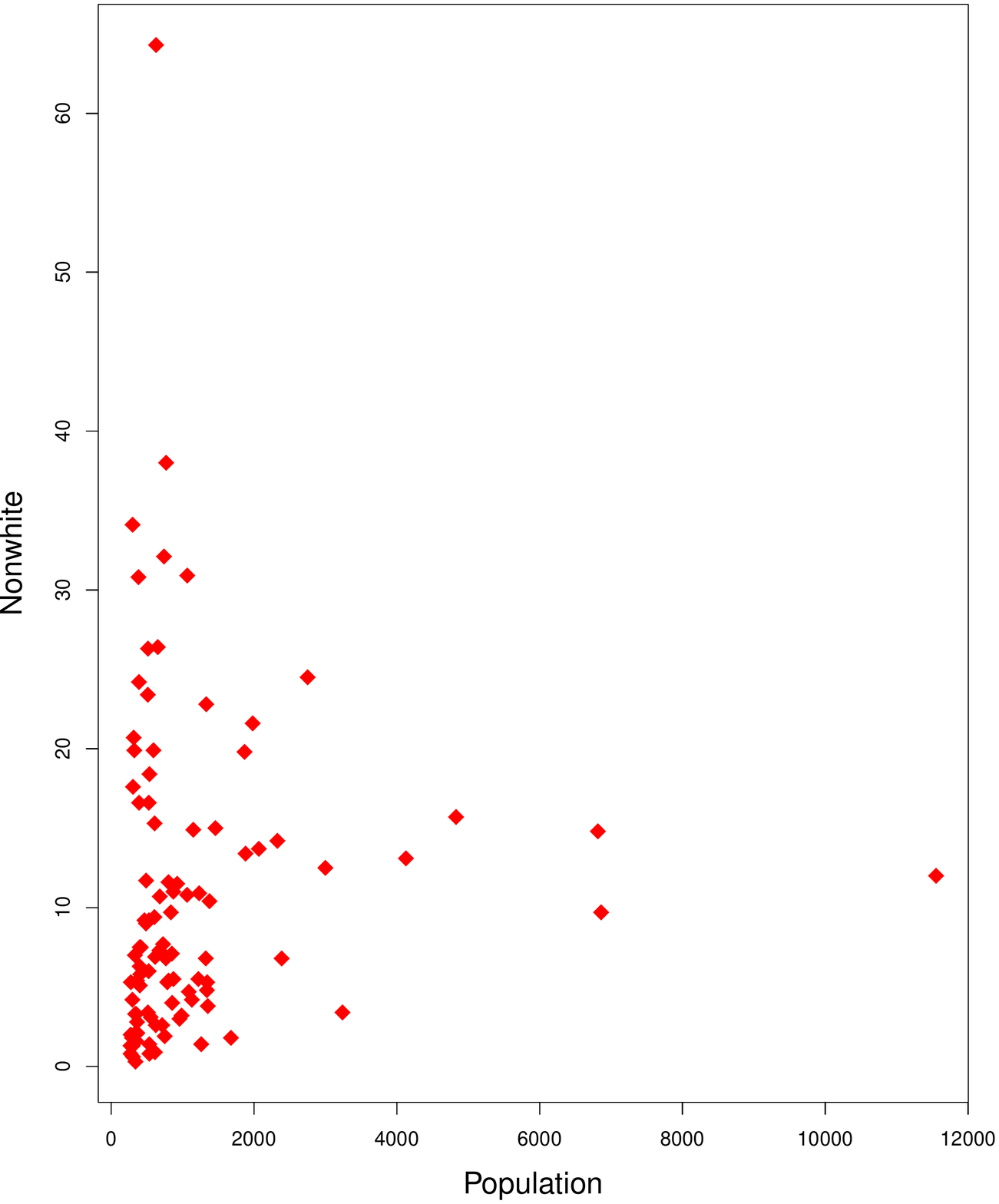}	
     			\includegraphics[height=6.5cm,width=7.5cm]{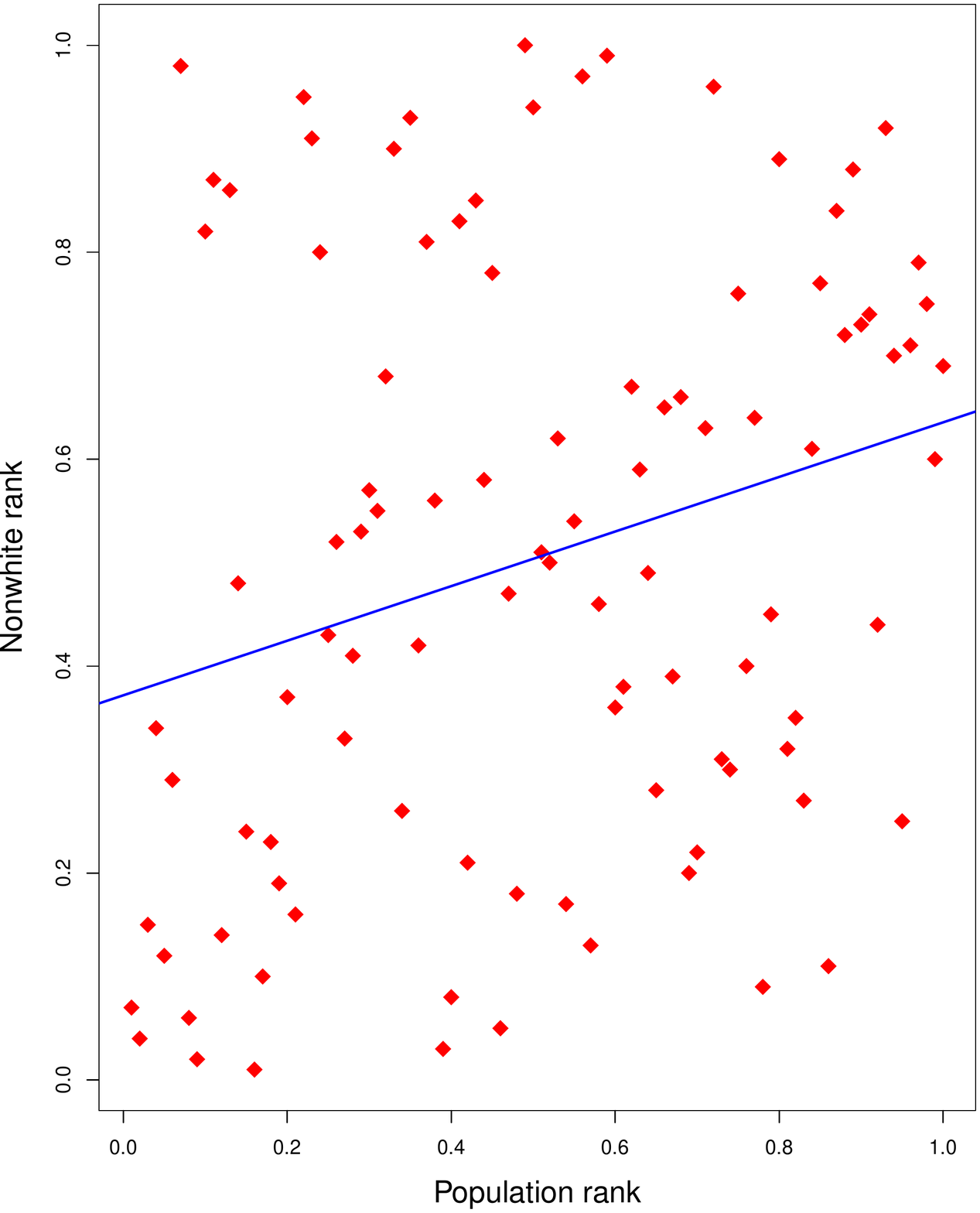}
     			\includegraphics[height=6.5cm,width=7.5cm]{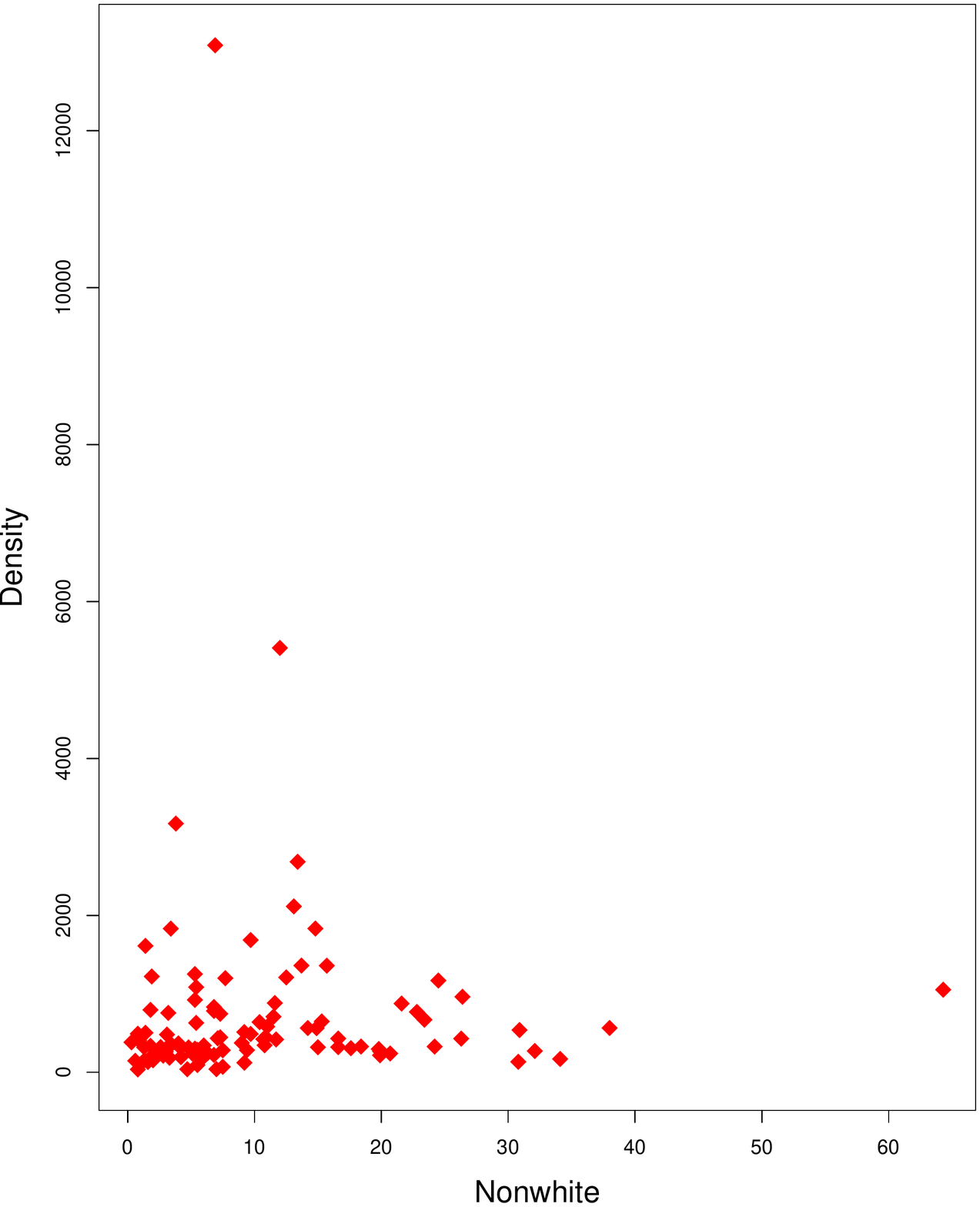}	
     			\includegraphics[height=6.5cm,width=7.5cm]{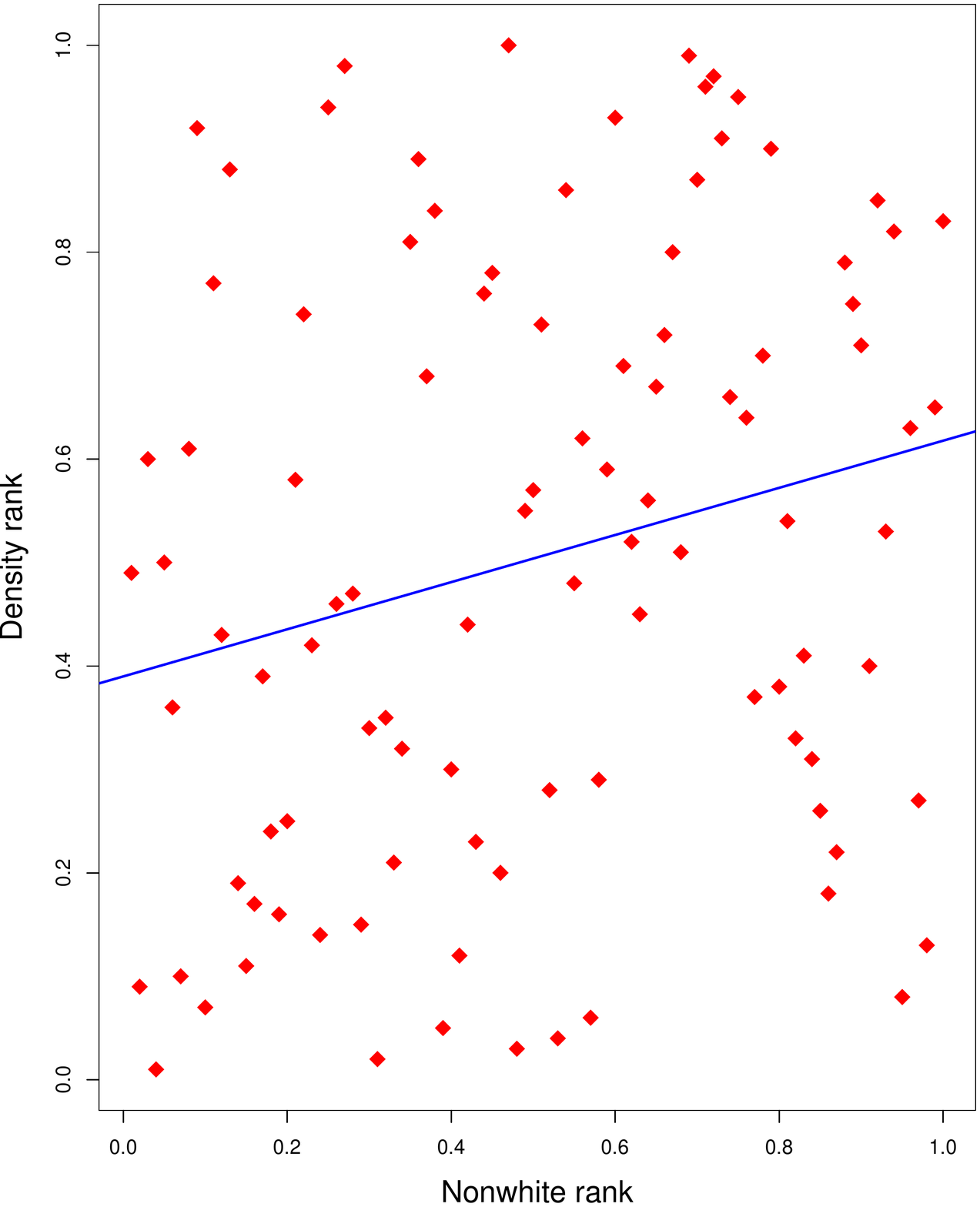}
     			\caption{The top left panel shows the scatter plot between total population $(i)$ and nonwhite population $(ii)$. The top right panel plots the population ranks from the $100$ metropolitan areas versus the corresponding nonwhite population ranks. The bottom left panel shows the scatter plot between nonwhite population $(ii)$ and population density $(iii)$. The bottom right panel plots the nonwhite population ranks from the $100$ metropolitan areas versus the corresponding population density ranks.}
     			\label{fig:UScrime}
     		\end{center}
     	\end{figure}
     	
     	\noindent Next, let us focus on the association between nonwhite population and population density. Firstly, it is natural to expect an association in this case, given that total and nonwhite population were associated (as argued in the previous paragraph). Once again, the corresponding scatter plot (see the bottom left panel of~\cref{fig:UScrime}) reveals no linear relationship, whereas the corresponding plot of ranks (see the bottom right panel in~\cref{fig:UScrime}) shows some linear relationship. This is supported by standard linear regression analysis as before:
     	\begin{verbatim}
     	> summary(fit)
     	Call: lm(formula = var2[, 2] ~ var1[, 2])
     	Coefficients:
     	Estimate Std. Error t value Pr(>|t|)    
     	(Intercept)          0.38992    0.05721   6.815 7.71e-10 ***
     	var1[, 2]  0.22789    0.09836   2.317   0.0226 *  
     	Signif. codes:  0 ‘***’ 0.001 ‘**’ 0.01 ‘*’ 0.05 ‘.’ 0.1 ‘ ’ 1
     	\end{verbatim}
     	Once again, there is evidence of association, but in this case, only our proposed rank distance covariance test detects this association whereas Pearson's correlation and distance correlation support the null hypothesis of independence. We believe that this is because of the presence of outliers, as evidenced by the scatter plot in the bottom left panel of~\cref{fig:UScrime}. As distance correlation requires finite moment assumptions for consistency, its performance can be affected adversely by the presence of outliers, which however, rank distance covariance is successfully robust to. This provides a real practical example when rank distance covariance may be more useful than usual distance covariance in exploratory association studies.
     	\begin{remark}\label{rem:ascutdata}
     		In the above data analysis, the performance of rank distance covariance stays the same irrespective of whether we use cutoffs from the universal distribution for fixed $n$ or the universal limit distribution (see~\cref{theo:indepasdistn}). We believe that the convergence of the fixed $n$ universal distribution to its asymptotic limit happens rather quickly when $d_1=d_2=1$. We will discuss more on this in~\cref{sec:univcutas}.
     	\end{remark}
     \end{example}
     \begin{example}[SONAR data set]\label{ex:sonar}
     	In this example, we look at a multivariate two-sample equality of distributions testing problem based on the benchmark data set \texttt{Sonar} (see~\cite{Sonardata}) available in the \texttt{R} package \texttt{mlbench}  The data comprises patterns obtained by reflecting sonar signals off a metal cylinder or a roughly cylindrical rock. Each pattern is a $60$ dimensional vector, with each entry between $0$ and $1$. Each number represents the energy of the signal within specific frequency bands, integrated over time. There are $111$ patterns from signals bounced off metal cylinders ({\it M-signals}) and $97$ patterns from signals bounced off rocks ({\it R-signals}). As in~\cite{munmun2014}, we are interested in testing whether the patterns arising out of {\it M-signals} and {\it R-signals} have the same distribution. \par 
     	Since this $60$ dimensional data is difficult to visualize, we first obtained the principal components corresponding to the {\it M-signals} and then projected the patterns corresponding to both kinds of signals along those $60$ directions. This results in $60$ one-dimensional projections corresponding to patterns from both {\it M-signals} and {\it R-signals}. In~\cref{fig:Sonar}, we show the QQ-plots between {\it M-signals} and {\it R-signals} corresponding to $3$ of these $60$ projections. It is clear in the plots that none of the distributions of the $3$ one-dimensional projections are the same. In fact, in $21$ out of these $60$ projections we get similar QQ-plots and corresponding $p$-values (from a two-sample Kolmogorov-Smirnov test) are all below $0.001$.\par 
     	The above discussion highlights that it is reasonable to expect that the null hypothesis from the two-sample equality of distributions test of interest should be rejected. At the standard $5\%$ level of significance, our proposed rank energy test (RE), Rosenbaum's crossmatch test (RC), the usual energy test (EN) and the Heller-Heller-Gorfine test (HHG) reject the null hypothesis, whereas the maximum mean discrepancy test (MMD) with Gaussian kernel fails to reject the null.
     	\begin{figure}[H]
     		\begin{center}
     			\includegraphics[height=6.5cm,width=5cm]{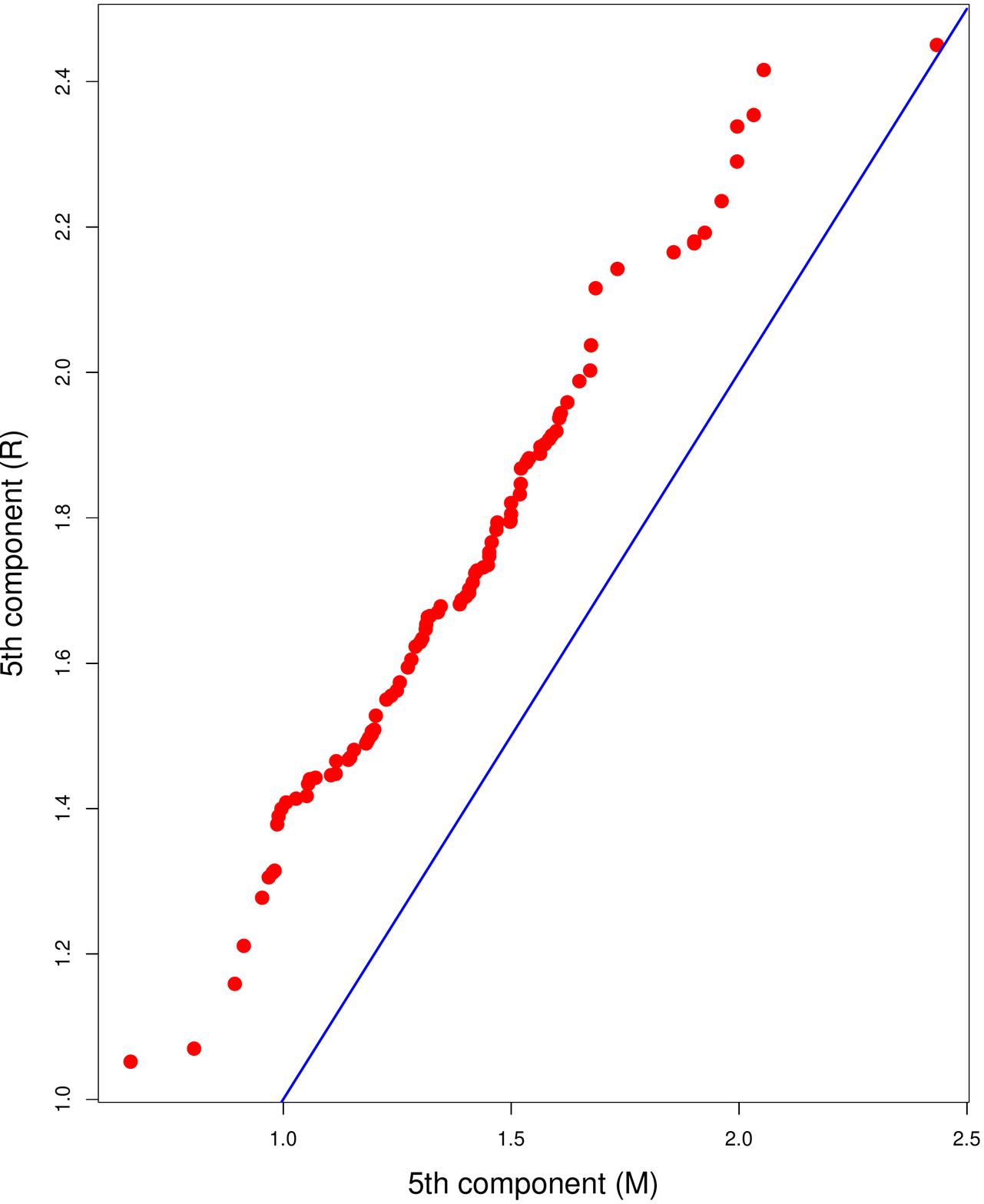}
     			\includegraphics[height=6.5cm,width=5cm]{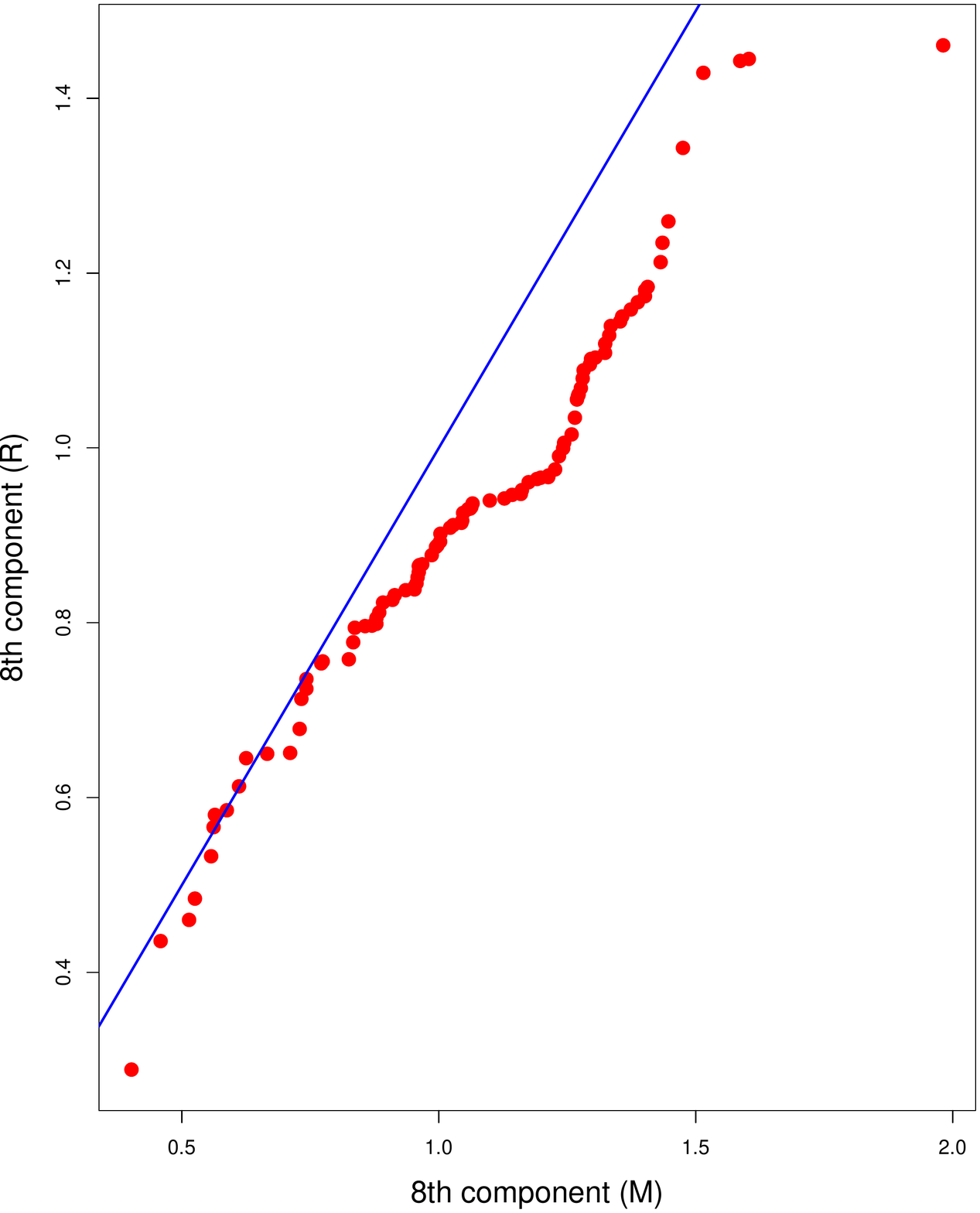}
     			\includegraphics[height=6.5cm,width=5cm]{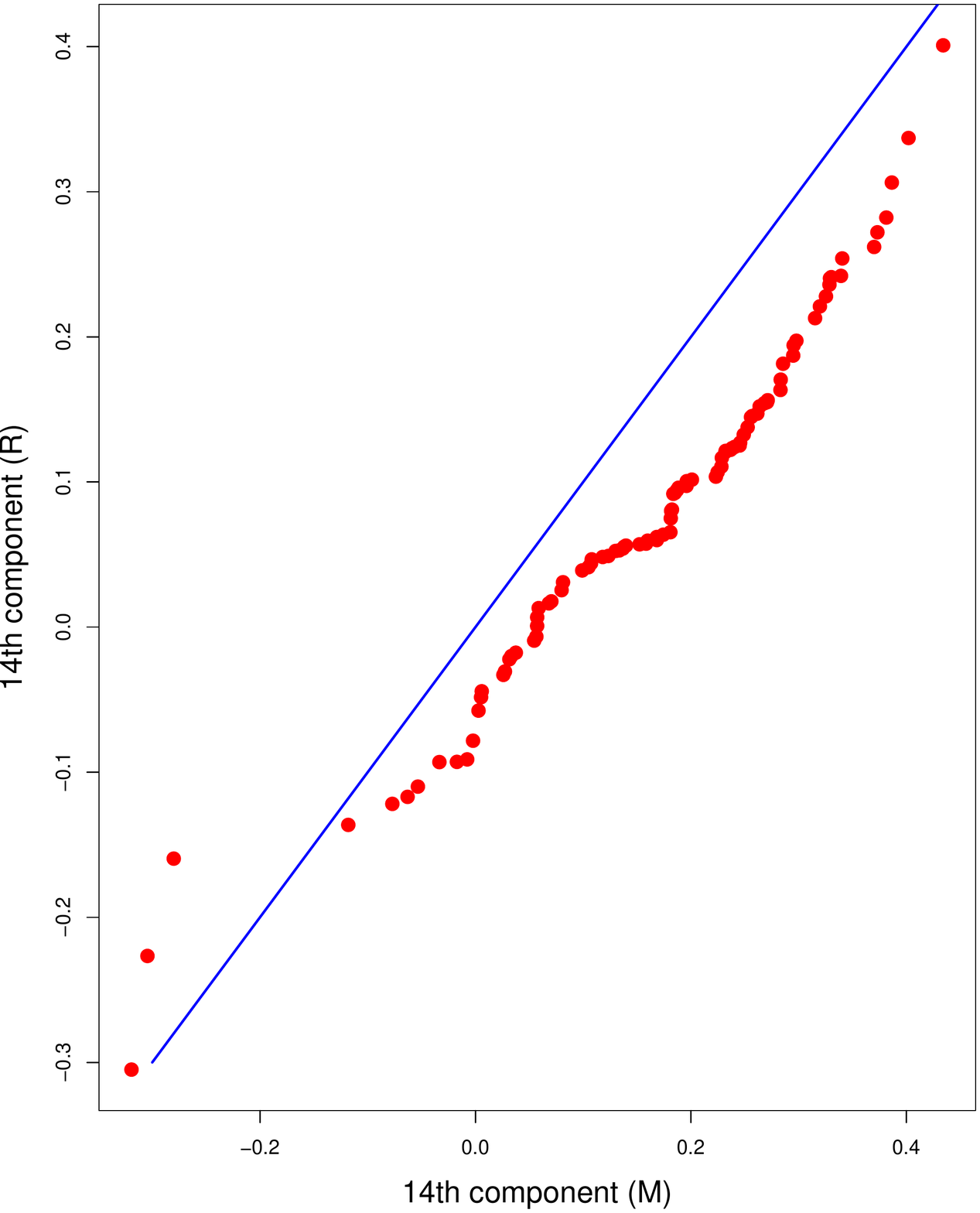}
     			\caption{QQ-plots corresponding to the one-dimensional projections of the {\it M-signals} and {\it R-signals} along the $5$'th, $8$'th and $14$'th principal components obtained from the data matrix corresponding to the {\it M-signals}.}
     			\label{fig:Sonar}		
     		\end{center}
     	\end{figure}
     \end{example} 
     \section{Computation of the test statistics}\label{sec:compute}
     We begin by introducing the~\emph{assignment problem} and illustrating its connection to the computation of our proposed multivariate ranks. Suppose that $n$ tasks are to be divided between $n$ agents. Any agent can be assigned to perform any task, incurring some cost that may vary depending on the agent-task assignment. It is however required that all agents perform one and only one task. Under this constraint, the assignment problem seeks to find the agent-task allotment which minimizes the overall cost. Suppose that the list of agents is denoted by  $\{p_1,p_2,\ldots ,p_n\}$ and the list of tasks by $\{t_1,t_2,\ldots ,t_n\}$. Also let $C(p_i,t_j)$ denote the cost of assigning task $t_j$ to person $p_i$ and finally, let $\mathcal{S}_n$ denote the set of all bijective functions from the set $\{p_1,\ldots ,p_n\}$ to $\{t_1,\ldots ,t_n\}$. Then the above problem may be stated as:
     \begin{equation}\label{eq:assign1}
     \min_{f\in \mathcal{S}_n} \sum_{i=1}^n C(p_i,f(p_i)).
     \end{equation}
     This problem has been studied extensively in the combinatorial optimization literature, see e.g.,~\cite{munkres1957,bertsekas1988}. One of the more efficient methods to solve~\eqref{eq:assign1} above is to use the Hungarian algorithm (see e.g.,~\cite{jonker1987}) which has worst case computational complexity $\mathcal{O}(n^3)$.
     
     It is easy to see the connection between the assignment problem as in~\eqref{eq:assign1} and our empirical transport problem as in~\eqref{eq:empopt}. As a result, by using the \emph{Hungarian algorithm}, we can obtain the vector of multivariate ranks in at most $\mathcal{O}(n^3)$ steps. Recall that distance covariance can be computed in $\mathcal{O}(n^2 d)$ steps. Therefore, our proposed multivariate rank-based test statistic can be computed in $\mathcal{O}(n^3+n^2 d)$ steps. A similar argument shows that our proposed multivariate rank energy statistic can be computed in $\mathcal{O}(m^3+n^3+mn d)$ steps.
     \par 
     It is important to note that the implementation of our proposed tests is extremely simple using, for example, \texttt{R}. First, we may generate the standard Halton sequence using the R-package \texttt{randtoolbox}, then we obtain the multivariate ranks using the R-package \texttt{clue} and finally we can calculate distance covariance (or the energy statistic) among the multivariate ranks using the R-package \texttt{energy}. All the methods described in this paper have been implemented using the \texttt{R} software. The relevant codes, including our simulation experiments, are available in the first author's \texttt{GitHub} page. 
     
     \section{Some additional simulations}\label{sec:addsim}
     Here we continue the discussion on the empirical performance of our proposed tests from~\cref{sec:sim}. Among other things, we give a detailed account on the empirical performance of $\Ren_{m,n}^2$, along with the universal asymptotic cutoffs for our proposed exactly distribution-free tests (see Tables~\ref{table:thresholdind} and~\ref{table:twosamascutoff}).	
     \subsection{Synthetic data experiments for two-sample goodness-of-fit  testing}\label{sec:syngof}
     We illustrate the empirical performance of $\Ren_{m,n}^2$ (hereafter referred to as $\Ren$) based on synthetic data. Throughout our simulation settings, we fix $m=n=200$, and $d=3$. We will compare our multivariate rank-based test ($\Ren$) to Rosenbaum's cross matching test (hereafter referred to as CMT; see~\cite{rosenbaum2005}) which is the only other computationally feasible distribution-free test for two-sample goodness-of-fit testing. This test has been implemented using the \texttt{R} package \texttt{crossmatch}. In addition, we also use the usual energy test (\cite{Gabor2013}; hereafter referred to as EN) from the \texttt{energy} package, and the Heller-Heller-Gorfine test (\cite{HHGtwosam}; hereafter called HHG) from the package \texttt{HHG}, as benchmarks. We also implemented the maximum mean discrepancy (MMD) based test (see~\cite{gretton2012}) using the \texttt{kmmd} function in the \texttt{R} package \texttt{kernlab} (see~\cite{kernlab}). However its performance was convincingly poorer than the other chosen methods for our proposed simulation settings, and so we refrain from providing those details. Of course, none of HHG, EN or MMD are exactly distribution-free. \par 
     As a general rule, we will write the three-dimensional independent vectors $\mathbf{X}$ and $\mathbf{Y}$ as $(X_1,X_2,X_3)$ and $(Y_1,Y_2,Y_3)$. We generate $200$ i.i.d.~copies of such random vectors and carry out a two-sample goodness-of-fit test based on these observations. Below we list our simulation settings:
     \begin{itemize}
     	\item[(V1)] $X_1,X_2,X_3,Y_1\overset{i.i.d.}{\sim}$ Cauchy$(0,1)$, and $Y_2,Y_3\overset{i.i.d.}{\sim}$ Cauchy$(0.2,1)$.
     	\item[(V2)] $X_1,Y_1$ are i.i.d. $\mathcal{U}^1$, $X_k=0.25+0.35\times X_{k-1}+U_k$ and $Y_k=0.25+0.5\times Y_{k-1}+V_k$ for $ k=2,3$. Here $U_2,U_3,V_2 ,V_3$ are i.i.d.~$\mathcal{U}^1$.
     	\item[(V3)] $\mathbf{X}\sim \mathcal{N}_3(\mathbf{0},\Sigma_1)$ and $\mathbf{Y}\sim\mathcal{N}_3(\mathbf{0},\Sigma_2)$ where $\Sigma_1(i,j)=0.35^{|i-j|}$ and $\Sigma_2(i,j)=0.65^{|i-j|}$ for $1\leq i,j\leq 3$.
     	\item[(V4)] $\mathbf{X}\sim \mathcal{N}_3(\mathbf{0},\Sigma_1)$ and $\mathbf{Y}\sim\mathcal{N}_3(\mathbf{0},\Sigma_2)$ where $\Sigma_1(i,j)=0.2$ for $i\neq j$ and $\Sigma_1(i,i)=1$, $\Sigma_2(i,j)=0.5$ for $i\neq j$ and $\Sigma_2(i,i)=1$, for $1\leq i,j\leq 3$.
     	\item[(V5)] $\mathbf{V}\sim \mathcal{N}_3(\mathbf{0},\Sigma_1)$ and $\mathbf{W}\sim\mathcal{N}_3(\mathbf{0},\Sigma_2)$ where $\Sigma_1(i,j)=0.35^{|i-j|}$ and $\Sigma_2(i,j)=0.75^{|i-j|}$ for $1\leq i,j\leq 3$. Set $X_i=\exp{(V_i)}$ and $Y_i=\exp{(W_i)}$ for $i=1,2,3$.
     	\item[(V6)] $\mathbf{V}\sim \mathcal{N}_3(\mathbf{0},\Sigma_1)$ and $\mathbf{W}\sim\mathcal{N}_3(\mathbf{0},\Sigma_2)$ where $\Sigma_1(i,j)=0.25$ for $i\neq j$ and $\Sigma_1(i,i)=1$, $\Sigma_2(i,j)=0.75$ for $i\neq j$ and $\Sigma_2(i,i)=1$, for $1\leq i,j\leq 3$. Finally, let $X_i=\exp{(V_i)}$ and $Y_i=\exp{(W_i)}$ for $i=1,2,3$.
     	\item[(V7)] $\mathbf{X}\sim \mathcal{N}_3(\mathbf{\mu}_1,3\mathbf{I})$ and $\mathbf{Y}\sim\mathcal{N}_3(\mathbf{\mu}_2,3\mathbf{I})$ where $\mathbf{\mu}_1=(0,0,0)$ and $\mathbf{\mu}_2=(0.25,0.25,0.25)$.
     	\item[(V8)] $\mathbf{V}\sim \mathcal{N}_3(\mathbf{\mu}_1,3\mathbf{I})$ and $\mathbf{W}\sim\mathcal{N}_3(\mathbf{\mu}_2,3\mathbf{I})$ where $\mathbf{\mu}_1=(0,0,0)$ and $\mathbf{\mu}_2=(0.25,0.25,0.25)$. Finally, let $X_i=V_i$ and $Y_i=W_i$ for $i=1,2,3$.
     	\item[(V9)] $X_1,X_2,X_3,V_1,V_2,V_3$ are i.i.d.~Gamma$(2,0.1)$, and $W_1,W_2,W_3$ are i.i.d.~with the same law as $\exp{(\exp(Z))}$ where $Z\sim\mathcal{N}(0,1)$. Finally, set $Y_i=W_iV_i$ for $i=1,2,3$.
     	\item[(V10)] $\mathbf{Z}_1,\mathbf{Z}_2$ are i.i.d.~$\mathcal{N}_3(\mathbf{1},\mathbf{I})$ and $A \sim$ Ber$(0.8)$. Let $\mathbf{W}=(W_1,W_2,W_3)$ be such that $W_i$'s are i.i.d.~$\mathcal{U}(10,11)$. Set $\mathbf{X}\coloneqq\mathbf{Z}_1$ and $\mathbf{Y}\coloneqq A\mathbf{Z}_2+(1-A)\mathbf{W}$.
     	\item[(V11)] $\mathbf{Z}_1,\mathbf{Z}_2$ are i.i.d.~$\mathcal{N}_3(\mathbf{1},\mathbf{I})$ and $A \sim$ Ber$(0.8)$. Let $\mathbf{W}=(W_1,W_2,,W_3)$ be such that $W_i$'s are i.i.d.~$\mathcal{N}(10,0.1)$. Set $\mathbf{X}\coloneqq\mathbf{Z}_1$ and $\mathbf{Y}\coloneqq A\mathbf{Z}_2+(1-A)\mathbf{W}$.
     \end{itemize}
     Many of the simulation settings above are similar to those considered in e.g.,~\cite{bbbm2019,munmun2014}. For instance, setting (V2) has been adopted from~\cite[Section 4]{munmun2014}. All the settings (V3)-(V8) are slightly modified versions of similar settings from~\cite[Section 3.3 and Appendix D]{bbbm2019}. These minor tweaks were made to make sure that the competing procedures have non-trivial power for the prescribed values of $n$ and $d$. Settings (V1) and (V9) deal with scenarios where the associated distributions do not have finite first moments. In settings (V10) and (V11), we look at settings featuring mixture distributions where there is a small proportion of noise added to one of the two otherwise identical distributions.
     In~\cref{table:twosamcompare}, we present our findings. The two columns corresponding to each methods represent the rejection probabilities (estimated using $1000$ independent replications) at levels $0.05$ and $0.1$.\par 
     \begin{table}[h]
     	\centering
     	\vspace{0.2in}
     	\begin{tabular}{| *{9}{Sc|}}
     		\hline
     		& \multicolumn{2}{c|}{(CMT)} & \multicolumn{2}{c|}{(HHG)} & \multicolumn{2}{c|}{(EN)} & \multicolumn{2}{c|}{($\Ren$)} \\ 
     		\hline
     		V1 & 0.08 & 0.13 & 0.08 & 0.15 & 0.07 & 0.13 & \textbf{0.23} & \textbf{0.34} \\ 
     		\hline
     		V2 & 0.24 & 0.34 & \textbf{0.91} & \textbf{0.94} & \textbf{0.90} & \textbf{0.94} & 0.84 & 0.89  \\ 
     		\hline
     		V3 & \textbf{0.29} & 0.41 & 0.19 & 0.34 & 0.17 & 0.34 & \textbf{0.26} & \textbf{0.46} \\ 
     		\hline
     		V4 & \textbf{0.24} & \textbf{0.34} & 0.16 & 0.31 & 0.17 & 0.33 & 0.18 & 0.32 \\ 
     		\hline
     		V5 & 0.61 & 0.73 & 0.54 & 0.70 & 0.35 & 0.56 & \textbf{0.77} & \textbf{0.93} \\ 
     		\hline
     		V6 & 0.84 & 0.90 & 0.77 & 0.88 & 0.59 & 0.82 & \textbf{0.96} & \textbf{0.99} \\ 
     		\hline
     		V7 & 0.08 & 0.13 & 0.38 & 0.51 & \textbf{0.52} & \textbf{0.65} & 0.49 & 0.63 \\ 
     		\hline
     		V8 & 0.07 & 0.11 & 0.27 & 0.39 & 0.25 & 0.35 & \textbf{0.29} & \textbf{0.43} \\ 
     		\hline
     		V9 & 0.06 & 0.06 & \textbf{1.00} & \textbf{1.00} & 0.97 & 0.97 & \textbf{1.00} & \textbf{1.00}\\ 
     		\hline
     		V10 & 0.50 & 0.61 & \textbf{1.00} & \textbf{1.00} & \textbf{1.00} & \textbf{1.00} & 0.97 & 0.99 \\ 
     		\hline 
     		V11 & 0.47 & 0.60 & \textbf{1.00} & \textbf{1.00} & \textbf{1.00} & \textbf{1.00} & 0.95 & 0.98 \\
     		\hline
     	\end{tabular}
     	\vspace{0.1in}
     	\caption{Proportion of times the null hypothesis was rejected across $11$ settings. Here $m=n=200$ and $d=3$.}
     	\label{table:twosamcompare}
     \end{table}
     \noindent\textbf{(CMT):} ~\cref{table:twosamcompare} reveals a number of interesting points of comparison between CMT and $\Ren$. Let us start with settings (V10) and (V11). In these settings $\mathbf{X}$ is multivariate Gaussian, whereas $\mathbf{Y}$ is the same multivariate Gaussian with a small fractions of Uniform and Gaussian noise ($20\%$) respectively. In both these settings $\Ren$ outperforms CMT. This shows that $\Ren$ is perhaps more robust to outliers than CMT (as mentioned in the Introduction). Next, let us look at the heavy-tailed settings (V1) and (V9) (no finite first moments). Here again $\Ren$ convincingly outperforms CMT, once again reinforcing the superior performance of our proposed multivariate rank-based tests in the absence of finite moments (as mentioned in the Introduction). Across settings (V3)-(V8) (based on multivariate normals and log-normals), once again $\Ren$ largely outperforms CMT, except perhaps in settings (V3) and (V4). These settings feature different orders of decaying correlation and equicorrelation respectively, in a multivariate Gaussian setting. The performances of $\Ren$ and CMT are almost indistinguishable in these two scenarios; in fact, CMT perhaps marginally outperforms $\Ren$ in (V4). As multivariate normals and log-normals are popular modeling choices in practice, we will look into such simulation settings in greater detail in~\cref{sec:multgausslog}. Finally, in (V2), which has the same flavor as a first order autoregressive model, once again $\Ren$ significantly outperforms CMT. Overall, $\Ren$ has a much superior performance than CMT in the proposed simulation settings. \par 
     \noindent\textbf{(EN):} The most striking feature from~\cref{table:twosamcompare} in terms of comparison between $\Ren$ and EN is perhaps that $\Ren$ largely outperforms EN in settings (V3), (V5), (V6) and (V8), whereas in settings (V4) and (V7) the performances are mostly comparable, EN being marginally better. This has been a recurrent observation in our simulations that the EN test (or equivalently DCoV for independence testing) loses out to $\Ren$ (equivalently $\Rdcov$) when dealing with multivariate log-normals [(V4), (V6), (V8)] whereas it performs comparably when dealing with multivariate Gaussian (or mixtures thereof). This observation has been studied in more details in~\cref{sec:multgausslog} using multivariate normal and log-normal location alternatives. In the heavy-tailed settings, as expected, $\Ren$ outperforms EN; convincingly in (V1) and marginally in (V9). In settings (V2), (V10) and (V11), EN and $\Ren$ seem to have almost identical performance; EN being marginally superior. This perhaps shows that the performance of EN is also somewhat robust to small proportions of noise in data. \par 
     \noindent\textbf{(HHG):} Across settings (V3)-(V8) all of which are based on multivariate normal and log-normal based alternatives, $\Ren$ convincingly and consistently outperforms HHG. A similar observation can be seen in the heavy-tailed setting (V1). In all the other settings, the two tests have comparable performance. In particular, being based on the ranks of pairwise distances, HHG too is perhaps somewhat robust to small proportions of outliers (as indicated by its performance in settings (V10) and (V11)).
     \subsection{Multivariate normal and log-normal settings}\label{sec:multgausslog}
     Multivariate normal and log-normal settings have been used in the context of comparing nonparametric testing procedures, e.g.,~\cite{bbbm2019}. In~\cref{sec:sim}, we used multivariate normal and log-normal settings for comparing $\Rdcov$ with competing procedures for multivariate independence testing. A similar exercise can be carried out for the two-sample goodness-of-fit testing problem as well. We will consider $m=n=200$, $d=3$ and consider the following two settings:\par 
     \noindent (TG) $(X_1,X_2,X_3)\sim\mathcal{N}_3((\mu,\mu,\mu),3\mathbf{I})$ where the mean parameter $\mu$ varies in $[-1,1]$, whereas the distribution of $(Y_1,Y_2,Y_3)\sim\mathcal{N}(\mathbf{0},3\mathbf{I})$ stays fixed. \par 
     \noindent (TGL) $(X_1,X_2,X_3)=(\exp{(Z_1)},\exp{(Z_2)},\exp{(Z_3)})$ where $(Z_1,Z_2,Z_3)\sim\mathcal{N}_3((\mu,\mu,\mu),3\mathbf{I})$ and the mean parameter $\mu$ varies in $[-1,1]$, whereas the distribution of $(Y_1,Y_2,Y_3)=(\exp{(W_1)},\exp{(W_2)},\exp{(W_3)})$, with $(W_1,W_2,W_3)\sim\mathcal{N}(\mathbf{0},3\mathbf{I})$, stays fixed.\par  
     \begin{table}[h]
     	\centering
     	\begin{tabular}{| *{9}{Sc|}}
     		\hline
     		Mean & \multicolumn{2}{c|}{(CMT)} & \multicolumn{2}{c|}{\textbf{(EN)}} & \multicolumn{2}{c|}{(HHG)} & \multicolumn{2}{c|}{$(\Ren)$} \\ 
     		\hline 
     		-0.60 & 0.33 & 0.44 & \textbf{1.00} & \textbf{1.00} & 1.00 & 1.00 & 1.00 & 1.00 \\ 
     		\hline
     		-0.52 & 0.23 & 0.32 & \textbf{1.00} & \textbf{1.00} & 0.98 & 0.99 & 0.99 & 1.00 \\ 
     		\hline
     		-0.45 & 0.18 & 0.26 & \textbf{0.98} & \textbf{0.99} & 0.92 & 0.96 & 0.97 & 0.98 \\ 
     		\hline
     		-0.37 & 0.14 & 0.20 & \textbf{0.89} & \textbf{0.93} & 0.78 & 0.86 & 0.85 & 0.91 \\ 
     		\hline
     		-0.30 & 0.10 & 0.17 & \textbf{0.70} & \textbf{0.79} & 0.52 & 0.66 & 0.65 & 0.76 \\ 
     		\hline
     		-0.22 & 0.07 & 0.14 & \textbf{0.43} & \textbf{0.57} & 0.31 & 0.43 & 0.39 & 0.52 \\ 
     		\hline
     		-0.15 & 0.06 & 0.10 & \textbf{0.20} & \textbf{0.32} & 0.15 & 0.23 & 0.19 & 0.28 \\ 
     		\hline
     		-0.07 & 0.05 & 0.09 & \textbf{0.09} & \textbf{0.15} & 0.08 & 0.14 & 0.08 & 0.15 \\ 
     		\hline
     		0.07 & 0.06 & 0.10 & \textbf{0.08} & \textbf{0.15} & 0.07 & 0.14 & 0.08 & 0.15 \\ 
     		\hline
     		0.15 & 0.06 & 0.11 & \textbf{0.21} & \textbf{0.33} & 0.15 & 0.25 & 0.19 & 0.30 \\ 
     		\hline
     		0.22 & 0.08 & 0.12 & \textbf{0.41} & \textbf{0.54} & 0.28 & 0.41 & 0.38 & 0.50 \\ 
     		\hline
     		0.30 & 0.09 & 0.14 & \textbf{0.71} & \textbf{0.80} & 0.50 & 0.66 & 0.64 & 0.76 \\ 
     		\hline
     		0.37 & 0.14 & 0.22 & \textbf{0.89} & \textbf{0.94} & 0.73 & 0.84 & 0.84 & 0.92 \\ 
     		\hline
     		0.45 & 0.17 & 0.26 & \textbf{0.97} & \textbf{0.99} & 0.90 & 0.95 & 0.96 & 0.98 \\ 
     		\hline
     		0.52 & 0.25 & 0.36 & \textbf{0.99} & \textbf{1.00} & 0.98 & 0.99 & 0.99 & 1.00 \\ 
     		\hline
     		0.60 & 0.33 & 0.45 & \textbf{1.00} & \textbf{1.00} & 0.99 & 1.00 & 1.00 & 1.00 \\ 
     		\hline
     	\end{tabular}
     	\vspace{0.1in}
     	\caption{Proportion of times the null hypothesis was rejected across different values of correlations for the multivariate log-normal setting (TG). In this case, EN (in bold) has the best empirical performance.}
     	\label{table:twosamnormal}
     \end{table}
     \begin{table}[h]
     	\centering
     	\begin{tabular}{| *{9}{Sc|}}
     		\hline
     		Mean & \multicolumn{2}{c|}{(CMT)} & \multicolumn{2}{c|}{(EN)} & \multicolumn{2}{c|}{(HHG)} & \multicolumn{2}{c|}{\textbf{(REN)}} \\ 
     		\hline 
     		-0.80 & 0.53 & 0.65 & 1.00 & 1.00 & 1.00 & 1.00 & \textbf{1.00} & \textbf{1.00} \\ 
     		\hline
     		-0.6 & 0.30 & 0.39 & 0.98 & 0.99 & 0.98 & 0.99 & \textbf{0.99} & \textbf{1.00} \\ 
     		\hline
     		-0.52 & 0.18 & 0.27 & 0.93 & 0.96 & 0.95 & 0.98 & \textbf{0.99} & \textbf{1.00} \\ 
     		\hline
     		-0.45 & 0.16 & 0.25 & 0.85 & 0.91 & 0.89 & 0.94 & \textbf{0.95} & \textbf{0.98} \\ 
     		\hline 
     		-0.37 & 0.11 & 0.18 & 0.69 & 0.79 & 0.74 & 0.84 & \textbf{0.83} & \textbf{0.90} \\ 
     		\hline
     		-0.30 & 0.10 & 0.16 & 0.49 & 0.59 & 0.53 & 0.64 & \textbf{0.61} & \textbf{0.73} \\ 
     		\hline
     		-0.22 & 0.08 & 0.12 & 0.31 & 0.41 & 0.32 & 0.44 & \textbf{0.35} & \textbf{0.48} \\ 
     		\hline
     		-0.15 & 0.06 & 0.11 & 0.15 & 0.24 & 0.17 & 0.26 & \textbf{0.17} & \textbf{0.26} \\ 
     		\hline
     		-0.07 & 0.07 & 0.11 & 0.07 & 0.15 & 0.09 & 0.15 & \textbf{0.08} & \textbf{0.16} \\ 
     		\hline
     		0.07 & 0.05 & 0.09 & 0.07 & 0.13 & 0.07 & 0.14 & \textbf{0.08} & \textbf{0.14} \\ 
     		\hline
     		0.15 & 0.06 & 0.10 & 0.16 & 0.24 & 0.16 & 0.26 & \textbf{0.17} & \textbf{0.28} \\ 
     		\hline
     		0.22 & 0.07 & 0.12 & 0.28 & 0.41 & 0.30 & 0.44 & \textbf{0.36} & \textbf{0.49} \\ 
     		\hline
     		0.30 & 0.09 & 0.14 & 0.50 & 0.62 & 0.53 & 0.66 & \textbf{0.60} & \textbf{0.72} \\ 
     		\hline
     		0.37 & 0.12 & 0.17 & 0.70 & 0.80 & 0.74 & 0.83 & \textbf{0.82} & \textbf{0.88} \\ 
     		\hline
     		0.45 & 0.16 & 0.25 & 0.86 & 0.92 & 0.88 & 0.93 & \textbf{0.95} & \textbf{0.98} \\ 
     		\hline
     		0.52 & 0.22 & 0.29 & 0.95 & 0.97 & 0.96 & 0.98 & \textbf{0.98} & \textbf{0.99} \\ 
     		\hline
     		0.60 & 0.27 & 0.38 & 0.98 & 0.99 & 0.99 & 0.99 & \textbf{1.00} & \textbf{1.00} \\ 
     		\hline
     	\end{tabular}
     	\vspace{0.1in}
     	\caption{Proportion of times the null hypothesis was rejected across different values of correlations for the multivariate log-normal setting (TGL). In this case, REN (in bold) has the best empirical performance.}
     	\label{table:twosamnormallog}
     \end{table}
     In Tables~\ref{table:twosamnormal} and~\ref{table:twosamnormallog}, we present the estimated powers (using $1000$ independent replicates) for CMT, HHG, EN and RE, at levels $0.05$ and $0.1$. Note that, for both settings and, $\Ren$ convincingly outperforms Rosenbaum's crossmatch test, which is the only other exactly distribution-free test for multivariate two-sample goodness-of-fit testing. In~\cref{table:twosamnormal}, the energy test has the best empirical performance, closely followed by $\Ren$. In this setting, $\Ren$ consistently outperforms HHG, rather convincingly for small values of $|\mu|$. Next, for setting $\Ren$ has the best performance as it outperforms all its competitors. Note that is a somewhat heavy-tailed setting, although the associated distribution has all moments finite (the exponential moment becomes infinite). The superior performance of $\Ren$ in this case further strengthens our belief that $\Ren$ can provide a better multivariate goodness-of-fit test in heavy-tailed settings.
     \subsection{Connection with distance covariance for bivariate normal distribution}\label{sec:bivdcovcon}
     In~\cite[Theorem 7]{Gabor2007}, the authors explicitly calculate the population distance correlation when $(X,Y)$ follows a bivariate Gaussian distribution with mean vector $\mathbf{0}$, variances $1$ and correlation $\rho$. Their result interestingly shows that the population distance correlation is a strictly increasing function of $|\rho|$. In this subsection, the goal is to inspect if analogous properties hold for the rank distance correlation, defined below.
     \begin{definition}[Rank distance correlation]\label{def:rdcorr}
     	The rank distance correlation $(\Rdcorr)$ between $\mathbf{Z}_1$ and $\mathbf{Z}_2$ is defined as the usual distance correlation (see~\cite[Equation 2.7]{Gabor2007}) between $R_1(\mathbf{Z}_1)$ and $R_2(\mathbf{Z}_2)$ (where $R_1(\cdot)$ and $R_2(\cdot)$ are as in~\cref{def:rdcov}). In other words,
     	\begin{equation}\label{eq:rdcorr1}
     	\Rdcorr^2(\mathbf{Z}_1,\mathbf{Z}_2)\coloneqq\frac{\Rdcov^2(\mathbf{Z}_1,\mathbf{Z}_2)}{\Rdcov(\mathbf{Z}_1,\mathbf{Z}_1)\Rdcov(\mathbf{Z}_2,\mathbf{Z}_2)}.
     	\end{equation}
     \end{definition}
     Note that~\eqref{eq:rdcorr1} is well-defined by~\cref{lem:proprdcov} (part (c)).  By~\cite[Theorem 3]{Gabor2007}, it follows directly that $\Rdcorr(\mathbf{Z}_1,\mathbf{Z}_2)\in [0,1]$. For our proposed rank distance correlation measure (population version) above, however, a closed form expression is rather difficult to obtain when $(X,Y)$ follows a bivariate normal distribution. But, note that the rank maps corresponding to the marginals $X$ and $Y$ are exactly the Gaussian cumulative distribution functions (as discussed in~\cref{sec:revOT}). Therefore, we can use Monte-Carlo approximations with this known rank map and obtain the population rank correlation measure. The subsequent plot reveals something interesting: 
     \begin{figure}[H]
     	\begin{center}
     		\includegraphics[height=7cm,width=10cm]{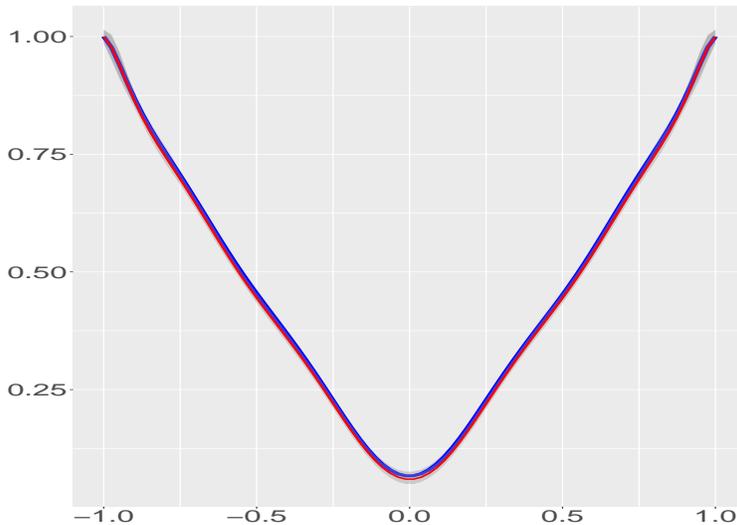}	
     		\caption{The plot depicts the population distance correlation and the population rank distance correlation  ($y$-axis) in blue and red respectively, as a function of $\rho$ ($x$-axis) in the bivariate Gaussian setting.}
     		\label{fig:dcovvsrankdcov}
     	\end{center}
     \end{figure}	
     Note that~\cref{fig:dcovvsrankdcov} shows (numerically) that the rank distance correlation is also an increasing function of $|\rho|$. The most striking feature about~\cref{fig:dcovvsrankdcov} is the proximity between the usual distance correlation and the rank distance correlation in the bivariate Gaussian case. As of now, we do not have any explanation for this striking feature. We believe that exploring the connections between rank and usual distance correlation could be an interesting area for future research.
     
     \subsection{Universal asymptotic cutoffs}\label{sec:univcutas}
     Recall that the nonparametric tests proposed in~\cref{sec:mtp} have exact distribution-freeness. Therefore, it is natural to ask: ``under the null hypothesis, how fast does this fixed sequence of finite sample distributions converge to its asymptotic limit?". This question is of immense practical interest. This would help practitioners avoid using different cutoffs for different sample sizes provided the sample size is ``large enough". As of now we do not have a theoretical answer to this question. In this section therefore, we attempt to answer this question through numerical experiments.\par 
     Let us start with $n\Rdcov_n^2$. Below we present the $p$-values from a two-sample Kolmogorov-Smirnov test for equality of distributions between $n\Rdcov_n^2$ and $1000\Rdcov_{1000}^2$ as $n$ varies, for $d_1=d_2=2$ and $d_1=d_2=8$ respectively.
     \begin{table}[h]
     	\centering
     	\begin{tabular}{| *{6}{Sc|}}
     		\hline
     		& \multicolumn{1}{c|}{(100)} & \multicolumn{1}{c|}{(300)} & \multicolumn{1}{c|}{(500)} & \multicolumn{1}{c|}{(700)} & \multicolumn{1}{c|}{(900)} \\ 
     		\hline
     		\mbox{$p$-value} & 0.09 & 0.98 & 0.78 & 0.14 & 0.94 \\ 
     		\hline
     	\end{tabular}
     	\vspace{0.1in}
     	\caption{$P$-values of two-sample Kolmogorov-Smirnov tests when comparing $n=100,300,500,700,900$ with $n=1000$, for $d_1=d_2=2$.}
     	\label{table:KS1}
     \end{table}
     
     \vspace{0.1in}
     
     \begin{table}[h]
     	\centering
     	\begin{tabular}{| *{6}{Sc|}}
     		\hline
     		& \multicolumn{1}{c|}{(100)} & \multicolumn{1}{c|}{(300)} & \multicolumn{1}{c|}{(500)} & \multicolumn{1}{c|}{(700)} & \multicolumn{1}{c|}{(900)} \\ 
     		\hline
     		\mbox{$p$-value} & 0.00 & 0.00 & 0.00 & 0.13 & 0.62 \\ 
     		\hline
     	\end{tabular}
     	\vspace{0.1in}
     	\caption{$P$-values of two-sample Kolmogorov-Smirnov tests when comparing $n=100,300,500,700,900$ with $n=1000$, for $d_1=d_2=8$.}
     	\label{table:KS2}
     \end{table}
     Tables~\ref{table:KS1} and~\ref{table:KS2} show that for $d_1=d_2=2$, the distribution of $n\Rdcov_n^2$ stabilizes by $n=300$ whereas for $d_1=d_2=8$, the distribution of $n\Rdcov_n^2$, as expected, takes longer to stabilize, potentially at around $n=700$. Note that, in practice, it may be more useful to see when the critical values at $5\%$ and $10\%$ levels stabilize as that is the question practitioners may be more interested in. Tables~\ref{table:threshold1} and~\ref{table:threshold2} provide the $5\%$ and $10\%$ cutoffs for the distribution of $n\Rdcov_n^2$ as $n$ varies, for $d_1=d_2=2$ and $d_1=d_2=8$ respectively.
     \begin{table}[h]
     	\centering
     	\begin{tabular}{| *{6}{Sc|}}
     		\hline
     		& \multicolumn{1}{c|}{(100)} & \multicolumn{1}{c|}{(300)} & \multicolumn{1}{c|}{(500)} & \multicolumn{1}{c|}{(700)} & \multicolumn{1}{c|}{(900)} \\ 
     		\hline
     		0.05 & 0.39 & 0.40 & 0.39 & 0.40 & 0.40 \\ 
     		\hline
     		0.1 & 0.36 & 0.36 & 0.36 & 0.36 & 0.36 \\ 
     		\hline
     	\end{tabular}
     	\vspace{0.1in}
     	\caption{Thresholds for $\alpha=0.05$, $0.1$ and $n=100,300,500,700,900$.}
     	\label{table:threshold1}
     \end{table}
     \begin{table}[h]
     	\centering
     	\begin{tabular}{| *{6}{Sc|}}
     		\hline
     		& \multicolumn{1}{c|}{(100)} & \multicolumn{1}{c|}{(300)} & \multicolumn{1}{c|}{(500)} & \multicolumn{1}{c|}{(700)} & \multicolumn{1}{c|}{(900)} \\ 
     		\hline
     		0.05 & 1.37 & 1.38 & 1.38 & 1.38 & 1.38 \\ 
     		\hline
     		0.1 & 1.34 & 1.35 & 1.35 & 1.35 & 1.35 \\ 
     		\hline
     	\end{tabular}
     	\vspace{0.1in}
     	\caption{Thresholds for $\alpha=0.05$, $0.1$ and $n=100,300,500,700,900$.}
     	\label{table:threshold2}
     \end{table}
     These tables show that the $5\%$ and $10\%$ quantiles probably stabilize much faster, by sample size $n=300$ even when $d_1=d_2=8$. The same observation recurs for other choices of $d_1,d_2$ if both are less then or equal  to 8. This allows us to provide universal asymptotic cutoffs as long as $d_1,d_2\leq 8$ (see~\cref{table:thresholdind}). One can of course make such statements for other choices of $d_1$ and $d_2$.\par 
     \begin{table}[h]
     	\centering
     	\begin{tabular}{|r|r|r|r|r|r|r|r|r|}
     		\hline
     		& $d_2=1$ & $d_2=2$ & $d_2=3$ & $d_2=4$ & $d_2=5$ & $d_2=6$ & $d_2=7$ & $d_2=8$ \\ 
     		\hline
     		$d_1=1$ & 0.23 & 0.30 & 0.35 & 0.38 & 0.41 & 0.44 & 0.47 & 0.49 \\ 
     		\hline
     		$d_1=2$ & 0.30 & 0.40 & 0.47 & 0.53 & 0.58 & 0.63 & 0.66 & 0.70 \\ 
     		\hline
     		$d_1=3$ & 0.35 & 0.47 & 0.56 & 0.63 & 0.70 & 0.76 & 0.81 & 0.86 \\ 
     		\hline
     		$d_1=4$ & 0.38 & 0.53 & 0.63 & 0.72 & 0.80 & 0.87 & 0.93 & 0.99 \\ 
     		\hline
     		$d_1=5$ & 0.41 & 0.58 & 0.70 & 0.80 & 0.89 & 0.97 & 1.04 & 1.10 \\ 
     		\hline
     		$d_1=6$ & 0.44 & 0.63 & 0.76 & 0.87 & 0.97 & 1.05 & 1.13 & 1.20 \\ 
     		\hline
     		$d_1=7$ & 0.47 & 0.66 & 0.81 & 0.93 & 1.04 & 1.13 & 1.21 & 1.30 \\ 
     		\hline
     		$d_1=8$ & 0.49 & 0.70 & 0.86 & 0.99 & 1.10 & 1.20 & 1.30 & 1.38 \\ 
     		\hline
     	\end{tabular}
     	\vspace{0.1in}
     	\caption{Asymptotic thresholds for $n\Rdcov_n^2$ when $\alpha=0.05$ and $d_1,d_2\leq 8$.}
     	\label{table:thresholdind}
     \end{table}
     We observe an exactly similar phenomenon for $mn(m+n)^{-1}\Ren_{m,n}^2$. We provide the asymptotic cutoffs corresponding to our proposed test for $1\leq d\leq 8$ in~\cref{table:twosamascutoff}.
     \begin{table}[h]
     	\centering
     	\begin{tabular}{rrrrrrrrr}
     		\hline
     		& 1 & 2 & 3 & 4 & 5 & 6 & 7 & 8 \\ 
     		\hline
     		1 & 0.94 & 1.12 & 1.26 & 1.37 & 1.45 & 1.54 & 1.61 & 1.67 \\ 
     		\hline
     	\end{tabular}
     	\vspace{0.1in}
     	\caption{Asymptotic thresholds for $mn(m+n)^{-1}\Ren_{m,n}^2$ when $\alpha=0.05$ and $d\leq 8$.}
     	\label{table:twosamascutoff}
     \end{table}
     \section{Some additional discussion}\label{sec:addis}
     We now provide some natural extensions of the proposed testing procedures to multivariate and multi-sample settings. We will also compare our proposals for defining ranks and associated test statistics with some competing proposals.
     \subsection{Extension to multivariate multi-sample ($\geq 2$) independence and equality of distributions testing}\label{sec:multext}
     In this section we discuss how one can construct distribution-free procedures for testing: \textbf{(I)} the mutual independence of two or more random vectors, and \textbf{(II)} the equality of two or more multivariate distributions.
     
     \noindent \textbf{(i)} \textbf{Testing for mutual independence of $K$ random vectors}: Suppose that we have observed i.i.d. data $\mathbf{X}_1,\ldots ,\mathbf{X}_n$ from some distribution $\mu$ on $\mathbb{R}^d$. For each $1\leq i\leq n$, suppose that $\mathbf{X}_i=(\mathbf{X}_i^1,\ldots ,\mathbf{X}_i^K)$ where $\mathbf{X}_1^j\sim\mu_j\in\mathcal{P}_{ac}(\mathbb{R}^{d_j})$, where $K\geq 2$. We are interested in testing the following hypothesis (of mutual independence):
     \begin{equation}\label{eq:multextindmain}
     \textrm{H}_0:\mu=\mu_1\otimes \mu_2\otimes \ldots \otimes \mu_K \qquad\qquad \mbox{versus}\qquad\qquad  \textrm{H}_1: \textrm{ Not H}_0. 
     \end{equation}
     Before proposing the test, let us start with some notation. Fix $1\leq j\leq K$, $1\leq i\leq n$, and define $\mathbf{X}^{j}:=\{\mathbf{X}_1^j,\ldots ,\mathbf{X}_n^j\}$ --- the $n$ data points for the $j$'th marginal, $\mathbf{X}_i^{j+}\coloneqq (\mathbf{X}_i^{j+1},\ldots ,\mathbf{X}_i^K)$ --- the $i$'th data point from the $(j+ 1)$'th marginal variable, and $\mathbf{X}^{j+}:=\{\mathbf{X}_1^{j+},\ldots ,\mathbf{X}_n^{j+}\}$ --- the collection of all sub-vectors from the $j$'th marginal. As before, define $\mathcal{H}_n^{d_j}$ to be the (fixed) sample multivariate ranks, $j=1,2,\ldots ,K$. Also let $\hat{R}_n^j(\cdot)$ denote the empirical rank map which transports the empirical distribution on $\mathbf{X}^j$ to that on $\mathcal{H}_n^{d_j}$ (see~\eqref{eq:empopt}), for $j=1,\ldots, K$. Finally, set $\hat{R}_n(\mathbf{X}_i^{j+}):=\big(\hat{R}_n^{j+1}(\mathbf{X}_i^{j+1}),\ldots, \hat{R}_n^{K}(\mathbf{X}_i^{K})\big)$ which is a vector of dimension $d_{j+1}+d_{j+2}+\ldots +d_K$.
     
     Now, we are in a position to define our test statistic $\Rdcov_n^2$ which is given by:
     \begin{equation*}
     \Rdcov_n^2\coloneqq\sum_{j=1}^{K-1} \Rdcov_{n,j}^2
     \end{equation*}
     where 
     \begin{align}\label{eq:multextind1}
     \Rdcov_{n,j}^2&\coloneqq \left(\frac{1}{n^2}\sum_{k,l=1}^n \lVert \hat{R}_n^{j}(\mathbf{X}_k^j)-\hat{R}_n^{j}(\mathbf{X}_l^j)\rVert\right)\times\left(\frac{1}{n^2}\sum_{k,l=1}^n \lVert \hat{R}_n^{j+}(\mathbf{X}_k^{j+})-\hat{R}_n^{j+}(\mathbf{X}_l^{j+})\rVert\right)\nonumber \\&\qquad \qquad \qquad +
     \frac{1}{n^2}\sum_{k,l=1}^n \lVert \hat{R}_n^{j}(\mathbf{X}_k^j)-\hat{R}_n^{j}(\mathbf{X}_l^j)\rVert\lVert \hat{R}_n^{j+}(\mathbf{X}_k^{j+})-\hat{R}_n^{j+}(\mathbf{X}_l^{j+})\rVert\nonumber \\ &\qquad\qquad-\frac{2}{n^3}\sum_{k,l,i=1}^n \lVert \hat{R}_n^{j}(\mathbf{X}_k^j)-\hat{R}_n^{j}(\mathbf{X}_l^j)\rVert\lVert \hat{R}_n^{j+}(\mathbf{X}_k^{j+})-\hat{R}_n^{j+}(\mathbf{X}_i^{j+})\rVert.
     \end{align}
     Note that the right hand side of~\eqref{eq:multextind1} can be interpreted as the usual distance covariance of the transformed data points $\big\{(\hat{R}_n^j(\mathbf{X}_i^j),\hat{R}_n^{j+}(\mathbf{X}_i^{j+})): 1\leq i\leq n\big\}$, for $j=1,\ldots, K-1$.
     
     We \emph{reject} the null hypothesis in~\eqref{eq:multextindmain} if $n\Rdcov_n^2$ is larger than $c_n$ where $$c_n\coloneqq \inf\{c:\mathbb{P}_{\textrm{H}_0}(n\Rdcov_n^2\geq c)\leq \alpha\};$$ here $\alpha$ is the usual prespecified type-I error level.~\cref{prop:multextind} below shows that $c_n$ can be obtained in a distribution-free fashion just as in the $K=2$ case (see~\cref{sec:Indep}); it also proves the consistency of the proposed test (see~\cref{proof:multextind} for the proof).
     \begin{prop}\label{prop:multextind}
     	If the marginals $\mu_j$'s are Lebesgue absolutely continuous, then, under $\textrm{H}_0$, $\Rdcov_n^2$ is a distribution-free statistic. Additionally, if the empirical distributions on $\mathcal{H}_n^{d_j}$ converge weakly to $\mathcal{U}^{d_j}$, for $j=1,2,\ldots ,K$, then the above test which rejects $\textrm{H}_0$ when $n\Rdcov_n^2\geq c_n$ is consistent.
     \end{prop}
     This extension from the case $K=2$ (in~\cref{sec:Indep}) to $K\geq 2$ (in this section) is based on the same framework as laid out in~\cite{matteson2014}.

     \noindent \textbf{(ii)} \textbf{Testing for $K$-sample equality of distributions}: Suppose that we have observed independent data $\{\mathbf{X}_i^j: 1\leq j\leq K, 1\leq i\leq n_j\}$, where $\mathbf{X}_1^j,\ldots ,\mathbf{X}_{n_j}^j\sim \mu_j\in\mathcal{P}_{ac}(\mathbb{R}^d)$, $\sum_{j=1}^K n_j=n$, and $K \ge 2$. In this setting, we are interested in testing the following hypothesis:
     \begin{equation}\label{eq:multextgofmain}
     \textrm{H}_0:\mu_1= \mu_2=\ldots = \mu_K\qquad\qquad \mbox{versus}\qquad\qquad  \textrm{H}_1: \textrm{ Not H}_0. 
     \end{equation}
     Let $\hat{R}_n(\cdot)$ denote the empirical rank map that transports the empirical distribution on the pooled sample (all $\mathbf{X}_i^j$'s taken together) to that on $\mathcal{H}_n^d$ (fixed sample multivariate ranks). We then construct our test statistic as:
     \begin{align*}
     \Ren_{1:K,n}^2\coloneqq \sum\limits_{j=1}^{K-1} \Ren_{j:(j+1),n}^2 
     \end{align*}
     where 
     \begin{align}\label{eq:multextgof1}
     \Ren_{j:(j+1),n}^2&\coloneqq \frac{2}{n_j n_{j+1}}\sum\limits_{m=1}^{n_j} \sum\limits_{l=1}^{n_{j+1}} \lVert \hat{R}_n(\mathbf{X}_m^{j})-\hat{R}_{n}(\mathbf{X}_l^{j+1})\rVert-\frac{1}{n_j^2}\sum\limits_{m,l=1}^{n_j} \lVert \hat{R}_n(\mathbf{X}_m^{j})-\hat{R}_n(\mathbf{X}_l^j)\rVert\nonumber \\&\qquad\qquad-\frac{1}{n_{j+1}^2}\sum\limits_{m,l=1}^{n_{j+1}} \lVert\hat{R}_n(\mathbf{X}_m^{j+1})-\hat{R}_n(\mathbf{X}_l^{j+1})\rVert.
     \end{align}
     Note that the right hand side of~\eqref{eq:multextgof1} can be interpreted as the usual energy statistic between the transformed data points $\big\{(\hat{R}_n(\mathbf{X}_i^j): {1\leq i\leq n_j}\big\}$ and $\big\{\hat{R}_n(\mathbf{X}_i^{j+1})): 1\leq i\leq n_{j+1}\big\}$, for $j = 1,\ldots, K-1$.
     
     We \emph{reject} the null hypothesis in~\eqref{eq:multextgofmain} if $n\Ren_{1:K,n}^2$ is larger than $c_n$ where $c_n\coloneqq \inf\{c:\mathbb{P}_{\textrm{H}_0}(n\Ren_{1:K,n}^2\geq c)\leq \alpha\}$.~\cref{prop:multextgof} below shows that $c_n$ can be obtained in a distribution-free fashion and demonstrates the consistency of the proposed test (see~\cref{proof:multextgof} for the proof).
     \begin{prop}\label{prop:multextgof}
     	If the $\mu_j$'s are Lebesgue absolutely continuous, then under $\textrm{H}_0$, $\Ren_{1:K,n}^2$ is a distribution-free statistic. Additionally, if $n_j/n\to\lambda_j\in (0,1)$ such that $\sum_{j=1}^K \lambda_j=1$, and the empirical distribution on $\mathcal{H}_n^d$ converges weakly to $\mathcal{U}^d$, then the above test which rejects $\textrm{H}_0$ when $n\Ren_{1:K,n}^2\geq c_n$, is consistent.
     \end{prop}
     \subsection{Comparison with competing rank-based approaches}\label{sec:compare}
     This section is devoted to the comparisons between the approach proposed in this paper with those in other papers that define multivariate ranks using measure transportation.\par 
     We start off with~\cite{del2018center}, which in fact, was the main motivation behind our approach. Recall that we define population ranks using measure transportation from the population distribution (say $\mu$) to the uniform distribution on $[0,1]^d$ (see~\cref{sec:defrank}). In~\cite{del2018center}, the authors suggest using a ``special uniform distribution on the unit sphere" instead of our proposed uniform distribution on $[0,1]^d$. For our notion of population ranks, the following interesting property holds: given a ``nice" measure $\mu=(\mu_1,\mu_2)$ with $\mu_1$ and $\mu_2$ supported on $\mathbb{R}^{d_1}$ and $\mathbb{R}^{d_2}$ respectively, suppose $F(\cdot),F_1(\cdot),F_2(\cdot)$ denote the corresponding population rank maps, then we must have $F=(F_1\; F_2)$ if and only if $\mu_1$ and $\mu_2$ are independent. In other words, the joint measure splitting into product of marginals is equivalent to the joint rank map splitting component-wise into marginal rank maps. This property is not true if we replace the uniform distribution on $[0,1]^d$ with the distribution proposed in~\cite{del2018center}. This interesting observation allows us to obtain an elegant solution to the mutual independence testing problem, as can be seen, for example, in the proof of~\cref{lem:proprdcov}. \par In addition, the definition of empirical ranks  in~\cite{del2018center} seems to have a distinct subjective element which makes it difficult to use in real-life testing problems. In particular, it relies on a certain factorization of the sample size, $n=n_R n_S+n_0$ where roughly, for large $n$, one should choose $n_R$ and $n_S$ as large as possible and $n_0$ as small as possible. 
     Consequently, the analogue of~\cref{prop:finproperties} (part (ii)) in~\cite{del2018center} provides a fixed sample universal distribution for their rank vector which now depends on $n_0$. Therefore, all nonparametric tests based on these ranks will also have null distribution depending on the choice of this subjective parameter.

     Another notion of multivariate rank was proposed in~\cite{boeckel2018multivariate} which is quite similar to our proposal. However the authors construct the empirical ranks by replacing the $\mathbf{h}_i$'s (fixed grid Halton sequence) in~\eqref{eq:empopt} with a random draw of $n$ i.i.d. uniforms. The problem here is that any test obtained using these ranks are now random variables, given the data (due to the external randomization), i.e., the test statistic will not be deterministically determined by the observed data (although the distribution of the test statistic will be). In~\cite{boeckel2018multivariate}, the authors propose a two-sample equality of distributions test where they use the Wasserstein distance between the first and second sample ranks from the pooled dataset, instead of the energy statistic that we use. Further, their theoretical results require much stronger assumptions on the underlying distributions, e.g., the assumption that the data generating distribution be compactly supported in addition to being absolutely continuous.
     
     Two other notions of multivariate ranks have been proposed in~\cite{chernozhukov2017monge,ghosal2019multivariate} based on what is called as the \emph{semi-discrete optimal transportation problem}. Due to space constraints, we do not describe this approach in detail here. However, the joint distribution of their proposed rank vectors is not distribution-free for fixed sample sizes. The theoretical results in~\cite{chernozhukov2017monge} (including uniform convergence of rank maps) assume compactly supported distributions, whereas~\cite{ghosal2019multivariate} proves a more general result under the condition that the population rank map is a homeomorphism --- still a relatively strong assumption. 
     One of the important observations in~\cite{ghosal2019multivariate} is that under much weaker conditions, we can get a weaker (than uniform convergence) notion of convergence of empirical ranks, which is sufficient to yield consistency of the proposed testing procedures herein. 
     Note that it is also very difficult to feasibly compute the notion of multivariate ranks proposed in~\cite{chernozhukov2017monge, ghosal2019multivariate} beyond $d\geq 4$.
     
     \subsection{Quasi-Monte Carlo methods}\label{sec:revQMC}
     Define $\mathfrak{F} :=\{f:[0,1]^d\mapsto\mathbb{R} | \int_{[0,1]^d} f^2(\mathbf{x})\,d\mathbf{x}<\infty, f\in C^2([0,1]^d), f\mbox{ has bounded HKV}\}$; here HKV stands for \textit{Hardy-Krause} variation, see~\cite{Morokoff1995} for detailed definitions. We begin with the following problem:
     \begin{equation*}
     \textbf{(P)}:\qquad \mbox{Calculate/Approximate }\int_{[0,1]^d} f(\mathbf{x})\,d\mathbf{x},\qquad \mbox{ for }f\in\mathfrak{F}.
     \end{equation*}
     Problem \textbf{(P)} has several applications in engineering, astronomy, mathematical and computational finance, and has consequently been studied extensively across various disciplines (for a book length treatment see~\cite{Philip1986} and the extensive list of references therein). In most applications, difficulties in solving \textbf{(P)} arise from two main sources --- (a) the antiderivative of $f(\cdot)$ may be hard to compute or it may not exist, and, (b) even a closed form expression for $f(\cdot)$ may not be available and the only information accessible to the analyst might be through evaluations of $f(\cdot)$ at certain query points of choice. Accordingly, the standard approach towards solving \textbf{(P)} is by choosing a set of query points $\textbf{x}^{(n)}=\{\textbf{x}_1,\ldots ,\textbf{x}_n\}$, $\mathbf{x}_i\in [0,1]^d$ and combining the evaluations $\{f(\mathbf{x}_i)\}_{1\leq i\leq n}$ by simple or weighted averaging. For the subsequent discussion, given any approximation $T_n\big(f;\mathbf{x}^{(n)}\big)$ to \textbf{(P)} based on the set of query points $\textbf{x}^{(n)}$, we define the \textit{approximation error} ($\epsilon_{n,f}$) as:
     \begin{equation*}
     \epsilon_{n,f}\coloneqq \Bigg|\int_{[0,1]^d} f(\mathbf{x})\,d\mathbf{x}-T_n\big(f;\mathbf{x}^{(n)}\big)\Bigg|.
     \end{equation*}
     \par One of the earliest approaches in literature was through the use of product numerical integration techniques. The principle idea is as follows: Suppose that $n=(m+1)^d$ for some $m\in\mathbb{N}\cup \{0\}$ and divide the $[0,1]$ interval along each of the coordinate axes into $m$ subintervals (of equal width). Consider the $d$ sets comprising of the endpoints from these intervals and form $\mathbf{x}^{(n)}$ by taking the Cartesian product of these $d$ sets. The collection of evaluations $\{f(\mathbf{x}_i)\}_{1\leq i\leq n}$ is then combined using a suitably chosen weighted averaging scheme. This idea forms the foundation of popular techniques such as the Simpson's rule, the Trapezoidal rule, etc. and leads to an \textit{approximation error} of $O\big(n^{-2/d}\big)$ in general (see~\cite{Atkinson1989,Uribe2002}). Note that these \textit{approximation errors} suffer from the ``curse of dimensionality" and become less useful for even moderately large $d$. \par
     A significant step in resolving this ``curse of dimensionality" problem was achieved through the development of \textit{Monte Carlo} methods (see~\cite{robert2013monte,hastings1970monte} and the references therein). The main idea here is to construct $\mathbf{x}^{(n)}$ using a random sample of $n$ vectors according to the uniform  distribution on $[0,1]^d$. The evaluations are then combined to form $T_n(f;\mathbf{x}^{(n)})$ as $n^{-1}\sum_{i=1}^n f(\mathbf{x}_i)$. In this case, $\epsilon_{n,f}$ (now a random variable) is $O_p\big(n^{-1/2}\big)$ (which follows from simple standard moment calculations). This dimension-free bound comes at a price. The error rates being random, there is no way to analyze the quality of approximation for a particular realization of $\mathbf{x}^{(n)}$. In certain sensitive practical problems, this can be a rather serious issue. \par
     In subsequent research, the utility of \textit{Monte Carlo} methods triggered this interesting idea: if a random sample corresponds to an average \textit{approximation error} bound of $O_p\big(n^{-1/2}\big)$, then there must be deterministic sequences which lead to error bounds of at most $O\big(n^{-1/2}\big)$. This served as the main motivation behind the long line of work, now referred to as \textit{Quasi-Monte Carlo} methods, which aims at finding fixed sequences which perform at least as well as \textit{Monte Carlo} methods for solving \textbf{(P)} (uniformly in $f(\cdot)$). Before delving further into this area, let us bring in some relevant notation and results. Define $\mathbb{J}\coloneqq\{\prod_{j=1}^d [0,u_j): (u_1,\ldots ,u_d)\in [0,1)^d\}$ and the \textit{star discrepancy} of a set of vectors $\mathbf{x}^{(n)}$ as:
     \begin{equation*}
     \mathcal{D}_n^*\big(\mathbf{x}^{(n)}\big)\coloneqq \sup_{J\in\mathbb{J}}\Bigg|\lambda_d(J)-\frac{1}{N}\cdot\#\{i:\mathbf{x}_i\in J\}\Bigg|
     \end{equation*}
     where $\lambda_d(J)$ denotes the Lebesgue measure of the set $J$. The following is a classical result in Quasi-Monte Carlo literature.
     \begin{prop}[Koksma-Hlawka inequality, see~\cite{Hlawka1961}]\label{prop:Hlawka}
     	Suppose $f\in\mathfrak{F}$ and let $V(f)$ denote the HKV of $f(\cdot)$. Then, for any $\mathbf{x}^{(n)}\coloneqq \{\mathbf{x}_1,\ldots ,\mathbf{x}_n\}$ with $\mathbf{x}_i\in [0,1)^d$, we have:
     	\begin{equation*}
     	\Bigg|\int_{[0,1]^d} f(\mathbf{x})\,d\mathbf{x}-\frac{1}{N}\sum_{i=1}^n f(\mathbf{x}_i)\Bigg|\leq V(f)\mathcal{D}_n^*\big(\mathbf{x}^{(n)}\big).
     	\end{equation*}
     \end{prop}
     \cref{prop:Hlawka} above implies that an upper bound to the star discrepancy of a sequence yields an upper bound to the \textit{approximation error}. As a result, constructing  sequences with small values of the star discrepancy has become a popular approach to addressing \textbf{(P)} over the past $50$ years or so. Some examples include \textit{Kronecker sequences}, \textit{Halton sequences} and digital sequences of \textit{Niederreiter type} (see~\cite{halton1964algorithm,Kuipers1974,Niederreiter1992}). For the time being, we will focus our attention on \textit{Halton sequences}.
     
     \noindent {\bf Halton sequences}: Given $n,d$, for $1\leq k\leq n$, $1\leq j\leq d$, form $x(k,j)$ as follows: let $p_i$ denote the $i^{th}$ prime and express $k$ in base $p_i$. Invert this $p_i$-ary expansion and place it after the decimal point. For example, in the case of $x(6,1)$, note that $6=(110)_2$, which makes $x(6,1)=(0.011)_2=3/8$. Finally, the \textit{Halton sequence} is formed by setting $\mathbf{x}^{(n)}\coloneqq \{\mathbf{x}_1,\ldots ,\mathbf{x}_n\}$ where $\mathbf{x}_i\coloneqq (x(i,1),\ldots ,x(i,d))\in [0,1]^d$. The \textit{star discrepancy} for these sequences is $O\big(\log^d(n)/n\big)$; see e.g.,~\cite{Hofer2009,Hofer2010}.

     \begin{corollary}\label{cor:weakcon}
     	Let $\mathbb{P}_n$ denote the empirical measure on the $d$-dimensional \textit{Halton sequence} of $n$ vectors. Then $\mathbb{P}_n\overset{w}{\longrightarrow}\mathcal{U}^d$. A similar result holds for the regular griding scheme used in product numerical integration or the random sampling scheme used in \textit{Monte Carlo} methods.
     \end{corollary}
     The above discussion is a very brief overview of this active research area. For the purposes of this paper it is perhaps instructive to filter out some key desirable properties of  Halton sequences which will be useful to us.
     \begin{table}[h]
     	\begin{tabular}{|p{40mm}|p{35mm}|p{35mm}|p{35mm}|} 
     		\hline
     		{\bf Properties}    & Regular grid & Random sampling & Halton sequences
     		\\
     		\hline
     		Sequence is deterministic & \textcolor{blue}{Yes} & No & \textcolor{blue}{Yes} \\
     		\hline
     		Can be constructed for all $n,d$   &   No, only when $n$ is a perfect power of $d$  &   \textcolor{blue}{Yes}  &   \textcolor{blue}{Yes}   \\
     		\hline
     		Rate of star discrepancy & $n^{-2/d}$ & $n^{-1/2}$ & \textcolor{blue}{$\log^d(n)/n$} \\
     		\hline
     		For fixed $d$, whole sequence must be recomputed if $n$ increases by $1$ & Yes & \textcolor{blue}{No, only one vector needs to be computed} & \textcolor{blue}{No, only one vector needs to be computed}\\ 
     		\hline
     	\end{tabular}
     	\newline
     	\caption{Properties of sequences discussed above. The most desirable properties are highlighted in \textcolor{blue}{blue}.}
     	\label{table:prop}
     \end{table}
     Although we restricted ourselves to the function class $\mathfrak{F}$ in this section, it is possible to consider more general function classes. However these technical details play no role in this paper, and hence we do not elaborated on these issues.

     \section{Proofs}\label{sec:appen}
     This section contains proofs of our main results.
     \subsection{Proof of~\cref{prop:finproperties}}\label{proof:finproperties}
     \begin{proof}
     	\textbf{(i)} This part of the proof is verbatim identical to~\cite[Proposition 1.6.1]{del2018center}.
     	
     	\noindent\textbf{(ii)} First, recall the definition of $\hat{\sigma}_n$ from~\eqref{eq:empopt}. Next, note that $\mathbf{X}_1,\mathbf{X}_2,\ldots ,\mathbf{X}_n$ are exchangeable (as they are i.i.d.). Moreover, $\mathbf{X}_i$'s are absolutely continuous. This implies that, given any two permutations $\sigma_1$, $\sigma_2$ in $S_n$, we have:
     	\begin{equation}\label{eq:finproperties1}
     	\sum_{i=1}^n \lVert \mathbf{X}_i-\mathbf{h}_{\sigma_1(i)}\rVert^2\neq \sum_{i=1}^n \lVert \mathbf{X}_i-\mathbf{h}_{\sigma_2(i)}\rVert^2 \qquad \mu\mbox{-a.e.}
     	\end{equation}
     	Moreover, the following set of $n!$ random variables,
     	$\left\{\sum_{i=1}^n \lVert \mathbf{X}_i-\mathbf{h}_{\sigma(i)}\rVert^2\right\}_{\sigma\in S_n}$ forms an exchangeable collection. Coupled with~\eqref{eq:finproperties1} the above observation yields that $\mathbb{P}(\hat{\sigma}_n=\sigma)=(n!)^{-1}$ for all $\sigma\in S_n$. This completes the proof.\par 
     	\noindent\textbf{(iii)} The independence between the multivariate ranks and the multivariate order statistics is now a direct consequence of the observation that the order statistics are complete sufficient (from part \textbf{(i)}), the multivariate ranks are ancillary for $\mu$ (as their distribution is free of $\mu$ from part \textbf{(ii)}), followed by an application of Basu's Theorem. 	\end{proof}
     
     \subsection{Proof of~\cref{lem:L2}}\label{proof:lemL-2}
     \begin{proof}
     	This proof requires a number of existing results from convex analysis which we present in~\cref{sec:auxi} for the convenience of the reader.\par 
     	Recall that $\mu_n$ is defined as the empirical measure on $\mathbf{\mathcal{D}}_n=\{\mathbf{X}_1,\ldots ,\mathbf{X}_n\}$ (see~\eqref{eq:mu_nu_n}) where $\mathbf{X}_1,\ldots ,\mathbf{X}_n\overset{i.i.d.}{\sim}\mu$. Also let all the random variables, i.e., $\mathbf{X}_i$'s be defined on the probability space $(\tilde{\Omega},\mathcal{A},\mathbb{P})$. We know that $\mu_n\overset{w}{\longrightarrow}\mu$ a.e. (see~\cite[Theorem 1]{Varadarajan1958}). In other words, there exists $\Omega \subset \tilde{\Omega}$, such that $\mathbb{P}(\Omega)=1$ and, for all $\omega\in\Omega$, $\mu_n({\omega})\overset{w}{\longrightarrow}\mu$. 
     	
     	Fix  $\omega\in\Omega$. Now consider the sequence of random measures $(\mbox{identity},\hat{R}_n)\#\mu_n({\omega})$; the randomness comes the random draw according to $\mu_n({\omega})$. As $\hat{R}_n({\omega})$ takes values in a compact set and $\mu_n({\omega})$ converges weakly, we have that $(\mbox{identity},\hat{R}_n)\#\mu_n({\omega})$ is asymptotically tight. Consequently, by Prokhorov's theorem, every subsequence of $(\mbox{identity},\hat{R}_n)\#\mu_n({\omega})$ has a further subsequence, say $(\mbox{identity},\hat{R}_{n_k})\#\mu_{n_k}({\omega})$ (by a relabeling if necessary) such that $(\mbox{identity},\hat{R}_{n_k})\#\mu_{n_k}({\omega})$ converges weakly on $\mathbb{R}^d\times\mathbb{R}^d$. Let us call this limit $\gamma(\omega)$ (which could also depend on the subsequence $\{n_k\}_{k\ge 1}$). Next, we will show that each of these subsequences have the same weak limit (which for now may depend on $\omega$).

     	Let $\Gamma(\mu,\mathcal{U}^d)$ be the family of probability distributions on $\mathbb{R}^d\times\mathbb{R}^d$ that have marginals (on the first and the last $d$-coordinates) $\mu$ and $\mathcal{U}^d$ respectively. Also, note that, by~\eqref{eq:empopt},~\eqref{eq:R_n_h} and~\cref{def:cycmon}, $(\mbox{identity},\hat{R}_{n_k})\#\mu_{n_k}({\omega})$ has cyclically monotone support. As $\mu_{n_k}({\omega})\overset{w}{\longrightarrow}\mu$ and $\hat{R}_{n_k}\#\mu_{n_k}({\omega})$ (which is simply the empirical measure on $\mathbf{\mathcal{H}}_{n_k}^d$) converges weakly to $\mathcal{U}^d$ (by assumption), we can conclude that $\gamma({\omega})$ has cyclically monotone support and $\gamma({\omega})\in\Gamma(\mu,\mathcal{U}^d)$ (by~\cite[Lemma 9]{Mccann1995}). Moreover, by~\cite[Corollary 14]{Mccann1995} and the fact that $\mu\in\mathcal{P}_{ac}(\mathbb{R}^d)$, there exists only one measure with cyclically monotone support in $\Gamma(\mu,\mathcal{U}^d)$. Therefore, irrespective of the subsequence $\{n_k\}_{k\geq 1}$, $(\mbox{identity},\hat{R}_{n_k})\#\mu_{n_k}({\omega})$ converges weakly to the same limit.

     	The above also shows that the weak limit is the same for all $\omega\in\Omega$; let us call it $\gamma$. Finally, by~\cite[Proposition 10]{Mccann1995} and the definition of $R(\cdot)$, we get that $\gamma=(\mbox{identity},R)\#\mu$. Therefore, we have proved that 
     	\begin{align}\label{eq:L2-1}
     	(\mbox{identity},\hat{R}_n)\#\mu_n({\omega})\overset{w}{\longrightarrow} (\mbox{identity},R)\#\mu\;\;\;\;\; \mbox{ for all}\;\; \omega\in\Omega.
     	\end{align}
     	Now, let $M_n$ be sampled according to the measure $\mu_n(\omega)$ and $M$ be sampled according to the measure $\mu$. Then~\eqref{eq:L2-1} can be restated as:
     	\begin{equation}\label{eq:L2-4}
     	(M_n,\hat{R}_n(M_n))\overset{w}{\longrightarrow} (M,R(M))\qquad \mbox{ for all}\;\; \omega\in\Omega.
     	\end{equation}
     	Let $g:\mathbb{R}^d\times \mathbb{R}^d\to [0,\infty)$ be defined as $g(x,y)\coloneqq\lVert y-R(x)\rVert$. Note that by Alexandroff theorem (see e.g.,~\cite{Alexandroff1939}), $R(\cdot)$ is continuous Lebesgue a.e., and consequently $\mu$-a.e.~(by the absolute continuity of $\mu$). Therefore the function  $g(\cdot)$ is discontinuous on a set which has measure 0 with respect to $(\mbox{identity},R)\#\mu$. Consequently, by applying the continuous mapping theorem with $g(\cdot)$ on~\eqref{eq:L2-4}, we get, for all $\omega\in\Omega$, 
     	\begin{align*}
     	g(M_n,\hat{R}_n(M_n)) = \lVert \hat{R}_n(M_n)-R(M_n)\rVert\overset{w}{\longrightarrow} g(M,R(M)) = 0.
     	\end{align*}
     	Finally, as $\hat{R}_n$ and $R$ are uniformly bounded, the dominated convergence theorem implies,
     	\begin{align*}
     	\int \lVert \hat{R}_n(\cdot)-R(\cdot)\rVert\,d\mu_n(\omega)=\frac{1}{n}\sum_{i=1}^n \lVert \hat{R}_n(\mathbf{X}_i)-R(\mathbf{X}_i)\rVert (\omega)\longrightarrow 0 \;\;\;\;\mbox{for all} \;\;\omega\in\Omega.
     	\end{align*}
     	This completes the proof.
     \end{proof}
     \subsection{Proof of~\cref{prop:srcc}}\label{proof:srcc}
     \begin{proof}
     	For the solution of this problem, as we are interested in limits under $|t|/|s|\to c$. Thus, let us assume without loss of generality that $|t|/|s|\in (c/2,2c)$. Note that, by a first order Taylor series expansion of the exponential function (as $G_1(Z_1)$ and $G_2(Z_2)$ are uniformly bounded),
     	\begin{align}\label{eq:srcc1}
     	F(t,s)&\coloneqq \mathbb{E}\left[\exp(itG_1(Z_1)+isG_2(Z_2))\right]\nonumber \\ &=1+it\mathbb{E}[G_1(Z_1)]+is\mathbb{E}[G_2(Z_2)]-\frac{t^2}{2}\mathbb{E}[G_1(Z_1)^2]-\frac{s^2}{2}\mathbb{E}[G_2(Z_2)^2]-st\mathbb{E}[G_1(Z_1)G_2(Z_2)]\nonumber \\ &\qquad\qquad+\mathcal{O}(\max\{|s|^3,|t|^3\}).
     	\end{align}
     	Now plugging in $t=0$, followed by $s=0$, alternatively in~\eqref{eq:srcc1} and multiplying the results, we get:
     	\begin{align*}
     	F(t,0)F(0,s)&=\mathbb{E}\left[\exp(itG_1(Z_1)\right]\,\mathbb{E}\left[isG_2(Z_2)\right]\nonumber \\ &=1+it\mathbb{E}[G_1(Z_1)]+is\mathbb{E}[G_2(Z_2)]-\frac{t^2}{2}\mathbb{E}[G_1(Z_1)^2]-\frac{s^2}{2}\mathbb{E}[G_2(Z_2)^2]\nonumber \\&\qquad\qquad-st\mathbb{E}[G_1(Z_1)]\mathbb{E}[G_2(Z_2)]+\mathcal{O}(\max\{|s|^3,|t|^3\}).
     	\end{align*}
     	This implies, following the notation from the problem statement,
     	\begin{align}\label{eq:srcc3}
     	f_{Z_1,Z_2}(t,s)&=F(t,s)-F(t,0)F(0,s)\nonumber \\ &=st\mathrm{Cov}(G_1(Z_1),G_2(Z_2))+\mathcal{O}(\max\{|s|^3,|t|^3\}).
     	\end{align}
     	In the same vein as~\eqref{eq:srcc3}, we further get the following estimates: 
     	\begin{align}\label{eq:srcc4}
     	f_{Z_1,Z_1}(t,t)=t^2\mathrm{Var}[G_1(Z_1)]+\mathcal{O}(|t|^3)\;\mbox{,}\; f_{Z_2,Z_2}(s,s)=s^2\mathrm{Var}[G_2(Z_2)]+\mathcal{O}(|s|^3).
     	\end{align}
     	By combining these observations from~\eqref{eq:srcc3} and~\eqref{eq:srcc4}, we get:
     	\begin{align*}
     	&\;\;\;\;\lim\limits_{t,s\to 0,|t|/|s|\to c}\frac{|f_{Z_1,Z_2}(t,s)|^2}{|f_{Z_1,Z_1}(t,t)||f_{Z_2,Z_2}(s,s)|}\nonumber \\ &=\lim\limits_{t,s\to 0,|t|/|s|\to c} \frac{(st\mathrm{Cov}(G_1(Z_1),G_2(Z_2))+\mathcal{O}(\max\{|t|^3,|s|^3\}))^2}{(t^2\mathrm{Var}[G_1(Z_1)]+\mathcal{O}(|t|^3))(s^2\mathrm{Var}[G_1(Z_1)]+\mathcal{O}(|s|^3))}\nonumber \\&=\rho^2(G_1(Z_1),G_2(Z_2)).
     	\end{align*}
     	This completes the proof.
     \end{proof}
     \subsection{Proof of~\cref{lem:proprdcov}}\label{proof:proprdcov}
     \begin{proof}
     	{(a)} This follows directly from the properties of distance covariance in~\cite[Remark 3]{Gabor2007} by noting that rank distance covariance is just the usual distance covariance between the multivariate ranks $R_1(\mathbf{Z}_1)$ and $R_2(\mathbf{Z}_2)$.\par 
     	\noindent {(b)} (if part). Note that if $\mathbf{Z}_1$ and $\mathbf{Z}_2$ are independent, then so are $\mathbf{R}_1(\mathbf{Z}_1)$ and $\mathbf{R}_2(\mathbf{Z}_2)$. Then,~\eqref{eq:rdcov1} clearly implies $\Rdcov^2(\mathbf{Z}_1,\mathbf{Z}_2)=0$.\par 
     	(only if part). By~\cite{Gabor2007} (or~\eqref{eq:rdcov1} coupled with continuity of characteristic functions), $\Rdcov^2(\mathbf{Z}_1,\mathbf{Z}_2)=0$ implies $R_1(\mathbf{Z}_1)$ and $R_2(\mathbf{Z}_2)$ are independent. Now note that, by \cref{prop:Mccan}, there exists maps $Q_1:[0,1]^{d_1}\to\mathbb{R}^{d_1}$ and $Q_2: [0,1]^{d_2}\to\mathbb{R}^{d_2}$  such that $$Q_1(R_1(\mathbf{Z}_1))=\mathbf{Z}_1\;\; a.e.\;\mu_1\qquad  \mbox{and}\qquad Q_2(R_2(\mathbf{Z}_2))=\mathbf{Z}_2\;\;a.e.\;\mu_2.$$
     	By the above display, $(\mathbf{Z}_1,\mathbf{Z}_2)\overset{d}{=}(Q_1(R_1(\mathbf{Z}_1)),Q_2(R_2(\mathbf{Z}_2)))$. As $R_1(\mathbf{Z}_1)$ and $R_2(\mathbf{Z}_2)$ are independent, so are $Q_1(R_1(\mathbf{Z}_1))$ and $Q_2(R_2(\mathbf{Z}_2))$, and consequently, so are $\mathbf{Z}_1$ and $\mathbf{Z}_2$. 	\par 
     	\noindent {(c)} From~\cite[Theorem 4]{Gabor2007}, we know that $\Rdcov(\mathbf{Z}_1,\mathbf{Z}_1)\geq 0$ with equality holding if and only if $R_1(\mathbf{Z}_1)$ has a degenerate distribution. However, $R_1(\mathbf{Z}_1)$ has a non-degenerate $d$-dimensional uniform distribution. So, $\Rdcov(\mathbf{Z}_1,\mathbf{Z}_1)>0$.
     	
     	\noindent {(d)} Set $\mathbf{Y}_1:=\mathbf{a}_1+b\mathbf{Z}_1$ where $\mathbf{a}_1\in\mathbb{R}^{d_1}$ and $b>0$. By~\cref{prop:Mccan}, we know that the (unique) rank map is the gradient of a convex function. Let the convex function corresponding to $\mathbf{Z}_1$ be $\phi_1(\cdot)$. This implies that $b\phi_1(b^{-1}(\mathbf{x}-\mathbf{a}_1))$ is also a convex function and its gradient is given by $R_1(b^{-1}(\mathbf{x}-\mathbf{a}_1))$. Set $\tilde{R}(\mathbf{y}_1)=R_1(b^{-1}(\mathbf{y}_1-\mathbf{a}_1))$ and note that 
     	\begin{equation}\label{eq:proprdcovpf1}
     	\tilde{R}(\mathbf{Y}_1)=R_1(b^{-1}(\mathbf{a}_1+b\mathbf{Z}_1-\mathbf{a}_1))=R_1(\mathbf{Z}_1)\overset{d}{=}\mathcal{U}^{d_1}.
     	\end{equation}
     	Therefore $\tilde{R}(\cdot)$ pushes the measure induced by $\mathbf{Y}_1$ to the $d$-dimensional uniform distribution, and it is the gradient of a convex function. As a result,~\cref{prop:Mccan} implies $\tilde{R}(\cdot)$ is the population tank map corresponding to $\mathbf{Y}_1$ is $\tilde{R}(\cdot)$. \par
     	The same argument as in~\eqref{eq:proprdcovpf1} also yields that $\bar{R}_2(\mathbf{Y}_2)=R_2(\mathbf{Z}_2)$ where $\mathbf{Y}_2 :=\mathbf{a}_2+b\mathbf{Z}_2$ and $\bar{R}(\cdot)$ denotes the population rank map corresponding to $\mathbf{Y}_2$. Therefore, $\Rdcov(\mathbf{Y}_1,\mathbf{Y}_2)=\Rdcov(\mathbf{Z}_1,\mathbf{Z}_2)$.
     	\par 	
     	\noindent{(e)} Let $R_1^n$ and $R_2^n$ denote the rank maps corresponding to $\mathbf{Z}_1^n$ and $\mathbf{Z}_2^n$ respectively. By repeating the same argument as in the proof of~\cref{lem:L2}, we get: $$(\mathbf{Z}_1^n,R_1^n(\mathbf{Z}_1^n))\overset{w}{\longrightarrow} (\mathbf{Z}_1,R_1(\mathbf{Z}_1))\qquad\mbox{and}\qquad (\mathbf{Z}_2^n,R_2^n(\mathbf{Z}_2^n))\overset{w}{\longrightarrow} (\mathbf{Z}_2,R_2(\mathbf{Z}_2))$$
     	which then by the continuous mapping theorem applied to the function $g(\mathbf{y},\mathbf{x}) :=\lVert \mathbf{y}-R_1(\mathbf{x})\rVert$ (equivalently, $g(\mathbf{y},\mathbf{x})=\lVert \mathbf{y}-R_2(\mathbf{x})\rVert$), as in the proof of~\cref{lem:L2}, implies, 
     	\begin{equation}\label{eq:proprdcovpf2}
     	\lVert R_1^n(\mathbf{Z}_1^{n})-R_1(\mathbf{Z}_1^{n})\rVert\overset{P}{\longrightarrow}0\qquad \mbox{and}\qquad \lVert R_2^n(\mathbf{Z}_2^{n})-R_2(\mathbf{Z}_2^{n})\rVert\overset{P}{\longrightarrow}0.
     	\end{equation}
     	Let $(\mathbf{Z}_1^{1,n},\mathbf{Z}_2^{1,n}),(\mathbf{Z}_1^{2,n},\mathbf{Z}_2^{2,n}),(\mathbf{Z}_1^{3,n},\mathbf{Z}_2^{3,n})$ be  i.i.d.~according to the same distribution as $(\mathbf{Z}_1^n,\mathbf{Z}_2^n)$ and let $(\mathbf{Z}_1^1,\mathbf{Z}_2^1),(\mathbf{Z}_1^2,\mathbf{Z}_2^2),(\mathbf{Z}_1^3,\mathbf{Z}_2^3)$ be i.i.d.~copies of $(\mathbf{Z}_1,\mathbf{Z}_2)$. Next, observe that, by the triangle inequality,
     	\begin{align*}
     	&\;\;\;\;\Big|\lVert R_1^n(\mathbf{Z}_1^{1,n})-R_1^n(\mathbf{Z}_1^{2,n})\rVert \lVert R_2^n(\mathbf{Z}_2^{1,n})-R_2^n(\mathbf{Z}_2^{2,n})\rVert-\lVert R_1(\mathbf{Z}_1^{1,n})-R_1(\mathbf{Z}_1^{2,n})\rVert \lVert R_2(\mathbf{Z}_2^{1,n})-R_2(\mathbf{Z}_2^{1,n})\rVert\Big|\nonumber \\&\leq \lVert R_1^n(\mathbf{Z}_1^{1,n})-R_1(\mathbf{Z}_1^{1,n})\rVert\lVert R_2^n(\mathbf{Z}_2^{1,n})-R_2^n(\mathbf{Z}_2^{2,n})\rVert+\lVert R_1^n(\mathbf{Z}_1^{2,n})-R_1(\mathbf{Z}_1^{2,n})\rVert\lVert R_2^n(\mathbf{Z}_2^{1,n})-R_2^n(\mathbf{Z}_2^{2,n})\rVert \nonumber \\& \hspace{1in} +\lVert R_1(\mathbf{Z}_1^{1,n})-R_1(\mathbf{Z}_1^{2,n})\rVert\lVert R_2^n(\mathbf{Z}_2^{1,n})-R_2(\mathbf{Z}_2^{1,n})\rVert\nonumber \\&\hspace{1in}+\lVert R_1(\mathbf{Z}_1^{1,n})-R_1(\mathbf{Z}_1^{2,n})\rVert\lVert R_2^n(\mathbf{Z}_2^{2,n})-R_2(\mathbf{Z}_2^{2,n})\rVert\overset{P}{\longrightarrow}0.
     	\end{align*}
     	Finally note that, by using the continuous mapping theorem on the joint weak convergence of $(\mathbf{Z}_1^n,\mathbf{Z}_2^n)$ to $(\mathbf{Z}_1,\mathbf{Z}_2)$, we further get: 
     	$$\lVert R_1(\mathbf{Z}_1^{1,n})-R_1(\mathbf{Z}_1^{2,n})\rVert\lVert R_2(\mathbf{Z}_2^{1,n})-R_2(\mathbf{Z}_2^{2,n})\rVert\overset{w}{\longrightarrow}\lVert R_1(\mathbf{Z}_1^{1})-R_1(\mathbf{Z}_1^{2})\rVert\lVert R_1(\mathbf{Z}_2^{1})-R_1(\mathbf{Z}_2^{2})\rVert.$$
     	Next, by applying the dominated convergence theorem, we get:
     	\begin{align}\label{eq:proprdcovpf3}
     	&\;\;\;\;\;\;\mathbb{E}\big[\lVert R_1(\mathbf{Z}_1^{1,n})-R_1(\mathbf{Z}_1^{2,n})\rVert\lVert R_2(\mathbf{Z}_2^{1,n})-R_2(\mathbf{Z}_2^{2,n})\rVert\big]\nonumber \\&\;\;\;\;\overset{n\to\infty}{\longrightarrow}\mathbb{E}\big[\lVert R_1(\mathbf{Z}_1^{1})-R_1(\mathbf{Z}_1^{2})\rVert\lVert R_1(\mathbf{Z}_2^{1})-R_1(\mathbf{Z}_2^{2})\rVert\big].
     	\end{align}
     	Combining~\eqref{eq:proprdcovpf3},~\eqref{eq:proprdcovpf2} with the dominated convergence theorem yields:
     	\begin{align}\label{eq:proprdcovpf4}
     	&\;\;\;\;\;\;\mathbb{E}\big[\lVert R_1^n(\mathbf{Z}_1^{1,n})-R_1^n(\mathbf{Z}_1^{2,n})\rVert \lVert R_2^n(\mathbf{Z}_2^{1,n})-R_2^n(\mathbf{Z}_2^{2,n})\rVert\big]\nonumber \\&\;\;\;\;\overset{n\to\infty}{\longrightarrow}\mathbb{E}\big[\lVert R_1(\mathbf{Z}_1^{1})-R_1(\mathbf{Z}_1^{2})\rVert\lVert R_1(\mathbf{Z}_2^{1})-R_1(\mathbf{Z}_2^{2})\rVert\big].
     	\end{align}
     	By using the same arguments as above, we can similarly show the following:
     	\begin{align}\label{eq:proprdcovpf5}
     	&\;\;\;\;\;\;\mathbb{E}\big[\lVert R_1^n(\mathbf{Z}_1^{1,n})-R_1^n(\mathbf{Z}_1^{2,n})\rVert\big]\mathbb{E}\big[ \lVert R_2^n(\mathbf{Z}_2^{1,n})-R_2^n(\mathbf{Z}_2^{2,n})\rVert\big]\nonumber \\&\;\;\;\;\overset{n\to\infty}{\longrightarrow}\mathbb{E}\big[\lVert R_1(\mathbf{Z}_1^{1})-R_1(\mathbf{Z}_1^{2})\rVert\big]\mathbb{E}\big[\lVert R_1(\mathbf{Z}_2^{1})-R_1(\mathbf{Z}_2^{2})\rVert\big]
     	\end{align}
     	and
     	\begin{align}\label{eq:proprdcovpf6}
     	&\;\;\;\;\;\;\mathbb{E}\big[\lVert R_1^n(\mathbf{Z}_1^{1,n})-R_1^n(\mathbf{Z}_1^{2,n})\rVert \lVert R_2^n(\mathbf{Z}_2^{1,n})-R_2^n(\mathbf{Z}_2^{3,n})\rVert\big]\nonumber \\&\;\;\;\;\overset{n\to\infty}{\longrightarrow}\mathbb{E}\big[\lVert R_1(\mathbf{Z}_1^{1})-R_1(\mathbf{Z}_1^{2})\rVert\lVert R_1(\mathbf{Z}_2^{1})-R_1(\mathbf{Z}_2^{3})\rVert\big].
     	\end{align}
     	Combining~\eqref{eq:proprdcovpf4},~\eqref{eq:proprdcovpf5},~\eqref{eq:proprdcovpf6} and~\cref{lem:proprdcov} (part (a)), we finally get:
     	$$\Rdcov^2(\mathbf{Z}_1^n,\mathbf{Z}_2^n)\overset{n\to\infty}{\longrightarrow}\Rdcov^2(\mathbf{Z}_1,\mathbf{Z}_2)$$
     	which completes the proof.
     \end{proof}
     \subsection{Proof of~\cref{lem:propren}}\label{proof:propren}
     \begin{proof}
     	{(a)} This is a direct consequence of \cite[Lemmas 2.2 and 2.3]{Baringhaus2004}.
     	\par 
     	\noindent{(b)} (if part). Assuming $\mathbf{Z}_1\overset{d}{=}\mathbf{Z}_2$, we have $\mathbb{P}(\mathbf{a}^\top\mathbf{Z}_1\leq t)=\mathbb{P}(\mathbf{a}^\top\mathbf{Z}_2\leq t)$ for all $\mathbf{a}\in \mathcal{S}^{d-1}$ and $t\in\mathbb{R}$. Therefore, by~\eqref{eq:ren1}, $\Ren_{\lambda}^2(\mathbf{Z}_1,\mathbf{Z}_2)=0$.\par 
     	(only if part). By~\cite[Theorem 2.1]{Baringhaus2004}, if $\Ren_{\lambda}^2(\mathbf{Z}_1,\mathbf{Z}_2)=0$, then we have:
     	\begin{align}\label{eq:proprenpf1}
     	R_{\lambda}(\mathbf{Z}_1)\overset{d}{=}R_{\lambda}(\mathbf{Z}_2).
     	\end{align}
     	Next, by~\cref{prop:Mccan}, there exists $Q_{\lambda}:\mathbb{R}^d\to\mathbb{R}^d$ such that
     	\begin{align}\label{eq:proprenpf2}
     	Q_{\lambda}(R_{\lambda}(\mathbf{Z}_1))=\mathbf{Z}_1\qquad \mbox{a.e. }\mu_1\qquad\mbox{and}\qquad Q_{\lambda}(R_{\lambda}(\mathbf{Z}_2))=\mathbf{Z}_2\qquad \mbox{a.e. }\mu_2.
     	\end{align}
     	Finally,~\eqref{eq:proprenpf2} combined with~\eqref{eq:proprenpf1} yields $\mathbf{Z}_1\overset{d}{=}\mathbf{Z}_2$ and completes the proof.\par 
     	\noindent{(c)} The proof is verbatim similar to that of~\cref{lem:proprdcov} (part (c)) in~\cref{proof:proprdcov}.\par 
     	\noindent{(d)} Note that, by part (a), $\Ren_{\lambda}^2(\mathbf{Z}_1,\mathbf{Z}_2)$ can be written as expectations of Euclidean distances between bounded random vectors. Therefore, the proof is exactly similar to that of~\cref{lem:proprdcov} (part (d)) in~\cref{proof:proprdcov}. We leave the details to the reader.
     \end{proof}
     \subsection{Proof of~\cref{rem:indepdistfree}}\label{proof:indepdistfree}
     \begin{proof}
     	Note that $\Rdcov_n^2$, as defined in~\eqref{eq:indepstat}, is a function of $(\hat{R}_n^{\mathbf{X}}(\mathbf{X}_1),\ldots ,\hat{R}_n^{\mathbf{X}}(\mathbf{X}_n),$ $\hat{R}_n^{\mathbf{Y}}(\mathbf{Y}_1),\ldots ,\hat{R}_n^{\mathbf{Y}}(\mathbf{Y}_n))$. Under $\mathrm{H}_0$, the distribution of the above vector further splits into the product of the marginal distributions of $(\hat{R}_n^{\mathbf{X}}(\mathbf{X}_1),\ldots ,\hat{R}_n^{\mathbf{X}}(\mathbf{X}_n))$ and $(\hat{R}_n^{\mathbf{Y}}(\mathbf{Y}_1),\ldots ,\hat{R}_n^{\mathbf{Y}}(\mathbf{Y}_n))$. By~\cref{prop:finproperties} (part (ii)), each of the marginals are distribution-free, i.e., their distribution does not depend on $\mu_{\mathbf{X}}$ and $\mu_{\mathbf{Y}}$. In fact, they are distributed uniformly over all $n!$ permutations of each of the fixed grids $\mathcal{H}_n^{d_1}$ and $\mathcal{H}_n^{d_2}$ respectively. This results in the distribution-free property of the statistic $\Rdcov_n^2$ (under $\mbox{H}_0$).
     \end{proof}
     \subsection{Proof of~\cref{lem:dcovequiv}}\label{proof:dcovequiv}
     \begin{proof}
     	Let us first prove~\eqref{eq:dcovequivmain1}. Suppose that $(X_1,Y_1),(X_2,Y_2)$ and $(X_3,Y_3)$ are random samples drawn from the joint distribution of $(X,Y)$. Let $W_i :=F^X(X_i)$ and $Z_i :=F^Y(Y_i)$ for $i=1,2,3$. Further, we will use $F^{W,Z}(\cdot)$, $F^W(\cdot)$ and $F^Z(\cdot)$ to denote the joint and marginal DFs of $W$ and $Z$ respectively. Note that,
     	\begin{equation}\label{eq:dcovequivpf1}
     	\Rdcov^2(X,Y)=\mathbb{E}[|W_1-W_2||Z_1-Z_2|]+\mathbb{E}|W_1-W_2|\mathbb{E}|Z_1-Z_2|-2\mathbb{E}|W_1-W_2||Z_1-Z_3|.
     	\end{equation}
     	Further, we can write the following simple algebraic identity:
     	\begin{equation}\label{eq:dcovequivpf2}
     	|W_1-W_2|=\int_{-\infty}^{\infty} \big[\mathbf{1}(W_1\leq u\leq W_2)+\mathbf{1}(W_2\leq u\leq W_1)\big]\,du.
     	\end{equation}
     	We can write a similar result for $|Z_1-Z_2|$ (as in~\eqref{eq:dcovequivpf2}), which, on multiplying with the right hand side of~\eqref{eq:dcovequivpf2}, yields:
     	\begin{align}\label{eq:dcovequivpf3}
     	&\;\;\;\;|W_1-W_2||Z_1-Z_2|\nonumber \\&=\int_{-\infty}^{\infty}\int_{-\infty}^{\infty} \bigg[\mathbf{1}(W_1\leq u\leq W_2)\mathbf{1}(Z_1\leq v\leq Z_2)+\mathbf{1}(W_1\leq u\leq W_2)\mathbf{1}(Z_2\leq v\leq Z_1)\nonumber \\&\qquad +\mathbf{1}(W_2\leq u\leq W_1)\mathbf{1}(Z_1\leq v\leq Z_2)+\mathbf{1}(W_2\leq u\leq W_1)\mathbf{1}(Z_2\leq v\leq Z_1)\,du\,dv\bigg]
     	\end{align}
     	Next, by applying Fubini's theorem on~\eqref{eq:dcovequivpf3}, we get:
     	\begin{align}\label{eq:dcovequivpf4}
     	\mathbb{E}\left[|W_1-W_2||Z_1-Z_2|\right]&=2\int_{-\infty}^{\infty} \Big[F^{W,Z}(u,v)+2(F^{W,Z}(u,v))^2-2F^W(u)F^{W,Z}(u,v)\nonumber \\&\qquad -2F^{W,Z}(u,v)F^Z(v)+F^W(u)F^Z(v)\Big]\,du\,dv.
     	\end{align}
     	Similar calculations also result in the following:
     	\begin{align}\label{eq:dcovequivpf5}
     	\mathbb{E}|W_1-W_2|\,\mathbb{E}|Z_1-Z_2|&=4\int_{-\infty}^{\infty}\int_{-\infty}^{\infty} \Big[F^W(u)F^Z(v)-(F^W(u))^2F^Z(v)\nonumber \\&\qquad-F^W(u)(F^Z(v))^2+(F^W(u)F^Z(v))^2\Big]\,du\,dv
     	\end{align}
     	and,
     	\begin{align}\label{eq:dcovequivpf6}
     	\mathbb{E}|W_1-W_2||Z_1-Z_3|&=\int_{-\infty}^{\infty}\int_{-\infty}^{\infty} \Big[3F^W(u)F^Z(v)+F^{W,Z}(u,v)+4F^{W,Z}(u,v)F^W(u)F^Z(v)\nonumber \\&\qquad -2F^W(u)F^{W,Z}(u,v)-2F^{W,Z}(u,v)F^Z(v)-2(F^W(u))^2F^Z(v)\nonumber \\&\qquad-2F^Z(u)(F^W(v))^2\Big]\,du\,dv.
     	\end{align}
     	Therefore, by combining~\eqref{eq:dcovequivpf1},~\eqref{eq:dcovequivpf4},~\eqref{eq:dcovequivpf5} and~\eqref{eq:dcovequivpf6}, we get:
     	\begin{align*}
     	\Rdcov^2(X,Y)&=4\int_{-\infty}^{\infty}\int_{-\infty}^{\infty} \Big[(F^{W,Z}(u,v))^2-2F^{W}(u)F^Z(v)F^{W,Z}(u,v)+(F^W(u)F^Z(v))^2\Big]\,du\,dv\nonumber \\&=4\int_{-\infty}^{\infty}\int_{-\infty}^{\infty} \Big(F^{W,Z}(u,v)-F^W(u)F^Z(v)\Big)^2\,du\,dv.
     	\end{align*}
     	Next, note that $F^{W,Z}(u,v)=F^{X,Y}\big((F^X(u))^{-1},(F^{Y}(v))^{-1}\big)$, $F^W(u)=u$ and $F^Z(v)=v$. Using this observation, coupled with a standard change of variable formula for Riemann integrals, we get:
     	\begin{align*}
     	\frac{1}{4}\Rdcov^2(X,Y)=\int_{-\infty}^{\infty}\int_{-\infty}^{\infty} \big(F^{X,Y}(x,y)-F^{X}(x)F^{Y}(y)\big)^2\,dF^X(x)\,dF^{Y}(y).
     	\end{align*}
     	This completes the proof of~\eqref{eq:dcovequivmain1}. \par 
     	In order to prove~\eqref{eq:dcovequivmain2}, let us start with some notation:
     	\begin{align*}
     	T_1\coloneqq\frac{1}{n^4}\sum_{i,j,k,l}\mathbf{1}(X_l<X_j,Y_k<Y_i)\; ,\; T_2\coloneqq \frac{1}{n^5}\sum_{i,j,k,l,m} \mathbf{1}(X_l<X_j,Y_k<Y_m, Y_i<Y_m)
     	\end{align*}
     	\vspace{-0.2in}
     	\begin{align*}
     	T_3\coloneqq \frac{1}{n^5}\sum_{i,j,k,l,m} \mathbf{1}(Y_i<Y_j,X_k<X_m,X_l<X_m)\; ,\; T_4\coloneqq \frac{1}{n^3}\sum_{i,j,k} \mathbf{1}(X_k<X_i,Y_k<Y_j),
     	\end{align*}
     	\vspace{-0.2in}
     	\begin{align*}
     	T_5\coloneqq \frac{1}{n^4}\sum_{i,j,k,l} \mathbf{1}(X_l<X_j)\mathbf{1}(Y_k<Y_i,Y_l<Y_i), 
     	\end{align*}
     	\vspace{-0.2in}
     	\begin{align*}
     	T_6\coloneqq \frac{1}{n^4}\sum_{i,j,k,l}\mathbf{1}(X_k<X_j,X_l<X_j)\mathbf{1}(Y_k<Y_i),
     	\end{align*}
     	\vspace{-0.2in}
     	\begin{align*}
     	T_7\coloneqq \frac{1}{n^5}\sum_{i,j,k,l,m}\mathbf{1}(X_k<X_j,X_l<X_j)\mathbf{1}(Y_k<Y_i,Y_m<Y_i),
     	\end{align*}
     	\vspace{-0.2in}
     	\begin{align*}
     	T_8\coloneqq \frac{1}{n^4}\sum_{i,j,k,l}\mathbf{1}(X_k<X_j,X_l<X_j)\mathbf{1}(Y_k<Y_i,Y_l<Y_i),
     	\end{align*}
     	\vspace{-0.2in}
     	\begin{align*}
     	T_9\coloneqq \frac{1}{n^6}\sum_{i,j,k,l,m,p}\mathbf{1}(X_k<X_j,X_l<X_j)\mathbf{1}(Y_m<Y_i,Y_p<Y_i).
     	\end{align*}
     	Next, note that:
     	\begin{align}\label{eq:dcovequivpf28}
     	S_1&\coloneqq \frac{1}{n^2}\sum_{k,l} \big|\hat{R}_n^X(X_k)-\hat{R}_n^X(X_l)\big|\big|\hat{R}_n^{Y}(Y_k)-\hat{R}_n^Y(Y_l)\big|\nonumber \\&=\frac{1}{n^4}\sum_{k,l,i,j}\big[\mathbf{1}(X_k<X_j<X_l)+\mathbf{1}(X_l<X_j<X_k)\big]\big[\mathbf{1}(Y_k<Y_i<Y_l)+\mathbf{1}(Y_l<Y_i<Y_k)\big]\nonumber \\&=4T_1-4T_2-4T_3+4T_9.
     	\end{align}
     	Similar calculations reveal that,
     	\begin{align}\label{eq:dcovequivpf29}
     	S_2&\coloneqq\left(\frac{1}{n^2}\sum_{k,l=1}^n \big| \hat{R}_n^{\mathbf{X}}(\mathbf{X}_k)-\hat{R}_n^{\mathbf{X}}(\mathbf{X}_l)\big|\right)\times\left(\frac{1}{n^2}\sum_{k,l=1}^n \big| \hat{R}_n^{\mathbf{Y}}(\mathbf{Y}_k)-\hat{R}_n^{\mathbf{Y}}(\mathbf{Y}_l)\big|\right)\nonumber \\&=2T_1+2T_4-4T_5-4T_6+4T_8,
     	\end{align}
     	and
     	\begin{align}\label{eq:dcovequivpf210}
     	S_3&\coloneqq \frac{1}{n^3}\sum_{k,l,m=1}^n \big| \hat{R}_n^{\mathbf{X}}(\mathbf{X}_k)-\hat{R}_n^{\mathbf{X}}(\mathbf{X}_l)\big|\big| \hat{R}_n^{\mathbf{Y}}(\mathbf{Y}_k)-\hat{R}_n^{\mathbf{Y}}(\mathbf{Y}_m)\big|\nonumber \\&=3T_1-2T_2-2T_3+T_4-2T_5-2T_6+4T_7.
     	\end{align}
     	Recall that $\Rdcov_n^2=S_1+S_2-2S_3$. Therefore, by~\eqref{eq:dcovequivpf28},~\eqref{eq:dcovequivpf29} and~\eqref{eq:dcovequivpf210}, we get:
     	\begin{align*}
     	\Rdcov_n^2&=4(T_7-2T_8+T_9)\nonumber \\ &=\frac{4}{n^2}\sum_{i,j}\left(F^{X,Y}_n(X_i,Y_j)-F^X_n(X_i)F^Y_n(Y_j)\right)^2\nonumber \\ &=4\int \int  \left(F^{X,Y}_n(x,y)-F^X_n(x)F^Y_n(y)\right)^2\,dF^X_n(x)\,dF^Y_n(y).
     	\end{align*}
     	This completes the proof of~\eqref{eq:dcovequivmain2}.
     \end{proof}
     \subsection{Proof of~\cref{theo:indepasdistn}}\label{proof:indepasdistn}
     \begin{proof}
     	Let us start with some notation. Let $\mathcal{H}_n^{d_1}=\{\mathbf{h}_1^{d_1},\mathbf{h}_2^{d_1},\ldots ,\mathbf{h}_n^{d_1}\}$ and $\mathcal{H}_n^{d_2}=\{\mathbf{h}_1^{d_2},\mathbf{h}_2^{d_2},\ldots ,\mathbf{h}_n^{d_2}\}$. Next, define:
     	\begin{align*}
     	w(\mathbf{t},\mathbf{s})\coloneqq\left(\frac{\pi^{1+(d_1+d_2)/2}}{\Gamma((1+d_1)/2)\Gamma((1+d_2)/2)}\lVert \mathbf{t}\rVert^{1+d_1}\lVert \mathbf{s}\rVert^{1+d_2}\right)^{-1},\\
     	f_{\mathbf{X},\mathbf{Y}}^n(\mathbf{t},\mathbf{s})\coloneqq \frac{1}{n}\sum_{j=1}^n \exp\left(i\mathbf{t}^{\top}\hat{R}_n^{\mathbf{X}}(\mathbf{X}_j)+i\mathbf{s}^{\top}\hat{R}_n^{\mathbf{Y}}(\mathbf{Y}_j)\right),\\
     	f_{\mathbf{X}}^n(\mathbf{t})\coloneqq \frac{1}{n}\sum_{j=1}^n \exp\left(i\mathbf{t}^{\top}\mathbf{h}_j^{d_1}\right)\qquad\mbox{and}\qquad f_{\mathbf{Y}}^n(\mathbf{s})\coloneqq \frac{1}{n}\sum_{j=1}^n \exp\left(i\mathbf{s}^{\top}\mathbf{h}_j^{d_2}\right)
     	\end{align*}
     	for $\mathbf{t}\in\mathbb{R}^{d_1},\; \mathbf{s}\in\mathbb{R}^{d_2}$ and $i=\sqrt{-1}$. Note that $f_{\mathbf{X}}^n(\cdot)$ and $f_{\mathbf{Y}}^n(\cdot)$ are deterministic quantities. Recall that, under $\mathrm{H}_0$, $(\hat{R}_n^{\mathbf{X}}(\mathbf{X}_1),\ldots ,\hat{R}_n^{\mathbf{X}}(\mathbf{X}_n))$ and $(\hat{R}_n^{\mathbf{Y}}(\mathbf{Y}_1),\ldots ,\hat{R}_n^{\mathbf{Y}}(\mathbf{Y}_n))$ are independent and distributed uniformly over all $n!$ permutations of the sets $\mathcal{H}_n^{d_1}$ and $\mathcal{H}_n^{d_2}$ respectively (see~\cref{rem:indepdistfree}). Let $\sigma_1$ and $\sigma_2$ be two independent random permutations of the set $\{1,2,\ldots ,n\}$. Then note that:
     	\begin{align}\label{eq:equivdistn}
     	f_{\mathbf{X},\mathbf{Y}}^n(\mathbf{t},\mathbf{s})&\overset{d}{=}\frac{1}{n}\sum_{j=1}^n \exp\left(i\mathbf{t}^{\top}\mathbf{h}^{d_1}_{\sigma_1(j)}+i\mathbf{s}^{\top}\mathbf{h}^{d_2}_{\sigma_2(j)}\right)\nonumber \\&\overset{d}{=}\frac{1}{n}\sum_{j=1}^n \exp\left(i\mathbf{t}^{\top}\mathbf{h}^{d_1}_j+i\mathbf{s}^{\top}\mathbf{h}^{d_2}_{\sigma_2(j)}\right).
     	\end{align}
     	Using~\eqref{eq:equivdistn}, along with~\cite[Theorem 1]{Gabor2007}, we get:
     	\begin{align}\label{eq:equivdistnmain}
     	\Rdcov_n^2\overset{d}{=}\mathop{\mathlarger{\mathlarger{\int}}} \Bigg|\frac{1}{n}\sum_{j=1}^n \exp\left(i\mathbf{t}^{\top}\mathbf{h}^{d_1}_j+i\mathbf{s}^{\top}\mathbf{h}^{d_2}_{\sigma_2(j)}\right)-f_{\mathbf{X}}^n(\mathbf{t})f_{\mathbf{Y}}^n(\mathbf{s})\Bigg|^2w(\mathbf{t},\mathbf{s})\,d\mathbf{t}\,d\mathbf{s}.
     	\end{align}
     	\cref{lem:rankdcov} (see~\cref{proof:rankdcov} for a proof) below deals with the asymptotic behavior of permutation statistics of the same form as in the right hand side of~\eqref{eq:equivdistnmain}. 
     	\begin{lemma}\label{lem:rankdcov}
     		Consider two infinite deterministic sequences, $\{\mathbf{U}_i\}_{i\geq 1}$ and $\{\mathbf{V}_i\}_{i\geq 1}$, $\mathbf{U}_i\in [0,1]^{d_1}$ and $\mathbf{V}_i\in [0,1]^{d_2}$, such that $\mathbb{P}_n=(1/n)\sum_{i=1}^n \delta_{\mathbf{U}_i}\overset{w}{\longrightarrow}\mathcal{U}^{d_1}$ and $\mathbb{Q}_n=(1/n)\sum_{i=1}^n \delta_{\mathbf{V}_i}\overset{w}{\longrightarrow}\mathcal{U}^{d_2}$, where $\mathcal{U}^{d_1}$ (and $\mathcal{U}^{d_2}$) are the standard Lebesgue measures on $[0,1]^{d_1}$ (and $[0,1]^{d_2}$). Further, let $S_n$ denote the set of all permutations of $\{1,2,\ldots ,n\}$ and $\mathbf{\pi}_n$ be a random permutation drawn uniformly from $S_n$. Define the following:
     		\begin{align*}
     		\xi_{n}(\mathbf{t},\mathbf{s})\coloneqq &\sqrt{n}\Bigg(\frac{1}{n}\sum_{k=1}^n \exp\big(i\mathbf{t}^{\top}\mathbf{U}_k+i\mathbf{s}^{\top}\mathbf{V}_{\pi_n(k)}\big)-\Bigg(\frac{1}{n}\sum_{k=1}^n \exp\big(i\mathbf{t}^{\top}\mathbf{U}_k\big)\Bigg)\nonumber \\&\qquad\times \Bigg(\frac{1}{n}\sum_{k=1}^n \exp\big(i\mathbf{s}^{\top}\mathbf{V}_{k}\big)\Bigg)\Bigg)
     		\end{align*}
     		for $\mathbf{t}\in\mathbb{R}^{d_1}$ and $\mathbf{s}\in \mathbb{R}^{d_2}$. Then, we have,
     		\begin{align}\label{eq:fluctuation}
     		D_n\coloneqq \int_{\mathbb{R}^{d_1}\times\mathbb{R}^{d_2}}\big|\xi_n(\mathbf{t},\mathbf{s})\big|^2w(\mathbf{t},\mathbf{s})\,d\mathbf{t}\,d\mathbf{s}\overset{w}{\longrightarrow}\int_{\mathbb{R}^{d_1}\times\mathbb{R}^{d_2}}\big|\xi(\mathbf{t},\mathbf{s})\big|^2w(\mathbf{t},\mathbf{s})\,d\mathbf{t}\,d\mathbf{s}\coloneqq D.
     		\end{align}
     		Here $\xi(\cdot,\cdot)$ is a complex-valued Gaussian process with mean $0$ and covariance kernel 
     		\begin{align*}
     		R((\mathbf{t}_1,\mathbf{s}_1),(\mathbf{t}_0,\mathbf{s}_0))&\coloneqq\left(f_{d_1}(\mathbf{t}_1-\mathbf{t}_0)-f_{d_1}(\mathbf{t}_1)\overline{f_{d_1}(\mathbf{t}_0)}\right)\left(f_{d_2}(\mathbf{s}_1-\mathbf{s}_0)-f_{d_2}(\mathbf{s}_1)\overline{f_{d_2}(\mathbf{s}_0)}\right)
     		\end{align*}
     		where $f_{d_1}(\cdot)$ and $f_{d_2}(\cdot)$ denote the characteristic functions of a $\mathcal{U}^{d_1}$ and a $\mathcal{U}^{d_2}$ random variable respectively.
     	\end{lemma}
     	With the above lemma in mind, note that the empirical distributions on $\mathcal{H}_n^{d_1}$ and $\mathcal{H}_n^{d_2}$ converge weakly to $\mathcal{U}^{d_1}$ and $\mathcal{U}^{d_2}$ respectively (by assumption \textbf{(AP2)}). By setting $\{\mathbf{U}_j\}_{j\geq 1}$ as $\{\mathbf{h}_j^{d_1}\}_{j\geq 1}$ and $\{\mathbf{V}_j\}_{j\geq 1}$ as $\{\mathbf{h}_j^{d_2}\}_{j\geq 1}$, and applying~\cref{lem:rankdcov} to the right side of~\eqref{eq:equivdistnmain}, we get that, as $n\to\infty$:
     	\begin{align}\label{eq:indepasdistnmain}
     	n\Rdcov_n^2\overset{w}{\longrightarrow} D\overset{d}{=}\sum_{j\geq 1}\lambda_jZ^2_j
     	\end{align}
     	where $D$ is defined in~\eqref{eq:fluctuation}, $\lambda_j$'s are fixed positive constants, and $Z_j$'s are independent standard normals. The last equivalence in~\eqref{eq:indepasdistnmain} follows from~\cite[Chapter 1, Section 2]{Kuo1975} using standard Karhunen-Lo\`{e}ve type expansions for Gaussian processes. This completes the proof.
     \end{proof}
     \subsection{Proof of~\cref{theo:indepconsis}}\label{proof:indepconsis}
     \begin{proof}
     	Recall that $\Rdcov_n^2=S_1+S_2-2S_3$
     	where 
     	\begin{align*}
     	S_1&\coloneqq \frac{1}{n^2}\sum_{k,l=1}^n \lVert \hat{R}_n^{\mathbf{X}}(\mathbf{X}_k)-\hat{R}_n^{\mathbf{X}}(\mathbf{X}_l)\rVert\lVert \hat{R}_n^{\mathbf{Y}}(\mathbf{Y}_k)-\hat{R}_n^{\mathbf{Y}}(\mathbf{Y}_l)\rVert,\\ S_2&\coloneqq\left(\frac{1}{n^2}\sum_{k,l=1}^n \lVert \hat{R}_n^{\mathbf{X}}(\mathbf{X}_k)-\hat{R}_n^{\mathbf{X}}(\mathbf{X}_l)\rVert\right)\times\left(\frac{1}{n^2}\sum_{k,l=1}^n \lVert \hat{R}_n^{\mathbf{Y}}(\mathbf{Y}_k)-\hat{R}_n^{\mathbf{Y}}(\mathbf{Y}_l)\rVert\right),\\ S_3&\coloneqq\frac{1}{n^3}\sum_{k,l,m=1}^n \lVert \hat{R}_n^{\mathbf{X}}(\mathbf{X}_k)-\hat{R}_n^{\mathbf{X}}(\mathbf{X}_l)\rVert\lVert \hat{R}_n^{\mathbf{Y}}(\mathbf{Y}_k)-\hat{R}_n^{\mathbf{Y}}(\mathbf{Y}_m)\rVert.
     	\end{align*}
     	Let us focus on $S_1$. Observe that by the triangle inequality, for any $k,l\in\{1,2,\ldots ,n\}$:  
     	\begin{equation}\label{eq:indepconsispf2}
     	\lVert \hat{R}_n^{\mathbf{X}}(\mathbf{X}_k)-\hat{R}_n^{\mathbf{X}}(\mathbf{X}_l)\rVert \geq \lVert R^{\mathbf{X}}(\mathbf{X}_k)-R^{\mathbf{X}}(\mathbf{X}_l)\rVert - \lVert \hat{R}_n^{\mathbf{X}}(\mathbf{X}_k)-R^{\mathbf{X}}(\mathbf{X}_k)\rVert-\lVert \hat{R}_n^{\mathbf{X}}(\mathbf{X}_l)-R^{\mathbf{X}}(\mathbf{X}_l)\rVert.
     	\end{equation}
     	Plugging in~\eqref{eq:indepconsispf2} in $S_1$, we get:
     	\begin{align}\label{eq:indepconsispf3}
     	\liminf\limits_{n\to\infty} S_1&\geq \liminf\limits_{n\to\infty}\left(\frac{1}{n^2}\sum_{k,l=1}^n \lVert R^{\mathbf{X}}(\mathbf{X}_k)-R^{\mathbf{X}}(\mathbf{X}_l)\rVert\lVert \hat{R}_n^{\mathbf{Y}}(\mathbf{Y}_k)-\hat{R}_n^{\mathbf{Y}}(\mathbf{Y}_l)\rVert\right)\nonumber \\ &\;\;\;\;-\limsup\limits_{n\to\infty} \left(\frac{1}{n^2}\sum_{k,l=1}^n \lVert R^{\mathbf{X}}(\mathbf{X}_k)-\hat{R}_n^{\mathbf{X}}(\mathbf{X}_k)\rVert\lVert \hat{R}_n^{\mathbf{Y}}(\mathbf{Y}_k)-\hat{R}_n^{\mathbf{Y}}(\mathbf{Y}_l)\rVert\right)\nonumber \\ &\;\;\;\;-\limsup\limits_{n\to\infty} \left(\frac{1}{n^2}\sum_{k,l=1}^n \lVert R^{\mathbf{X}}(\mathbf{X}_l)-\hat{R}_n^{\mathbf{X}}(\mathbf{X}_l)\rVert\lVert \hat{R}_n^{\mathbf{Y}}(\mathbf{Y}_k)-\hat{R}_n^{\mathbf{Y}}(\mathbf{Y}_l)\rVert\right).
     	\end{align}
     	Now, by~\cref{lem:L2} the last two terms on the right side of~\eqref{eq:indepconsispf3} equal $0$ a.s. Therefore, 
     	\begin{align}\label{eq:indepconsispf4}
     	\liminf\limits_{n\to\infty} S_1\geq \liminf\limits_{n\to\infty}\left(\frac{1}{n^2}\sum_{k,l=1}^n \lVert R^{\mathbf{X}}(\mathbf{X}_k)-R^{\mathbf{X}}(\mathbf{X}_l)\rVert\lVert \hat{R}_n^{\mathbf{Y}}(\mathbf{Y}_k)-\hat{R}_n^{\mathbf{Y}}(\mathbf{Y}_l)\rVert\right) \;\;\mbox{a.s}.
     	\end{align}
     	Next, starting from the right side of~\eqref{eq:indepconsispf4} and repeating the same argument as above on the $\mathbf{Y}$'s instead of $\mathbf{X}$'s, we get:
     	\begin{align}\label{eq:indepconsispf5}
     	\liminf\limits_{n\to\infty} S_1\geq \liminf\limits_{n\to\infty}\left(\frac{1}{n^2}\sum_{k,l=1}^n \lVert R^{\mathbf{X}}(\mathbf{X}_k)-R^{\mathbf{X}}(\mathbf{X}_l)\rVert\lVert R^{\mathbf{Y}}(\mathbf{Y}_k)-R^{\mathbf{Y}}(\mathbf{Y}_l)\rVert\right)\;\;\mbox{a.s}.
     	\end{align}
     	Note that $\{R^{\mathbf{X}}(\mathbf{X}_i),R^{\mathbf{Y}}(\mathbf{Y}_i)\}_{1\leq i\leq n}$ are i.i.d.~random vectors, which implies that the right side of~\eqref{eq:indepconsispf5} is a standard V-statistic. Consequently, by invoking the strong law of large numbers for V-statistics, we have:
     	\begin{align}\label{eq:indepconsispf6}
     	\liminf\limits_{n\to\infty} S_1\geq \mathbb{E}\big[\lVert \mathbf{R}^{\mathbf{X}}(\mathbf{Z}_1^1)-R^{\mathbf{X}}(\mathbf{Z}_1^2)\rVert\lVert R^{\mathbf{Y}}(\mathbf{Z}_2^1)-R^{\mathbf{Y}}(\mathbf{Z}_2^2)\rVert\big]\qquad\qquad \;\;\mbox{a.s}.
     	\end{align}
     	where $(\mathbf{Z}_1^1,\mathbf{Z}^1_2),(\mathbf{Z}^2_1,\mathbf{Z}_2^2)$ are independent observations having the same distribution as $(\mathbf{X}_1,\mathbf{Y}_1)$. Also, analogous to~\eqref{eq:indepconsispf2}, an application of the triangle inequality also yields the following:
     	\begin{align}\label{eq:indepconsispf7}
     	&\;\;\;\;\;\lVert \hat{R}_n^{\mathbf{X}}(\mathbf{X}_k)-\hat{R}_n^{\mathbf{X}}(\mathbf{X}_l)\rVert\nonumber\\ &\leq \lVert R^{\mathbf{X}}(\mathbf{X}_k)-R^{\mathbf{X}}(\mathbf{X}_l)\rVert + \lVert \hat{R}_n^{\mathbf{X}}(\mathbf{X}_k)-R^{\mathbf{X}}(\mathbf{X}_k)\rVert+\lVert \hat{R}_n^{\mathbf{X}}(\mathbf{X}_l)-R^{\mathbf{X}}(\mathbf{X}_l)\rVert.
     	\end{align}
     	Starting from~\eqref{eq:indepconsispf7}, and repeating the same arguments as in~\eqref{eq:indepconsispf3},~\eqref{eq:indepconsispf4},~\eqref{eq:indepconsispf5} and~\eqref{eq:indepconsispf6}, we get:
     	\begin{align}\label{eq:indepconsispf8}
     	S_1\overset{n\to\infty}{\longrightarrow}\mathbb{E}\big[\lVert \mathbf{R}^{\mathbf{X}}(\mathbf{Z}_1^1)-R^{\mathbf{X}}(\mathbf{Z}_1^2)\rVert\lVert R^{\mathbf{Y}}(\mathbf{Z}_2^1)-R^{\mathbf{Y}}(\mathbf{Z}_2^2)\rVert\big]\qquad\qquad \;\;\mbox{a.s}.
     	\end{align}
     	Similar arguments applied to $S_2$ and $S_3$ yield the following:
     	\begin{align}\label{eq:indepconsispf9}
     	S_2\overset{n\to\infty}{\longrightarrow}\mathbb{E}\big[\lVert R^{\mathbf{X}}(\mathbf{Z}_1^1)-R^{\mathbf{X}}(\mathbf{Z}_1^2)\rVert\big] \mathbb{E}\big[\lVert R^{\mathbf{Y}}(\mathbf{Z}_2^1)-R^{\mathbf{Y}}(\mathbf{Z}_2^2)\rVert\big]\qquad \;\;\mbox{a.s}.
     	\end{align}
     	and,
     	\begin{align}\label{eq:indepconsispf10}
     	S_3\overset{n\to\infty}{\longrightarrow}\mathbb{E}\big[\lVert R^{\mathbf{X}}(\mathbf{Z}_1^1)-R^{\mathbf{X}}(\mathbf{Z}_1^2)\rVert\lVert R^{\mathbf{Y}}(\mathbf{Z}^1_2)-R^{\mathbf{Y}}(\mathbf{Z}^3_2)\rVert\big]\qquad \;\;\mbox{a.s}.
     	\end{align}
     	where $(\mathbf{Z}_1^1,\mathbf{Z}^1_2),(\mathbf{Z}^2_1,\mathbf{Z}_2^2),(\mathbf{Z}_1^3,\mathbf{Z}_2^3)$ are independent observations having the same distribution as $(\mathbf{X}_1,\mathbf{Y}_1)$. Finally, by combining ~\eqref{eq:indepconsispf8},~\eqref{eq:indepconsispf9},~\eqref{eq:indepconsispf10} and~\eqref{eq:proprdcovequiv}, we get:
     	\begin{equation}\label{eq:indepconsispf11}
     	\Rdcov_n^2=S_1+S_2-2S_3\overset{n\to\infty}{\longrightarrow}\Rdcov^2(\mathbf{X},\mathbf{Y}) \qquad \;\;\mbox{a.s}.
     	\end{equation}
     	which completes the proof of the first part.\par 
     	For the second part, note that,~\cref{theo:indepasdistn} implies $c_n=O(1)$. Also,~\cref{lem:proprdcov} (part (ii)) and~\eqref{eq:indepconsispf11} imply that whenever $\mu\neq \mu_{\mathbf{X}}\otimes \mu_{\mathbf{Y}}$, we will have:
     	\begin{align*}
     	n\Rdcov_n^2\overset{n\to\infty}{\longrightarrow}\infty \qquad \;\;\mbox{a.s}.
     	\end{align*}
     	which in turn, yields $\mathbb{P}(n\Rdcov_n^2>c_n)\overset{n\to\infty}{\longrightarrow}1$ if $\mu\neq\mu_{\mathbf{X}}\otimes\mu_{\mathbf{Y}}$, thereby completing the proof.
     \end{proof}
     \subsection{Proof of~\cref{rem:twosamdistfree}}\label{proof:twosamdistfree} 
     \begin{proof}
     	$\Ren_{m,n}^2$, as defined in~\eqref{eq:twogofp}, is a function of $(\hat{R}_{m,n}^{\mathbf{X},\mathbf{Y}}(\mathbf{X}_1),\ldots ,\hat{R}_{m,n}^{\mathbf{X},\mathbf{Y}}(\mathbf{X}_m),\hat{R}_{m,n}^{\mathbf{X},\mathbf{Y}}(\mathbf{Y}_1),\ldots ,\hat{R}_{m,n}^{\mathbf{X},\mathbf{Y}}(\mathbf{Y}_n))$. Note that the above vector is uniformly distributed over the set of all $(m+n)!$ permutations of the fixed grid $\mathcal{H}_{m+n}^d$ under $\mbox{H}_0$ and \textbf{(AP3)} (see~\cref{prop:finproperties}, part (ii)). This results in the distribution-free property of the statistic $\Ren_{m,n}^2$ (under $\mbox{H}_0$).
     \end{proof}
     \subsection{Proof of~\cref{lem:gofeqcvm}}\label{proof:gofeqcvm}
     \begin{proof}
     	Let us first prove \eqref{eq:gofeqcvmmain2}. Note that $\mathbb{P}\big(H^{X,Y}_{\lambda}(X)\leq t\big)=F^X\big((H_{\lambda}^{X,Y})^{-1}(t)\big)$ and $\mathbb{P}\big(H^{X,Y}_{\lambda}(Y)\leq t\big)=G^Y\big((H_{\lambda}^{X,Y})^{-1}(t)\big)$ for all $t$ in the support of $H_{\lambda}^{X,Y}$. By~\eqref{eq:ren1}, we have:
     	\begin{align*}
     	\frac{1}{2}\Ren_{\lambda}^2(X,Y)&=\int_0^1 \bigg(\mathbb{P}(H^{X,Y}_{\lambda}(X)\leq t)-\mathbb{P}(H^{X,Y}_{\lambda}(Y)\leq t)\bigg)^2\,dt\nonumber \\&=\int_0^1 \bigg(F^X\big((H_{\lambda}^{X,Y})^{-1}(t)\big)-G^Y\big((H_{\lambda}^{X,Y})^{-1}(t)\big)\bigg)^2\,dt\nonumber \\&=\int_{-\infty}^{\infty} \bigg(F^X(t)-G^Y(t)\bigg)^2\,dH^{X,Y}_{\lambda}(t)
     	\end{align*}
     	where the last line follows by a simple change of variable argument. This  proves~\eqref{eq:gofeqcvmmain2}.\par 
     	In order to prove~\eqref{eq:gofeqcvmmain1}, note that, by using~\cite[Equation 5]{Baringhaus2004}, we get:
     	\begin{equation}\label{eq:gofeqcvmpf2}
     	\frac{1}{2}\Ren_{m,n}^2=\int_{-\infty}^{\infty}\left(\frac{1}{m}\sum_{j=1}^m \mathbf{1}\big(\hat{R}_{m,n}^{X,Y}(X_j)\leq t)-\frac{1}{n}\sum_{j=1}^n \mathbf{1}\big(\hat{R}_{m,n}^{X,Y}(Y_j)\leq t)\right)^2\,dt.
     	\end{equation}
     	Clearly, for $t<(m+n)^{-1}$ or $t>1$, the right hand side of~\eqref{eq:gofeqcvmpf2} equals $0$. For $t\in [k(m+n)^{-1},(k+1)(m+n)^{-1})$, $k\in \{1,2,\ldots ,m+n-1\}$, observe that:
     	\begin{equation}\label{eq:gofeqcvmpf3}
     	\frac{1}{m}\sum_{j=1}^m \mathbf{1}\big(\hat{R}_{m,n}^{X,Y}(X_j)\leq t)=\frac{1}{m}\sum_{j=1}^m \mathbf{1}(X_j\leq Z_{(k)})
     	\end{equation}
     	where $Z_{(k)}$ denotes the $k$'th order statistic of the pooled sample $\{X_1,\ldots ,X_m,Y_1,\ldots ,Y_n\}$. A similar observation as in~\eqref{eq:gofeqcvmpf3} also holds for the $Y_j$'s. Consequently, by adding up over $k\in \{1,2,\ldots ,m+n\}$, we get that the right hand side of~\eqref{eq:gofeqcvmpf2} equals
     	\begin{align*}
     	&\;\;\;\;\frac{1}{m+n}\sum_{i=1}^{m+n}\left(\frac{1}{m}\sum_{j=1}^m \mathbf{1}(X_j\leq Z_{(i)})-\frac{1}{n}\sum_{j=1}^n \mathbf{1}(Y_j\leq Z_{(i)})\right)^2\nonumber \\&=\int_{-\infty}^{\infty} (F_m^X(t)-G_n^Y(t))^2\,dH_{m+n}(t)
     	\end{align*}
     	which completes the proof of~\eqref{eq:gofeqcvmmain1}.
     \end{proof}
     \subsection{Proof of~\cref{theo:twosamasdistn}}\label{proof:twosamasdistn}
     \begin{proof}
     	First we start with some notation: set $\mathcal{H}_{m+n}^d\coloneqq \{\mathbf{h}_1^d,\mathbf{h}_2^d,\ldots ,\mathbf{h}_{m+n}^d\}$ and
     	\begin{equation}\label{eq:maintwosamdef}
     	F_m^{\mathbf{a}}(r)\coloneqq\frac{1}{m}\sum_{j=1}^m \mathbf{1}\left(\mathbf{a}^{\top}\hat{R}_{m,n}^{\mathbf{X},\mathbf{Y}}(\mathbf{X}_j)\leq r\right)\quad \mbox{and}\quad G_n^{\mathbf{a}}(r)\coloneqq\frac{1}{n}\sum_{j=1}^n \mathbf{1}\left(\mathbf{a}^{\top}\hat{R}_{m,n}^{\mathbf{X},\mathbf{Y}}(\mathbf{Y}_j)\leq r\right)
     	\end{equation}
     	where $\mathbf{a}\in\mathcal{S}^{d-1}=\{\mathbf{z}:\lVert \mathbf{z}\rVert=1\}$, $r\in\mathbb{R}$ and $\mathbf{1}(\cdot)$ denotes the standard indicator function. Recall that $\kappa(\cdot)$ is the uniform measure on $\mathcal{S}^{d-1}$ and $\gamma_d=(2\Gamma(d/2))^{-1}\sqrt{\pi}(d-1)\Gamma((d-1)/2)$ for $d>1$ and $\gamma_d=1$ for $d=1$.\par 
     	Next, recall that $(\hat{R}_{m,n}^{\mathbf{X},\mathbf{Y}}(\mathbf{X}_1),\ldots ,\hat{R}_{m,n}^{\mathbf{X},\mathbf{Y}}(\mathbf{X}_m),\hat{R}_{m,n}^{\mathbf{X},\mathbf{Y}}(\mathbf{Y}_1),\ldots ,\hat{R}_{m,n}^{\mathbf{X},\mathbf{Y}}(\mathbf{Y}_n))$ is uniformly distributed over the set of all $(m+n)!$ permutations of the set $\mathcal{H}_{m+n}^d$ (see the proof of~\cref{rem:twosamdistfree} in~\cref{proof:twosamdistfree}). From~\eqref{eq:maintwosamdef}, this implies the following:
     	\begin{equation*}
     	\left(F_m^{\mathbf{a}}(r),G_n^{\mathbf{a}}(r)\right)\overset{d}{=}\left(\frac{1}{m}\sum_{j=1}^m \mathbf{1}\big(\mathbf{a}^{\top}\mathbf{h}^d_{\sigma_1(j)}\leq r\big),\frac{1}{n}\sum_{j=m+1}^{m+n}\mathbf{1}\big(\mathbf{a}^{\top}\mathbf{h}^d_{\sigma_1(j)}\leq r\big)\right)
     	\end{equation*}
     	for a random permutation $\sigma_1$ of the set $\{1,2,\ldots ,m+n\}$. Now, by~\cite[Lemma 2.3, Equation 5]{Baringhaus2004}, we further get:
     	\begin{align}\label{eq:twosamdistequivmain}
     	&\;\;\;\;\frac{mn}{m+n}\Ren_{m,n}^2\nonumber \\&\overset{d}{=}\frac{ mn}{m+n}\cdot\gamma_d\int\limits_{\mathcal{S}^{d-1}}\int_{-\infty}^{\infty}\left(\frac{1}{m}\sum_{j=1}^m \mathbf{1}\big(\mathbf{a}^{\top}\mathbf{h}^d_{\sigma_1(j)}\leq r\big)-\frac{1}{n}\sum_{j=m+1}^{m+n}\mathbf{1}\big(\mathbf{a}^{\top}\mathbf{h}^d_{\sigma_1(j)}\leq r\big)\right)^2 dr\;d\kappa(\mathbf{a}).
     	\end{align}
     	Let us take a look at~\cref{lem:rankenergy} below (see~\cref{proof:rankenergy} for a proof) which provides a general result on permutation statistics of the same form as in the right side of~\eqref{eq:twosamdistequivmain}.
     	\begin{lemma}\label{lem:rankenergy}
     		Consider an infinite sequence $\{\mathbf{U}_i\}_{i\geq 1}$, $\mathbf{U}_i\in\mathbb{R}^{d}$, such that $\mathbb{P}_{m+n}=(m+n)^{-1}\sum_{i=1}^{m+n}\delta_{\mathbf{U}_i} \overset{w}{\longrightarrow}\mathcal{U}^d$ (uniform distribution on $[0,1]^d$), as $\min{(m,n)}\to\infty$. Further, let $S_{N}$ denote the set of all permutations of $\{1,2,\ldots ,N\}$ and $\pi_N$ be a random permutation drawn uniformly from $S_N$. Define $\mathcal{S}^{d-1}=\{\mathbf{x}\in\mathbb{R}^d:\lVert \mathbf{x}\rVert =1\}$ and $\Theta_{m,n}:\mathcal{S}^{d-1}\times \mathbb{R}\mapsto\mathbb{R}$ by,
     		\begin{align*}
     		\Theta_{m,n}(\mathbf{a},r)=\sqrt{\frac{mn}{m+n}}\cdot \left(\frac{1}{m}\sum_{i=1}^m \mathbf{1}(\mathbf{a}^{\top}\mathbf{U}_{\pi_N(i)}\leq r)-\frac{1}{n}\sum_{j=m+1}^{m+n}\mathbf{1}(\mathbf{a}^{\top}\mathbf{U}_{\pi_N(j)}\leq r)\right).
     		\end{align*}
     		Then we have,
     		\begin{align*}
     		E_{m,n}\coloneqq \int_{\mathcal{S}^{d-1}}\int_{\mathbb{R}}\Theta_{m,n}^2(\mathbf{a},r)\,dr\,d\kappa(\mathbf{a})\overset{w}{\longrightarrow}\int_{\mathcal{S}^{d-1}}\int_{\mathbb{R}}\Theta^2(\mathbf{a},r)\,dr\,d\kappa(\mathbf{a})\coloneqq E
     		\end{align*}
     		as $\min{(m,n)}\to\infty$. Here $\kappa(\cdot)$ denotes the uniform measure on $\mathcal{S}^{d-1}$ and $\Theta(\cdot,\cdot)$ is a mean $0$ Gaussian process with covariance kernel given by, 
     		\begin{align}\label{eq:covenergy}
     		C\big((\mathbf{a}_1,r_1),(\mathbf{a}_2,r_2)\big)\coloneqq \mathbb{P}(\mathbf{a}_1^{\top}\tilde{\mathbf{U}}\leq r_1, \mathbf{a}_2^{\top}\tilde{\mathbf{U}}\leq r_2)-\mathbb{P}(\mathbf{a}_1^{\top}\tilde{\mathbf{U}}\leq r_1)\cdot \mathbb{P}(\mathbf{a}_2^{\top}\tilde{\mathbf{U}}\leq r_2),
     		\end{align} 
     		for $(\mathbf{a}_1,r_1)$, $(\mathbf{a}_2,r_2)\in \mathcal{S}^{d-1}\times\mathbb{R}$, where $\tilde{\mathbf{U}}$ has the same distribution as $\mathcal{U}^d$.
     	\end{lemma}
     	With the above lemma in mind, note that the empirical distribution on $\mathcal{H}_n^d$ converges weakly to $\mathcal{U}^{d}$ under assumption \textbf{(AP4)}. By setting $\{\mathbf{U}_j\}_{j\geq 1}$ as $\{\mathbf{h}_j^d\}_{j\geq 1}$, and applying~\cref{lem:rankenergy} to the right hand side of~\eqref{eq:twosamdistequivmain}, we get that, as $\min{(m,n)}\to\infty$, 
     	\begin{equation}\label{eq:twosamasdistmain}
     	\frac{mn}{m+n}\Ren_{m,n}^2\overset{w}{\longrightarrow}\gamma_d E\overset{d}{=} \sum_{j=1}^{\infty}\tau_jZ_j^2
     	\end{equation}
     	where $\tau_j$'s are fixed nonnegative constants and $Z_j$'s are i.i.d.~standard normals. The last equality in~\eqref{eq:twosamasdistmain} follows from~\cite[Chapter 1, Section 2]{Kuo1975} using standard Karhunen-Lo\`eve type expansions of Gaussian processes. This completes the proof.
     \end{proof}
     
     \subsection{Proof of~\cref{theo:twosamconsis}}\label{proof:twosamconsis}
     \begin{proof}
     	Note that $\Ren_{m,n}^2$ is the average of products of Euclidean distances between bounded random vectors, in the same spirit as $\Rdcov_n^2$. Therefore, The proof is exactly similar to the proof of~\cref{theo:indepconsis} in~\cref{proof:indepconsis}. We leave the details to the reader.
     \end{proof}
     \subsection{Proof of~\cref{prop:multextind}}\label{proof:multextind}
     \begin{proof}
     	By the same argument as in~\cref{proof:finproperties}, we get that the rank vectors $(\hat{R}_n^j(\mathbf{X}^j_1),\ldots,\hat{R}_n^j(\mathbf{X}^j_n))$ are uniformly distributed over all $n!$ possible permutations of the set $\mathcal{H}_n^{d_j}$. Also these rank vectors are independent over $j\in \{1,2,\ldots ,K\}$ under $\mathrm{H}_0$. This implies that our test statistic, which is a function of $(\hat{R}_n^j(\mathbf{X}^j_1),\ldots,\hat{R}_n^j(\mathbf{X}^j_n))_{1\leq j\leq K}$, is distribution-free under $\mathrm{H}_0$. Consequently, the threshold $c_n$, is also distribution-free under $\mathrm{H}_0$.\par 
     	Finally, by~\cref{theo:indepconsis}, we see that $\Rdcov_{n,j}^2$ converges a.s. to the following quantity:
     	\begin{equation}\label{eq:multextindpf1}
     	\Rdcov_n^2\overset{a.s.}{\longrightarrow}\sum_{j=1}^{K-1}\Rdcov^2_{*}(\mathbf{X}^j,\mathbf{X}^{j+})
     	\end{equation}
     	where $\Rdcov^2_{*}(\mathbf{X}^j,\mathbf{X}^{j+})$ is the usual distance covariance between $R_j(\mathbf{X}^j)$ and $(R_{j+1}(\mathbf{X}^{j+1}),\ldots ,R_k(\mathbf{X}^k))$ and $R_j(\cdot)$ is the population rank map for $\mu_j$, $1\leq j\leq K$.\par 
     	Next, note that under $\mathrm{H}_0$, $n\Rdcov_{n,j}^2=\mathcal{O}_p(1)$ for $1\leq j\leq K-1$ by~\cref{theo:indepasdistn} (the same second moment computations as in~\eqref{eq:covfnver} and~\eqref{eq:secmoment}). Therefore, $n\Rdcov_n^2=\mathcal{O}_p(1)$ which implies $c_n=\mathcal{O}(1)$.\par
     	Now, under $\mathrm{H}_1$, we will show that the right side of~\eqref{eq:multextindpf1} is strictly positive. In this direction, let us proceed by contradiction. Suppose that the right side of~\eqref{eq:multextindpf1} equals $0$ under $\mathrm{H}_1$. This implies $\Rdcov^2_{*}(\mathbf{X}^j,\mathbf{X}^{j+})=0$ for all $1\leq j\leq K-1$. Therefore, by the same argument as in~\cref{lem:proprdcov} (part (b)), $\mathbf{X}^{j}$ and $\mathbf{X}^{j+}$ are independent for $1\leq j\leq K-1$. An alternate way of stating this would be as follows:
     	\begin{align}\label{eq:multextindpf2}
     	\mathbb{E}\left[\exp\left(i\sum_{l=j}^K \mathbf{t}_l^{\top}\mathbf{X}^l\right)\right]=\mathbb{E}\left[\exp\left(i\mathbf{t}_j^{\top}\mathbf{X}^j\right)\right]\times\mathbb{E}\left[\exp\left(i\sum_{l=j+1}^K \mathbf{t}_l^{\top}\mathbf{X}^l\right)\right]
     	\end{align}
     	for all $(\mathbf{t}^1,\ldots ,\mathbf{t}_K)\in \mathbb{R}^{d_1}\times \ldots \times \mathbb{R}^{d_k}$ and $1\leq j\leq K-1$. Next, by~\cite{matteson2014}, we have:
     	\begin{align}\label{eq:multextindpf3}
     	&\;\;\;\;\Bigg|\mathbb{E}\left[\exp\left(i\sum_{l=1}^K \mathbf{t}_l^{\top}\mathbf{X}^l\right)\right]-\prod\limits_{l=1}^K \mathbb{E}\left[\exp\left(i\mathbf{t}_l^{\top}\mathbf{X}^l\right)\right]\Bigg|\nonumber \\&\leq \sum_{j=1}^{K-1}\Bigg|\mathbb{E}\left[\exp\left(i\sum_{l=j}^K \mathbf{t}_l^{\top}\mathbf{X}^l\right)\right]-\mathbb{E}\left[\exp\left(i\mathbf{t}_j^{\top}\mathbf{X}^j\right)\right]\times\mathbb{E}\left[\exp\left(i\sum_{l=j+1}^K \mathbf{t}_l^{\top}\mathbf{X}^l\right)\right]\Bigg|.
     	\end{align}
     	The right hand side of~\eqref{eq:multextindpf3} equals $0$ by~\eqref{eq:multextindpf2}, and consequently the left hand side of~\eqref{eq:multextindpf3} equals $0$. Therefore $\mu=\mu_1\otimes \ldots \otimes \mu_K$, which contradicts $\mathrm{H}_1$ and consequently proves the claim. As a result we also have $n\Rdcov_n^2\overset{a.s.}{\longrightarrow}\infty$ under $\mathrm{H}_1$. As $c_n=\mathcal{O}(1)$, we have
     	\begin{align*}
     	\mathbb{P}(n\Rdcov_n^2\geq c_n)\overset{n\to\infty}{\longrightarrow}1
     	\end{align*}
     	which completes the proof of consistency.
     	\subsection{Proof of~\cref{prop:multextgof}}\label{proof:multextgof}
     	Note that the set of $n$ vectors $(\hat{R}_n(\mathbf{X}^j_1),\ldots ,\hat{R}_n(\mathbf{X}^j_{n_j}))_{1\leq j\leq K}$ are uniformly distributed over the $n!$ permutations of the set $\mathcal{H}_n^d$ (same argument as in~\cref{prop:finproperties}). This implies that the test statistic $\Ren_{1:K,n}^2$ is distribution-free and consequently, so is $c_n$.\par 
     	By the same argument as in the proof of~\cref{theo:twosamconsis}, we have:
     	\begin{align}\label{eq:multextgofpf1}
     	\Ren_{1:K,n}^2=\sum_{j=1}^{K-1}\Ren_{j:(j+1),n}^2\overset{a.s.}{\longrightarrow}\sum_{j=1}^{K-1}\Ren^2_{\lambda,*}(\mathbf{X}^j,\mathbf{X}^{j+1})
     	\end{align}
     	where $\Ren^2_{\lambda,*}(\mathbf{X}^j,\mathbf{X}^{j+1})$ is the usual energy distance between $R_{\lambda}(\mathbf{X}^j)$ and $R_{\lambda}(\mathbf{X}^{j+1})$, $R_{\lambda}(\cdot)$ denotes the population rank map corresponding to the measure $\sum_{j=1}^K \lambda_j\mu_j$.\par 
     	Under $\mathrm{H}_0$, by similar second moment calculations as in the proof of~\cref{theo:twosamasdistn}, we have $n\Ren^2_{1:K,n}=\mathcal{O}_p(1)$ which in turn implies $c_n=\mathcal{O}(1)$.\par 
     	Under $\mathrm{H}_1$, there exists $\tilde{j}\leq K-1$ such that $\mu_{\tilde{j}}\neq \mu_{\tilde{j}+1}$ which implies $\Ren^2_{\lambda,*}(\mathbf{X}^{\tilde{j}},\mathbf{X}^{\tilde{j}+1})>0$ and consequently the right hand side of~\eqref{eq:multextgofpf1} is strictly positive. Therefore, $n\Ren_{1:K,n}^2\overset{a.s.}{\longrightarrow}\infty$ and as $c_n=\mathcal{O}(1)$, we have $$\mathbb{P}(n\Ren^2_{1:K,n}\geq c_n)\overset{n\to\infty}{\longrightarrow}1$$ which completes the proof of consistency.
     \end{proof}
     \subsection{Proof of~\cref{lem:multsym}}\label{proof:multsym}
     
     \begin{proof}
     	Note that, by the independence of the $\mathbf{Z}_i$'s, $\{\tilde{R}_n(\mathbf{Z}_i)\}_{1\leq i\leq n}$ is distributed uniformly over the set of all $n!$ permutations of $\mathcal{H}_n^{2d}$. Moreover, under $\textrm{H}_0$, $\mathbf{Z}_1\overset{d}{=}-\mathbf{Z}_1$. Based on this observation, it is easy to check that:
     	\begin{equation*}
     	\mathbb{P}\big(\hat{R}_n(\mathbf{X}_1)=r_1,\hat{R}_n(-\mathbf{X}_1)=r_2,\ldots ,\hat{R}_n(\mathbf{X}_n)=r_{2n-1},\hat{R}_n(-\mathbf{X}_n)=r_{2n}\big)=\frac{1}{2^n n!}
     	\end{equation*}
     	for $(r_1,r_2,\ldots ,r_{2n-1},r_{2n})\in \mathcal{J}$ where
     	\begin{align*}
     	\mathcal{J}&\coloneqq\big\{(j_1,\ldots ,j_{2n}):j_i\in\mathbb{R}^d, (j_i, j_{i+1})=\tilde{h}_i\textrm{ or }(j_{i+1}, j_i)=\tilde{h}_i\textrm{ for }i\in \{1,3,\ldots ,2n-1\}, \nonumber \\&\qquad\qquad (\tilde{h}_1,\ldots ,\tilde{h}_n)\textrm{ is a permutation of the set }\mathcal{H}_n^{2d}\big\}.
     	\end{align*}
     \end{proof}
     \section{Some general results on permutation statistics}\label{sec:permappen}
     In this section, we will prove some general results on asymptotic distributions of certain permutation based statistics which were used in~\cref{sec:appen}. Since the distinction between vectors and scalars is pretty clear in this section, we will drop the boldface fonts (previously used to denote vectors) subsequently for notational convenience. 
     \subsection{Proof of~\cref{lem:rankdcov}}\label{proof:rankdcov}
     Recall the notation in~\cref{lem:rankdcov}. Note that $D$ (in~\eqref{eq:fluctuation}) is well-defined (see~\cite[Chapter 1, Section 2]{Kuo1975}) and has finite expectation. For any $c>0$, define the following:
     \begin{align*}
     W_{n,c}\coloneqq\int_{1/c\leq \lVert t\rVert,\lVert s\rVert\leq c}\big|\xi_n(t,s)\big|^2w(t,s)\,dt\,ds,\;\;\;\;W_c=\int_{1/c\leq \lVert t\rVert,\lVert s\rVert\leq c}\big|\xi(t,s)\big|^2w(t,s)\,dt\,ds.
     \end{align*}
     Our proof will proceed through the following three steps:
     \begin{lemma}\label{lem:doublimit}
     	For any $\delta>0$, $\lim_{c\to\infty}\limsup_{n\to\infty}\mathbb{P}(|W_{n,c}-D_n|>\delta)=0$.
     \end{lemma}
     \begin{prop}\label{lem:singlim}
     	For any $\delta>0$, $\lim_{c\to\infty}\mathbb{P}(|W_{c}-D|>\delta)=0$.
     \end{prop}
     \begin{lemma}\label{lem:weaklim}
     	For any $c>0$, $W_{n,c}\overset{w}{\longrightarrow}W_c$ as $n\to\infty$.
     \end{lemma}
     \noindent Combining~\cref{lem:singlim}, Lemmas~\ref{lem:doublimit} and~\ref{lem:weaklim} with~\cite[Lemma 2.5]{Sen2010}, yields $D_n\overset{w}{\longrightarrow}D$ as $n\to\infty$ and completes the proof.
     
     \begin{proof}[Proof of~\cref{lem:weaklim}]
     	For $(t,s)\in\mathbb{R}^{d_1}\times\mathbb{R}^{d_2}$, define:
     	\begin{align*}
     	f_U^n(t)\coloneqq\frac{1}{n}\sum_{k=1}^n\exp\left(it^{\top}U_k\right)\;\; ,\;\;f_V^n(t)\coloneqq\frac{1}{n}\sum_{k=1}^n\exp\left(it^{\top}V_k\right)
     	\end{align*}
     	In the following, $U_k$'s and $V_k$'s are fixed, so the expectations are taken with respect to the randomness arising out of the randomly drawn permutation $\pi_n\in S_n$. Note that, for each $(t_1,s_1), (t_0,s_0)$, we have $\mathbb{E}[\xi_n(t_1,s_1)]=0$ and 
     	\begin{align}\label{eq:covfnver}
     	\mathbb{E}\left[\xi_n(t_1,s_1)\overline{\xi_n(t_0,s_0)}\right]&=\frac{n}{n-1}\left(f_U^n(t_1-t_0)-f_U^n(t_1)\overline{f_U^n(t_0)}\right)\left(f_V^n(s_1-s_0)-f_V^n(s_1)\overline{f_V^n(s_0)}\right)\nonumber\\ &\overset{n\to\infty}{\longrightarrow}\left(f_{d_1}(t_1-t_0)-f_{d_1}(t_1)\overline{f_{d_1}(t_0)}\right)\left(f_{d_2}(s_1-s_0)-f_{d_2}(s_1)\overline{f_{d_2}(s_0)}\right)\nonumber \\ &=R((t_1,s_1),(t_0,s_0))
     	\end{align}
     	where the convergence follows from the assumption that $\mathbb{P}_n\overset{w}{\longrightarrow}\mathcal{U}^{d_1}$ and $\mathbb{Q}_n\overset{w}{\longrightarrow}\mathcal{U}^{d_2}$.
     	
     	Now fix $M\geq 1$, and two sequences of real numbers $(\alpha_1,\ldots ,\alpha_M)$ and $(\beta_1,\ldots,\beta_M)$. For $\{(t_m,s_m)\in\mathbb{R}^{d_1}\times\mathbb{R}^{d_2}\}_{m=1}^M$, define:
     	\begin{align*}
     	&\Lambda_n(t_m,s_m)\coloneqq\sqrt{n}\left[\frac{1}{n}\sum_{k=1}^n \cos\big(t_m^{\top}U_k+s_m^{\top}V_{\pi_n(k)}\big)-\frac{1}{n^2}\sum_{k,l}\cos\big(t_m^{\top}U_k+s_m^{\top}V_l\big)\right]\nonumber\\ &\Theta_n(t_m,s_m)\coloneqq\sqrt{n}\left[\frac{1}{n}\sum_{k=1}^n \sin\big(t_m^{\top}U_k+s_m^{\top}V_{\pi_n(k)}\big)-\frac{1}{n^2}\sum_{k,l}\sin\big(t_m^{\top}U_k+s_m^{\top}V_l\big)\right].
     	\end{align*}
     	Note that $\Lambda_n$ and $\Theta_n$ are centered random variables. Further, define two matrices $(A_{kl})_{n\times n}$ and $(C_{kl})_{n\times n}$ as:
     	\begin{align*}
     	&A_{kl}\coloneqq\frac{1}{\sqrt{n}}\sum_{m=1}^M\left[\alpha_m\cos\big(t_m^{\top}U_k+s_m^{\top}V_l\big)+\alpha_m\sin\big(t_m^{\top}U_k+s_m^{\top}V_l\big)\right],\\  &C_{kl}\coloneqq A_{kl}-\overline{A_{k.}}-\overline{A_{.l}}+A_{..},
     	\end{align*}
     	and note that $\sum_{k=1}^n C_{k\pi_n(k)}=\sum_{m=1}^M \big(\alpha_m\Lambda_n(t_m,s_m)+\beta_m\Theta_n(t_m,s_m)\big)$. By Hoeffding's Combinatorial Central Limit theorem (see e.g.,~\cite[Theorem 1.1]{Chen2015}), we will have $\sum_{k=1}^n C_{k\pi_n(k)}\overset{w}{\longrightarrow}\mathcal{N}(0,\sigma^2)$ if the following two conditions hold: (a) $\sum_{k,l} |C_{kl}|^3=o(n)$ and (b) $\mbox{Var}\big(\sum_{k=1}^n C_{k\pi_n(k)}\big)\overset{n\to\infty}{\longrightarrow}\sigma^2$. As sine and cosine functions are bounded, and $U_k$'s, $V_k$'s all lie in fixed (free of $n$) compact sets, it is easy to check that $\sum_{k,l} |C_{ij}|^3=O\big(\sqrt{n}\big)$, so condition (a) is satisfied. For verifying (b), let us first set up some prerequisites. Observe that:
     	\begin{align*}	&\;\;\;\;\;\mathbb{E}\left[\Lambda_n(t_{1},s_{1})\Theta_n(t_{0},s_{0})\right]\nonumber \\ &=\frac{1}{n}\sum_{k,l}\mbox{Cov}\bigg(\cos\big(t_{1}^{\top}U_k+s_{1}^{\top}V_{\pi_n(k)}\big),\sin\big(t_{0}^{\top}U_{l}+s_{0}^{\top}V_{\pi_n(l)}\big)\bigg)\nonumber\\ &=\frac{1}{n^2}\sum_{k,l} \cos\big(t_{1}^{\top}U_k+s_{1}^{\top}V_l\big)\sin\big(t_{0}^{\top}U_k+s_{0}^{\top}V_l\big)+\frac{1}{n^3(n-1)}\sum_{k\neq l,p\neq q}\cos\big(t_{1}^{\top}U_k+s_{1}^{\top}V_p\big)\nonumber \\&\qquad\qquad \times \sin\big(t_{0}^{\top}U_l+s_{0}^{\top}V_q\big)-\frac{1}{n^3}\sum_{p,k\neq l}\cos\big(t_1^{\top}U_k+s_1^{\top}V_p\big)\sin\big(t_0^{\top}U_l+s_0^{\top}V_p\big)\nonumber \\&\qquad\qquad -\frac{1}{n^3}\sum_{k,p\neq q}\cos\big(t_1^{\top}U_k+s_1^{\top}V_p\big)\sin\big(t_0^{\top}U_k+s_0^{\top}V_q\big)\nonumber\\ &\overset{n\to\infty}{\longrightarrow}\mathbb{E}\Big[\cos\big(t_1^{\top}\tilde{U}_1+s_1^{\top}\tilde{V}_1\big)\sin\big(t_0^{\top}\tilde{U}_1+s_0^{\top}\tilde{V}_1\big)+\cos\big(t_1^{\top}\tilde{U}_1+s_1^{\top}\tilde{V}_1\big)\sin\big(t_0^{\top}\tilde{U}_2+s_0^{\top}\tilde{V}_2\big)\nonumber\\ &\;\;-\cos\big(t_1^{\top}\tilde{U}_1+s_1^{\top}\tilde{V}_1\big)\sin\big(t_0^{\top}\tilde{U}_1+s_0^{\top}\tilde{V}_2\big)-\cos\big(t_1^{\top}\tilde{U}_1+s_1^{\top}\tilde{V}_1\big)\sin\big(t_0^{\top}\tilde{U}_2+s_1^{\top}\tilde{V}_1\big)\Big].
     	\end{align*}
     	Here $\tilde{U}_1,\tilde{U}_2$ are i.i.d.~$\mathcal{U}^{d_1}$ and $\tilde{V}_1,\tilde{V}_2$ are i.i.d.~$\mathcal{U}^{d_2}$. A similar calculation shows that, for any $m_1,m_2\in \{1,2,\ldots ,M\}$, $\mathbb{E}[\Lambda_n^2(t_{m_1},s_{m_1})]$, $\mathbb{E}[\Theta_n^2(t_{m_1},s_{m_1})]$, $\mathbb{E}[\Theta_n(t_{m_1},s_{m_1})\Theta_n(t_{m_2},s_{m_2})]$,  $\mathbb{E}[\Lambda_n(t_{m_1},s_{m_1})\Lambda_n(t_{m_2},s_{m_2})]$ and $\mathbb{E}[\Lambda_n(t_{m_1},s_{m_1})\Theta_n(t_{m_2},s_{m_2})]$ also converges. Denote the corresponding limits as $\Lambda_{m_1}^2$, $\Theta_{m_1}^2$, $\Theta_{m_1 m_2}$, $\Lambda_{m_1 m_2}$ and $\Lambda_{m_1}\Theta_{m_2}$ respectively. This implies that $\mbox{Var}\big(\sum_{k=1}^n C_{k\pi_n(k)}\big)$ equals,
     	\begin{align*}
     	&\;\;\mathbb{E}\Bigg[\sum_{m=1}^M \alpha_m^2\Lambda_n^2(t_m,s_m)+\sum_{m=1}^M \beta_m^2\Theta_n^2(t_m,s_m)+\sum_{m_1\neq m_2} \alpha_{m_1}\alpha_{m_2}\Lambda_n(t_{m_1},s_{m_1})\Lambda_n(t_{m_2},s_{m_2})\nonumber\\ &+\sum_{m_1\neq m_2} \beta_{m_1}\beta_{m_2}\Theta_n(t_{m_1},s_{m_1})\Theta_n(t_{m_2},s_{m_2})+\sum_{m_1,m_2}\alpha_{m_1}\beta_{m_2}\Lambda_n(t_{m_1},s_{m_1})\Theta_n(t_{m_2},s_{m_2})\Bigg]\nonumber\\ &\overset{n\to\infty}{\longrightarrow}\sum_{m=1}^M \alpha_m^2\Lambda_m^2+\sum_{m=1}^M \beta_m^2\Theta_m^2+\sum_{m_1\neq m_2} \alpha_{m_1}\alpha_{m_2}\Lambda_{m_1 m_2}+\sum_{m_1\neq m_2} \beta_{m_1}\beta_{m_2}\Theta_{m_1 m_2}\nonumber \\&\qquad\qquad+\sum_{m_1,m_2}\alpha_{m_1}\beta_{m_2}\Lambda_{m_1}\Theta_{m_2}.
     	\end{align*}
     	This completes the proof of (b) and therefore, $\sum_{k=1}^n C_{k\pi_n(k)}$ converges to a Gaussian limit. By the Cram\'er-Wold theorem, the vector $$\Gamma_n\coloneqq(\Lambda_n(t_1,s_1),\ldots ,\Lambda_n(t_M,s_M),\Theta_n(t_1,s_1),\ldots ,\Theta_n(t_M,s_M))$$ converges to a multivariate Gaussian distribution. Take $f:\mathbb{R}^{2M}\to\mathbb{C}^M$ such that $f(x_1,\ldots ,x_{2M})=(x_1+ix_{M+1},\ldots ,x_{M}+ix_{2M})$ where $i=\sqrt{-1}$. Then, by the continuous mapping theorem, $f(\Gamma_n)=(\xi_n(t_1,s_1),\ldots ,\xi_n(t_M,s_M))$ converges to a complex-valued Gaussian process with covariance kernel given by $R(\cdot,\cdot)$, as shown in~\eqref{eq:covfnver}. \newline 
     	For $c>0$, define $A_c=\{t\in\mathbb{R}^{d_1}:(1/c)\leq \lVert t\rVert\leq c\}$ and $B_c=\{s\in\mathbb{R}^{d_2}:(1/c)\leq \lVert s\rVert\leq c\}$. Note that $A_c\times B_c$ is compact. The preceding discussion yields the convergence for the finite dimensional distributions of the process $\xi_n(\cdot,\cdot)$ over $A_c\times B_c$. In order to show $\xi_n(\cdot,\cdot)\overset{w}{\longrightarrow}\xi(\cdot,\cdot)$ in $L^{\infty}[A_c\times B_c]$ (in the sense of~\cite{Hoffman1991}), we would need to show asymptotic equicontinuity (see~\cite{Hoffman1991}), i.e., given any $\epsilon>0$,
     	\begin{equation}\label{eq:asequicon}
     	\lim_{\delta\to 0}\limsup_{n\to\infty}\mathbb{P}^*\left(\sup_{\substack{s_1,s_0\in B_c,\ t_1,t_0\in A_c,\\ \sqrt{\lVert s_1-s_0\rVert^2+\lVert t_1-t_0\rVert^2}\leq \delta}}\big|\xi_n(t_1,s_1)-\xi_n(t_0,s_0)\big|>\epsilon\right)=0
     	\end{equation}
     	where $\mathbb{P}^*$ denotes the outer probability. Note that $\xi_n(t,s)=\sum_{k=1}^n Z_{nk}(t,s)$ where $Z_{nk}(t,s)\coloneqq n^{-1/2}\big[\exp\big(it^{\top}U_k+is^{\top}V_{\pi_n(k)}\big)-f_U^n(t)f_V^n(s)\big]$. Let $\rho((t_1,s_1),(t_0,s_0))\coloneqq \sqrt{\lVert t_1-t_0\rVert^2+\lVert s_1-s_0\rVert^2}$. From the proof of~\cite[Theorem 2.11.1]{vaart1996},~\eqref{eq:asequicon} follows if we can show the following:
     	\begin{enumerate}
     		\item[Step 1.] There exists a sequence $\eta_n\downarrow 0$ such that $\sup_{1\leq k\leq n}\sup_{(t,s)\in A_c\times B_c}\big|Z_{nk}(t,s)\big|\leq \eta_n$.
     		\item[Step 2.] For any sequence $\delta_n\downarrow 0$, we have:
     		\begin{align*}
     		\sup_{s_1,s_0\in B_c,\ t_1,t_0\in A_c, \rho((t_1,s_1),(t_0,s_0))\leq \delta_n} \mathbb{E}^*\bigg|\sum_{k=1}^n \big[Z_{nk}(t_1,s_1)-Z_{nk}(t_0,s_0)\big]\bigg|^2\overset{n\to\infty}{\longrightarrow}0.
     		\end{align*}
     		\item[Step 3.] For any sequence $\delta_n\downarrow 0$, $\int_0^{\delta_n} \sqrt{\log{N(\epsilon,A_c\times B_c,d_n(\cdot,\cdot))}}\,d\epsilon \overset{\mathbb{P}^*}{\longrightarrow}0$ where $N(\epsilon,A_c\times B_c,d_n(\cdot,\cdot))$ denotes the $\epsilon$ covering number of the set $A_c\times B_c$ based on the random metric $d_n(\cdot,\cdot)$ satisfying $d_n^2((t_1,s_1),(t_0,s_0))=\sum_{k=1}^n |Z_{nk}(t_1,s_1)-Z_{nk}(t_0,s_0)|^2$ (for related definitions in this context, see~\cite[Chapter 2.2]{vaart1996}).
     	\end{enumerate}
     	Note that $|Z_{nk}(\cdot,\cdot)|\leq 2n^{-1/2}$ a.s. and so \textbf{step 1} holds. Note that, for any $s_1,s_0\in B_c$, $t_1,t_0\in A_c$ and $\rho((t_1,s_1),(t_0,s_0))\leq \delta_n$, we have:
     	\begin{align}\label{eq:secmoment}
     	&\;\;\mathbb{E}^*\bigg|\sum_{k=1}^n \big[Z_{nk}(t_1,s_1)-Z_{nk}(t_0,s_0)\big]\bigg|^2\nonumber \\&=\frac{n}{n-1}\Bigg[\bigg(1-f_U^n(t_0-t_1)f_V^n(s_0-s_1)\bigg)+\bigg(1-f_U^n(t_1-t_0)f_V^n(s_1-s_0)\bigg)\nonumber \\&\qquad+\bigg(f_U^n(t_0)f_V^n(s_0)\big(\overline{f_U^n(t_0)}\overline{f_V^n(s_0)}-\overline{f_U^n(t_1)}\overline{f_V^n(s_1)}\big)\bigg)+\bigg(f_U^n(t_1)f_V^n(s_1)\big(\overline{f_U^n(t_1)}\overline{f_V^n(s_1)}\nonumber \\&\qquad-\overline{f_U^n(t_0)}\overline{f_V^n(s_0)}\big)\bigg)\nonumber -\bigg(f_U^n(t_1)\big(\overline{f_U^n(t_1)}-\overline{f_V^n(s_1-s_0)}\overline{f_U^n(t_0)}\big)\bigg)\nonumber \\&\qquad-\bigg(f_U^n(t_0)\big(\overline{f_U^n(t_0)}-\overline{f_V^n(s_0-s_1)}\overline{f_U^n(t_1)}\big)\bigg)-\bigg(f_V^n(s_1)\big(\overline{f_V^n(s_1)}-\overline{f_U^n(t_1-t_0)}\overline{f_V^n(s_0)}\big)\bigg)\nonumber \\&\qquad-\bigg(f_V^n(s_0)\big(\overline{f_V^n(s_0)}-\overline{f_U^n(t_0-t_1)}\overline{f_V^n(s_1)}\big)\bigg)\Bigg].
     	\end{align} 
     	As $\sin(\cdot)$ and $\cos(\cdot)$ are Lipschitz functions with Lipschitz norm bounded by $1$, each term within a parenthesis on the right hand side of~\eqref{eq:secmoment} may be bounded in modulus by $4\delta_n$. This completes the proof of \textbf{step 2}. Further, once again using the Lipschitz nature of $\cos(\cdot)$ and $\sin(\cdot)$, it is easy to check that:
     	\begin{align}\label{eq:covnum1}
     	d_n((t_1,s_1),(t_0,s_0))=\sqrt{\sum_{k=1}^n |Z_{nk}(t_1,s_1)-Z_{nk}(t_0,s_0)|^2}\leq 10\rho((t_1,s_1),(t_0,s_0))
     	\end{align}
     	where the last event happens with $\mathbb{P}^*$ outer probability $1$. By using $\lesssim_c$ to hide constants that only depend on $p,q$ and $c$, we get the following chain of inequalities:
     	\begin{align*}
     	\int_0^{\delta_n} \sqrt{\log{N(\epsilon,A_c\times B_c,d_n(\cdot,\cdot))}}\,d\epsilon &\overset{(a)}{\leq}\int_0^{\delta_n} \sqrt{\log{N(\epsilon/10,A_c\times B_c,\rho(\cdot,\cdot))}}\,d\epsilon\nonumber \\ &\overset{(b)}{\lesssim_c}\int_0^{\delta_n}\epsilon^{-1/2}\,d\epsilon\overset{n\to\infty}{\longrightarrow}0
     	\end{align*}
     	where (a) happens with $\mathbb{P}^*$ outer probability $1$ and follows from~\eqref{eq:covnum1}, (b) follows from a standard volumetric argument for estimating covering numbers, see e.g.,~\cite[Lemma 4.5]{vande2000}. This completes the proof of \textbf{step 3}. Therefore, by combining \textbf{steps 1, 2 and 3}, we get $\xi_n(\cdot,\cdot)\overset{w}{\longrightarrow}\xi(\cdot,\cdot)$ in $L^{\infty}[A_c\times B_c]$. Finally, note that $w(\cdot,\cdot)$ is bounded in $A_c\times B_c$. ~\cref{lem:weaklim} then follows by the continuous mapping theorem with the integral (over $A_c\times B_c$) operator.
     \end{proof}
     \begin{proof}[Proof of~\cref{lem:doublimit}]
     	Define $\mathcal{B}_{d_1}(1)=\{z\in\mathbb{R}^{d_1}:\lVert z\rVert \leq 1\}$ and a function $G:(0,\infty)\times\mathcal{B}_{d_1}(1)\to \mathbb{R}$ as,
     	\begin{align}\label{eq:defg}
     	G(y,w)\coloneqq \int_{\lVert z\rVert\leq y} \frac{1-\cos\langle w,z\rangle}{\lVert z\rVert ^{1+d_1}}\,dz.
     	\end{align}
     	By~\cite[Lemma 1]{Gabor2007}, $G(\cdot,\cdot)$ is uniformly bounded, and by an application of the dominated convergence theorem, $\lim_{\delta\downarrow 0}G(y,w)=0$, for each $w\in\mathcal{B}_{d_1}(1)$. Next note that, for any $c>1$,
     	\begin{align}\label{eq:parineq}
     	|W_{n,c}-D_n|&\leq \int_{\{\lVert t\rVert \leq 1/c\}\cup \{\lVert t\rVert \geq c\}} |\xi_n(t,s)|^2w(t,s)\,dt\,ds\nonumber \\&\qquad +\int_{\{\lVert s\rVert \leq 1/c\}\cup \{\lVert s\rVert \geq c\}} |\xi_n(t,s)|^2w(t,s)\,dt\,ds.
     	\end{align}
     	We will use $\lesssim$ to hide constants which depend only on $d_1$ and $d_2$. Therefore,
     	\begin{align}\label{eq:boundf1}
     	&\mathbb{E}\int_{\lVert t\rVert \leq 1/c} |\xi_n(t,s)|^2w(t,s)\,dt\,ds\nonumber \\ &\overset{(a)}{\lesssim} \frac{n}{n-1}\int_{\lVert t\rVert \leq 1/c} \frac{(1-|f_U^n(t)|^2)(1-|f_V^n(s)|^2)}{\lVert t\rVert^{1+d_1}\lVert s\rVert^{1+d_2}}\,dt\,ds\nonumber\\ &\overset{(b)}{=}\frac{n}{n-1}\cdot\frac{1}{n^4}\sum_{k,l,m,h}\int_{\lVert t\rVert \leq 1/c}\frac{1-\cos\langle t,U_k-U_l\rangle}{\lVert t\rVert^{1+d_1}}\cdot \frac{1-\cos\langle s,V_m-V_h\rangle}{\lVert s\rVert^{1+d_2}}\,dt\,ds\nonumber\\ &\overset{(c)}{\lesssim}\frac{n}{n-1}\cdot \frac{1}{n^4}\sum_{k,l,m,h} G\left(\frac{\lVert U_k-U_l\rVert}{c},\frac{U_k-U_l}{\lVert U_k-U_l\rVert}\right)\cdot \lVert V_m-V_h\rVert
     	\end{align} 
     	where (a) follows from Fubini's Theorem and the calculations from~\eqref{eq:covfnver}, (b) uses the fact that $\sin(\cdot)$ is an odd function and hence integrates to $0$ when integrated over symmetric sets, (c) uses the definition from~\eqref{eq:defg} and~\cite[Lemma 1]{Gabor2007}. The right hand side of~\eqref{eq:boundf1} converges to
     	\begin{align}\label{eq:firstlim}
     	\mathbb{E}\left[G\left(\frac{\lVert \tilde{U}_1-\tilde{U}_2\rVert}{c},\frac{\tilde{U}_1-\tilde{U}_2}{\lVert \tilde{U}_1-\tilde{U}_2\rVert}\right)\right]\cdot \mathbb{E}\bigg[\lVert \tilde{V}_1-\tilde{V}_2\rVert\bigg]
     	\end{align}
     	where $\tilde{U}_1,\tilde{U}_2\sim \mathcal{U}^{d_1}$ and $\tilde{V}_1,\tilde{V}_2\sim\mathcal{U}^{d_2}$ are four independent random variables. Finally, by an application of the dominated convergence theorem,~\eqref{eq:firstlim} converges to $0$ as $c\to\infty$. By the same calculation as in~\eqref{eq:firstlim}, we get:
     	\begin{align}\label{eq:boundf2}
     	\mathbb{E}\int_{\lVert t\rVert \geq c} |\xi_n(t,s)|^2w(t,s)\,dt\,ds &\lesssim \frac{n}{n-1}\cdot\frac{1}{n^4}\sum_{k,l,m,h}\lVert V_m-V_h\rVert\cdot \int_{\lVert t\rVert \geq c}\frac{\,dt}{\lVert t\rVert^{1+d_1}}\nonumber \\ &\lesssim \int_{\lVert t\rVert \geq c}\frac{\,dt}{\lVert t\rVert^{1+d_1}}. 
     	\end{align}
     	Clearly, the right hand side of~\eqref{eq:boundf2} converges to $0$ as limits are taken over $n\to\infty$ followed by $c\to\infty$. We can use the same arguments from~\eqref{eq:boundf1} and~\eqref{eq:boundf2} on the second term in the right hand side of~\eqref{eq:parineq} to get the same conclusion. Therefore, by an application of Markov's inequality, for any $\epsilon>0$,
     	\begin{align*}
     	\lim_{c\to\infty}\limsup_{n\to\infty}\mathbb{P}[|W_{n,c}-D_n|>\epsilon]\leq \lim_{c\to\infty}\limsup_{n\to\infty}\frac{1}{\epsilon}\cdot \mathbb{E}[|W_{n,c}-D_n|]=0.
     	\end{align*}  
     	This completes the proof.
     \end{proof}
     \begin{proof}[Proof of~\cref{lem:singlim}]
     	This proof is exactly the same as that of~\cref{lem:doublimit} and we leave the details to the reader. One can also use the tightness of $D$ (see~\eqref{eq:fluctuation}) as shown in~\cite[Chapter 1, Section 2]{Kuo1975}.
     \end{proof}
     
     \subsection{Proof of~\cref{lem:rankenergy}}\label{proof:rankenergy}
     Note that, for any $a\in\mathcal{S}^{d-1}$ and $U_i\in [0,1]^d$, $|a^{\top}U_i|\leq \lVert U_i\rVert\leq\sqrt{d}$. Therefore, 
     \begin{align*}
     E_{m,n}=\int_{\mathcal{S}^{d-1}}\int_{-\sqrt{d}}^{\sqrt{d}}\Theta_{m,n}^2(a,r)\,dr\,d\kappa(a)\;\;\mbox{and}\;\;E=\int_{\mathcal{S}^{d-1}}\int_{-\sqrt{d}}^{\sqrt{d}}\Theta^2(a,r)\,dr\,d\kappa(a).
     \end{align*}
     From~\cref{lem:steinweaklim} and~\cref{cor:applystein}, we have convergence (weakly and in second moments) of the finite dimensional distributions of the process $\Theta_{m,n}(a,r)$, $a\in\mathcal{S}^{d-1}$, $r\in [-\sqrt{d},\sqrt{d}]$. Next note that $\Theta_{m,n}(a,r)$ may be rewritten as,
     \begin{align*}
     \Theta_{m,n}(a,r)=\sqrt{\frac{n}{m}}\cdot\sqrt{m+n}\left(\frac{1}{N}\sum_{i=1}^N \mathbf{1}(a^{\top}U_i\leq t)-\frac{1}{n}\sum_{j=m+1}^{m+n}\mathbf{1}(a^{\top}U_{\pi_N(j)}\leq t)\right).
     \end{align*}
     Further observe that the set $\mathcal{F}\coloneqq \{\mathbf{1}(a^{\top}\cdot\leq t):(a,t)\in\mathcal{S}^{d-1}\times [-\sqrt{d},\sqrt{d}]\}$ of indicator functions on closed half-spaces is a VC class with index $d+1$ (see e.g.,~\cite{Dudley1978}) and consequently satisfies the uniform entropy condition, as in~\cite[Equation 2.5.1]{vaart1996} (see e.g.,~\cite[Theorem 2.6.7]{vaart1996}). The asymptotic equicontinuity of the process $\Theta_{m,n}(a,r)$ over $\mathcal{S}^{d-1}\times [-\sqrt{d},\sqrt{d}]$ then follows from the same proof as in~\cite[Theorem 2.5.2]{vaart1996} (as it uses similar empirical process tools as in the proof of~\cref{lem:rankdcov}, we leave the details to the interested reader). This then implies that $\Theta_{m,n}(\cdot,\cdot)$ converges weakly to $\Theta(\cdot,\cdot)$ in $L^{\infty}\big(\mathcal{S}^{d-1}\times [-\sqrt{d},\sqrt{d}]\big)$ (see~\cite{Hoffman1991}). The weak convergence of $E_{m,n}$ to $E$ then follows from a direct application of the continuous mapping theorem.
     
     \begin{lemma}\label{lem:steinweaklim}
     	Recall the notation and assumptions introduced in~\cref{lem:rankenergy}. Consider a $K$-tuple, $(a_1,r_1),\ldots ,(a_K,r_K)$, where $(a_i,r_i)\in \mathcal{S}^{d-1}\times\mathbb{R}$. Then the vector $\big(\Theta_{m,n}(a_1,r_1),\ldots ,\Theta_{m,n}(a_K,r_K)\big)$ converges weakly to a multivariate Gaussian distribution with mean $0$ and covariance matrix $\Sigma_{K\times K}$, where $\Sigma_{ij}=C\big((a_i,r_i),(a_j,r_j)\big)$ (given in~\eqref{eq:covenergy}), as $\min{(m,n)}\to\infty$. 
     \end{lemma}
     
     \begin{proof}
     	For the sake of simplicity, we will work with $K=2$. Set $\alpha^{\top}\coloneqq (\alpha_1,\alpha_2)\in\mathbb{R}^2$ and $\Theta_{m,n}^{\top}\coloneqq\big(\Theta_{m,n}(a_1,r_1),\Theta_{m,n}(a_2,r_2)\big)$. It suffices to show (by the Cram\'er-Wold Theorem), $\alpha^{\top}\Theta_{m,n}\overset{w}{\longrightarrow}\mathcal{N}(0,\alpha^{\top}\Sigma \alpha)$ as $\min{(m,n)}\to\infty$. Our proof proceeds using Stein's method of exchangeable pairs, see e.g.,~\cite{Cha2011}. For the reader's convenience, we also present this result in~\cref{prop:Steinexch}. Define $T_{m,n}\coloneqq\alpha^{\top}\Theta_{m,n}$. Draw two random indices $I$ and $J$, without replacement, from the set $\{1,2,\ldots ,N\}$. Construct a new permutation, $\tilde{\pi}_N$ as $\tilde{\pi}_N(I)=\pi_N(J)$, $\tilde{\pi}_N(J)=\pi_N(I)$ and $\tilde{\pi}_N(k)=\pi_N(k)$ for $k\neq I,J$. It is easy to check that $(\pi_N,\tilde{\pi}_N)$ forms an exchangeable pair of random vectors. Let $\tilde{T}_{m,n}\coloneqq \alpha^{\top}\tilde{\Theta}_{m,n}$ where $\tilde{\Theta}_{m,n}$ is calculated by replacing $\pi_N$ with $\tilde{\pi}_N$ in $\theta_{m,n}$. Note that,
     	\begin{align}\label{eq:exchange1}
     	&\mathbb{E}[T_{m,n}-\tilde{T}_{m,n}|\pi_N]\nonumber \\&=\mathbb{E}\bigg[\frac{2\alpha_1\sqrt{m+n}}{\sqrt{mn}}\bigg(\mathbf{1}(a_1^{\top}U_{\pi_N(I)}\leq r_1)-\mathbf{1}(a_1^{\top}U_{\pi_N(I)}\leq r_1)\bigg)\mathbf{1}(I\leq m,J\geq m+1)\nonumber\\ &+\frac{2\alpha_2\sqrt{m+n}}{\sqrt{mn}}\bigg(\mathbf{1}(a_2^{\top}U_{\pi_N(I)}\leq r_2)-\mathbf{1}(a_2^{\top}U_{\pi_N(I)}\leq r_2)\bigg)\mathbf{1}(I\leq m,J\geq m+1)\bigg|\pi_N\bigg]\nonumber\\ &=2(m+n-1)^{-1}T_{m,n}
     	\end{align}
     	which in turn implies, $\mathbb{E}[T_{m,n}-\tilde{T}_{m,n}|T_{m,n}]=2(m+n-1)^{-1}T_{m,n}$. Define $c_0\coloneqq(m+n-1)(2\alpha^{\top}\Sigma\alpha)^{-1}$. Note that $|T_{m,n}-\tilde{T}_{m,n}|\leq 2(|\alpha_1|+|\alpha_2|)(\min{(m,n)})^{-1/2}$. By~\cite[Theorem 1.2]{Cha2011}, our desired conclusion follows if we can show the following:
     	\begin{align}\label{eq:steincon}
     	\mathbb{E}\bigg|1-\frac{c_0}{2}\mathbb{E}\left[(T_{m,n}-\tilde{T}_{m,n})^2\big|T_{m,n}\right]\bigg|\overset{\min{(m,n)}\to\infty}{\longrightarrow}0.
     	\end{align}
     	Note that,
     	\begin{align}\label{eq:versec}
     	&\;\;\;\mathbb{E}\left[\frac{c_0}{2}\cdot \big(T_{m,n}-\tilde{T}_{m,n}\big)^2\big|\pi_N\right]\nonumber\\ &=\frac{(m+n)(m+n-1)}{(2\alpha^{\top}\Sigma\alpha)(mn)}\mathbb{E}\bigg[\bigg\{\alpha_1\bigg(\mathbf{1}(a_1^{\top}U_{\pi_N(I)}\leq r_1)-\mathbf{1}(a_1^{\top}U_{\pi_N(J)}\leq r_1)\bigg)\mathbf{1}(I\leq m,\nonumber \\&J\geq m+1)+\alpha_2\bigg(\mathbf{1}(a_2^{\top}U_{\pi_N(I)}\leq r_2)-\mathbf{1}(a_2^{\top}U_{\pi_N(J)}\leq r_2)\bigg)\mathbf{1}(I\leq m,J\geq m+1)\bigg\}^2\bigg|\pi_N\bigg]\nonumber\\ &=\frac{\alpha_1^2}{(mn)(2\alpha^{\top}\Sigma\alpha)}\bigg[n\sum_{i=1}^m \mathbf{1}(a_1^{\top}U_{\pi_N(i)}\leq r_1)+m\sum_{j=m+1}^{m+n} \mathbf{1}(a_1^{\top}U_{\pi_N(j)}\leq r_1)\nonumber \\ &-2\sum_{i\leq m,j\geq m+1} \mathbf{1}(a_1^{\top}U_{\pi_N(i)}\leq r_1, a_1^{\top}U_{\pi_N(j)}\leq r_1)\bigg]+\frac{\alpha_2^2}{(mn)(2\alpha^{\top}\Sigma\alpha)}\bigg[n\sum_{i=1}^m \mathbf{1}(a_2^{\top}U_{\pi_N(i)}\leq r_2)\nonumber \\&+m\sum_{j=m+1}^{m+n} \mathbf{1}(a_2^{\top}U_{\pi_N(j)}\leq r_2)-2\sum_{i\leq m,j\geq m+1} \mathbf{1}(a_2^{\top}U_{\pi_N(i)}\leq r_2,a_2^{\top}U_{\pi_N(j)}\leq r_2)\bigg]\nonumber \\&+\frac{\alpha_1\alpha_2}{(mn)(2\alpha^{\top}\Sigma\alpha)}\bigg[n\sum_{i=1}^m \mathbf{1}(a_1^{\top}U_{\pi_N(i)}\leq r_1,a_2^{\top}U_{\pi_N(i)}\leq r_2)+m\sum_{j=m+1}^{m+n}\mathbf{1}(a_1^{\top}U_{\pi_N(j)}\leq r_,1\nonumber \\&a_2^{\top}U_{\pi_N(j)}\leq r_2)-2\sum_{i\leq m,j\geq m+1} \mathbf{1}(a_1^{\top}U_{\pi_N(i)}\leq r_1,a_2^{\top}U_{\pi_N(j)}\leq r_2)\bigg].
     	\end{align}
     	Further $\mathbb{E}[m^{-1}\sum_{i\leq m}\mathbf{1}(a_1^{\top}U_{\pi_N(i)}\leq r_1)]=N^{-1}\sum_{i\leq N}\mathbf{1}(a_1^{\top}U_i\leq r_1)\to\mathbb{P}(a_1^{\top}\tilde{U}\leq r_1)$ and  $\mbox{Var}[m^{-1}\sum_{i\leq m}\mathbf{1}(a_1^{\top}U_{\pi_N(i)}\leq r_1)]=\mathcal{O}(m^{-1})$. Therefore, $m^{-1}\sum_{i\leq m}\mathbf{1}(a_1^{\top}U_{\pi_N(i)}\leq r_1)\overset{\mathbb{P}}{\longrightarrow}\mathbb{P}(a_1^{\top}\tilde{U}\leq r_1)$. Similar arguments may be used to prove that
     	\begin{align}\label{eq:morecon}
     	&\frac{1}{m}\sum_{i=1}^m \mathbf{1}(a_1^{\top}U_{\pi_N(i)}\leq r_1,a_2^{\top}U_{\pi_N(i)}\leq r_2)\to\mathbb{P}(a_1^{\top}\tilde{U}\leq r_1,a_2^{\top}\tilde{U}\leq r_2)\;\;\qquad\qquad\mbox{and}\nonumber \\&\frac{1}{mn}\sum_{i\leq m,j\geq m+1} \mathbf{1}(a_1^{\top}U_{\pi_N(i)}\leq r_1,a_2^{\top}U_{\pi_N(j)}\leq r_2)\to\mathbb{P}(a_1^{\top}\tilde{U}\leq r_1)\mathbf{P}(a_2^{\top}\tilde{U}\leq r_2).
     	\end{align}
     	Recall the definition of $C(\cdot,\cdot)$ from~\eqref{eq:covenergy}. Note that~\eqref{eq:morecon} implies~\eqref{eq:versec} converges in probability to
     	\begin{align}\label{eq:varmatch}
     	\frac{1}{\alpha^{\top}\Sigma\alpha}\left[\alpha_1^2C((a_1,r_1),(a_1,r_1))+2\alpha_1\alpha_2C((a_1,r_1),(a_2,r_2))+\alpha_2^2C((a_2,r_2),(a_2,r_2))\right]=1.
     	\end{align}
     	Finally, as $c_0(T_{m,n}-\tilde{T}_{m,n})^2$ is uniformly bounded,~\eqref{eq:varmatch} implies~\eqref{eq:steincon} by the dominated convergence theorem.
     \end{proof}
     \begin{corollary}\label{cor:applystein}
     	Recall the notation from the statement and proof of~\cref{lem:steinweaklim}. Then $\mathbb{E}[T_{m,n}^2]\to\alpha^{\top}\Sigma\alpha$ as $\min{(m,n)}\to\infty$.
     \end{corollary}
     \begin{proof}
     	In the proof of~\cref{lem:steinweaklim}, we showed that $(c_0/2)\mathbb{E}(T_{m,n}-\tilde{T}_{m,n})^2\to 1$. Considering all limits to be under $\min{(m,n)}\to\infty$, we get: 
     	\begin{align*}\label{eq:proofapplystein}
     	1=\lim (c_0/2)\mathbb{E}(T_{m,n}-\tilde{T}_{m,n})^2&=\lim c_0\mathbb{E}\big[T_{m,n}\mathbb{E}[T_{m,n}-\tilde{T}_{m,n}|T_{m,n}]\big]\nonumber \\&\overset{(a)}{=}\lim 2c_0(m+n-1)^{-1}\mathbb{E}[T_{m,n}^2]
     	\end{align*}
     	which completes the proof. Here, (a) follows from~\eqref{eq:exchange1}.
     \end{proof}
     \section{Auxiliary Results}\label{sec:auxi}
     \begin{lemma}[Alexandroff Theorem,~\citet{Alexandroff1939}]\label{lem:Alexandroff}
     	Let $f:U\to\mathbb{R}$ be a convex function, where $U$ is an open convex subset of $\mathbb{R}^n$. Then $f$ has a second derivative Lebesgue a.e.~in $U$.	
     \end{lemma}
     \begin{lemma}(Almost sure weak convergence of empirical measure,~\citet{Varadarajan1958})\label{lem:varada}
     	Let $(W,d)$ be a separable metric space and $\mu$ be a probability measure supported on $W$. Also, say $\mu_n$ denotes the empirical counterpart of $\mu$. Then $d_{W}(\mu_n,\mu)\overset{a.s.}{\longrightarrow}0$ where $d_W(\cdot,\cdot)$ is any metric on the space of probability measures on $(W,d)$ that equivalently characterizes weak convergence.
     \end{lemma}
     \begin{lemma}(\citet[Lemma 9]{Mccann1995})\label{lem:Mcclem9}
     	Suppose $\mu_n\in \mathcal{P}(\mathbb{R}^d\times\mathbb{R}^d)$ converges weakly to $\mu\in\mathcal{P}(\mathbb{R}^d\times\mathbb{R}^d)$. Then,
     	\begin{itemize}
     		\item[(i)] If $\mu_n$ has cyclically monotone support for all large $n$, then so does $\mu$.
     		\item[(ii)] Let $\Gamma(\nu_1,\nu_2)$ denote the subset of $\mathcal{P}(\mathbb{R}^d\times\mathbb{R}^d)$ with first and second marginals $\nu_1$ and $\nu_2$ respectively; $\nu_1,\nu_2\in\mathcal{P}(\mathbb{R}^d)$. If $\mu_n\in\Gamma(\nu_1^n,\nu_2^n)$ where $\nu_1^n\overset{w}{\longrightarrow}\nu_1$ and $\nu_2^n\overset{w}{\longrightarrow}\nu_2$, then $\mu\in\Gamma(\nu_1,\nu_2)$.
     	\end{itemize}
     \end{lemma}
     \begin{definition}[Cyclically monotone maps]\label{def:cycmon}
     	A subset $S$ of $\mathbb{R}^d\times\mathbb{R}^d$ is said to be {\it cyclically monotone} if, given any finite subset of ${S}$, say $\{(\mathbf{x}_1,\mathbf{y}_1),\ldots ,(\mathbf{x}_k,\mathbf{y}_k)\}$, we have:
     	\begin{equation*}
     	\langle \mathbf{y}_1,\mathbf{x}_2-\mathbf{x}_1\rangle+\ldots +\langle \mathbf{y}_{k-1},\mathbf{x}_k-\mathbf{x}_{k-1}\rangle+\langle \mathbf{y}_k,\mathbf{x}_1-\mathbf{x}_k\rangle\leq 0.
     	\end{equation*}
     	A multi-valued map $f:\mathbb{R}^d\to\mathbb{R}^d$ is said to be a {\it cyclically monotone map} if, given any finite subset $\{\mathbf{x}_1,\ldots ,\mathbf{x}_k\}$ of $\mathbb{R}^d$, the set $\{(\mathbf{x}_1,f(\mathbf{x}_1)),\ldots ,(\mathbf{x}_k,f(\mathbf{x}_k))\}$ is {\it cyclically monotone}. 
     \end{definition}
     \begin{definition}[Subdifferential of a convex function]\label{def:subdiffconv}
     	Let $f:\mathbb{R}^d\to\mathbb{R}$ be a proper, lower semicontinuous convex function. Then the {\it subdifferential} of $f(\cdot)$ at a point $\mathbf{x}\in\mathbb{R}^d$ is defined as:
     	\begin{equation*}
     	\partial f(\mathbf{x})\coloneqq \{\mathbf{z}:f(\mathbf{y})-f(\mathbf{x})\geq \langle \mathbf{z},\mathbf{y}-\mathbf{x}\rangle,\mbox{ }\mbox{ for all } \mathbf{y}\in\mathbb{R}^d\}.
     	\end{equation*}
     \end{definition}
     \begin{lemma}(Cyclic monotonicity and subdifferential of convex functions;~\citet[Theorem 1]{Rockafellar1966})\label{lem:Rockcyc}
     	The graph of the subdifferential $\partial f(\cdot)$ of a convex function $f:\mathbb{R}^d\to\mathbb{R}$ is a cyclically monotone subset of $\mathbb{R}^d\times\mathbb{R}^d$. Moreover, any cyclically monotone subset of $\mathbb{R}^d\times\mathbb{R}^d$ is contained in the graph of the subdifferential of a proper, lower semicontinuous convex function from $\mathbb{R}^d\to\mathbb{R}$.
     \end{lemma}
     \begin{lemma}(Uniqueness of measure preserving couplings, see~\citet[Corollary 14]{Mccann1995})\label{lem:Mccor14}
     	Let $\nu_1,\nu_2\in\mathcal{P}(\mathbb{R}^d)$, and suppose that one of these two measures in Lebesgue absolutely continuous. Then, there exists one and only one measure $\nu\in\Gamma(\nu_1,\nu_2)$ (see (ii) from~\cref{lem:Mcclem9}) with cyclically monotone support.
     \end{lemma}
     \begin{lemma}(Existence of measure transformation maps; see~\cite[Proposition 10]{Mccann1995})\label{lem:MccProp10}
     	Assume that $\nu\in\Gamma(\nu_1,\nu_2)$ (see (ii) from~\cref{lem:Mcclem9}) is supported on the graph of the subdifferential $\partial f(\cdot)$ of some proper, lower semicontinuous convex function $f:\mathbb{R}^d\to\mathbb{R}$ (i.e., the support of $\nu$ is a subset of the graph of $\partial f(\cdot)$). Further, suppose that $\nu_1\in\mathcal{P}_{ac}(\mathbb{R}^d)$. Then $\nabla f(\cdot)$ pushes $\nu_1$ to $\nu_2$, i.e., $\nu=(\mbox{identity}\times \nabla f)\#\nu_1$.
     \end{lemma}
     \begin{prop}(Hoeffding's Central Limit theorem; see~\citet[Theorem 1.1]{Chen2015})\label{prop:HCLT}
     	Suppose $\mathfrak{X}=\{X_{ij}:1\leq i,j\leq n\}$ be a $n\times n$ array of independent random variables where $n\geq 2$, $\mathbb{E}[X_{ij}]=c_{ij}$, $\mathrm{Var}(X_{ij})=\sigma^2_{ij}\geq 0$ and $\mathbb{E}|X_{ij}|^3<\infty$. Assume that,
     	\begin{equation*}
     	c_{i\cdot}\coloneqq \frac{1}{n}\sum_j c_{ij}=0\qquad \mbox{and}\qquad c_{\cdot j}\coloneqq \frac{1}{n}\sum_i c_{ij}=0.
     	\end{equation*}
     	Let $\pi_n$ be an uniform permutation drawn from $S_n$ (all permutations of $\{1,2,\ldots ,n\}$) independent of $\mathfrak{X}$. Let $W_n=\sum_i X_{i\pi_n(i)}$. Then,
     	\begin{equation*}
     	\mathrm{Var}(W_n)=\frac{1}{n}\sum_{i,j}\sigma^2_{ij}+\frac{1}{n-1}\sum_{i,j}c^2_{ij}.
     	\end{equation*}
     	Further, if $\mathrm{Var}(W_n)=1$, then we have:
     	\begin{equation*}
     	\sup_{z\in\mathbb{R}} \bigg|\mathbb{P}(W_n<z)-\Phi(z)\bigg|\leq \frac{451}{n}\sum_{i,j} \mathbb{E}|X_{ij}|^3
     	\end{equation*}
     	where $\Phi(\cdot)$ denotes the standard Gaussian distribution function.
     \end{prop}
     \begin{prop}(Stein's method of exchangeable pairs; see~\citet[Theorem 1.2]{Cha2011})\label{prop:Steinexch}
     	Let $(W,W')$ denote an exchangeable pair of random variables where $W$ has finite variance. Also, suppose that,
     	\begin{equation*}
     	\mathbb{E}(W-W'|W)=g(W)+r(W)\qquad \mathrm{and} \qquad |W-W'|\leq \delta,
     	\end{equation*}
     	where $\delta$ is a constant (non-random) and $g(\cdot)$ is a differentiable function on $\mathbb{R}$. Let $G(t)=\int_0^t g(s)\,ds$ and  $p(t)=c_1\exp(-c_0 G(t))$, where $c_1=(\exp(-c_0G(t)))^{-1}>0$. Here $c_0$ is some positive real number. Further, let us also assume the following conditions:
     	\begin{itemize}
     		\item[(i)] $g(\cdot)$ is nondecreasing, $g(t)\geq 0$ for $t\geq 0$ and $g(t)\leq 0$ for $t\leq 0$.
     		\item[(ii)] There exists $c_2<\infty$ such that for all $x\in\mathbb{R}$, 
     		\begin{equation*}
     		\min{\big(1/c_1,1/|c_0 g(x)|\big)}\big(|x|+3/c_1\big)c_0|g'(x)|\leq c_2.
     		\end{equation*}
     	\end{itemize}
     	Finally, let $Y$ be a random variable with density $p_1(\cdot)$. Under all the above conditions, the following holds:
     	\begin{align*}
     	\sup_{z\in\mathbb{R}} \big|\mathbb{P}(W\leq z)-\mathbb{P}(Y\leq z)\big|&\leq 3\big|1-(c_0/2)\mathbb{E}[(W-W')^2|W]\big|+c_1\max{(1,c_2)}\delta\nonumber \\&+2(c_0/c_1)\mathbb{E}\big|r(W)\big|+\delta^3 c_0\big\{(2+c_2)/2\mathbb{E}\big|c_0g(W)\big|+c_1c_2/2\big\}.
     	\end{align*}
     \end{prop}
\end{document}